\documentclass[a4paper, 10pt, abstracton]{scrartcl}

\usepackage{amsmath}
\usepackage{amssymb}
\usepackage{graphicx}
\usepackage{epstopdf}
\usepackage{inputenc}
\usepackage{geometry}
\usepackage{a4wide}

\usepackage[font=small]{caption}
\usepackage{caption,subcaption}
\usepackage{cite}
\usepackage{calc}
\usepackage{float}
\usepackage{grffile} 
\usepackage{enumitem}

\usepackage{color}
\usepackage{xcolor}

\usepackage{mathtools}

\usepackage{comment} 

\newcommand{\U}{\mathcal{U}}

\newcommand{\R}{\mathbb{R}}
\newcommand{\X}{\mathcal{X}}

\newcommand{\Np}{N}
\newcommand{\pb}{\pmb{x}_*}
\newcommand{\gb}{\hat{\pmb{x}}_*}

\usepackage{algpseudocode}
\usepackage{algorithm}

\usepackage{multirow}
\usepackage{rotating}

\usepackage{authblk}

\usepackage{lscape} 

\usepackage[final]{pdfpages}

\usepackage[normalem]{ulem} 

\algnewcommand{\IIf}[1]{\State\algorithmicif\ #1\ \algorithmicthen} 
\algnewcommand{\EndIIf}{\unskip\ \algorithmicend\ \algorithmicif} 

\usepackage{mathrsfs} 
\usepackage{eufrak} 

\begin{document}

\title{Automatic Generation of Algorithms for Black-Box Robust Optimisation Problems \thanks{Partially funded through EPSRC grants EP/L504804/1 and EP/M506369/1.}}
\author[1]{Martin Hughes\footnote{Corresponding author. Email: \texttt{m.hughes6@lancaster.ac.uk}} }
\author[2]{Marc Goerigk}
\author[1]{Trivikram Dokka}
\affil[1]{Department of Management Science, Lancaster University, United Kingdom}
\affil[2]{Network and Data Science Management, University of Siegen, Germany}

\date{} 

\maketitle


\abstract{We develop algorithms capable of tackling robust black-box optimisation problems, where the number of model runs is limited. When a desired solution cannot be implemented exactly the aim is to find a robust one, where the worst case in an uncertainty neighbourhood around a solution still performs well. This requires a local maximisation within a global minimisation.

To investigate improved optimisation methods for robust problems, and remove the need to manually determine an effective heuristic and parameter settings, we employ an automatic generation of algorithms approach: Grammar-Guided Genetic Programming. We develop algorithmic building blocks to be implemented in a Particle Swarm Optimisation framework, define the rules for constructing heuristics from these components, and evolve populations of search algorithms. Our algorithmic building blocks combine elements of existing techniques and new features, resulting in the investigation of a novel heuristic solution space.

As a result of this evolutionary process we obtain algorithms which improve upon the current state of the art. We also analyse the component level breakdowns of the populations of algorithms developed against their performance, to identify high-performing heuristic components for robust problems.
}



\section{Introduction}
\label{sec:introduction}

The use of optimisation search techniques to investigate a decision variable solution space and identify good solutions is common when using models to support informed decision making. However the search may be impacted by issues such as model run times, the size of the solution space, and uncertainty, see \cite{BenTalElGhaouiNemirovski2009, GoerigkSchobel2016}. In this work we are concerned with optimisation under implementation uncertainty, and where some budget on the number of model runs restricts the search.

If a model can take the form of a mathematical program, optimisation may be tackled efficiently and exactly. Here we assume this is not the case, and instead some model is employed which from an optimisation perspective can be considered a black-box where decision variable values are input and an objective extracted. In this case only an approximate global optimum is sought, and so in this work we consider metaheuristic techniques applicable to general, likely non-convex problems.

With implementation uncertainty an ideal solution cannot be achieved exactly, so solutions are sought where all points in the uncertainty neighbourhood around a candidate still perform well. When it is known how the uncertainty is distributed, the problem is one of stochastic optimisation, see \cite{PaenkeBrankeJin2006, HomemdeMelloBayraksan2014}. Instead we assume the uncertainty takes the form of some set containing all uncertainty scenarios such as an interval, making the problem one of robust optimisation. Specifically a classic robust setting is considered, where the worst (inner maximum) model value in the uncertainty region around a candidate solution is sought in the context of an overarching (outer) minimisation objective, \cite{BenTalNemirovski1998}.

Here our aim is to develop improved search techniques through the automatic generation of metaheuristics, actively seeking good heuristics and avoiding the need for the manual determination of the search algorithm and parameter settings. A hyper-heuristic approach is employed, genetic programming (GP) \cite{Nohejl2011}, an evolutionary process where each individual in a population is an algorithm -- here a metaheuristic for a robust problem. From the initial population some measure of fitness is determined for each heuristic, and a new generation established through typical evolutionary selection, combination and mutation processes. After multiple generations the fittest heuristic is chosen and applied to the problem at hand.

To facilitate the GP search, heuristic sub-components are generated. When combined correctly these algorithmic building blocks form a complete heuristic. The sub-components form a language, and the design rules specifying how they combine to create a heuristic represent a grammar. This is Grammar-Guided Genetic Programming (GGGP) \cite{Nohejl2011}.

As with any evolutionary approach, GP employs combination and mutation operations to generate improved (fitter) solutions. However integrating sub-algorithms (computer sub-programs in the more general GP sense), may not be straightforward when the intention is to form a coherent, executable higher level algorithm. A common GGGP approach uses a tree-based representations of the overarching algorithm \cite{MasciaLopezIbanezDuboisLacosteStutzle2014, ContrerasBoltonParada2015, MirandaPrudencio2016, MirandaPrudencio2017}. This approach-representation is adopted here, where we specify heuristic sub-components in terms of a context-free grammar (CFG) and use standard tree-based random combination and mutation operators \cite{MirandaPrudencio2016}.

\paragraph{Contributions and outline.} Improved global metaheuristics are developed for robust black-box problems under implementation uncertainty, for problems of 30 dimensions (30D) and 100D and assuming a budget of 2,000 model runs. A GGGP search of the solution space of heuristics for robust  problems is used to identify the best approaches. The previously uninvestigated heuristic solution space comprises algorithmic building blocks that combine to form a complete particle swarm based heuristic. A large number of sub-components are developed using existing approaches and novel implementations.

New algorithms are tested on a suite of problems, and improved heuristics for general robust problems are identified. The significance of individual algorithmic sub-components is also assessed against heuristic performance. The effectiveness of an inner maximisation by random sampling on a small number of points and using a particle level stopping condition, is established. For the outer minimisation a small swarm of particles performs well, as does communication via a Global typology. The preferred particle movement uses an inertia based velocity equation plus specialised capabilities drawn from the largest empty hypersphere \cite{HughesGoerigkWright2019, HughesGoerigkDokka2020a} and descent directions \cite{BertsimasNohadaniTeo2007, BertsimasNohadaniTeo2010nonconvex, BertsimasNohadaniTeo2010, HughesGoerigkDokka2020a} heuristics.   

In Section~\ref{sec:robustOpti} we outline the optimisation problem of concern here, and current approaches for addressing it. We include descriptions of the heuristics that form the basis for the building blocks in the GP analysis. Section~\ref{sec:autoGenAlgos} gives an overview on the automatic generation of algorithms, and in Section~\ref{sec:autoRob} GP is discussed in detail. This includes sub-component descriptions, the design rules for constructing complete heuristics, and our GP approach including tree-based representation and operators. Section~\ref{sec:experiments} describes the experimental analysis, results for the best heuristics identified, and an analysis of heuristic sub-component performance. Section~\ref{sec:concusionsFurtherWork} provides conclusions and possible directions for future work.


\section{Robust optimisation}
\label{sec:robustOpti}

\subsection{Problem description}
\label{sec:minMax}

A general optimisation problem without consideration of uncertainty takes the form:

\begin{align*}
	\quad \min\ & f(\pmb{x}) \\
	\text{s.t. } & \pmb{x} \in \X
\end{align*}
\vspace{1mm}

\begin{figure}[htb]
	\centering
	\includegraphics[width=0.6\textwidth]{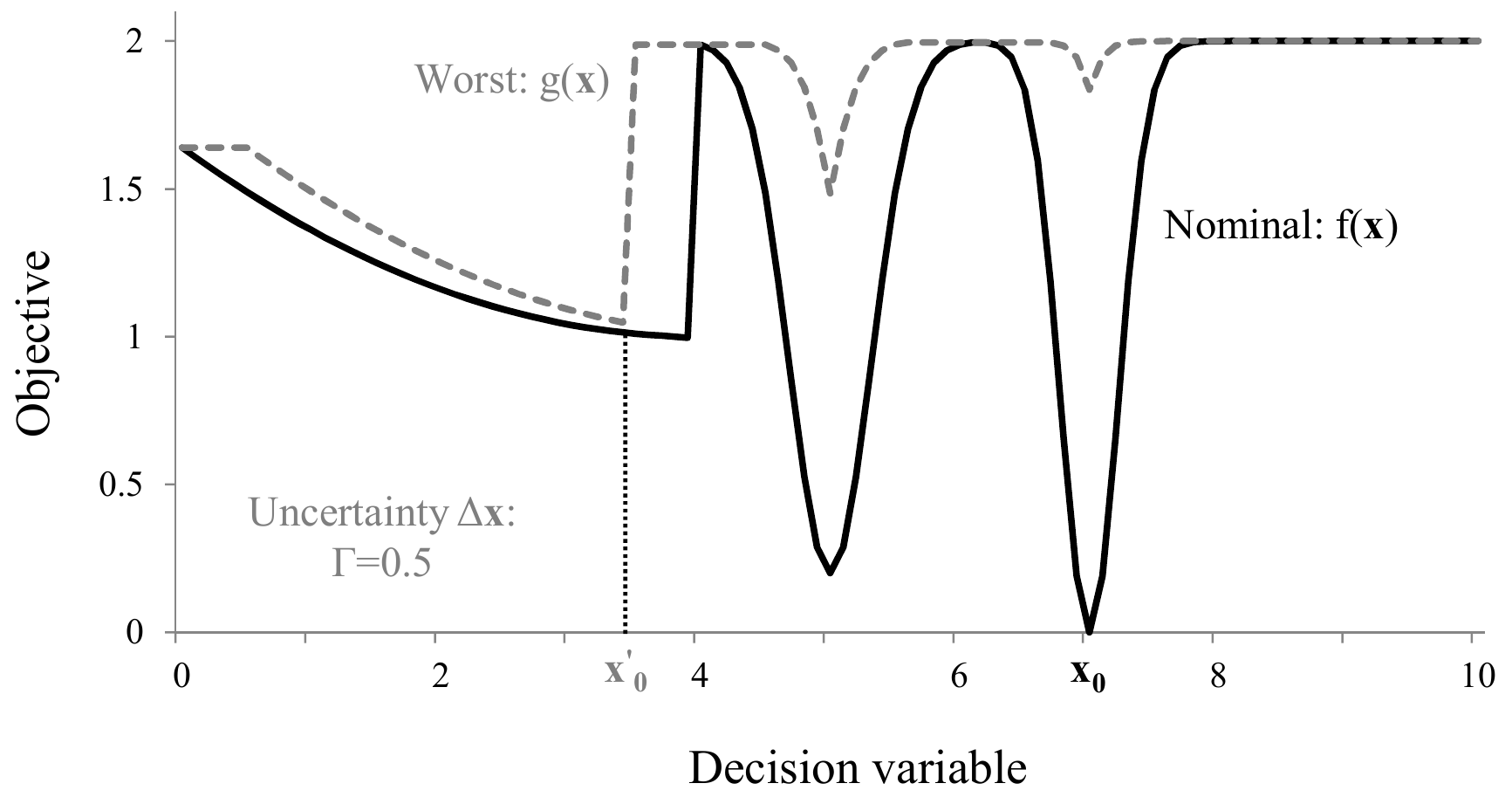} 
	\caption{The worst case cost curve (dashed grey) is generated by determining the maximum objective value in the uncertainty neighbourhood around all points $\pmb{x}$ on the nominal (solid black) curve. Due to the uncertainty the global optimum shifts to $\pmb{x}'_0$.}
\label{fig:Kruisselbrink}
\end{figure}

This is the nominal problem. The objective $f: \R^n \to \R$ operates on the $n$-dimensional vector of decision variables $\pmb{x}=(x_{1}, x_{2}, \ldots, x_{n})^T$ in the feasible region $\X \subseteq \R^n$. Here we assume box constraints $\X = \prod_{i\in[n]} [l_i,u_i]$. Any other feasibility constraint is assumed to be dealt with by a penalty in the objective. The notation $[n]:=\{1,\ldots,n\}$ is used. Consider the problem due to \cite{Kruisselbrink2012} in Figure~\ref{fig:Kruisselbrink}, where $\X \subseteq \R^1$, $l_1=0$ and $u_1=10$. The nominal problem is the black curve.

\vspace{4mm}
\vspace{4mm}

Introducing uncertainty $\Delta\pmb{x}$ around the intended solution $\pmb{x}$, makes only a solution $\tilde{\pmb{x}}=\pmb{x}+\Delta\pmb{x}$ achievable. If it is assumed that the uncertainty neighbourhood around a candidate is completely defined by a radius $\Gamma > 0$, the uncertainty set is \cite{BertsimasNohadaniTeo2010}:
\[
	\U:=\{ \Delta \pmb{x}\in\R^n \mid \| \Delta \pmb{x} \| \leq \Gamma \}
\]
where $\|\cdot\|$ represents the Euclidean norm. Using a local maximisation to find a robust solution $\pmb{x}$, the worst case value $g(\pmb{x})$ is optimised for any $\tilde{\pmb{x}}$ in the uncertainty neighbourhood of $\pmb{x}$:
\[
	g(\pmb{x}):=\max_{\Delta \pmb{x} \in \U} f(\pmb{x} + \Delta \pmb{x})
\]
In Figure~\ref{fig:Kruisselbrink}, $\Gamma=0.5$ and each point on the worst case cost (dashed grey) curve $g(\pmb{x})$ is generated by finding the maximum value on the nominal curve within a range of $-0.5$ to $+0.5$ of the desired solution $\pmb{x}$. 

The complete min max robust problem then involves finding the outer minimum objective in $\X$, where that objective is itself an inner maximisation in the uncertainty neighbourhood around each solution $\pmb{x}\in\X$ for the nominal objective function:
\[
	\min_{\pmb{x}\in\X} g(\pmb{x}) = \min_{\pmb{x}\in\X} \max_{\Delta\pmb{x} \in \U} f(\pmb{x} + \Delta \pmb{x}) \tag{MinMax}
\]
In the example this moves the global minimisation search from the black curve to the grey curve, where the global optimum shifts from $\pmb{x}_0$ to $\pmb{x}'_0$.

As $\pmb{x}+\Delta\pmb{x}$ may be outside of $\X$, here it is not assumed that $f$ is restricted to $\X$. If it is required that $\pmb{x}+\Delta\pmb{x}\in\X$ for all $\Delta\pmb{x}\in\U$, this could be achieved for example through the reduction of the original $\X$ by $\Gamma$.

\subsection{State of the art}
\label{sec:literature}

Since its initial formalisation \cite{KouvelisYu1997, BenTalNemirovski1998} robust optimisation has been heavily aligned with mathematical programming, see \cite{bertsimas2011theory,gabrel2014recent,GoerigkSchobel2016}. Where mathematical programming techniques cannot be applied, metaheuristics may be considered. However only limited consideration has been given to robust black-box optimisation under implementation uncertainty, see \cite{MarzatWalterPietLahanier2013, GoerigkSchobel2016, MarzatWalterPietLahanier2016}.

Whilst standard metaheuristics can be extended to the robust worst case through the brute force addition of an inner maximisation routine into an outer minimisation setting e.g. \cite{HughesGoerigkWright2019, HughesGoerigkDokka2020a}, more refined robust-specific methods may be desirable. Such techniques include co-evolutionary approaches \cite{Herrmann1999, ShiKrohling2002, Jensen2004, CramerSudhoffZivi2009, MasudaKuriharaAiyoshi2011}, robust evolutionary approaches \cite{TsutsuiGhosh1997, BeyerSendhoff2007}, and the use of emulation by surrogates or meta-models \cite{OngNairLum2006, BeyerSendhoff2007, ZhouZhang2010, VuDAmbrosioHamadiLiberti2016} including Kriging \cite{MarzatWalterPietLahanier2013, urRehmanLangelaarvanKeulen2014, MarzatWalterPietLahanier2016, urRehmanLangelaar2017} and Bayesian techniques \cite{ChenLucierSingerSyrgkanis2017, SandersEversonFieldsendRahat2019}. However specific assumptions or simplifications are typically required for such methods to be effective, or there are limitations on the problems that can be addressed e.g. due to dimensionality. 

Two existing general robust approaches requiring no further assumptions or simplifications are given specific attention here. They form the basis for some of the algorithmic building blocks which constitute the grammar in our GP analysis. These are the local descent directions (d.d.) approach \cite{BertsimasNohadaniTeo2007, BertsimasNohadaniTeo2010nonconvex, BertsimasNohadaniTeo2010} and the global largest empty hypersphere (LEH) method \cite{HughesGoerigkWright2019}. Both are single-solution techniques, although elements of these approaches have been incorporated into robust population-based approaches \cite{HughesGoerigkDokka2020a}. First, however, we consider the non-robust particle swarm optimisation (PSO) metaheuristic \cite{KennedyEberhart1995, KennedyEberhartShi2001, Ghazali2009}, which is the basis for all heuristics in our GP search. Constituent elements of a typical PSO algorithm are included as building blocks in our GP grammar.

\subsubsection{Particle swarm optimisation}
\label{sec:pso}

PSO is a population-based approach which moves a `swarm' of particles through points in the decision variable space, performing function evaluations and iterating particle positions through particle level `velocity' vectors. Velocities are based on particle histories, shared information from the swarm, scaling, and randomisation. The intention is for the behaviour of this complex systems of particles to approximate a global optimisation search. There are very many PSO formulations building on this general concept, see for example \cite{Kameyama2009, NickabadiMehdiEbadzadehSafabakhsh2011, ZhangWangJi2015, Kiran2017, SenguptaBasakPeters2018}. 

In a robust setting two-swarm co-evolutionary PSO techniques have been considered \cite{ShiKrohling2002, MasudaKuriharaAiyoshi2011}, whilst \cite{HughesGoerigkWright2019} employ a simple robust PSO (rPSO)  as a comparator heuristic, with inner maximisation by random sampling. The rPSO from \cite{HughesGoerigkWright2019} is extended in \cite{HughesGoerigkDokka2020a} through the addition of d.d.\ and LEH elements. Here the framework for each heuristic in the GP population is a basic PSO formulation \cite{KennedyEberhart1995, KennedyEberhartShi2001, Ghazali2009}, built upon through the availability of more complex algorithmic building blocks in our grammar. 

In a basic non-robust PSO formulation, the swarm (population) of $\Np$ particles start at iteration $t=0$ randomly located at points $\pmb{x}^j(0)$ in $\X$, where the function is evaluated; here $j=1, \ldots, \Np$. Each particle stores information on the best position it has visited in its history, $\pb^j$. Best refers to the lowest objective function value.

Information sharing is a key element of PSO, with each particle associated with a neighbourhood of other particles. Within a neighbourhood information is shared on the best point visited by any neighbourhood particle in their entire histories, $\gb$. A number of neighbourhood topologies are included as components in the grammar here, see Section~\ref{sec:networks}.

A particle is moved to a location $\pmb{x}^j(t)$ at iteration $t$, through the addition of that particle's current velocity vector $\pmb{v}^j$ to its previous position:
\[
	\pmb{x}^j(t)=\pmb{x}^j(t-1) \> + \> \pmb{v}^j(t)  \tag{PSOmove}
\]
There are a number of alternative velocity formulations. In the grammar in Section~\ref{sec:CFG} two of the most basic formulations are considered, including so-called inertia \cite{ShiEberhart1998, KennedyEberhartShi2001} and constriction \cite{ClercKennedy2002, KennedyEberhartShi2001} coefficients: 
\[
	\pmb{v}^j(t)=\omega \cdot \pmb{v}^j(t-1) \> + \> C_1\cdot\pmb{r}_1\cdot(\pb^j-\pmb{x}^j(t-1)) \> + \> C_2\cdot\pmb{r}_2\cdot(\gb-\pmb{x}^j(t-1))  \tag{Inertia}
\]
\[
	\pmb{v}^j(t)=\chi \cdot \Bigg( \pmb{v}^j(t-1) \> + \> C_1\cdot\pmb{r}_1\cdot(\pb^j-\pmb{x}^j(t-1)) \> + \> C_2\cdot\pmb{r}_2\cdot(\gb-\pmb{x}^j(t-1)) \Bigg)  \tag{Constriction}
\]
where
\[
	\chi= \frac{2}{\left| 2-\phi-\sqrt{ \phi^2-4\phi } \right|}
\]
with
\[
	\phi=C_1+C_2
\]
Here particle velocities $\pmb{v}^j(0)$ are initialised by uniform random sampling $\sim U(0\> , \>0.1)^n$ \cite{Engelbrecht2012}. Each component of the random vectors $\pmb{r}$ is typically randomly sampled individually, $\pmb{r}_1 \> , \> \pmb{r}_2 \sim U(0\> , \>1)^n$. Vector multiplication is component-wise. The scalar terms $C_1$, $C_2$  represent weightings that a particle puts on its $\pb^j$ ($C_1$) versus $\gb$ ($C_2$) location data. The inertia scalar $\omega$ moderates the significance of the preceding velocity, whilst the constriction scalar $\chi$ is used to avoid particles 'exploding' -- disappearing out of the feasible region. Here an invisible boundary condition is adopted \cite{RobinsonRahmatSamii2004}, with particles allowed to leave the feasible region but no function evaluations undertaken for particle locations outside of $\X$.

A non-robust PSO can be extended to a robust approach through the addition of an inner maximisation search component. This is the approach adopted here. The inner maximisation techniques available as grammar components are discussed in Section~\ref{sec:InnerMax}.

The (Inertia) and (Constriction) formulations represent an rPSO baseline movement capability. For any given heuristic in our GGGP the movement calculation can be extended through the addition of components, described in Section~\ref{sec:CFG}. This includes building blocks based on a series of metaheuristics for robust problems developed using rPSO as a framework in \cite{HughesGoerigkDokka2020a}, and novel features here. These developments are largely based around two robust search techniques, d.d.\ and LEH.

\subsubsection{Descent direction}
\label{sec:dd}

Descent directions is an exploitation-focussed, individual-based robust local search technique \cite{BertsimasNohadaniTeo2007,BertsimasNohadaniTeo2010,BertsimasNohadaniTeo2010nonconvex} for solving (MinMax), although it can easily be extended to approximate a global search through random re-starts each time a local search completes \cite{HughesGoerigkWright2019}. We briefly summarise the method outlined in \cite{BertsimasNohadaniTeo2007,BertsimasNohadaniTeo2010,BertsimasNohadaniTeo2010nonconvex}.

In d.d.\ at each candidate point $\pmb{x}$ in the decision variable space that the search moves to, an inner maximisation search is performed to assess that point's uncertainty neighbourhood $N(\pmb{x})=\{\pmb{x}+\Delta\pmb{x} \mid \Delta\pmb{x}\in\U\}$ and approximate the worst case cost $\tilde{g}(\pmb{x})\approx g(\pmb{x})$. Function evaluations are stored in a history set $H$, and at each candidate the local information is further exploited through the identification of poor `high cost' points (hcps), those with the greatest objective function value, in $H$ and within the $\Gamma$-radius uncertainty region. At a candidate point $\pmb{x}$ the  high cost set $H_\sigma(\pmb{x})$ is defined as: 
\[
	H_\sigma(\pmb{x}) := \{ \pmb{x}' \in H\cap N(\pmb{x}) \mid f(\pmb{x}')  \geq  \tilde{g}(\pmb{x}) -  \sigma \}
\]
$\sigma$ is the threshold value for determining what constitutes an hcp. 

The optimum (descent) direction originating at the current candidate $\pmb{x}(t)$ at step $t$, and pointing away from the hcps, is then calculated using mathematical programming. The angle $\theta$ between the vectors connecting the points in $H_\sigma(\pmb{x}(t))$ to $\pmb{x}(t)$, and $\pmb{d}$, is maximised: 

\begin{align*}
\min_{\pmb{d},\beta} \ & \beta \tag{Soc} \\
\text{s.t. } &\| \pmb{d}\| \le 1 \tag{Con1} \\
& \pmb{d}^T \Bigg( \frac{\pmb{h} - \pmb{x}(t)}{\| \pmb{h} - \pmb{x}(t)  \|} \Bigg) \le \beta & \forall \pmb{h} \in H_\sigma(\pmb{x}(t)) \tag{Con2}\\
& \beta \le -\varepsilon \tag{Con3}
\end{align*}

Here $\varepsilon$ is a small positive scalar, so from (Con3) $\beta$ is negative. The left hand side of constraint (Con2) is the multiplication of $\cos \theta$ and $\| \pmb{d}\|$, for all hcps in $H_\sigma(\pmb{x}(t))$ and a feasible direction $\pmb{d}$. (Con2) therefore relates $\beta$ to the maximum value for $\cos \theta$ across all hcps. As the objective (Soc) is to minimise $\beta$, and $\beta$ is negative, the angle $\theta$ will be greater than $90^o$ and maximised. Also minimising $\beta$ in combination with (Con1) normalises $\pmb{d}$. A standard solver such as CPLEX can be used to solve (Soc). When a feasible direction cannot be found, that is (Soc) cannot be solved, the algorithm stops: a robust local minimum has been reached.

The local search repeats at step $t$ by moving away from the current candidate $\pmb{x}(t)$, in this optimum direction $\pmb{d}$ with a step size $\rho(t)$ large enough that the points in $H_\sigma(\pmb{x}(t))$ are at a minimum on the boundary of the uncertainty neighbourhood of the next candidate at step $t+1$. Then $ \pmb{x}(t+1) = \pmb{x}(t) + \rho(t)\cdot\pmb{d} $, where:
\[
	\rho(t) = \min \left\{ \pmb{d}^T(\pmb{h} - \pmb{x}(t)) + \sqrt{(\pmb{d}^T(\pmb{h} - \pmb{x}(t)))^2 - {\| \pmb{h} - \pmb{x}(t)  \|}^2 + \Gamma^2} \mid \pmb{h} \in H_\sigma(\pmb{x}(t)) \right\}  \tag{Rho} 
\]
Steps are repeated until a local minimum is reached. 

In one of the rPSO variants described in \cite{HughesGoerigkDokka2020a}, given neighbourhood uncertainty information for each particle at each step, a descent direction vector is calculated. This vector $\pmb{d}^j(t-1)$ for particle $j$ at step $t$ is used for the calculation of an additional velocity component:
\[
	C_3 \cdot \pmb{r}_3 \cdot \pmb{d}^j(t-1) \tag{ddVel}
\] 
In \cite{HughesGoerigkDokka2020a} each component of $\pmb{r}_3$ is randomly sampled individually, $\pmb{r}_3 \sim U(0\> , \>1)^n$, vector multiplication is component wise, and the scalar term $C_3$ represents a weighting that a particle puts on its local descent direction vector. From step $t=1$ onwards a variant on the baseline PSO (Inertia) velocity formulation is then used in \cite{HughesGoerigkDokka2020a}:
\[
	\pmb{v}^j(t)=\omega \cdot \pmb{v}^j(t-1) \> + \> C_1 \cdot \pmb{r}_1 \cdot (\pb^j-\pmb{x}^j(t-1)) \> + \> C_2 \cdot \pmb{r}_2 \cdot (\gb-\pmb{x}^j(t-1)) \> + \> C_3 \cdot \pmb{r}_3 \cdot \pmb{d}^j(t-1) \tag{InertiaV2}
\]
Building blocks components based on this d.d.\ approach, and associated parameters, are considered in the grammar here. Details are given in Section~\ref{sec:CFG}.

\subsubsection{Largest empty hypersphere}
\label{sec:leh}

Largest empty hypersphere is an exploration-focussed individual-based robust global search technique \cite{HughesGoerigkWright2019} for solving (MinMax). LEH takes the d.d.\ concept of hcps to a global setting, identifying a high cost set $H_\tau$ of poor points from within the global history set $H$, and moving to the centre of the region completely devoid of all such points. $H_\tau$ contains those points in $H$  with nominal objective function value $f(\pmb{x})$ greater than a threshold $\tau$. In LEH $\tau$ equals the current estimated robust global minimum value.

The centre of the LEH, $\pmb{x}(t)\in\X$ at iteration $t$, is the estimated point furthest from all hcps in $H_\tau$, and is approximated using a genetic algorithm (GA). Movement from centre of LEH to centre of LEH repeats until no point $\pmb{x}(t)\in\X$ which is at least $\Gamma$ away from all hcps can be identified, or a defined budget of available objective function evaluations (model runs) is exhausted. The final estimate for the global robust minimum is accepted. 

A key feature of LEH is the early stopping of neighbourhood searches at any candidate where an improved estimated robust global optimum cannot be achieved. In theory an inner maximisation is performed at each candidate $\pmb{x}(t)$, however in LEH each objective function evaluation in a neighbourhood analysis $f(\pmb{x}(t) + \Delta\pmb{x}(t))$ is compared to $\tau$, with the inner search terminating if that value exceeds $\tau$. This recognises that the current point $\pmb{x}(t)$ won't improve on the estimated robust global optimum, and has the potential to afford considerable savings in local function evaluations and so enable a more efficient exploration of $\X$.

One of the rPSO variants in \cite{HughesGoerigkDokka2020a} is based around core elements of the LEH approach. Firstly the stopping condition is employed at a particle level for each particle in each iteration. For any particle $j$ an inner maximisation search may begin but is terminated early if an uncertainty neighbourhood point exceeds the particle best information $\pb^j$. In fact by first assessing the complete history set $H$ of all previous function evaluations, no inner maximisation may be necessary if it is determined that some historical value in the particle's uncertainty neighbourhood already exceeds the best information $\pb^j$.

Using a second novel LEH-based feature, particles are assessed for `dormancy', defined as a state where no function evaluations have been required by a particle for a specified number of iterations. This may be due to the repeated identification of existing neighbourhood points which exceed the particle's best information $\pb^j$, prior to undertaking an inner maximisation. Or it may be due to the particle repeatedly moving outside the feasible region, linked to the use of an invisible boundary condition \cite{RobinsonRahmatSamii2004}. In either case dormancy suggests that a particle has become 'stuck'. In \cite{HughesGoerigkDokka2020a} dormant particles are relocated to the centre of the largest empty hypersphere devoid of all hcps, using the current robust global minimum as the high cost threshold. 

More details of the LEH-based components and associated parameters available in the grammar here are given in Section~\ref{sec:CFG}.


\section{The automatic generation of heuristics}
\label{sec:autoGenAlgos}

In seeking to develop improved metaheuristics for robust problems an obvious question is what features should be included in the search methodology. This is a step beyond the issues of what existing search technique a decision maker might employ, or what parameter settings might be used for any given problem. These issues impact the effectiveness of any optimisation search.

Given a problem for which an optimisation search is to be undertaken, the field of hyper-heuristics encompasses techniques which employ a search methodology to automatically identify or generate heuristics for application to that problem. The hyper-heuristic itself does not search the problem solution space, but rather seeks a heuristic for application to the problem. A high level classification of hyper-heuristic approaches distinguishes between methods for selecting a heuristic, from a space of heuristics, and methods for generating a heuristic \cite{BurkeHydeKendallOchoaOzcanWoodward2009, BurkeGendreauHydeKendallOchoaOzcanQu2013}. Our interest is in the latter.

The automatic generation of a search heuristic is a specific application of the broader theme of the automatic generation of algorithms, or the automatic generation of computer programs. One technique which can be applied in the general case and to the specific issue of automatically generating a search heuristic is genetic programming (GP) \cite{Koza1992, BurkeHydeKendallOchoaOzcanWoodward2009}. This employs the well known high level concepts of selection, combination and mutation to evolve a population of computer programs, or in our case search heuristics for robust problems. 

When considering the applications of genetic programming to automatically generate optimisation search approaches \cite{BurkeGendreauHydeKendallOchoaOzcanQu2013}, the use of tree-based context-free grammar-guided GP \cite{McKayHoaiWhighamShanONeill2010} to the generation of PSO heuristics described in \cite{MirandaPrudencio2016, MirandaPrudencio2017} is of particular interest here. We adopt that approach and apply it to PSO based heuristics for robust problems. 

Relatively little work has been undertaken on the application of GP to optimisation search techniques for uncertain problems. One example from the field of stochastic optimisation is \cite{MeiZhang2018}, where GP is applied to a vehicle routing problem including uncertainty. In terms of optimisation for robust problems, \cite{GoerigkHughes2019} use a simple GP-based approach to evolve techniques for robust combinatorial optimisation problems. However, to the author's knowledge \cite{GoerigkHughes2019} is the only explicit use of a GP-based approach applied to a robust problem, and there is no application of GP to metaheuristics for black-box robust optimisation problems prior to the work outlined here.


\section{The genetic programming of metaheuristics for robust problems}
\label{sec:autoRob}

\subsection{Genetic programming}
\label{sec:treeGP}

Our aim is to develop improved metaheuristics for robust problems and remove the manual determination of feature-technique-parameter choices, through the automatic generation of algorithms by genetic programming \cite{Koza1992, McKayHoaiWhighamShanONeill2010, Nohejl2011}. GP is an evolutionary process, and here each individual in the GP population is a heuristic. For an initial population of heuristics, a measure of fitness is calculated for each individual. A new generation of heuristics is then determined through typical evolutionary algorithm fitness-based selection, combination and mutation processes. This repeats over multiple generations, at the end of which the fittest heuristic is chosen. 

Each heuristic in the GP process is made up of multiple algorithmic sub-components and their parameter settings, which when combined appropriately form executable search heuristics. So the GP solution space consists of sub-components and their parameters.

We define algorithmic sub-components along with the production rules which determine how they combine to form complete heuristics. Sub-components are designed to integrate effectively under those construction rules. This is our grammar \cite{Koza1992}, which is employed within an evolutionary framework. That framework must be capable of performing combination and mutation operations on heuristics constructed from building blocks. Here a tree-based GP evolutionary process is used to facilitate these processes, as described in Section~\ref{sec:tree}. Details of the individual sub-components are given now in Section~\ref{sec:CFG}.

\subsection{Grammar}
\label{sec:CFG}

\subsubsection{Structure}
\label{sec:CFGstructure}

The heuristics considered here consist of outer minimisation and inner maximisation searches, wrapped around a black-box model. Each model run generates a single objective function output corresponding to a point in the model decision variable space $\X$. We use a PSO frame for all heuristics, comprising a swarm of particles moving over a series of iterations. This constitutes the outer minimisation, with inner maximisations undertaken at the particle level to determine the robust objective function value at a point in $\X$.

Every member of the population in the GP analysis has the same basic algorithmic structure, described by Algorithm~\ref{BaselineRPSO}. We assume a limit on the number of function evaluations (model runs) available. The swarm is initialised randomly, and the defined form of inner maximisation undertaken to determine robust objective values for each particle. Particle movement is then controlled by the velocity equation formulation and the forms of topology and movement, i.e. how particle velocities are calculated, how particles share information, and how these elements are used. The swarm moves and particle level inner maximisations are undertaken again. This repeats until the budget is exhausted. On completion the current best estimate for robust global minimum is accepted. 

\vspace{2mm}
\begin{algorithm}[H] 
\caption{Overview of a robust particle swarm optimisation algorithm} \label{BaselineRPSO}
\vspace{2mm} 
\hspace*{\algorithmicindent} \textbf{Inputs:} Swarm $size$, extent of inner search ($inExt$), $budget$ of function evaluations \\
\hspace*{\algorithmicindent} \textbf{Parameters:} Form of $inner$, form of $topology$, form of $velocity$, form of $movement$ \\
\hspace*{\algorithmicindent} \textbf{Parameters:} $innerParams$, $topolParams$, $velParams$, $moveParams$
\vspace{2mm} 
 
\begin{algorithmic}[1]
			
	\State $t \gets 0$ 

	\While{($budget>0$)}

		\ForAll{($j$ in $1,\ldots,size$)}

				\If{($t = 0$)}  
				
					\State Randomly initialise particle $\pmb{x}^j(0) \in \X$
			
				\Else 
			
					\State Update particle velocity according to ($velocity, velParams$)
					\State Update particle position according to ($movement, moveParams, topology, topolParams$) 
			
				\EndIf 
							
				\If{($\pmb{x}^j(t) \in \X$)} 
		
					\State Perform inner maximisation: 	
					
					\ForAll{($k$ in $1,\ldots,inExt$)}	
							\State Select uncertainty neighbourhood point: ($inner, innerParams$)	
							\State Evaluate function (run model, generate objective)		
							\State $budget \gets budget-1$
					\IIf{($budget=0$)} \textbf{break}: goto end \EndIIf 				
					\EndFor 

				\EndIf 
		
		\EndFor 
		
		\State $t \gets t+1$ 
		
	\EndWhile  

	\State \Return Current estimate of robust global best
	
\end{algorithmic}
\end{algorithm}

\medskip

Generating metaheuristics for robust problems in a GP process requires the definition of a grammar: algorithmic sub-components and the rules for combining them. Here the high level outline of each heuristic, Algorithm~\ref{BaselineRPSO}, also forms the high level design criteria in the grammar: the outer minimisation layer as a swarm of particles, some movement formulation, a topology dictating particle information-sharing, and an inner maximisation layer.

\begin{figure}[htbp]
	\centering
	\vspace{-3mm}	
	\includegraphics[width=1.0\textwidth]{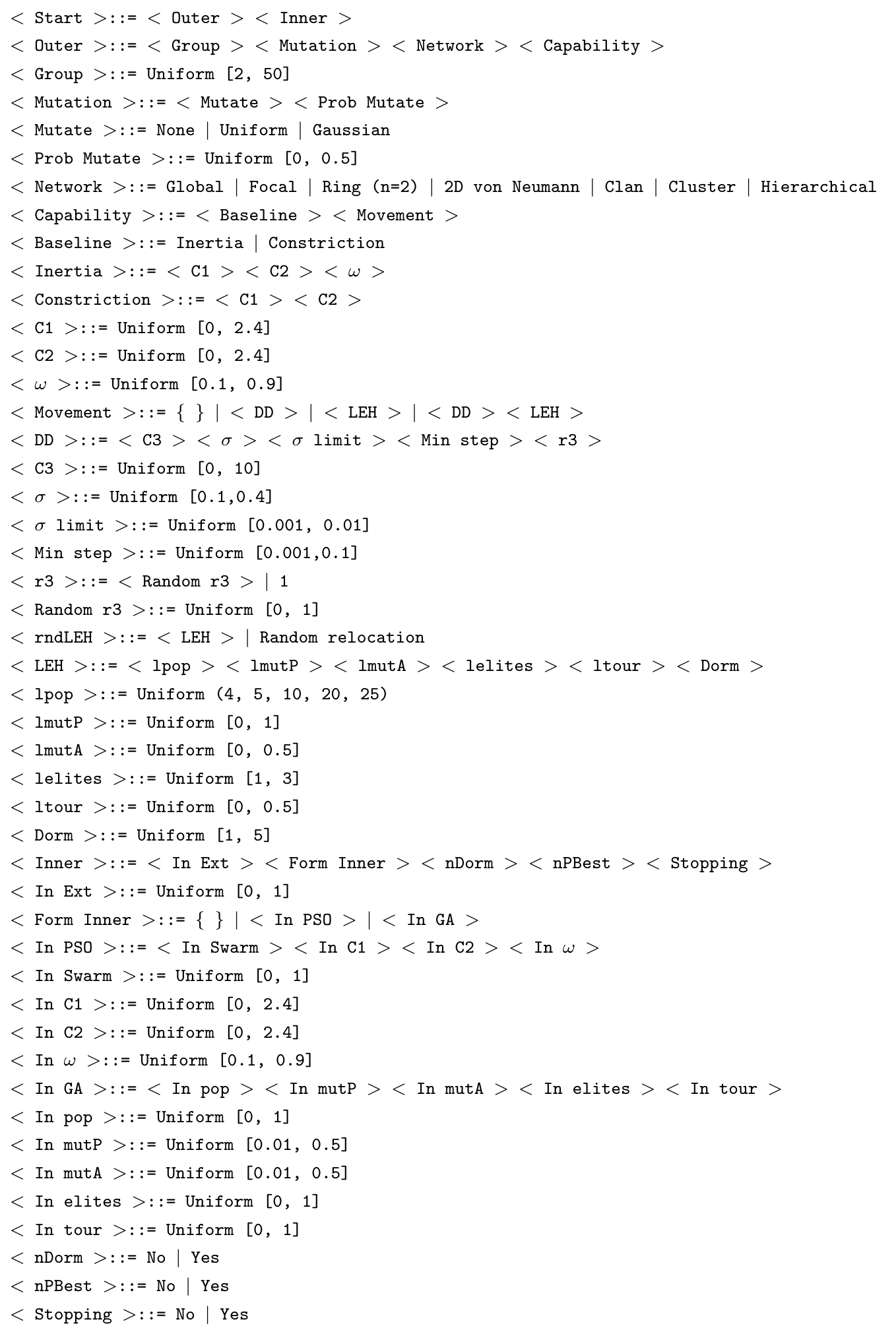} 
	\vspace{-3mm}	
	\caption{Context-free grammar employed here for the construction of metaheuristics for robust problems. }
	\label{fig:grammar}
\end{figure}

Our sub-components and the production rules for generating complete heuristics are defined in the grammar in Figure~\ref{fig:grammar}. The specific approach adopted here is known as context-free grammar genetic programming (CFG-GP). This uses a tree-based representation of algorithms and standard tree-based operators in the evolutionary process, see \cite{McKayHoaiWhighamShanONeill2010, MirandaPrudencio2016}.

The grammar includes non-terminal nodes, indicated by \texttt{< >}, terminal nodes, and the production rules (\texttt{::=}). The generation of a heuristic begins at the \texttt{<Start>} node, resulting in the generation of more nodes by following the rules in Figure~\ref{fig:grammar}. Each non-terminal node leads to the generation of further nodes according to the production rules, with each non-terminal node expanded upon until a terminal node is reached. The result is the generation of a series of nodes corresponding to elements of the heuristic. Non-terminal nodes do not result in the generation of further nodes, but instead in the determination of parameter settings: parameter values or choices of individual sub-components. On reaching a non-terminal node that portion of the heuristic is complete. The final heuristic is achieved when there are no more non-terminal nodes to expand upon in the sequence.

Sub-component details are given in Sections~\ref{sec:rPSO} to~\ref{sec:stopping}. The high level \texttt{<Outer>} and \texttt{<Inner>} elements are identifiable in Figure~\ref{fig:grammar}. \texttt{<Outer>} consists of \texttt{<Group>} (swarm size), \texttt{<Mutation>}, \texttt{<Network>} and \texttt{<Capability>} elements. \texttt{<Mutation>} refers to random variations applied to a particle's next location, \texttt{<Network>} specifies the rules for particle information sharing, and \texttt{<Capability>} covers a number of sub-components which combine to form the rules for particle movement. \texttt{<Capability>} breaks down into \texttt{<Baseline>} and \texttt{<Movement>}, where \texttt{<Baseline>} refers to core PSO velocity equations, and \texttt{<Movement>} refers to extended capabilities built around d.d. \cite{BertsimasNohadaniTeo2010} and LEH \cite{HughesGoerigkWright2019} techniques and their variants \cite{HughesGoerigkDokka2020a}. \texttt{<Inner>} is by random sampling, or a PSO or GA search, along with additional \texttt{<nDorm>}, \texttt{<nPBest>} and \texttt{<Stopping>} capabilities, based on features in \cite{HughesGoerigkDokka2020a} and explained here in Sections~\ref{sec:dormancy} to~\ref{sec:stopping}. There are also a number of sub-components and parameters associated with many of these elements, which in total constitutes our grammar.

\begin{figure}[htb]
	\centering
	\includegraphics[width=0.9\textwidth]{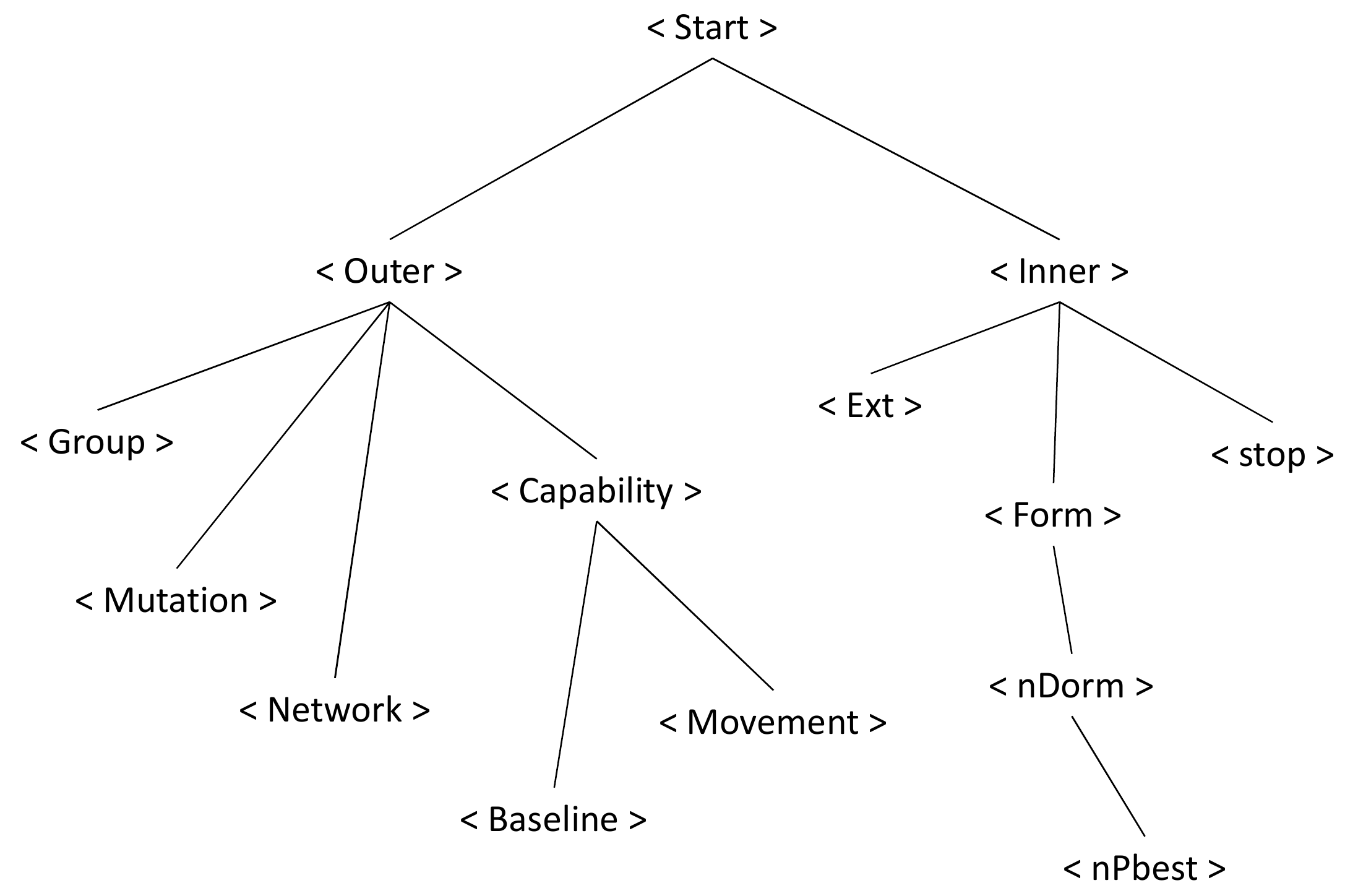} 
	\caption{Solution representation: high level tree-based representation of the heuristic generated by following the CFG-GP grammar production rules in Figure~\ref{fig:grammar}.}
	\label{fig:exampleTree}
\end{figure}

One way to visualise this process is in the form of a tree \cite{Koza1992, Whigham1995}, Figure~\ref{fig:exampleTree}, showing the high level structure of a heuristic generated by the production rules in Figure~\ref{fig:grammar}. \texttt{<Start>} produces the non-terminal nodes \texttt{<Outer>} and \texttt{<Inner>}. \texttt{<Outer>} is expanded upon, and when it is complete the \texttt{<Inner>} node is returned to. From the \texttt{<Outer>} node \texttt{<Group>}, \texttt{<Mutation>}, \texttt{<Network>} and \texttt{<Capability>} are generated one at a time, fully expanding on \texttt{<Group>} before moving to \texttt{<Mutation>} and so on. When \texttt{<Capability>} is complete \texttt{<Outer>} is complete.

\texttt{<Group>} generates a terminal node, a randomly sampled value between 2 and 50 for the number of particles (swarm size) in the heuristic. Having reached a terminal node, the next non-terminal node in the chain generated so far but not yet expanded upon, is \texttt{<Mutation>}. This refers to the mutating of individual particle positions. \texttt{<Mutation>} leads to the non-terminal combination of \texttt{<Mutate>} and \texttt{<Prob Mutate>}, both leading to terminal choices, respectively either \texttt{None} (no mutation), or \texttt{Uniform} or \texttt{Gaussian} mutation, and if mutation the probability of mutation from the range 0 to 0.5. The symbol \texttt{$\mid$} in the production rules designates a choice of one of the alternatives. On randomly choosing a mutation alternative and probability, if required, the next non-terminal node in the chain, \texttt{<Network>}, is returned to. And so on.

\subsubsection{Building blocks: Particle swarm framework}
\label{sec:rPSO}

A basic PSO extended by an inner maximisation forms the basis for all heuristics here. Core PSO elements are a \texttt{<Group>} (swarm) of particles, a \texttt{<Baseline>} velocity equation of either \texttt{Inertia} or \texttt{Constriction} forms described in Section~\ref{sec:pso}, and some system of particle information sharing. The latter, \texttt{<Network>}, is discussed in Section~\ref{sec:networks}. In all heuristics the \texttt{<Group>} size, \texttt{<C1>} and \texttt{<C2>} parameter values are sampled from the ranges defined in the grammar in Figure~\ref{fig:grammar}. The need for the \texttt{<$\omega$>} term depends on the choice of \texttt{<Baseline>}.

\subsubsection{Building blocks: Mutation}
\label{sec:mutation}

The use of the non-terminal node \texttt{<Mutation>} is considered at a particle level, after the candidate position in the next iteration has been determined, see \cite{MirandaPrudencio2017}. If used, mutation is considered separately for each particle and at the dimensional level, as a final stage in the movement calculation. For each particle it is determined whether or not to mutate by sampling against the probability \texttt{<Prob Mutate>}. If mutation is confirmed, any given dimension is mutated with probability randomly sampled from between 0 and $1/n$, so on average only one dimension is changed. The magnitude of change is sampled from either the Uniform or Gaussian distributions as appropriate, and related to the dimensional bounds.

\subsubsection{Building blocks: Networks}
\label{sec:networks}

The sharing of information throughout the swarm to inform movement at the individual particle level, is a core PSO element. Here that form of sharing between particles is determined by the \texttt{<Network>} component. Of the large number of networks available \cite{KennedyEberhart1995, KennedyMendes2002, MendesKennedyNeves2003, JansonMiddendorf2005, deCarvalhoBastosFilho2009, WangYangOrchard2016,  MirandaPrudencio2017}, we consider seven alternatives. Each particle is assigned to a network neighbourhood. At each iteration information on the best neighbourhood point visited by any particle in the network across all iterations, $\gb$, is shared in a manner defined by the baseline velocity equations (Inertia) or (Constriction). 

\begin{itemize}[leftmargin=*]	

	\item \textit{Global}: This is the most basic formulation, with all particles accessing the same neighbourhood information -- the current robust global minimum location \cite{KennedyEberhart1995, KennedyEberhartShi2001}.  
	\item \textit{Focal}: A singe particle is randomly selected as the focal. All particles access the same neighbourhood information, the focal particle's best  information, \cite{KennedyEberhartShi2001, KennedyMendes2002}. 
	\item \textit{Ring (size=2)}: In a network sense all particles may be randomly arranged into a ring formation. With this topology a particle has access to the best information from the adjoining particles in the ring. Here we set the size equal to two, so a particle has access to its two neighbour's (one either side in the ring formation) best information \cite{KennedyEberhartShi2001, KennedyMendes2002}.
	\item \textit{2D von Neumann}: In a network sense particles may be randomly arranged into a 2D grid, or more correctly the surface of a torus where the grid wraps around so that the top and bottom join, as do the left and right hand sides. Each particle has four neighbours, the nearest north, south, east and west particles in this grid formation, accessing the best information in this neighbourhood \cite{KennedyMendes2002, MirandaPrudencio2017}.
	\item \textit{Clan}: Each particle is randomly placed in a network sub-group, or clan. Each clan is linked to each other clan via a clan leader. The leader in each clan is the particle with the best performance, so leaders may change over iterations of the swarm. Each leader shares their information with all other clans \cite{deCarvalhoBastosFilho2009, MirandaPrudencio2017}.
	\item \textit{Cluster}: Each particle is randomly placed in a network sub-group, or cluster. Within each cluster a number of `informant' particles are randomly assigned. The number of informants is one less than the number of clusters, and within each cluster one informant is linked to one other cluster. Informants remain fixed. Within a cluster the best information is shared between all particles. Informant particles share their information with the single cluster they link to \cite{MendesKennedyNeves2003, MirandaPrudencio2017}. 
	\item \textit{Hierarchical}: All particles are randomly arranged in a tree formation, in a network sense. The depth and width of the tree is dependent on the number of particles (swarm size). Each particle communicates with the particle above it in the tree. At each iteration of the swarm the positions in the tree can shift: if a particle below another one in the tree performs better, the two swap positions. This applies to all particles in each iteration	\cite{JansonMiddendorf2005, MirandaPrudencio2017}.

\end{itemize}

\subsubsection{Building blocks: Additional movement capability}
\label{sec:movement}

The baseline PSO capability can be augmented by additional movement components based on the descent directions \cite{BertsimasNohadaniTeo2010} and largest empty hypersphere \cite{HughesGoerigkWright2019} approaches, as proposed in \cite{HughesGoerigkDokka2020a}. The additional formulations available in our grammar are none \texttt{\{ \}}, a d.d.\ based approach \texttt{<DD>}, an LEH based approach \texttt{<LEH>}, or a combined d.d.\ and LEH based approach \texttt{<DD>} \texttt{<LEH>}. In the case of \texttt{\{ \}} just the rPSO \texttt{<Baseline>} formulation is used. Otherwise the rPSO d.d.\ or LEH approaches, or both, augment the baseline rPSO formulation at the particle level as described in Sections~\ref{sec:dd} and ~\ref{sec:leh}, and \cite{HughesGoerigkDokka2020a}. 

Both \texttt{<DD>} and \texttt{<LEH>} require the determination of further non-terminal nodes. \texttt{<DD>} employs the nodes: \texttt{<C3>}, \texttt{<$\sigma$>}, \texttt{<$\sigma$ limit>}, \texttt{<Min step>} and \texttt{<r3>}. These are d.d.\ parameters whose descriptions can be found in \cite{BertsimasNohadaniTeo2010, HughesGoerigkWright2019, HughesGoerigkDokka2020a}. They all terminate once parameter values have been determined, with the exception of \texttt{<r3>} which relates to the additional $C_3$ component in the d.d.\ equations (ddVel) and (InertiaV2) in Section~\ref{sec:dd}. In the original formulation each element of $\pmb{r}_3$ is randomly sampled individually, $\pmb{r}_3 \> \sim U(0\> , \>1)^n$. This alternative is available in the component \texttt{<Random r3>} in the grammar, along with another where each element of $\pmb{r}_3$ is set to unity. The latter is a recognition that a locally calculated d.d.\ vector might be more effective without added random variation.

\texttt{<LEH>} relates to the relocation of a particle deemed 'dormant', and requires the determination of either the non-terminal node \texttt{<LEH relocation>} or the terminal selection of \texttt{Random relocation}. In the original formulation a particle is moved to the centre of the LEH devoid of all identified high cost points \cite{HughesGoerigkDokka2020a}, as described in Section~\ref{sec:leh}. Here this is designated by \texttt{<LEH relocation>}, and if selected the parameters \texttt{<lpop>}, \texttt{<lmutP>},  \texttt{<lmutA>}, \texttt{<lelites>}, \texttt{<ltour>} and \texttt{<Dorm>} must be determined. In the grammar all of these parameters terminate once their values have been generated; their descriptions can be found in \cite{HughesGoerigkWright2019, HughesGoerigkDokka2020a}. However an alternative is available here, \texttt{Random relocation}, which as the name suggests simply relocates a particle randomly in $\X$. This does not use any additional parameters.

\subsubsection{Building blocks: Inner maximisation}
\label{sec:InnerMax}

In theory an inner maximisation search is required to accurately estimate the worst objective function value in a candidate point's uncertainty neighbourhood. In practice, issues such as the run time for each function evaluation (model run) may be prohibitive. Here we assume a limit on the number of function evaluations that are possible. This will likely restrict the accuracy of any search, as it will cause some trade-off between the extent of an inner search (robustness) and the level of global exploration. Such a trade-off is not simple \cite{MirjaliliLewisMostaghim2015, DiazHandlXu2017}. In this context the choice of approach for, and the extent of, the inner maximisation is not obvious. Here the non-terminal \texttt{<Inner>} node generates several further non-terminal nodes: \texttt{<In Ext>}, \texttt{<Form Inner>}, \texttt{<nDorm>}, \texttt{<nPBest>} and \texttt{<Stopping>}.

\texttt{<nDorm>}, \texttt{<nPBest>} and \texttt{<Stopping>} are discussed in Sections~\ref{sec:dormancy} to~\ref{sec:stopping}. \texttt{<In Ext>} and \texttt{<Form Inner>} relate to the extent and form of inner search. \texttt{<In Ext>} is based on a randomly sampled value in the range 0 to 1. This value is related to the outer particle group size and budget of function evaluations, to determine a corresponding integer size of inner search. There are three alternatives for \texttt{<Form Inner>}: random sampling \texttt{\{ \}}, or inner PSO \texttt{<In PSO>} or genetic algorithm \texttt{<In GA>} searches. All apply to a candidate point's $\Gamma$-radius uncertainty neighbourhood, Section~\ref{sec:minMax}. If random sampling is used no additional parameters are required. Multiple parameter nodes are required for either \texttt{<In PSO>} or \texttt{<In GA>}. 

\texttt{<In PSO>} requires the determination of \texttt{<In Swarm>}, \texttt{<In C1>}, \texttt{<In C2>} and \texttt{<In $\omega$>}, the parameters for an (Inertia) form of PSO: an inner swarm size and settings for $C_1$, $C_2$ and $\omega$. For an inner PSO the (Inertia) formulation is fixed. \texttt{<In GA>} requires the determination of \texttt{<In pop>}, \texttt{<In mutP>}, \texttt{<In mutA>}, \texttt{<In elites>} and \texttt{<In tour>}, parameters for a standard form of GA \cite{Ghazali2009}: an inner population size, and settings for the probability of and amount of mutation, the number of elites, and a tournament size. If employed, \texttt{<In Swarm>} or \texttt{<In pop>} are initially determined in the range 0 to 1, and then related to \texttt{<In Ext>} to give a corresponding integer value. For an inner GA, \texttt{<In elites>} and \texttt{<In tour>} are initially determined in the range 0 to 1, then related to \texttt{<In pop>} to give integer values.

\subsubsection{Building blocks: Dormancy -- use of neighbourhood information}
\label{sec:dormancy}

The consideration of particle dormancy leading to its relocation \cite{HughesGoerigkDokka2020a}, is described in Section~\ref{sec:leh}. Dormancy refers to a particle becoming 'stuck'. Of concern here is when this might be due to the particle being in an already identified poor region of the solution space, and thereby repeatedly not requiring any function evaluations. In our grammar the determination of dormancy for each particle in each generation may (Yes) or may not (No) make use of the history set $H$ of all points evaluated, and specifically those points within a particle's uncertainty neighbourhood. The choice is represented in node \texttt{<nDorm>}.

\subsubsection{Building blocks: Supplement $\pb^j$ -- neighbourhood information}
\label{sec:supplementPbest}

A particle's robust value is based on an inner search, and leads to the determination of the particle's personal best location $\pb^j$ as employed in the (Inertia) or (Constriction) velocity formulation. Given a completed inner search, if relevant the identified robust value can be updated by the worst point already identified in the history set $H$ within the particle's current uncertainty neighbourhood. The choice of whether (Yes) or not (No) historic information is used in this way is included as a component here, represented by node \texttt{<nPBest>}.

\subsubsection{Building blocks: Stopping condition}
\label{sec:stopping}

The use of a stopping condition in an inner search for a given particle, if a point is identified with objective function value exceeding that particle's personal best information, has the potential to generate significant efficiencies in terms of function evaluations, see \cite{HughesGoerigkDokka2020a} and Section~\ref{sec:leh}. In our grammar the choice of whether (Yes) or not (No) to employ a stopping condition is included as a component, represented by node \texttt{<Stopping>}.

\subsection{Tree-based representation and evolutionary operators}
\label{sec:tree}

The evolutionary GP process begins with the random generation of a population of heuristics, constructed following the grammar production rules in Figure~\ref{fig:grammar}. Populating subsequent generations requires the selection, combination and mutation of heuristics. 

Fitness-based selection can be undertaken as in any standard evolutionary process. Here each heuristic in the population, in each generation, is applied to a single test problem or group of problems, as appropriate. A heuristic is run on any single problem multiple times to generate a sample. For each heuristic applied to each problem, the mean of the samples is used as a fitness measure. If only a single problem is under consideration the fitness across the population of heuristics can be determined directly from a comparison of the means. If multiple problems are considered, means must be calculated for each heuristic across multiple problems. The description of how these means are combined into a single fitness measure for each heuristic is given in the experimental analysis Section~\ref{sec:experiments}.

The calculated fitnesses are used in tournament selections to identify two parent heuristics per each individual in the following GP generation, see e.g. \cite{Eiben2012}. A number of elite, unchanged, heuristics are also retained from generation to generation.

\begin{figure}[htbp]
	\centering

	 \vspace*{3mm} 
	
	\begin{subfigure}{.5\textwidth}
		\centering
		\includegraphics[width=2.8in, height=2.4in]{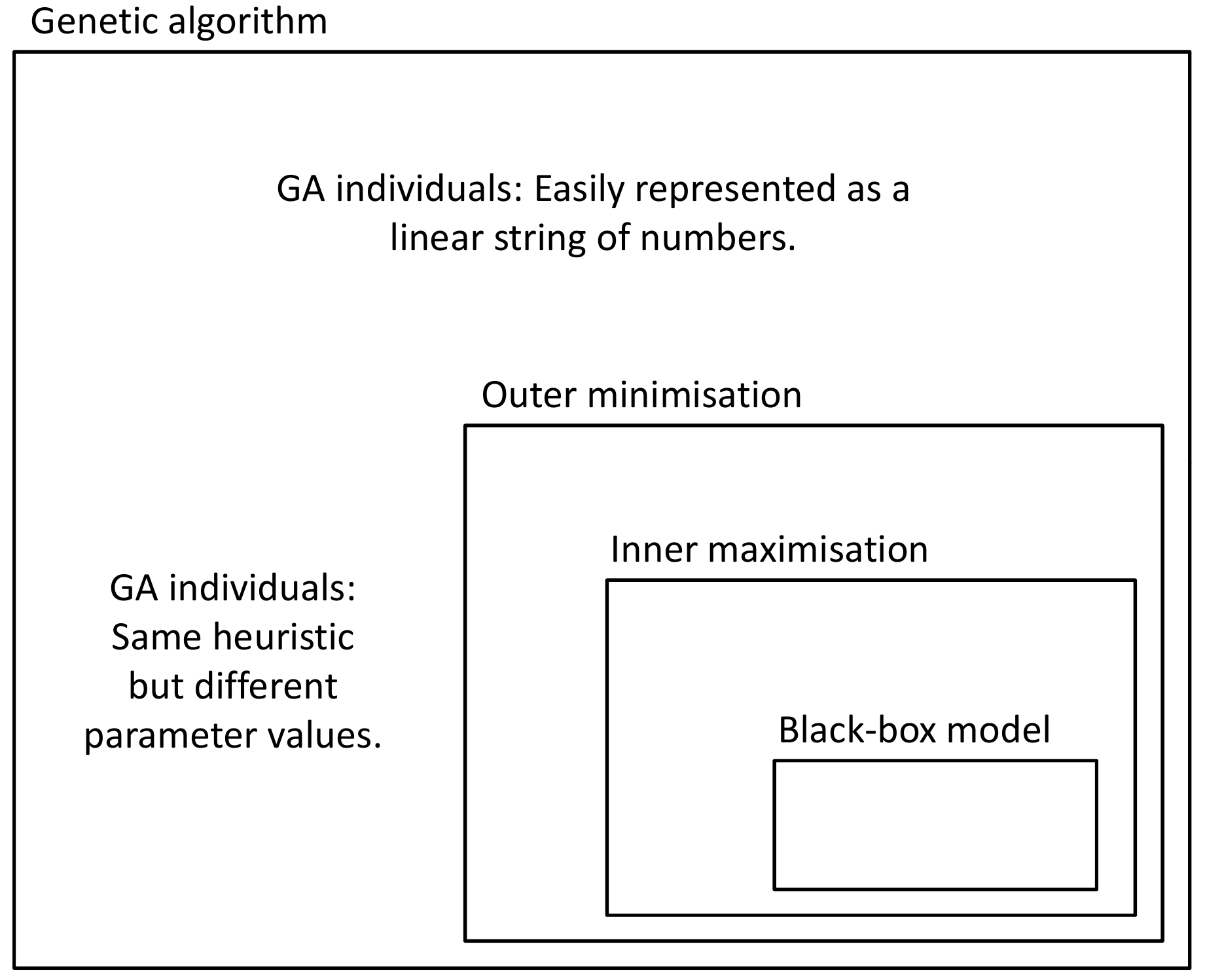} 	
		\caption{GA parameter tuning}
		\label{ParamTune}
	\end{subfigure}%
	\begin{subfigure}{.5\textwidth}
		\centering
		\includegraphics[width=2.8in, height=2.4in]{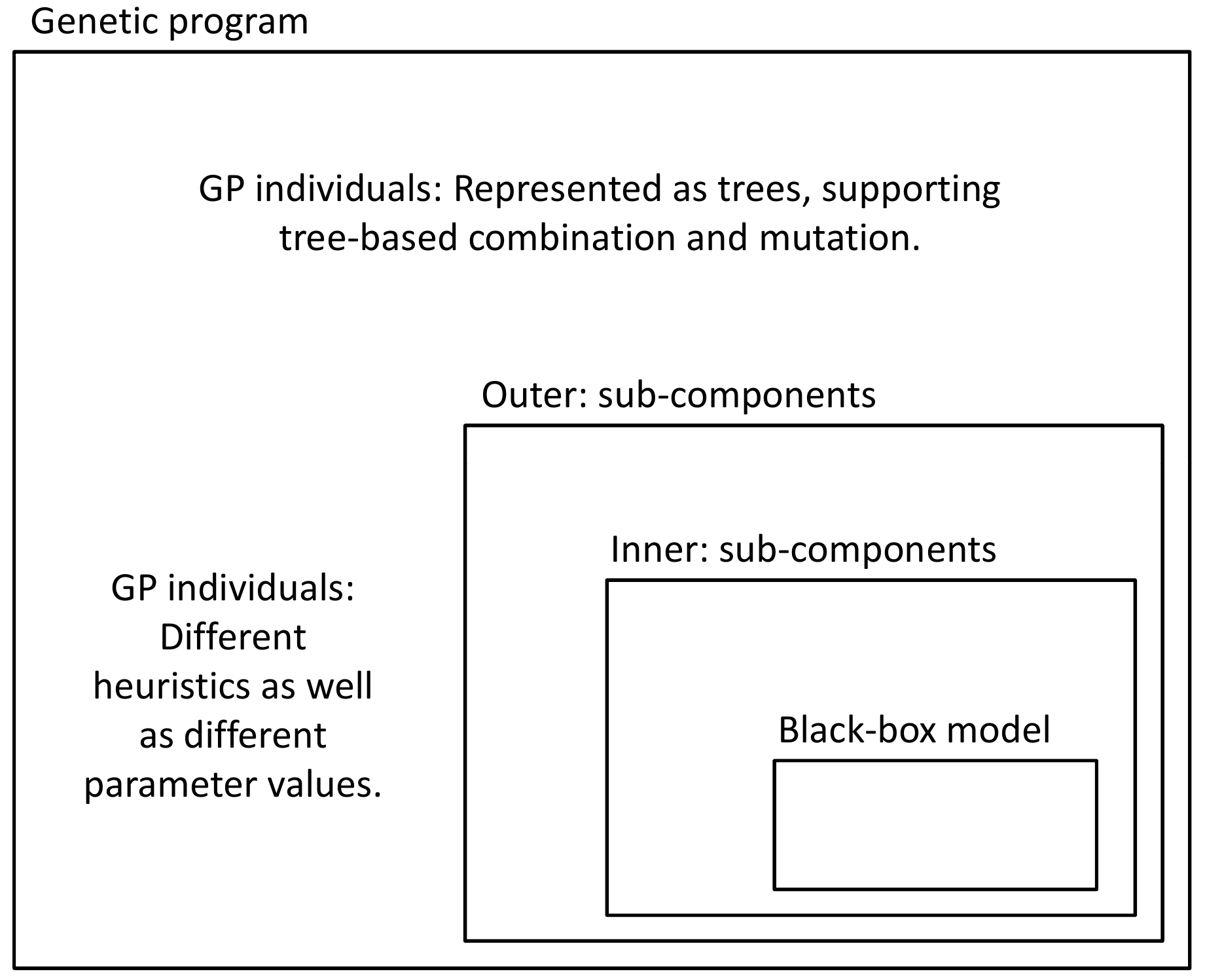}
    \caption{Tree-based GP}
		\label{TreeGP}
	\end{subfigure}
	
  \vspace*{5mm} 
		
	\caption{GA and GP applied to a metaheuristic for robust problems, consisting of an outer minimisation search and inner maximisation search operating on a black-box model.}
	\label{fig:evolAuto}
\end{figure}

For combination and mutation operations, less standard operators may be required. Consider, for example, the differences between the use of a GA to tune the parameters for a specific heuristic compared to the GP evolution of different heuristics, illustrated in Figure~\ref{fig:evolAuto}. All computational evolutionary processes require that an individual object (e.g. a heuristic) has a representative form for the evolution, and in particular combination and mutation operations, to be performed on. In the case of a GA tuning, the solution space consists of the parameters for a single heuristic, which can be represented as a simple linear string of values. So standard GA combination and mutation processes can be employed, see e.g. \cite{Ghazali2009}. Whereas in the GP, each heuristic may comprise different sets of sub-components, complicating a linear representation. For example two such strings would likely be of different lengths, with `corresponding' sections representing different sub-components and so different parameters. Combining and mutating these strings would introduce difficulties.

Fortunately GP offers an alternative heuristic representation, a tree. This is a common representation for computer programs and algorithms \cite{PoliDiChioLangdon2005, McKayHoaiWhighamShanONeill2010, Nohejl2011, MirandaPrudencio2017}, and lends itself to standard tree-based GP operators. The CFG-GP approach we employ uses this representative form for a heuristic generated by our grammar, Figures~\ref{fig:grammar} and~\ref{fig:exampleTree}, with standard random tree-based combination and mutation operators \cite{McKayHoaiWhighamShanONeill2010, MirandaPrudencio2016}.

\begin{figure}[htb]
	\centering
	\includegraphics[width=0.9\textwidth]{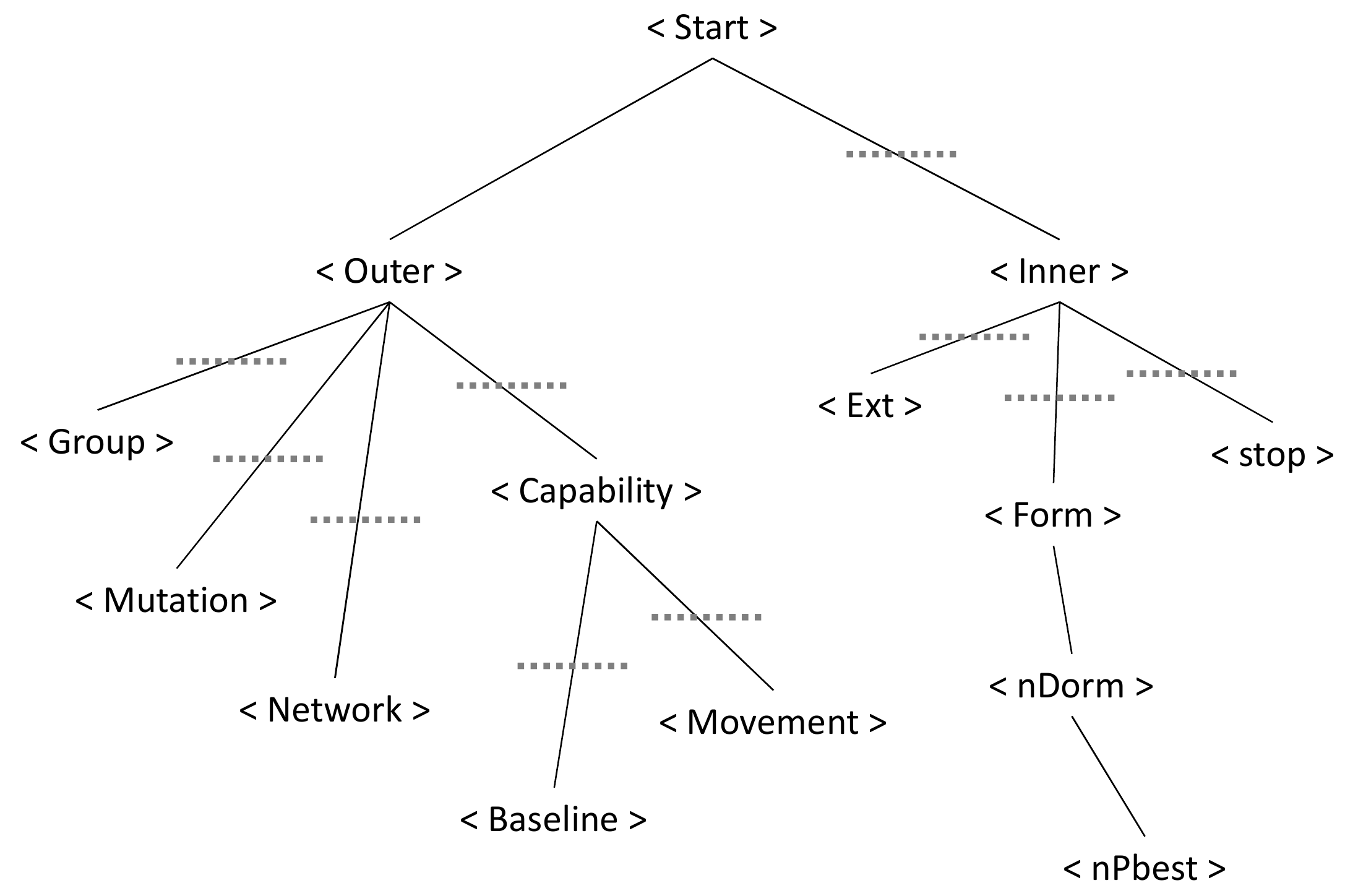} 
	\caption{High level tree-based representation of the heuristic generated by the grammar production rules in Figure~\ref{fig:grammar}, with cut points for combination and mutation operations.}
	\label{fig:treeCuts}
\end{figure}

Consider the high-level heuristic tree representation in Figure~\ref{fig:exampleTree} with the addition of `cut' points, Figure~\ref{fig:treeCuts}. Any two trees generated by our grammar have this overall structure, so the cut points will apply to all of our heuristics. Any two parent heuristics fitness-selected in the GP process, along with one randomly selected cut point, can be combined by merging the branches below the cut in parent tree 1 with the branches above the cut in parent tree 2. The resulting tree is an executable heuristic. This is the combination operation used here.

A newly combined heuristic, represented by a single tree, can be mutated by randomly selecting another cut point. Below the cut point completely new branches can be randomly generated by following the grammar in Figure~\ref{fig:grammar}, whilst retaining the existing branches from above the cut. This is the mutation operation used here, in conjunction with sampling against a probability of mutation to determine whether to mutate. Thus the requirement for selection, combination and mutation operators applicable directly to the heuristics generated by combining sub-components following our grammar, has been fulfilled.

Note that the generation of the heuristic from the corresponding tree is achieved simply by reading off the sub-components and associated parameter values from the terminal nodes (leaves) at the ends of each branch of the tree.

A final point should be made about the benefits of the CFG-GP approach. In a GP process it is not an absolute requirement to always generate executable algorithms, e.g. a fitness value of zero could be assigned to non-executable algorithms. However it can be appreciated that there is a considerable likelihood of  generating non-executable algorithms when randomly combing sub-algorithms. Not only is this very inefficient but it could result in any executable algorithm, whether effective or not, being deemed relatively fit and therefore propagating across many generations. A CFG-GP approach avoids such pitfalls, \cite{McKayHoaiWhighamShanONeill2010, MirandaPrudencio2016}.


\section{Computational experiments}
\label{sec:experiments}

\subsection{Experimental set up}
\label{sec:setUp}

The experimental analysis employs 10 established multi-dimensional robust test problems. The problems are listed in Table~\ref{tab:funcs} along with the feasible regions and $\Gamma$-radius uncertainty values used. Problem formulations and 2D representations are provided in Appendix~\ref{sec:testFunctionFormulae}. In our experiments 30D and 100D versions of these problems are considered. 


\begin{table}[h]
\begin{center}
\begin{tabular}{c|cc}
Name & $\X$ & $\Gamma$ \\
\hline
Rastrigin & $[14.88, 25.12]^n$ & $0.5$\\
Multipeak F1 & $[-5, -4]^n$ & $0.0625$\\
Multipeak F2 & $[10, 20]^n$ & $0.5$\\
Branke's Multipeak & $[-7, -3]^n$ & $0.5$\\
Pickelhaube & $[-40, -20]^n$ & $1$\\
Heaviside Sphere & $[-30, -10]^n$ & $1$\\
Sawtooth & $[-6, -4]^n$ & $0.2$\\
Ackley & $[17.232, 82.768]^n$ & $3$\\
Sphere & $[15, 25]^n$ & $1$\\
Rosenbrock & $[7.952, 12.048]^n$ & $0.25$
\end{tabular}
\caption{Test functions.}\label{tab:funcs}
\end{center}
\end{table}


A single GP run applies each heuristic to a test function or functions, in order to determine fitness and inform the evolutionary process. Here 22 GP runs are considered, once for each test problems individually (individual cases) and once for a combined run where each heuristic is applied to all 10 problems (general case). That is there are 10 individual case GP runs. Within each individual case GP run all heuristics are applied to the same single test problem. There is also one general case GP run, where all heuristics are applied to all 10 test problems. This is repeated for 30D and 100D. A budget of 2,000 function evaluations is assumed in each heuristic run. 

In the GP runs, when a  heuristic is applied to a problem this is repeated 20 times to generate a mean. For the individual case runs this is taken as the fitness of the heuristic. In a general case run the 10 separate means for each heuristic are used to determine 10 separate fitness rankings. For each heuristic the 10 rankings are averaged to give a combined initial ranking. This ranking is then refined using an elimination process. The worst performing heuristic is ranked lowest and removed. For the remaining heuristics the 10 rankings and combined ranking are recalculated, the new lowest performing heuristic is ranked second lowest overall and removed. This repeats until all heuristics have been ranked.

On completion of a GP run the best heuristic in the final population is accepted. To properly assess its performance 200 sample runs of the heuristic are undertaken, applied to the problem or problems on which it has been evolved. Each run generates an estimate of the location of the robust global optimum for the problem(s) at hand. The corresponding robust value at each global optimum location is re-estimated in a post-processing stage, as the worst value identified by randomly sampling a million points in the $\Gamma$-uncertainty neighbourhood of the optimum.

Algorithms are written in Java. For all d.d.\ calculations the solution of (Soc) includes a call to the IBM ILOG CPLEX Optimization Studio V12.6.3 software. 

We now report two analyses. The first, Section~\ref{sec:results}, considers the quality of the best solutions (heuristics) found in the GP runs. The second, Section~\ref{sec:analysis}, assesses the structure (component breakdowns) of the heuristics generated in the GP runs, against heuristic performance.

\subsection{Results for the best performing heuristics}
\label{sec:results}

\begin{figure}[H]
	\centering
	
	\vspace{5mm} 

	\begin{subfigure}[t]{.18\textwidth}
		\includegraphics[width=\textwidth]{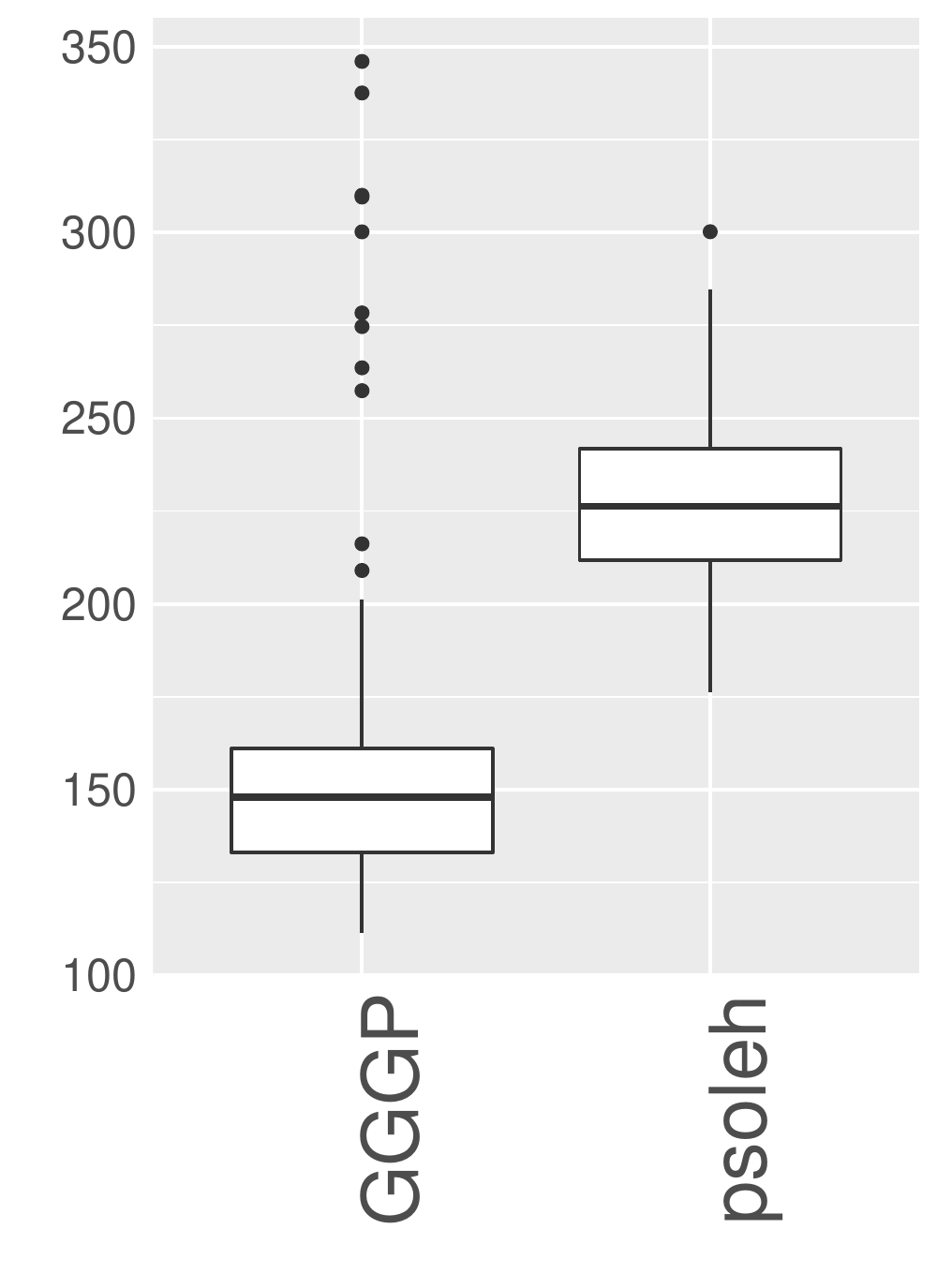}
		\vspace{-4mm} 
  	\caption{\scriptsize{Rastrigin}} \label{fig:box3051i}
	\end{subfigure}%
	\hspace{-1mm} 
	\begin{subfigure}[t]{.18\textwidth}
		\includegraphics[width=\textwidth]{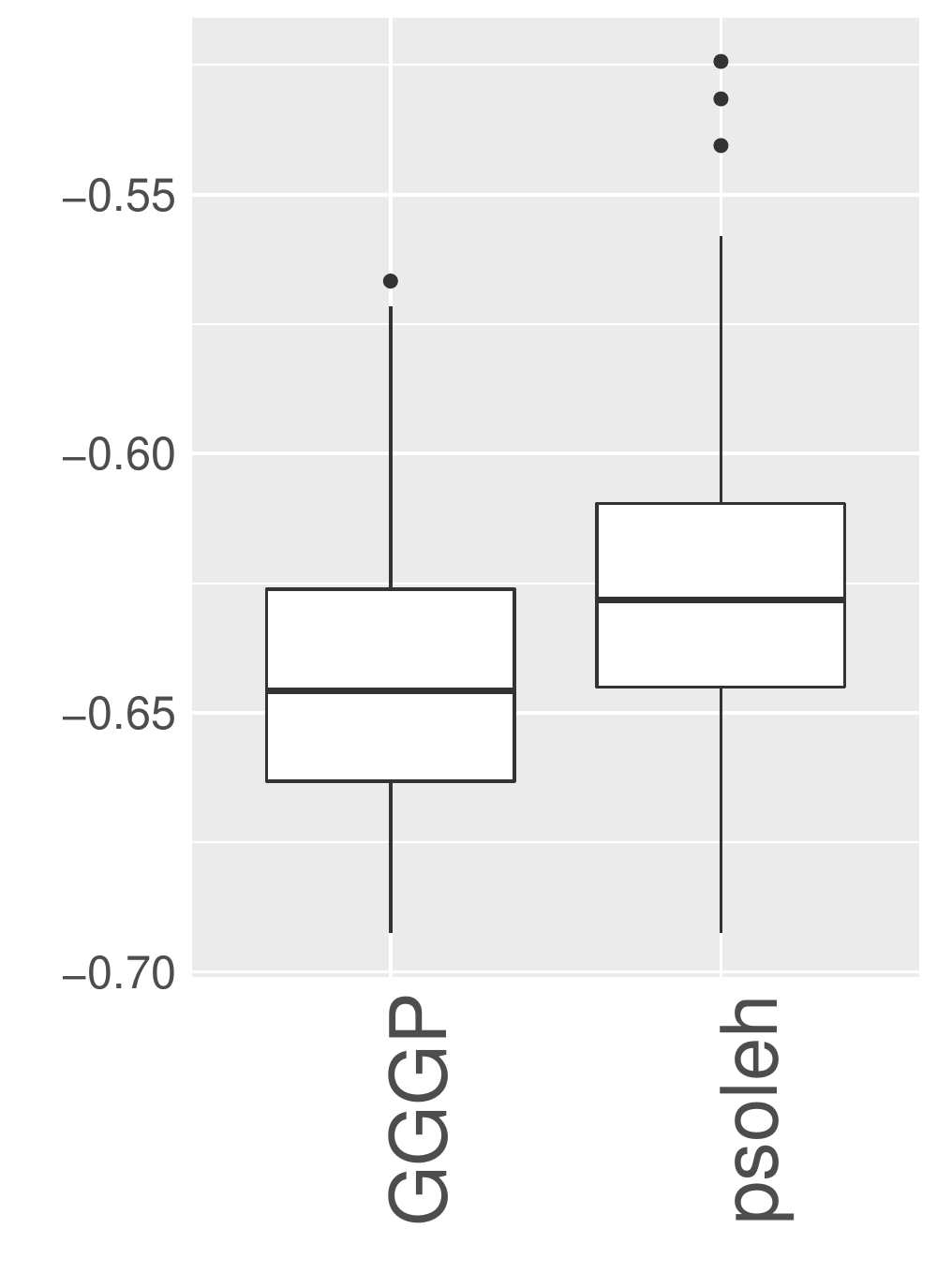}
		\vspace{-4mm} 
  	\caption{\scriptsize{Multipeak F1}} \label{fig:box3057i}
	\end{subfigure}%
	\hspace{-1mm} 
	\begin{subfigure}[t]{.18\textwidth}
		\includegraphics[width=\textwidth]{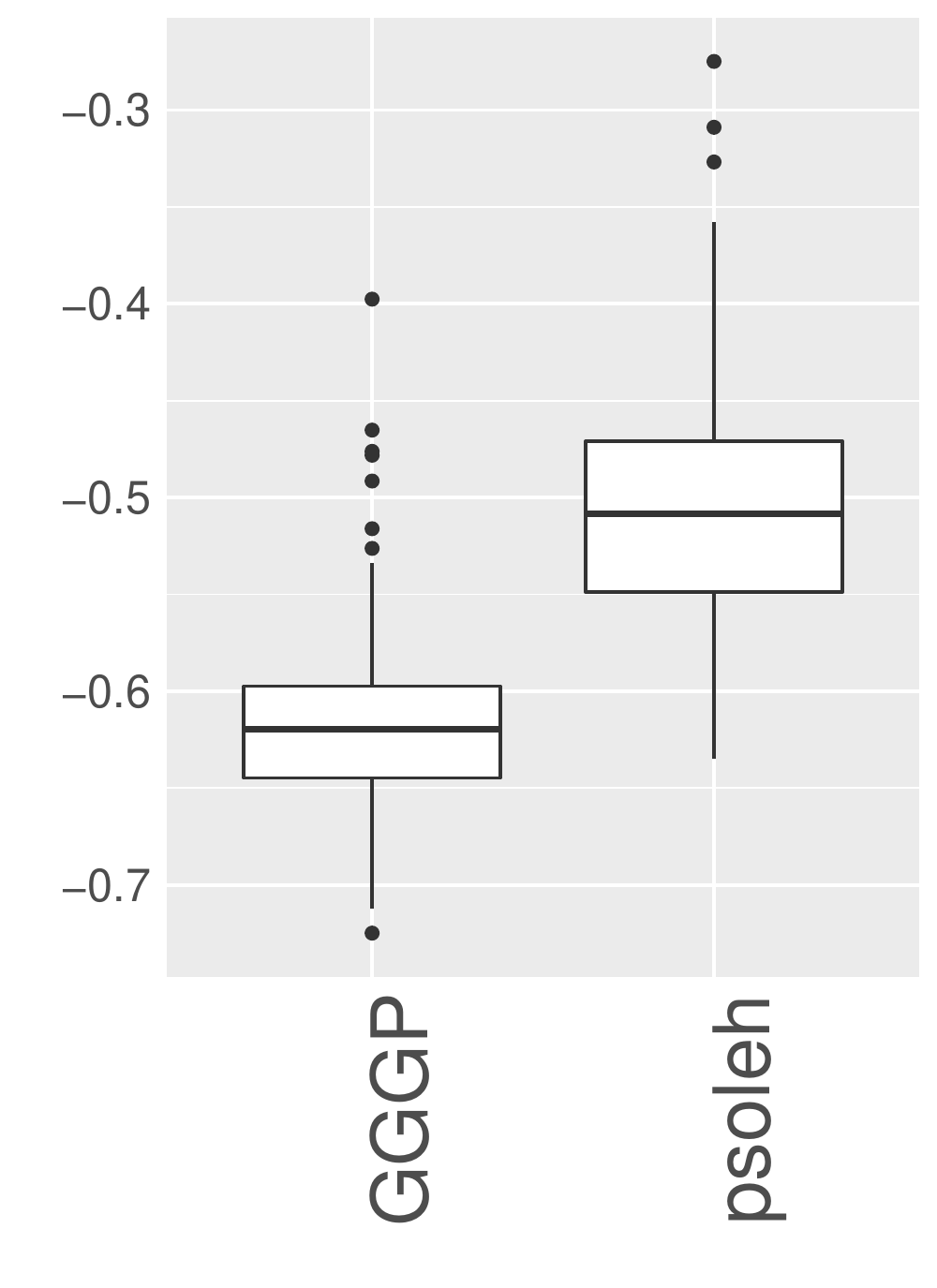}
		\vspace{-4mm} 
  	\caption{\scriptsize{Multipeak F2}} \label{fig:box3058i}
	\end{subfigure}%
	\hspace{-1mm} 
	\begin{subfigure}[t]{.18\textwidth}
		\includegraphics[width=\textwidth]{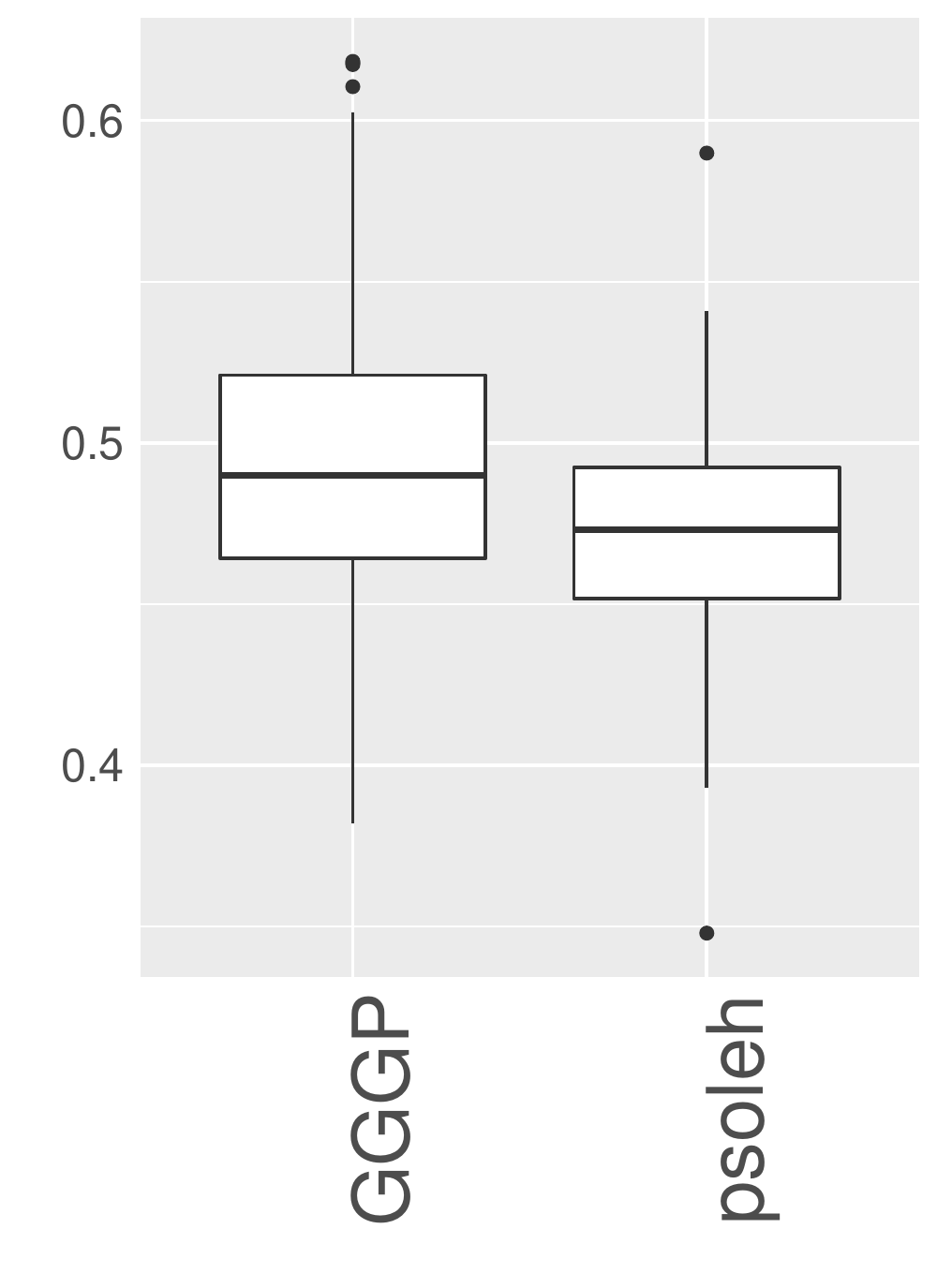}
		\vspace{-4mm} 
  	\caption{\scriptsize{Brankes}} \label{fig:box3067i}
	\end{subfigure}%
	\hspace{-1mm} 
	\begin{subfigure}[t]{.18\textwidth}
		\includegraphics[width=\textwidth]{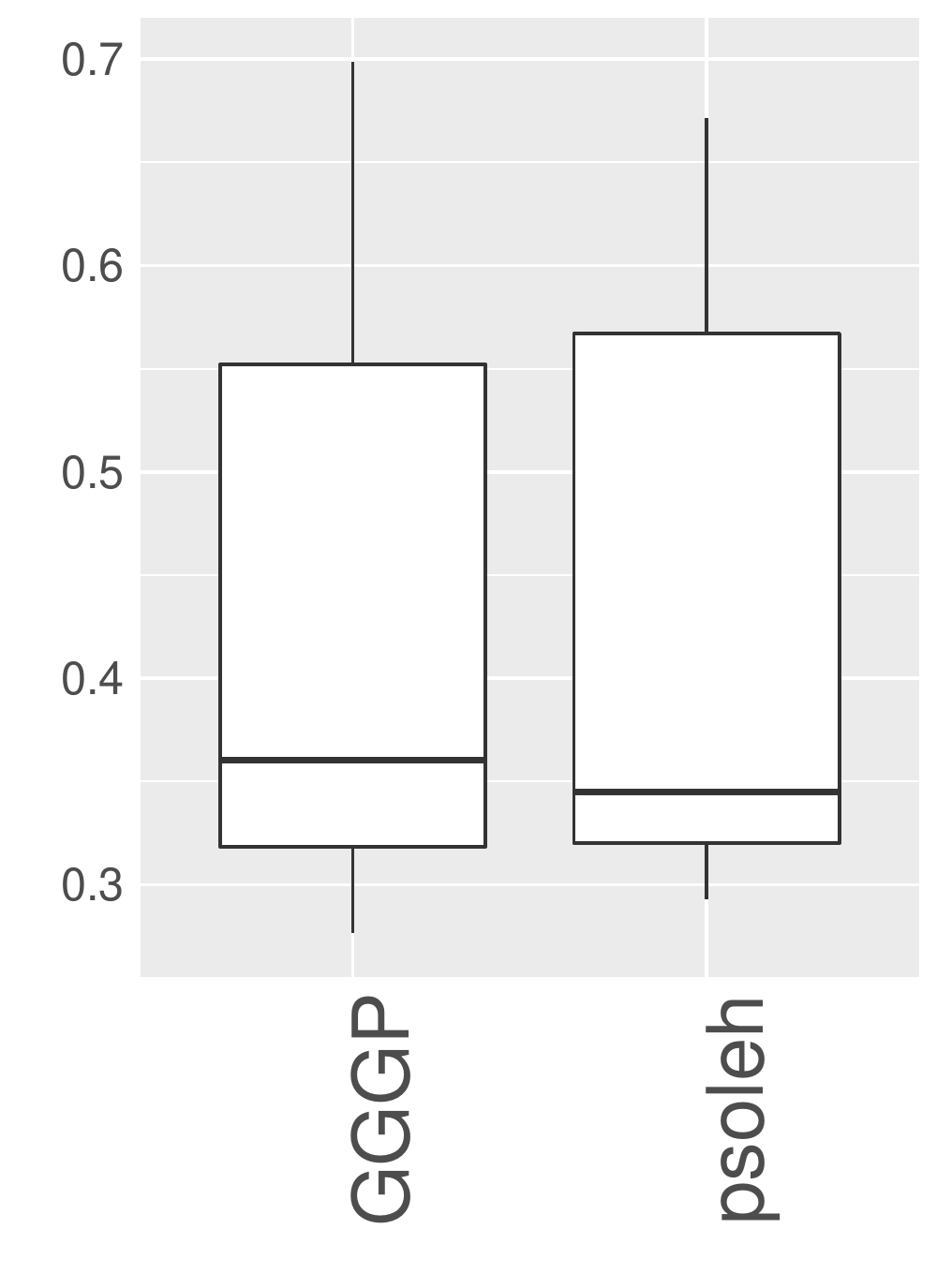}
		\vspace{-4mm} 
  	\caption{\scriptsize{Pickelhaube}} \label{fig:box3065i}
	\end{subfigure}
	
	\vspace{4mm} 
		
	\begin{subfigure}[t]{.18\textwidth}
		\includegraphics[width=\textwidth]{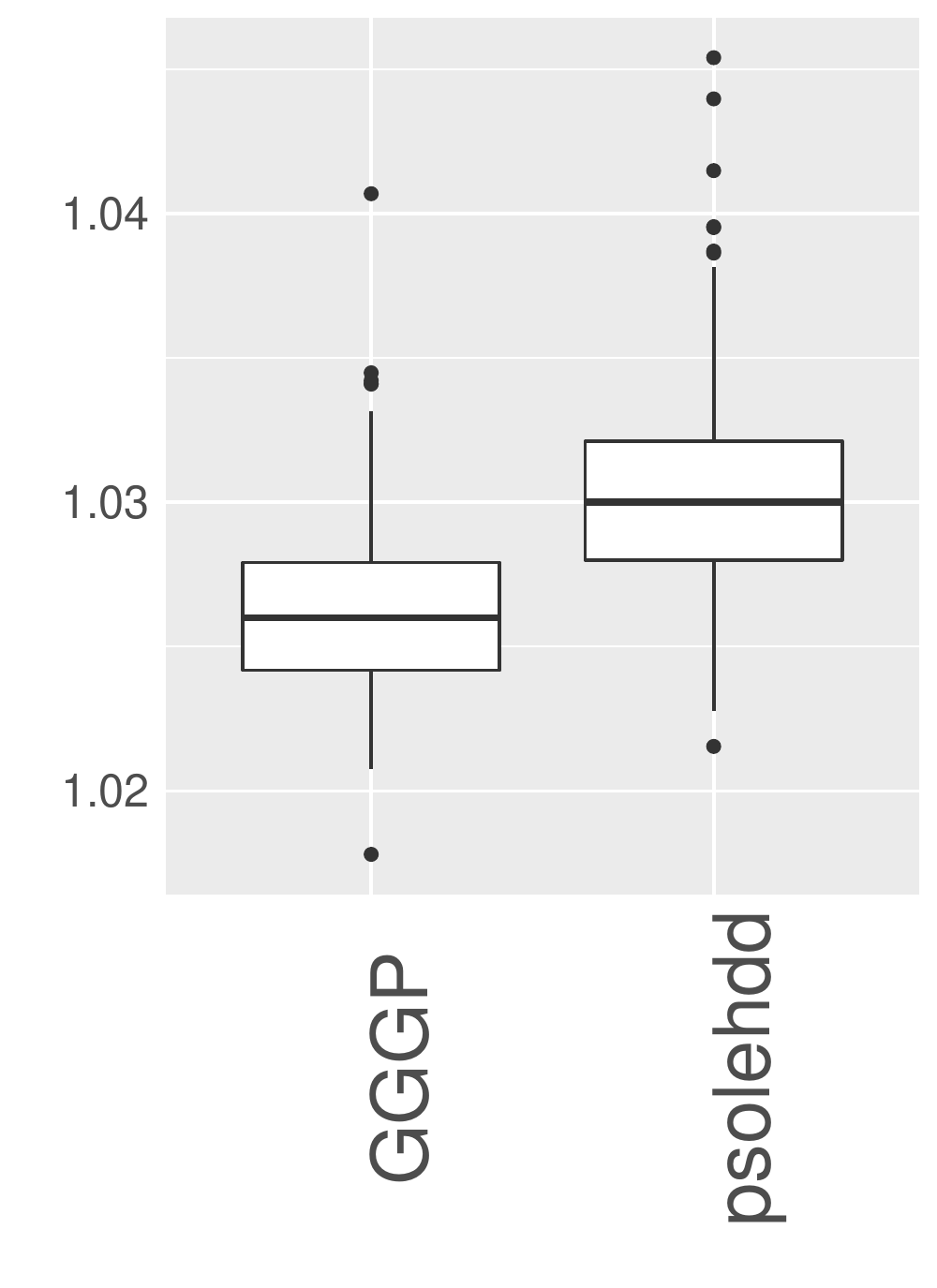}
		\vspace{-4mm} 
  	\caption{\scriptsize{Heaviside}} \label{fig:box3066i}
	\end{subfigure}%
	\hspace{-1mm} 
	\begin{subfigure}[t]{.18\textwidth}
		\includegraphics[width=\textwidth]{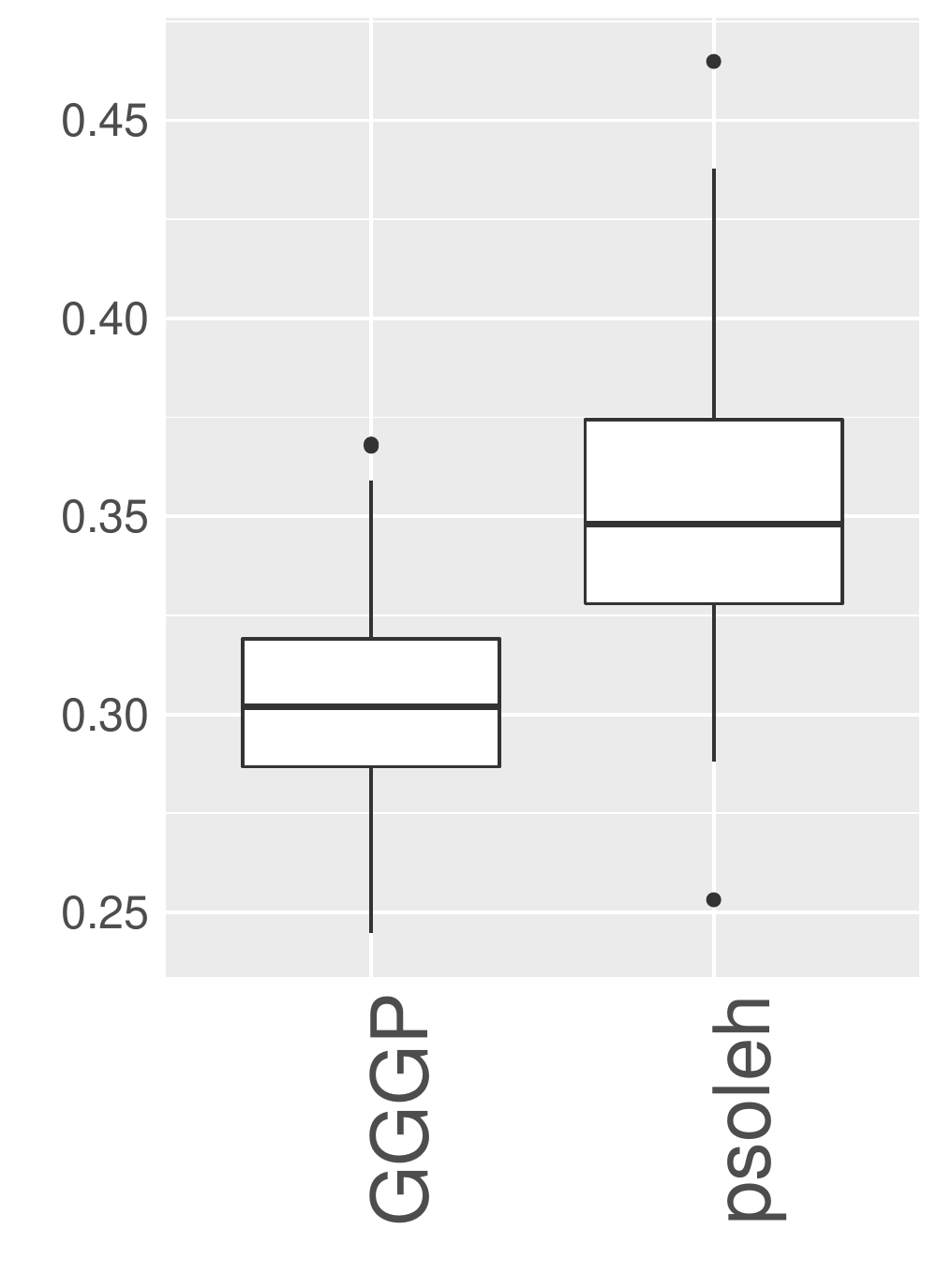}
		\vspace{-4mm} 
  	\caption{\scriptsize{Sawtooth}} \label{fig:box3055i}
	\end{subfigure}%
	\hspace{-1mm} 
	\begin{subfigure}[t]{.18\textwidth}
		\includegraphics[width=\textwidth]{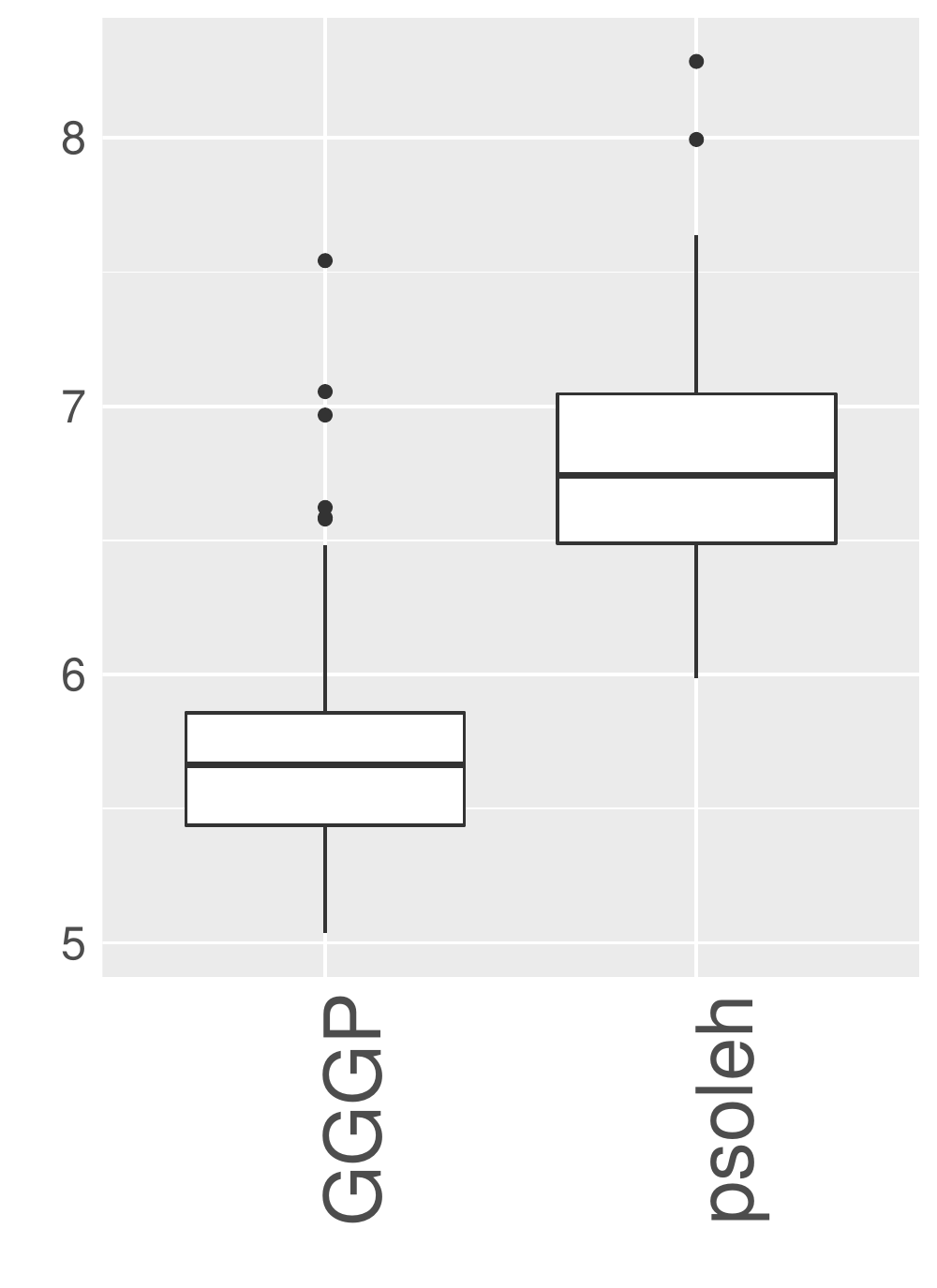}
		\vspace{-4mm} 
  	\caption{\scriptsize{Ackley}} \label{fig:box3052i}
	\end{subfigure}%
	\hspace{-1mm} 
	\begin{subfigure}[t]{.18\textwidth}
		\includegraphics[width=\textwidth]{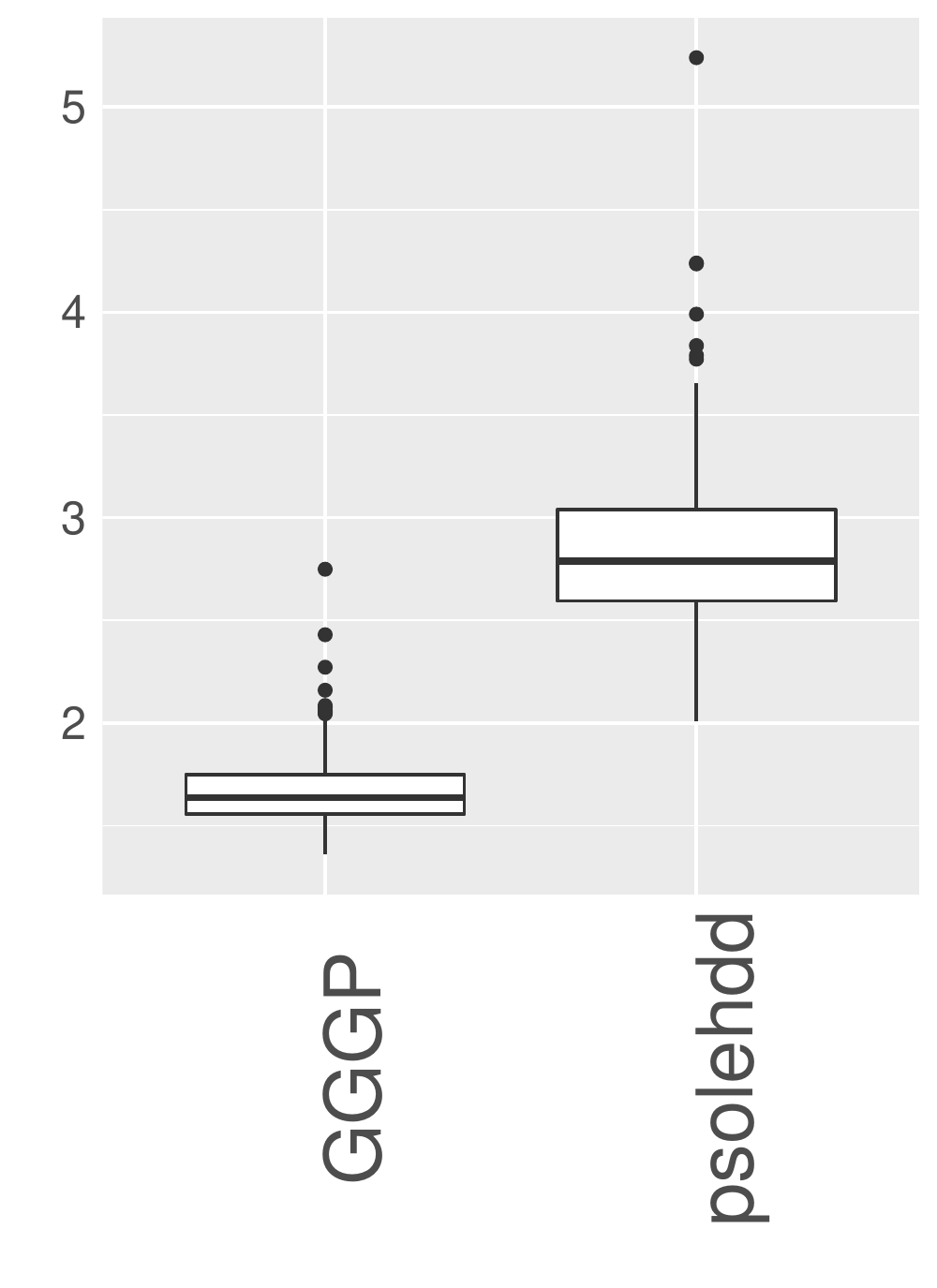}
		\vspace{-4mm} 
  	\caption{\scriptsize{Sphere}} \label{fig:box3054i}
	\end{subfigure}%
	\hspace{-1mm} 
	\begin{subfigure}[t]{.18\textwidth}
		\includegraphics[width=\textwidth]{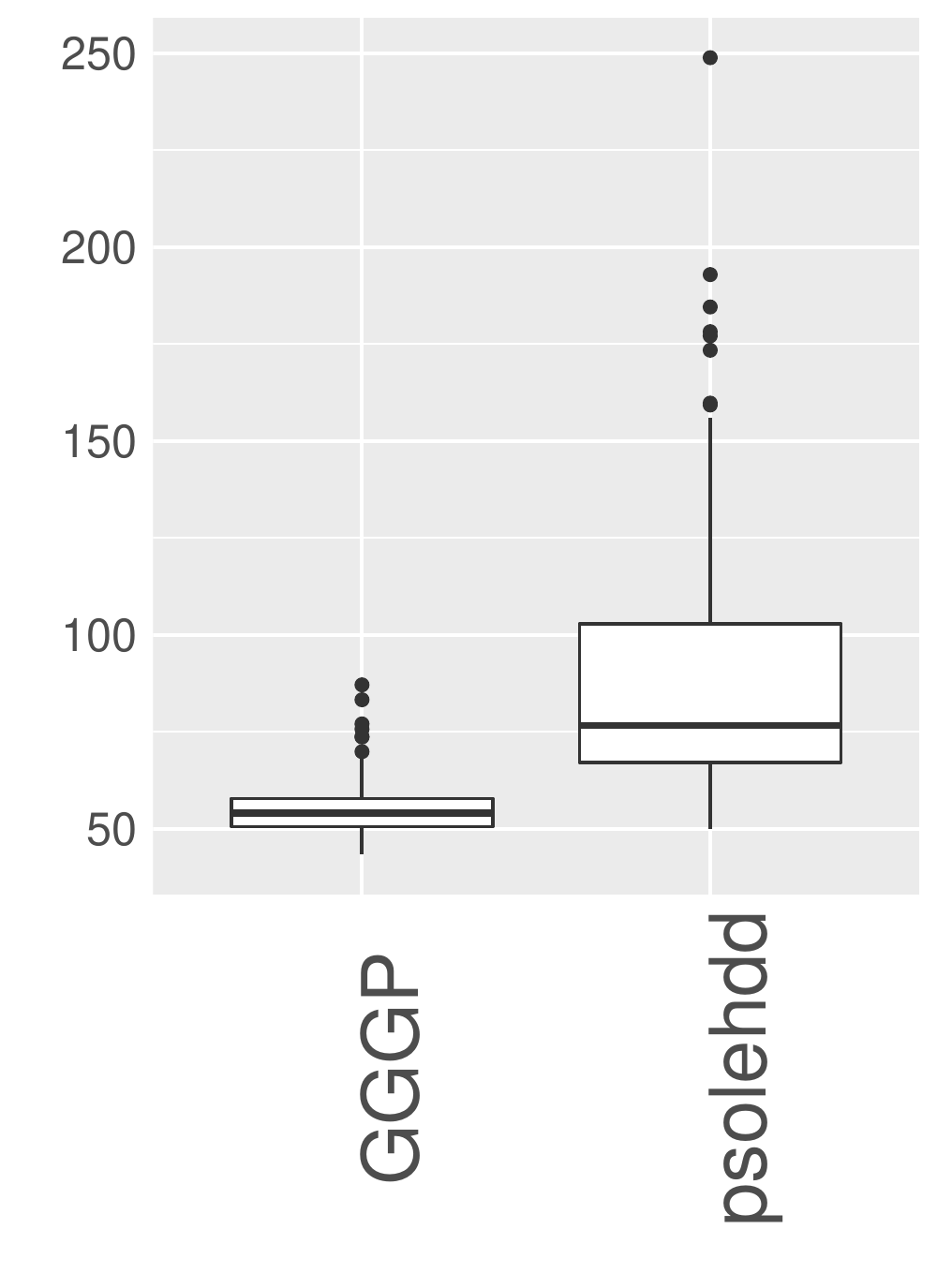}
		\vspace{-4mm} 
  	\caption{\scriptsize{Rosenbrock}} \label{fig:box3053i}
	\end{subfigure}
	\vspace{2mm} 
	\caption{30D individual bests box plots. 200 sample runs with a budget of 2,000 function evaluations. The comparators are taken from \cite{HughesGoerigkDokka2020a}, where the budget was 5,000 evaluations.}
	\label{fig:box30ind}
	
\end{figure}

\begin{figure}[H]
	\centering
	

	\begin{subfigure}[t]{.18\textwidth}
		\includegraphics[width=\textwidth]{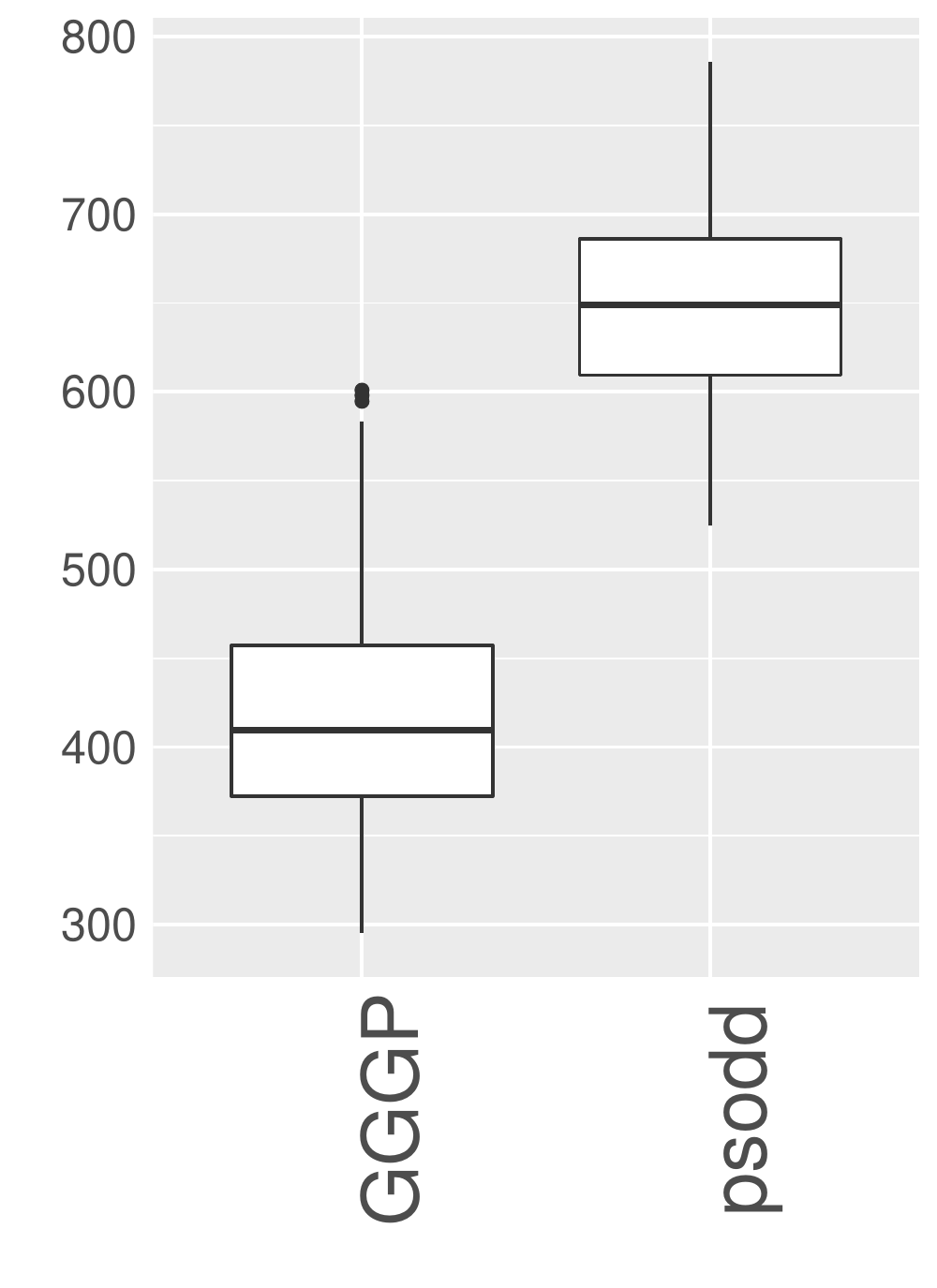}
		\vspace{-4mm} 
  	\caption{\scriptsize{Rastrigin}} \label{fig:box10051i}
	\end{subfigure}%
	\hspace{-1mm} 
	\begin{subfigure}[t]{.18\textwidth}
		\includegraphics[width=\textwidth]{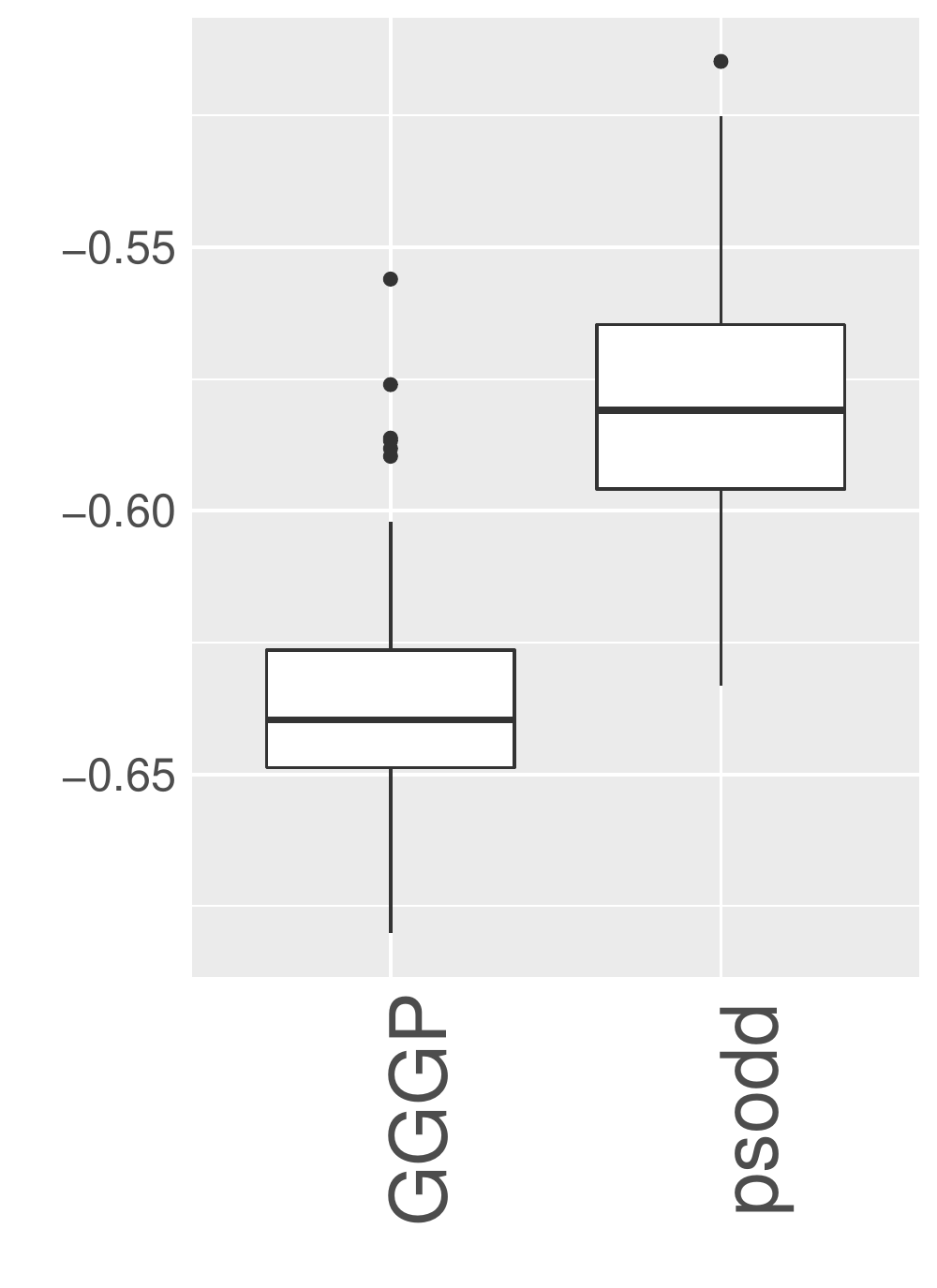}
		\vspace{-4mm} 
  	\caption{\scriptsize{Multipeak F1}} \label{fig:box10057i}
	\end{subfigure}%
	\hspace{-1mm} 
	\begin{subfigure}[t]{.18\textwidth}
		\includegraphics[width=\textwidth]{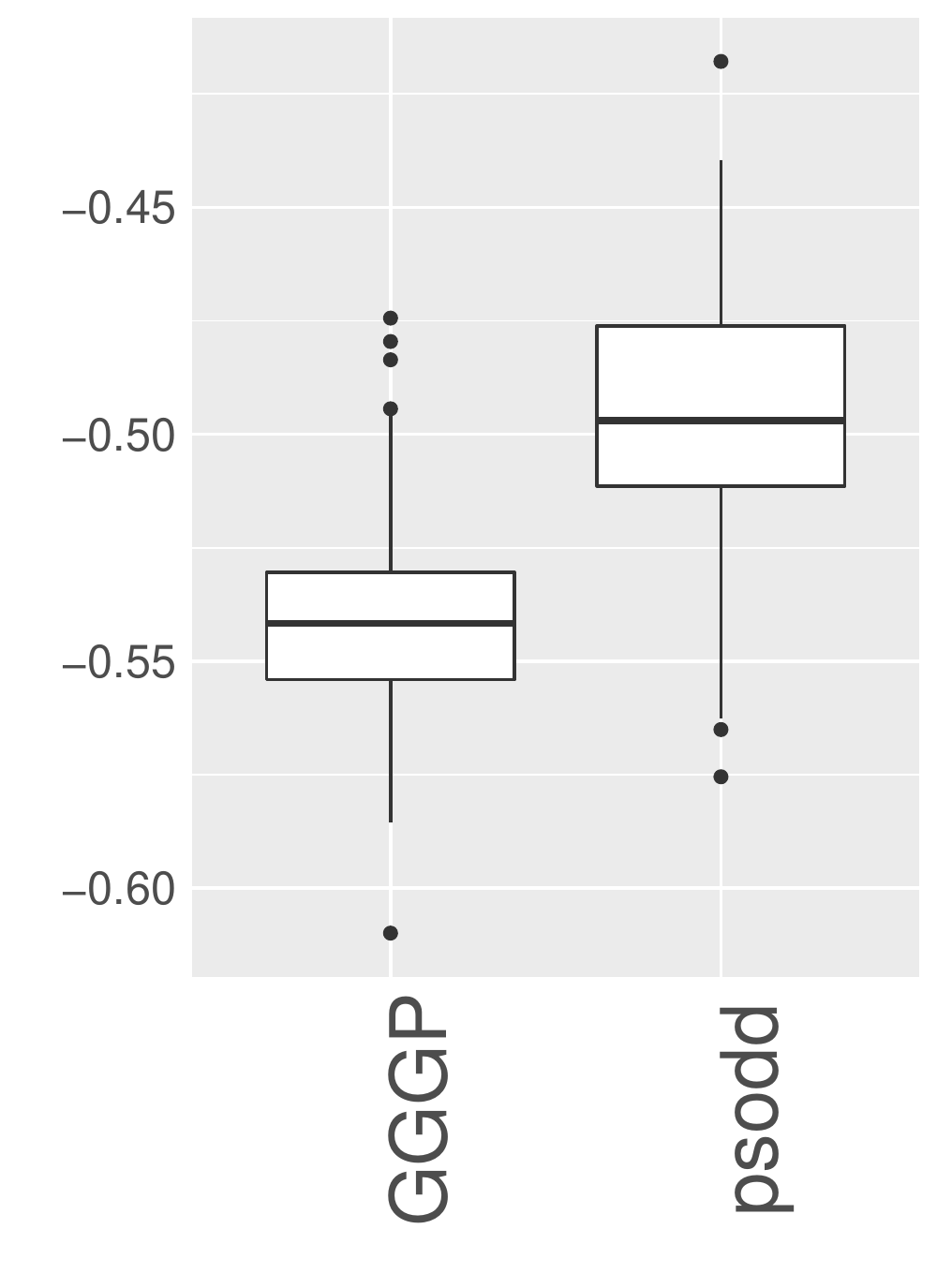}
		\vspace{-4mm} 
  	\caption{\scriptsize{Multipeak F2}} \label{fig:box10058i}
	\end{subfigure}%
	\hspace{-1mm} 
	\begin{subfigure}[t]{.18\textwidth}
		\includegraphics[width=\textwidth]{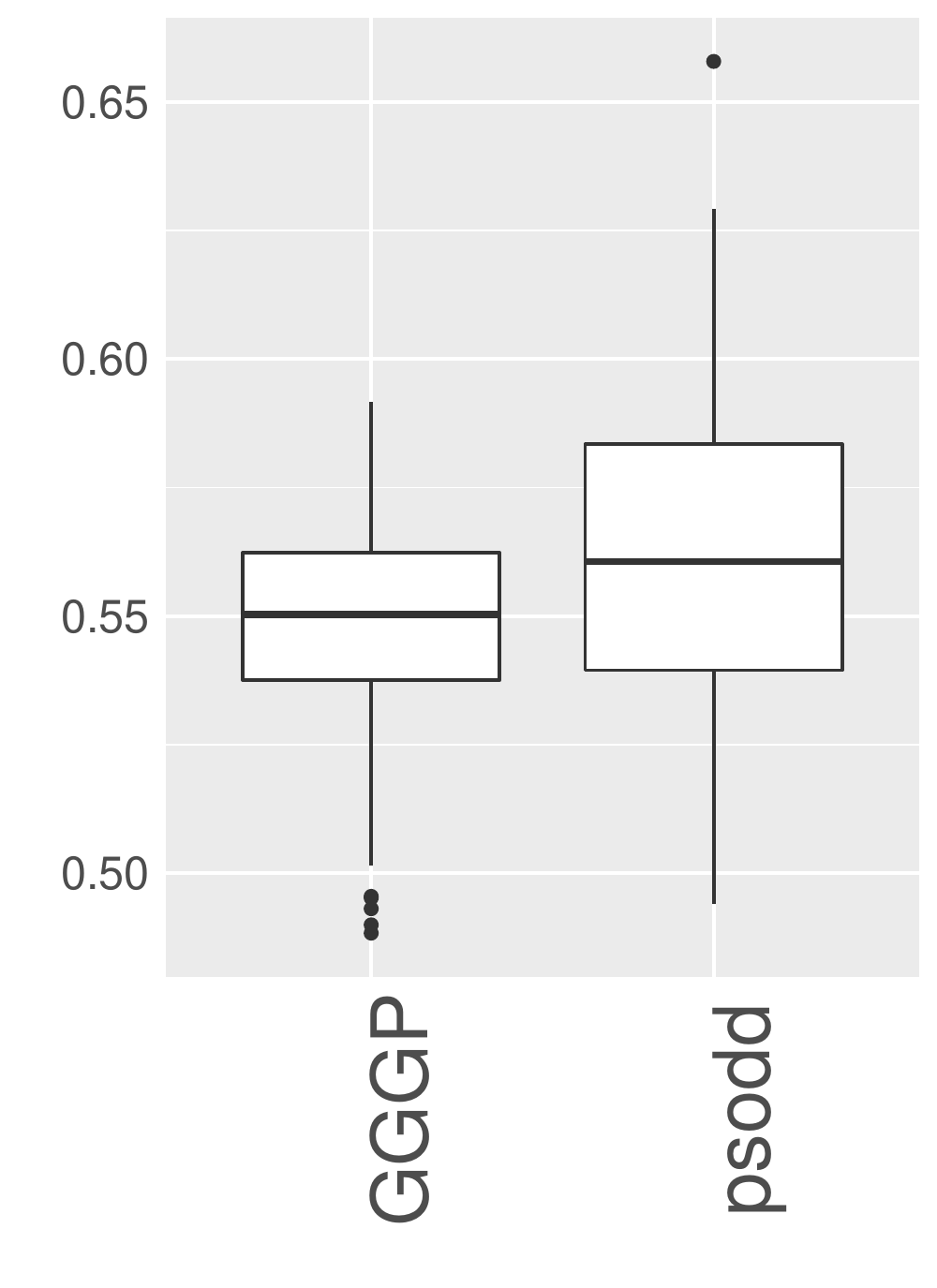}
		\vspace{-4mm} 
  	\caption{\scriptsize{Brankes}} \label{fig:box10067i}
	\end{subfigure}%
	\hspace{-1mm} 
	\begin{subfigure}[t]{.18\textwidth}
		\includegraphics[width=\textwidth]{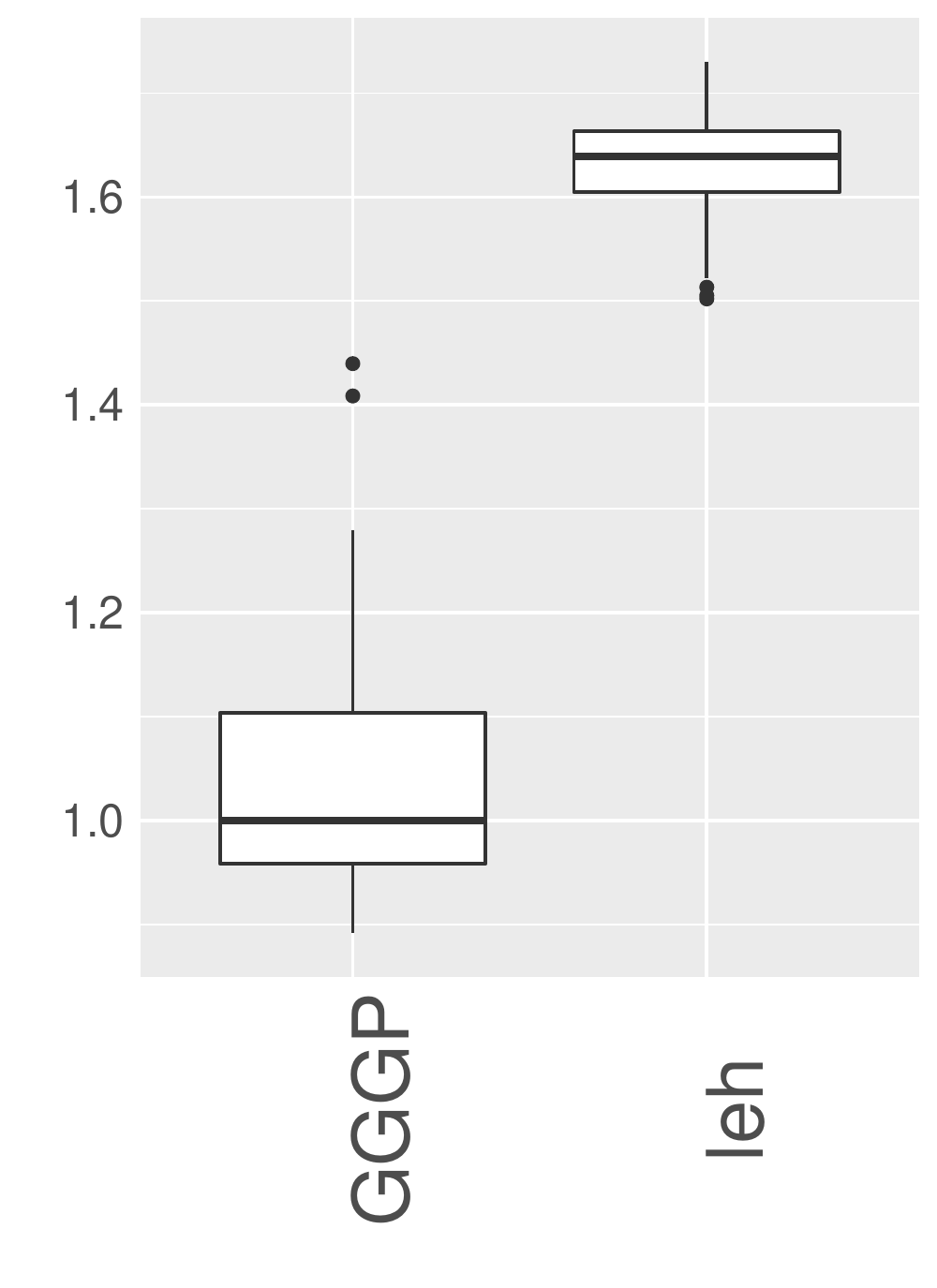}
		\vspace{-4mm} 
  	\caption{\scriptsize{Pickelhaube}} \label{fig:box10065i}
	\end{subfigure}
	
	\vspace{4mm} 
		
	\begin{subfigure}[t]{.18\textwidth}
		\includegraphics[width=\textwidth]{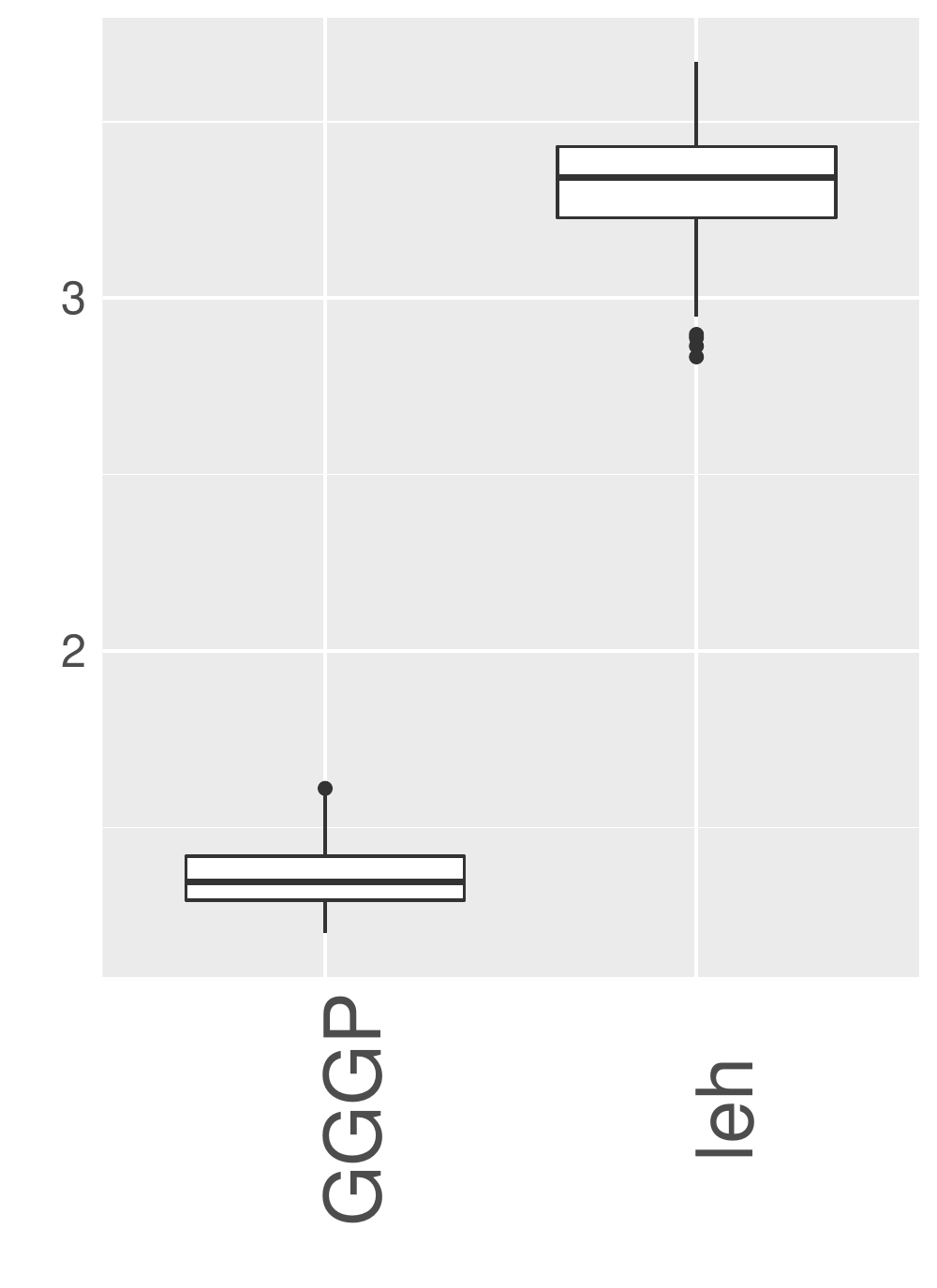}
		\vspace{-4mm} 
  	\caption{\scriptsize{Heaviside}} \label{fig:box10066i}
	\end{subfigure}%
	\hspace{-1mm} 
	\begin{subfigure}[t]{.18\textwidth}
		\includegraphics[width=\textwidth]{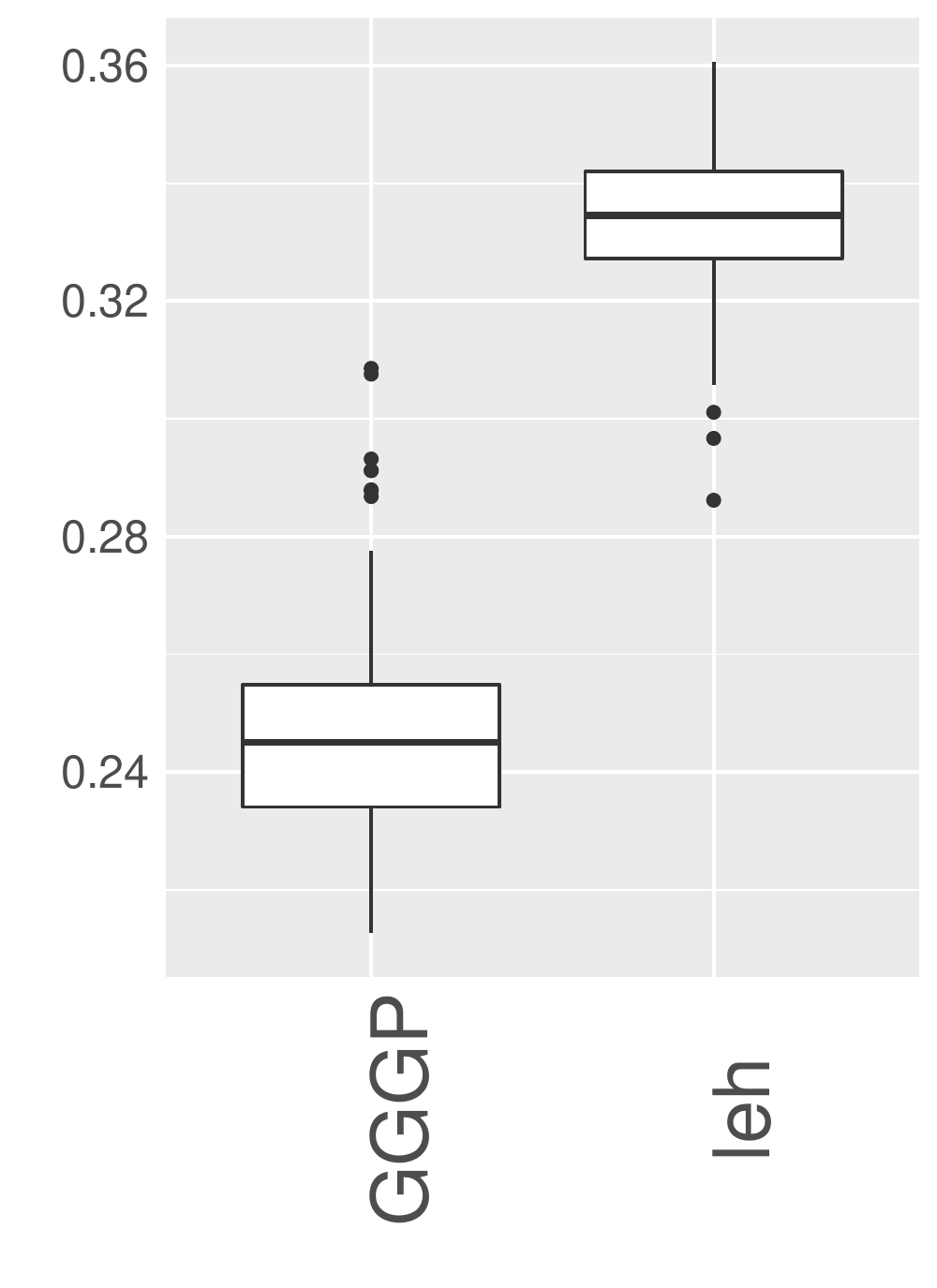}
		\vspace{-4mm} 
  	\caption{\scriptsize{Sawtooth}} \label{fig:box10055i}
	\end{subfigure}%
	\hspace{-1mm} 
	\begin{subfigure}[t]{.18\textwidth}
		\includegraphics[width=\textwidth]{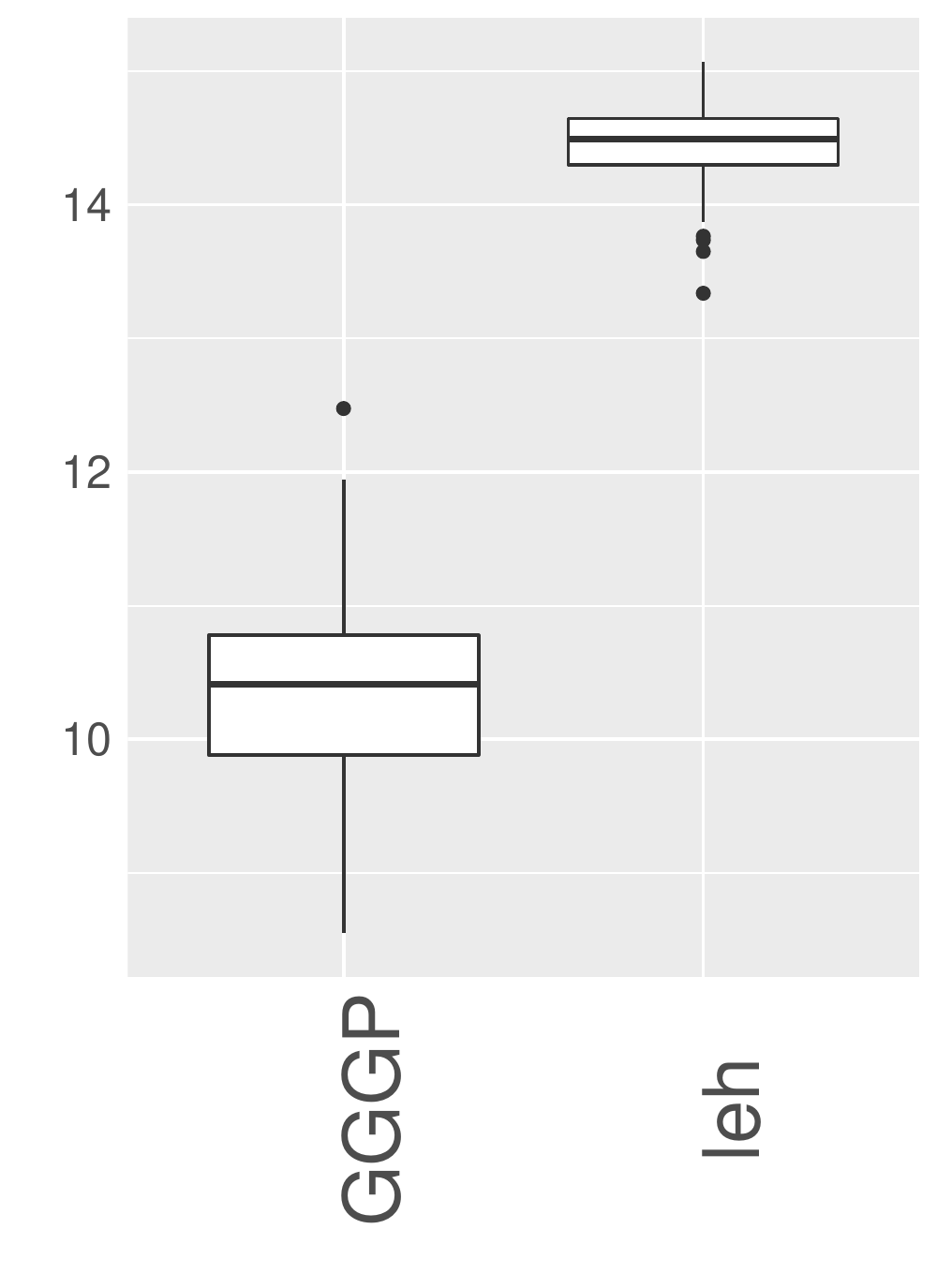}
		\vspace{-4mm} 
  	\caption{\scriptsize{Ackley}} \label{fig:box10052i}
	\end{subfigure}%
	\hspace{-1mm} 
	\begin{subfigure}[t]{.18\textwidth}
		\includegraphics[width=\textwidth]{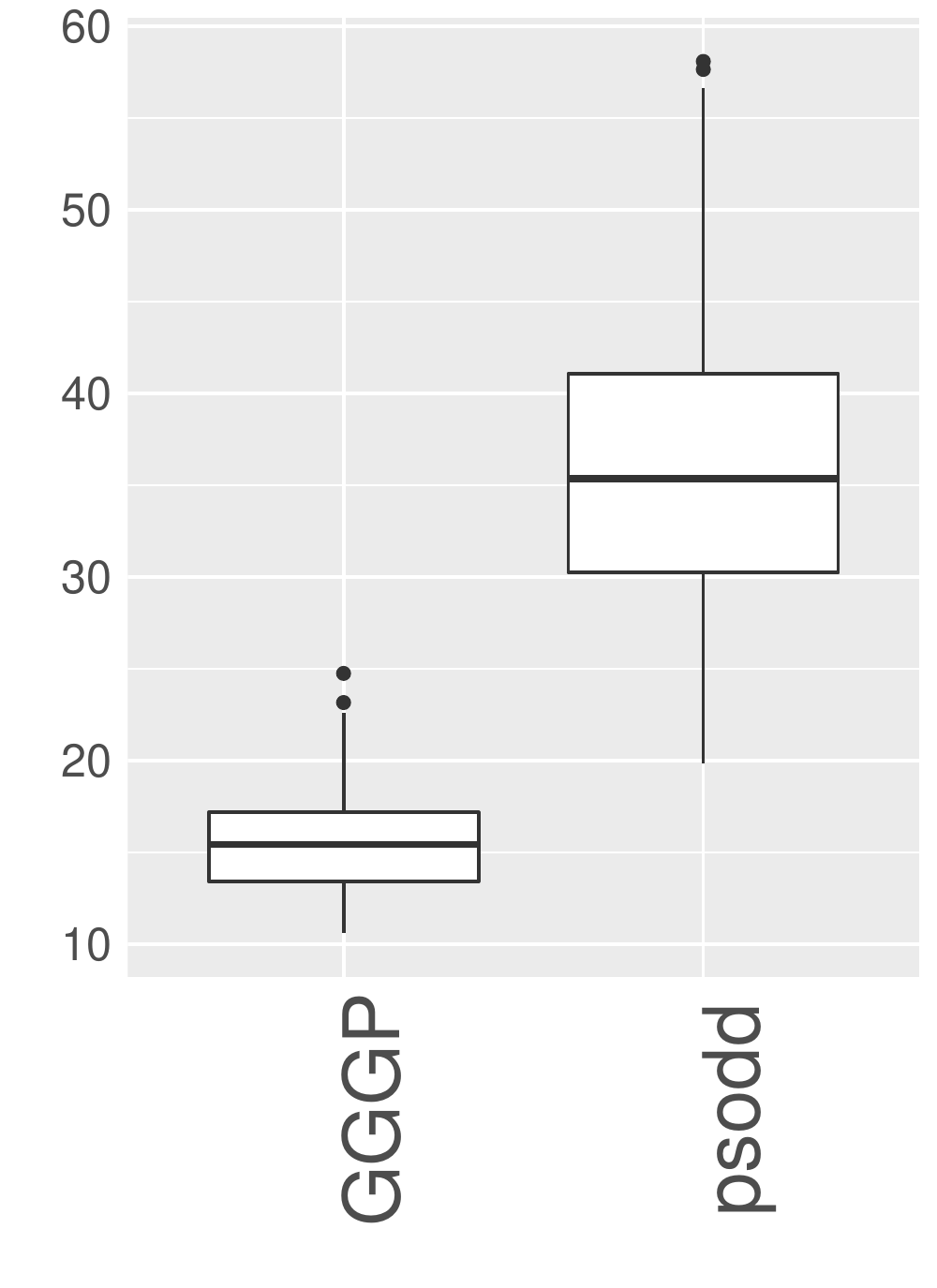}
		\vspace{-4mm} 
  	\caption{\scriptsize{Sphere}} \label{fig:box10054i}
	\end{subfigure}%
	\hspace{-1mm} 
	\begin{subfigure}[t]{.18\textwidth}
		\includegraphics[width=\textwidth]{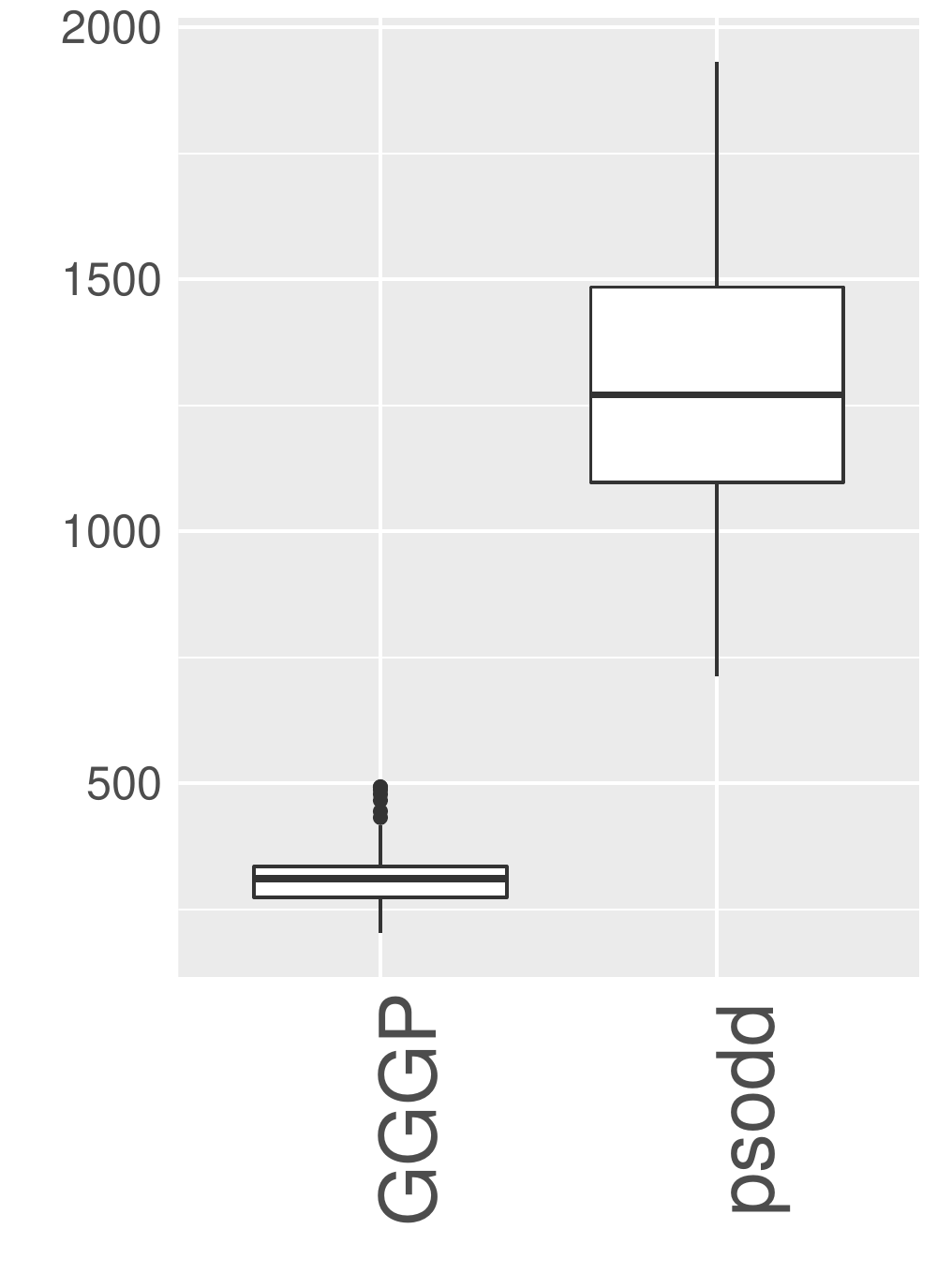}
		\vspace{-4mm} 
  	\caption{\scriptsize{Rosenbrock}} \label{fig:box10053i}
	\end{subfigure}
	\vspace{2mm} 
	\caption{100D individual bests box plots. 200 sample runs with a budget of 2,000 function evaluations. The comparators are taken from \cite{HughesGoerigkDokka2020a}, where the budget was 5,000 evaluations.}
	\label{fig:box100ind}
	
\end{figure}

\begin{figure}[H]
	\centering
	

	\begin{subfigure}[t]{.18\textwidth}
		\includegraphics[width=\textwidth]{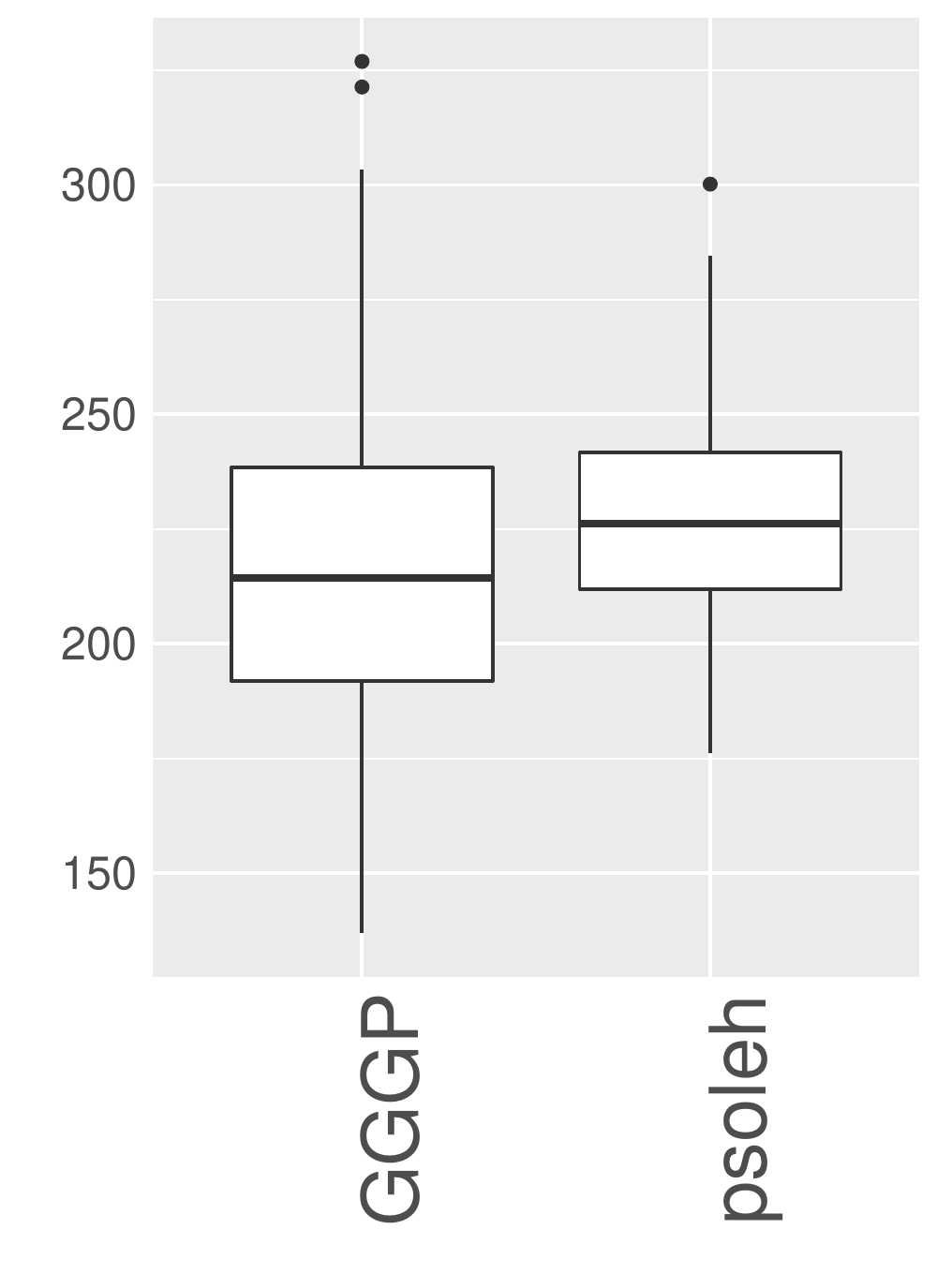}
		\vspace{-4mm} 
  	\caption{\scriptsize{Rastrigin}} \label{fig:box3051a}
	\end{subfigure}%
	\hspace{-1mm} 
	\begin{subfigure}[t]{.18\textwidth}
		\includegraphics[width=\textwidth]{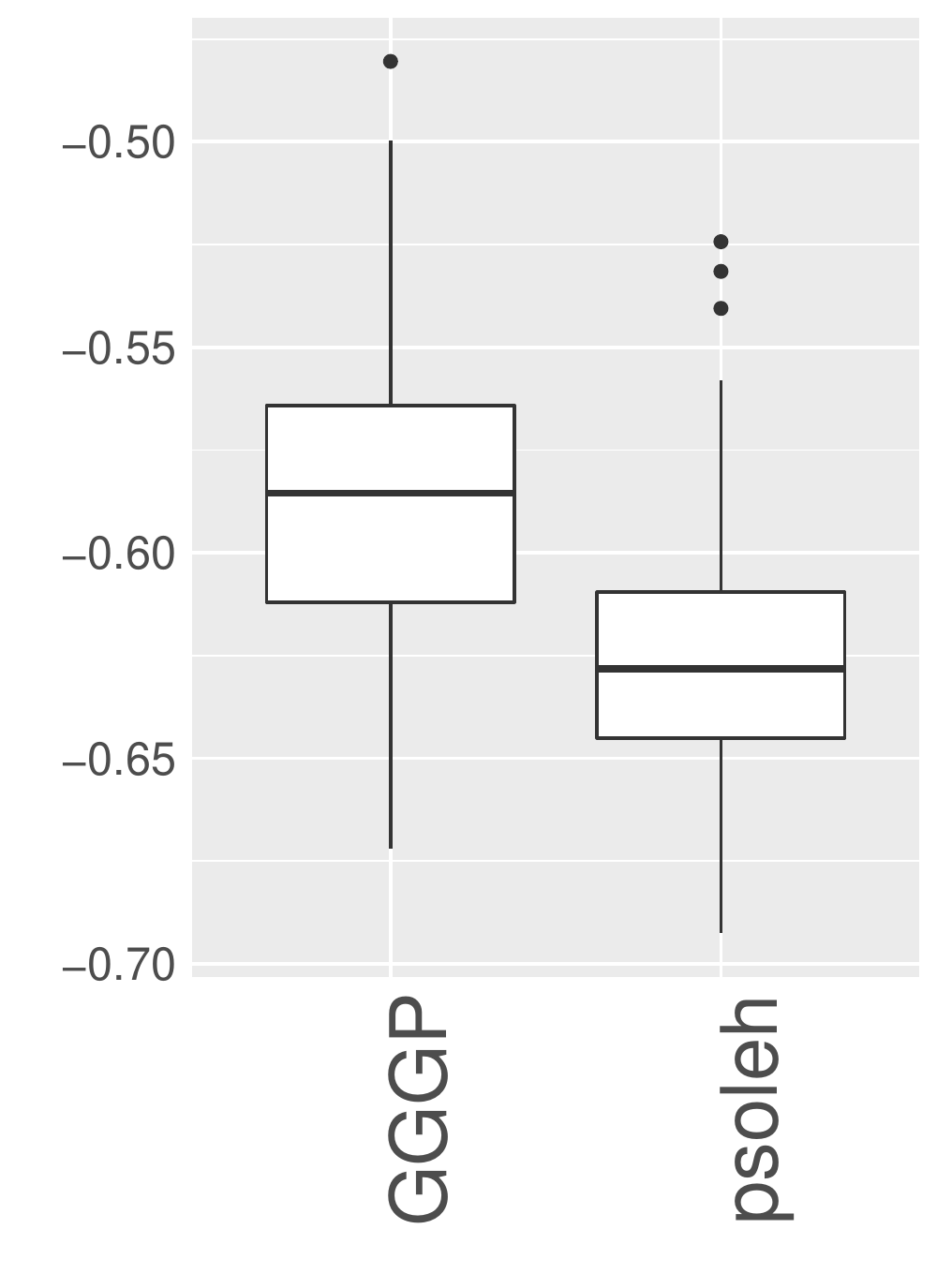}
		\vspace{-4mm} 
  	\caption{\scriptsize{Multipeak F1}} \label{fig:box3057a}
	\end{subfigure}%
	\hspace{-1mm} 
	\begin{subfigure}[t]{.18\textwidth}
		\includegraphics[width=\textwidth]{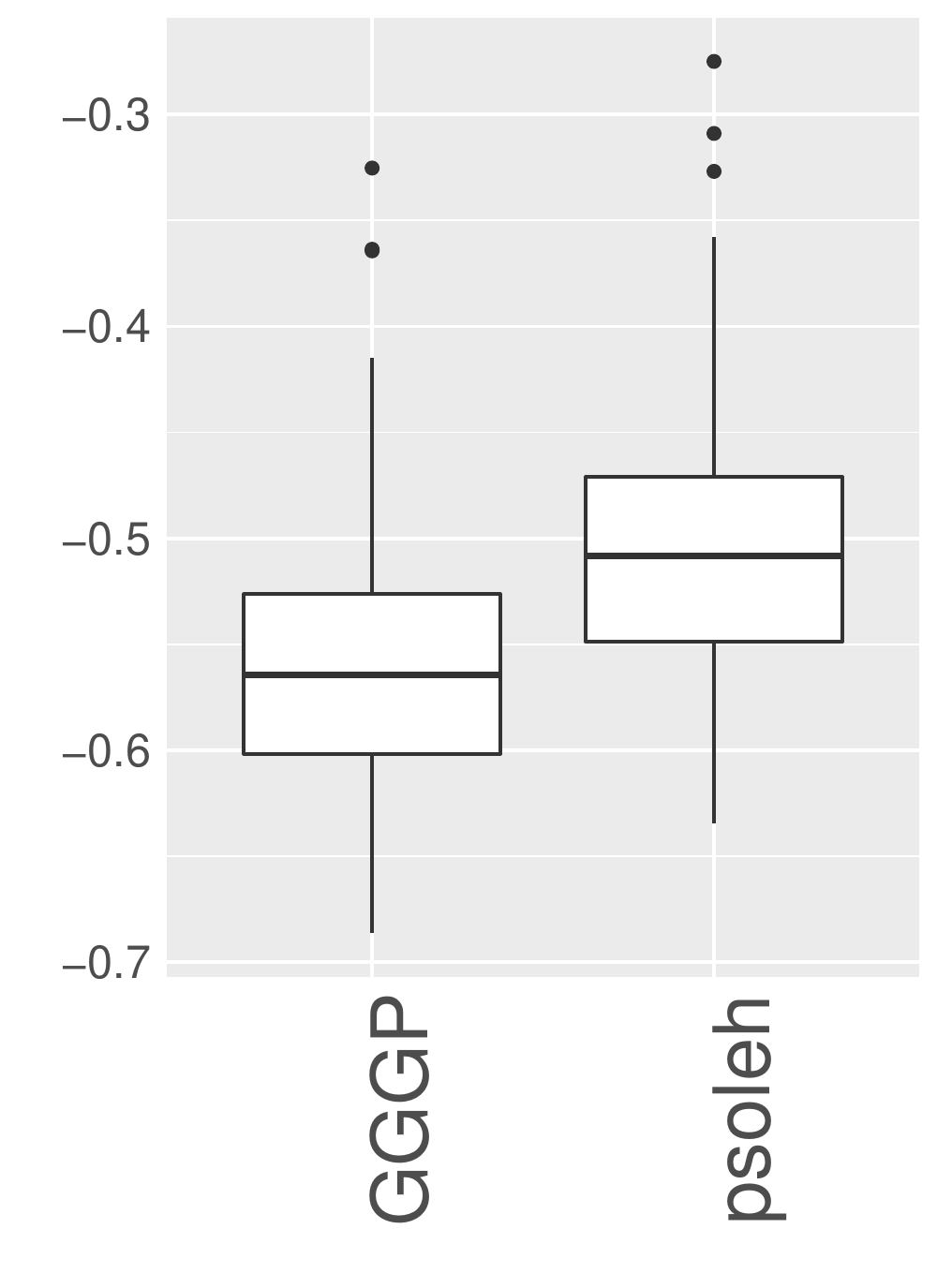}
		\vspace{-4mm} 
  	\caption{\scriptsize{Multipeak F2}} \label{fig:box3058a}
	\end{subfigure}%
	\hspace{-1mm} 
	\begin{subfigure}[t]{.18\textwidth}
		\includegraphics[width=\textwidth]{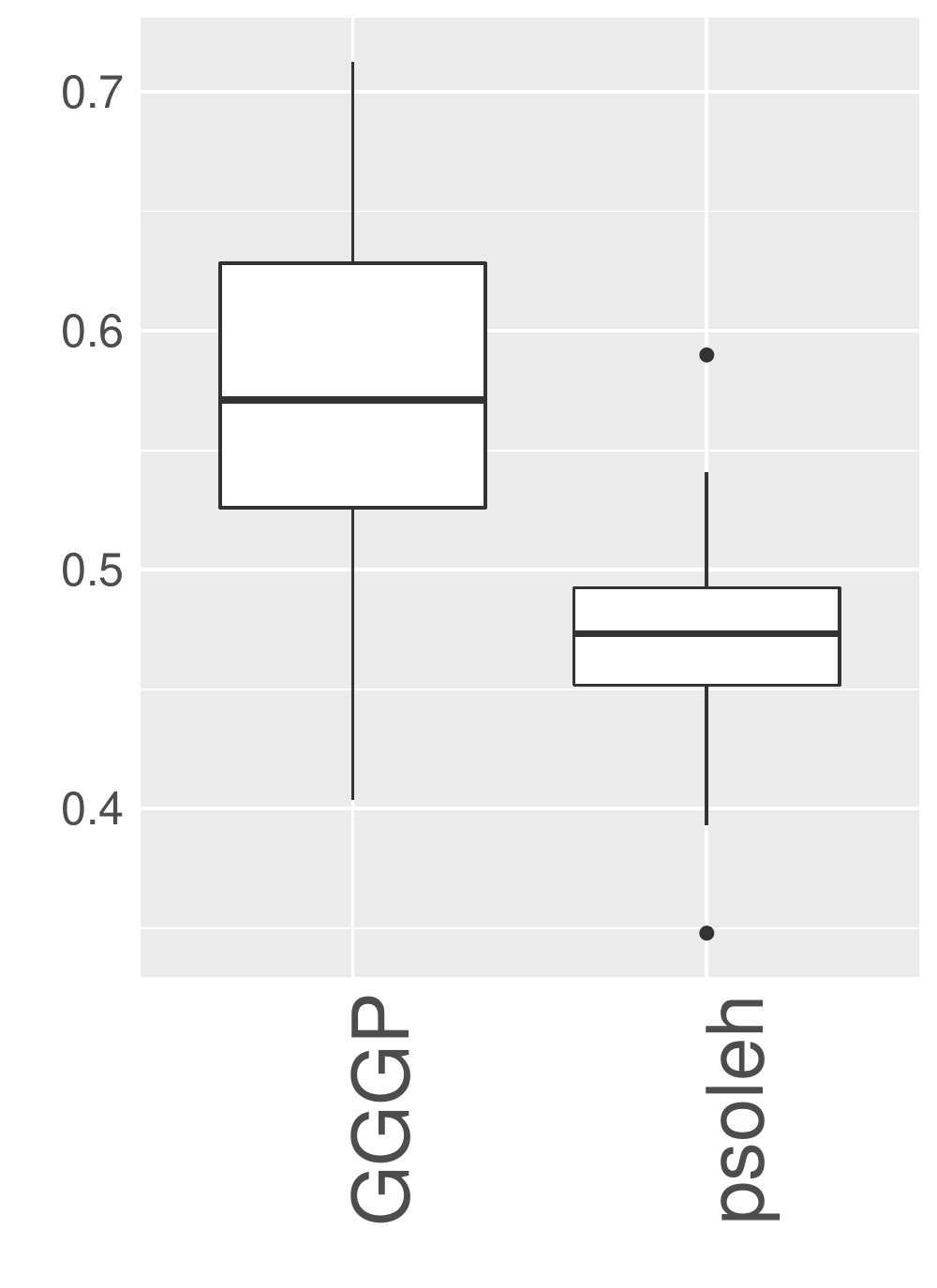}
		\vspace{-4mm} 
  	\caption{\scriptsize{Brankes}} \label{fig:box3067a}
	\end{subfigure}%
	\hspace{-1mm} 
	\begin{subfigure}[t]{.18\textwidth}
		\includegraphics[width=\textwidth]{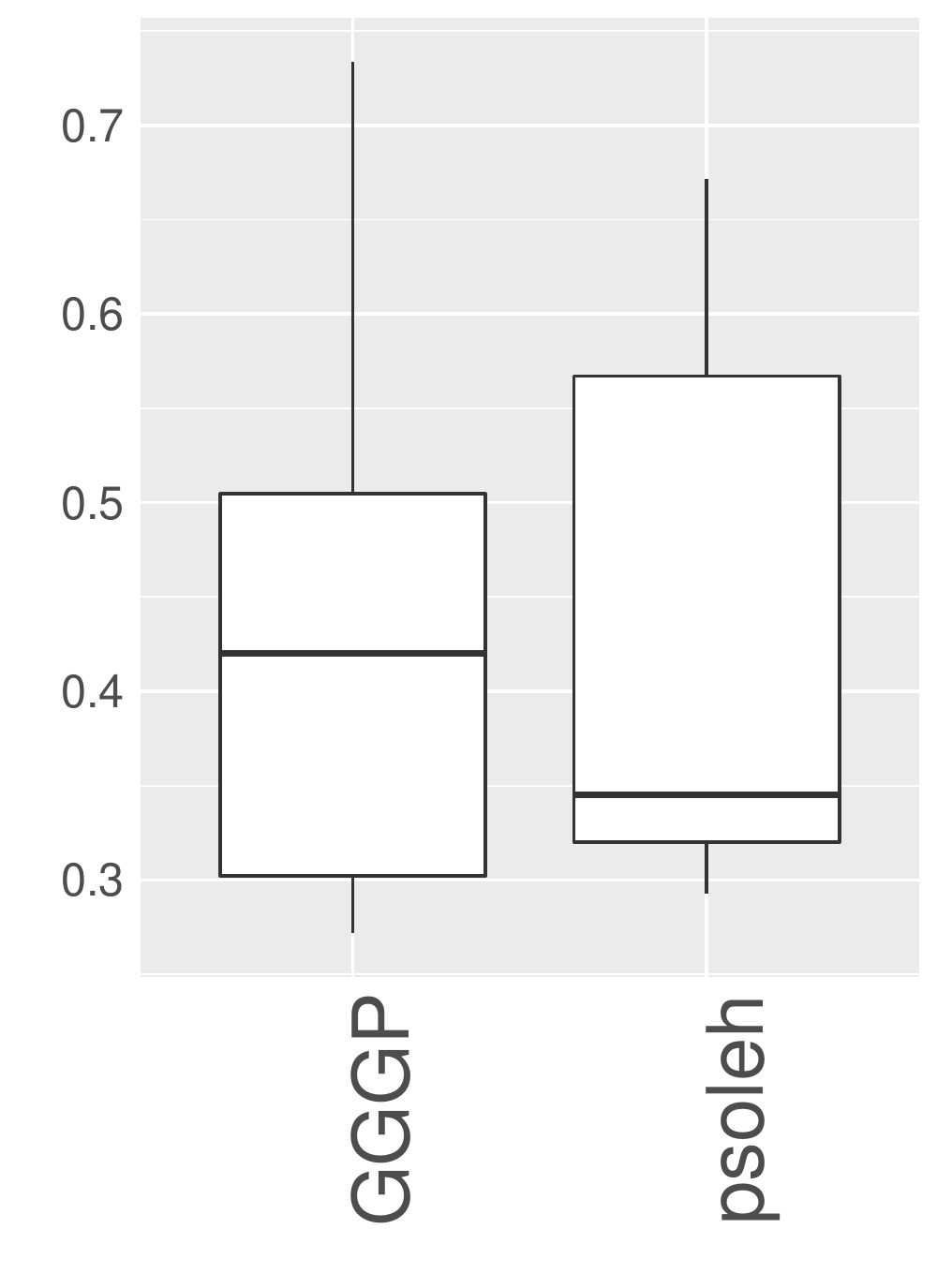}
		\vspace{-4mm} 
  	\caption{\scriptsize{Pickelhaube}} \label{fig:box3065a}
	\end{subfigure}
	
	\vspace{4mm} 
		
	\begin{subfigure}[t]{.18\textwidth}
		\includegraphics[width=\textwidth]{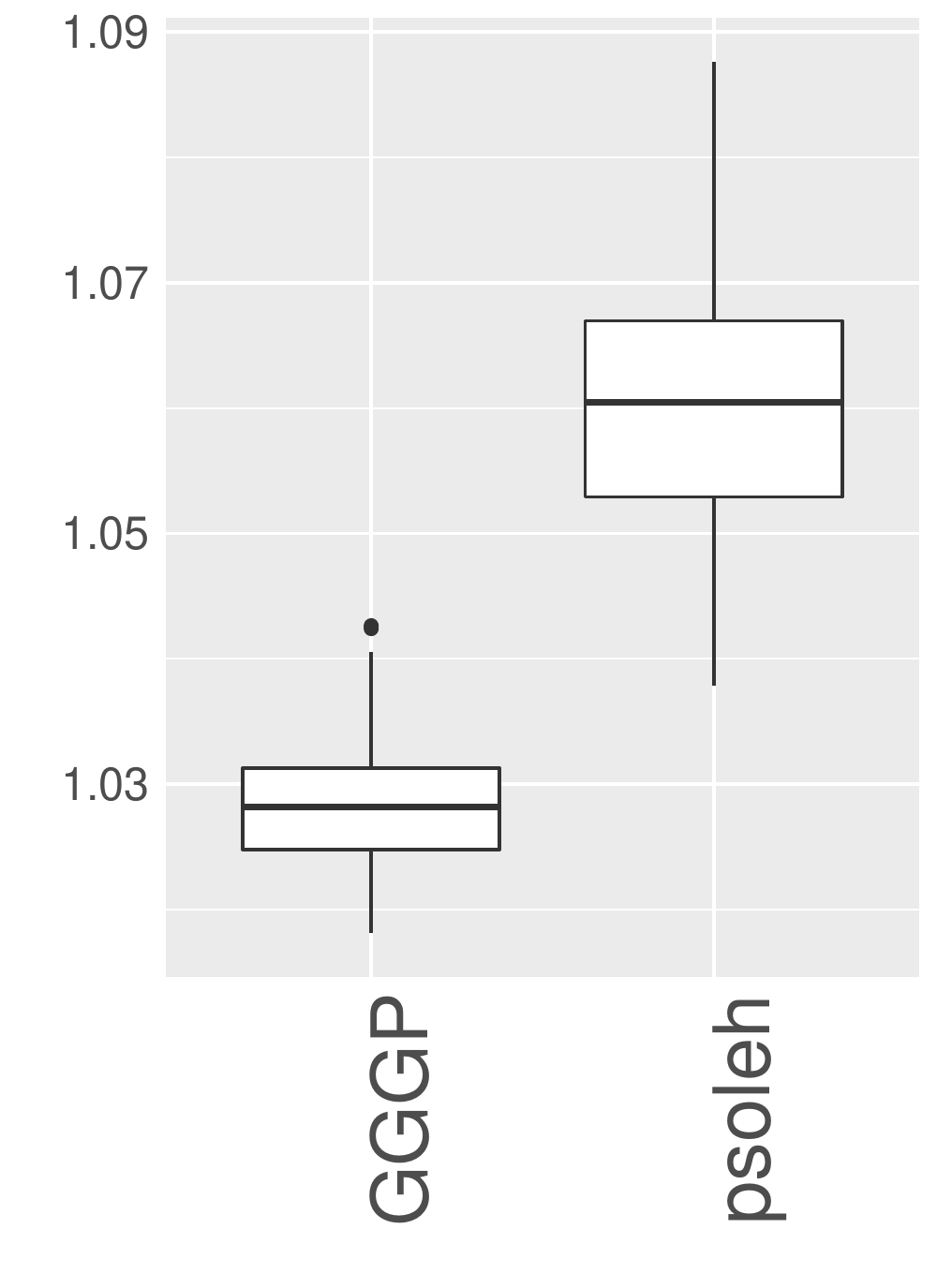}
		\vspace{-4mm} 
  	\caption{\scriptsize{Heaviside}} \label{fig:box3066a}
	\end{subfigure}%
	\hspace{-1mm} 
	\begin{subfigure}[t]{.18\textwidth}
		\includegraphics[width=\textwidth]{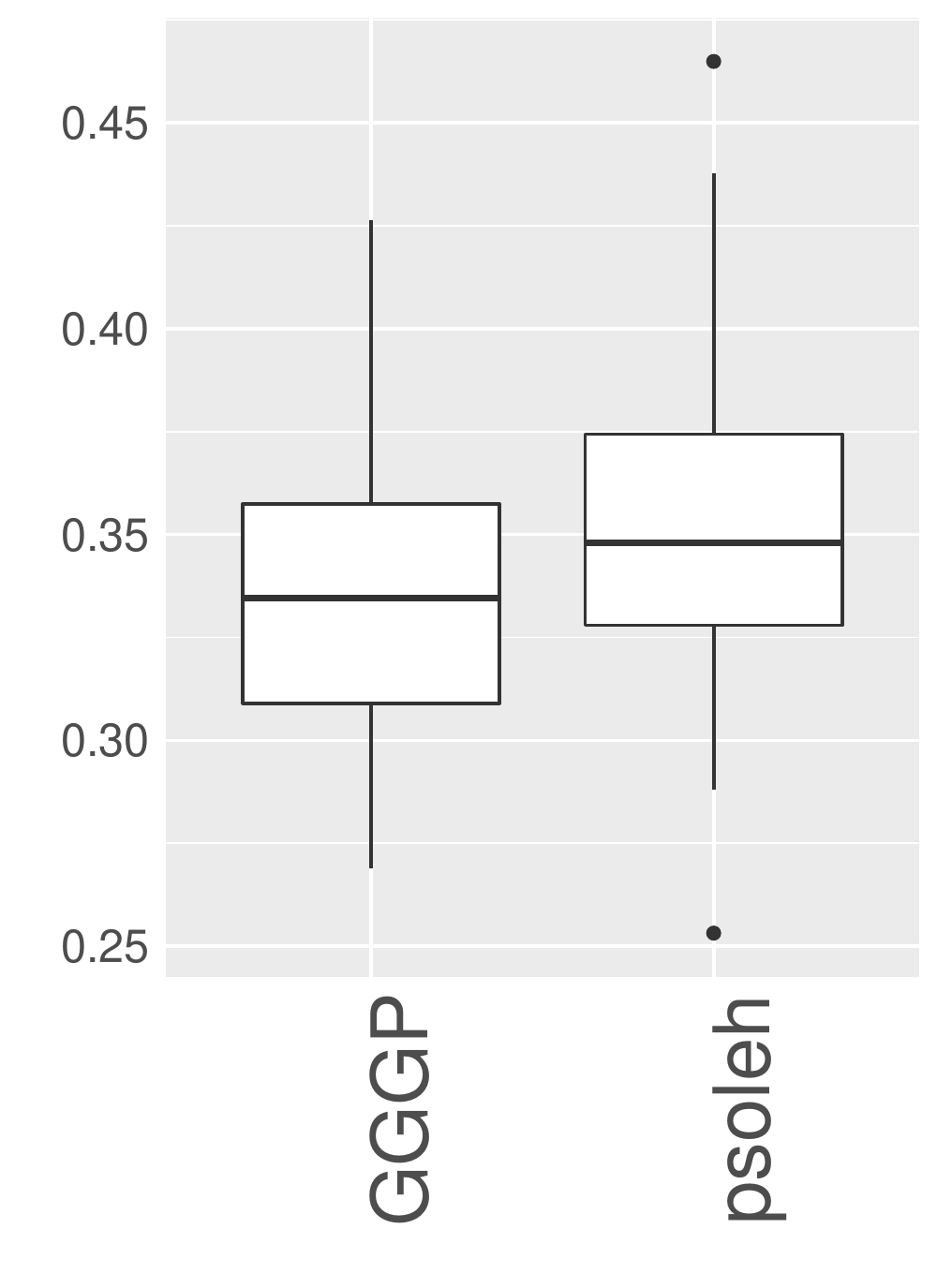}
		\vspace{-4mm} 
  	\caption{\scriptsize{Sawtooth}} \label{fig:box3055a}
	\end{subfigure}%
	\hspace{-1mm} 
	\begin{subfigure}[t]{.18\textwidth}
		\includegraphics[width=\textwidth]{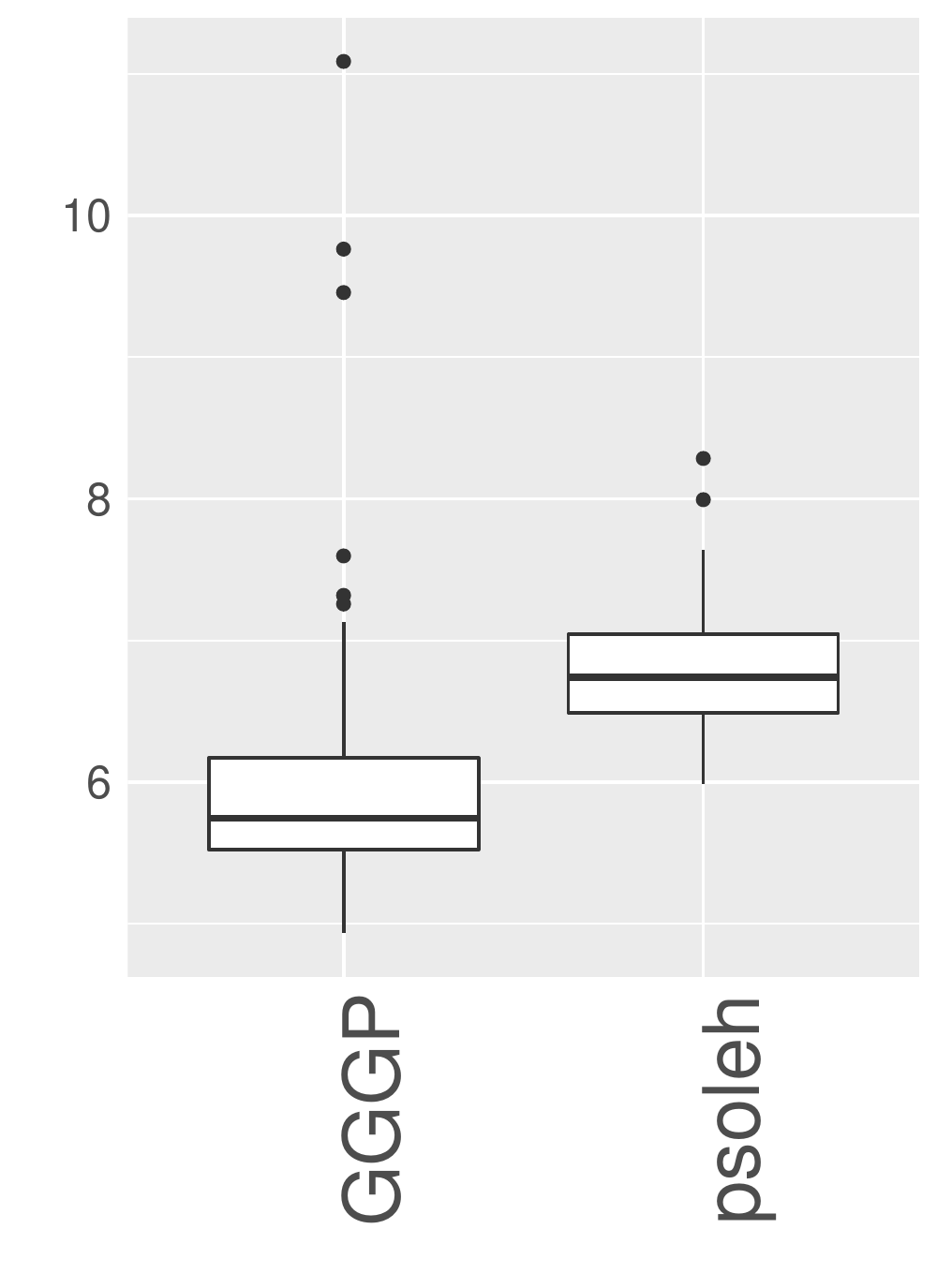}
		\vspace{-4mm} 
  	\caption{\scriptsize{Ackley}} \label{fig:box3052a}
	\end{subfigure}%
	\hspace{-1mm} 
	\begin{subfigure}[t]{.18\textwidth}
		\includegraphics[width=\textwidth]{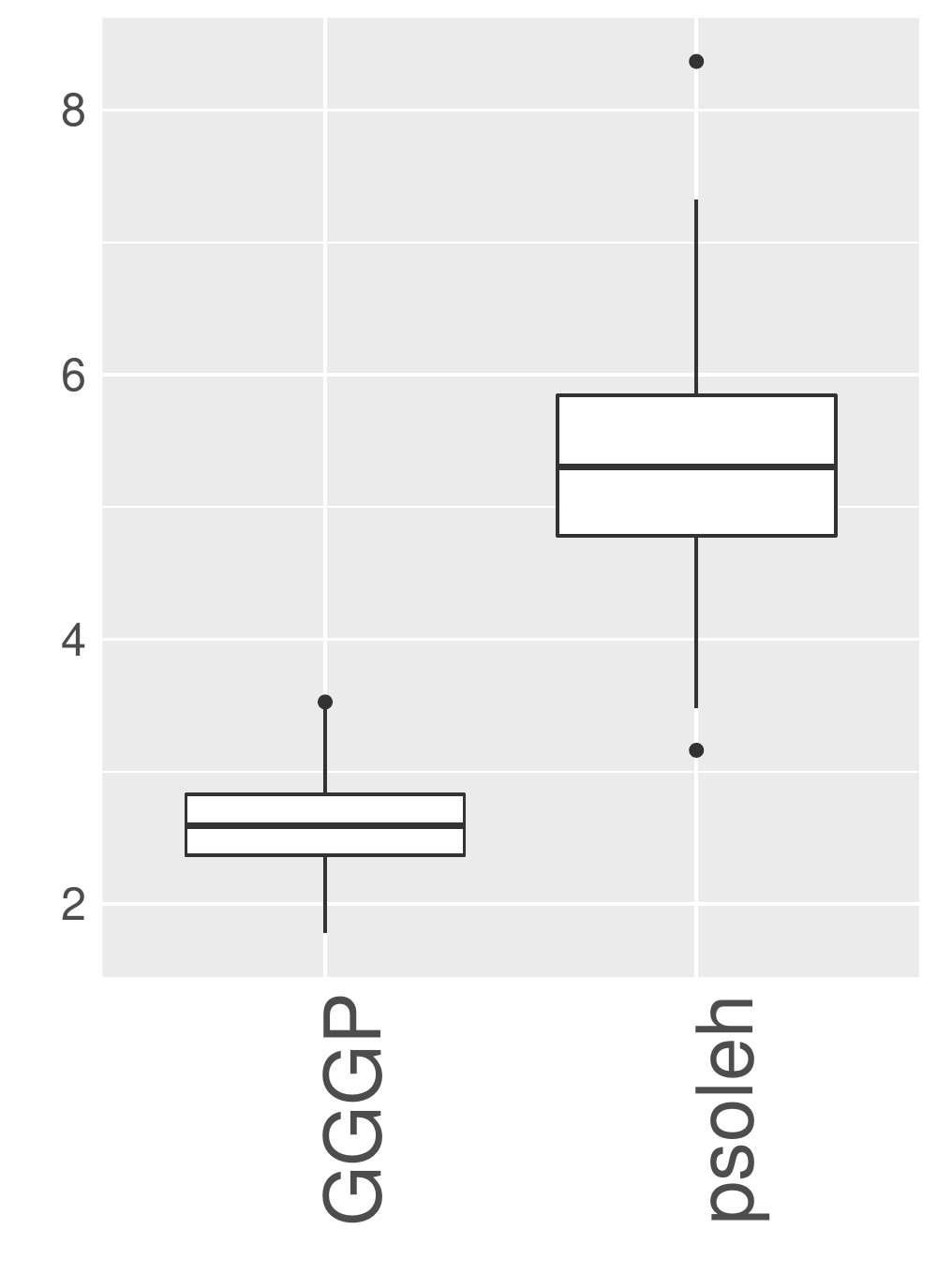}
		\vspace{-4mm} 
  	\caption{\scriptsize{Sphere}} \label{fig:box3054a}
	\end{subfigure}%
	\hspace{-1mm} 
	\begin{subfigure}[t]{.18\textwidth}
		\includegraphics[width=\textwidth]{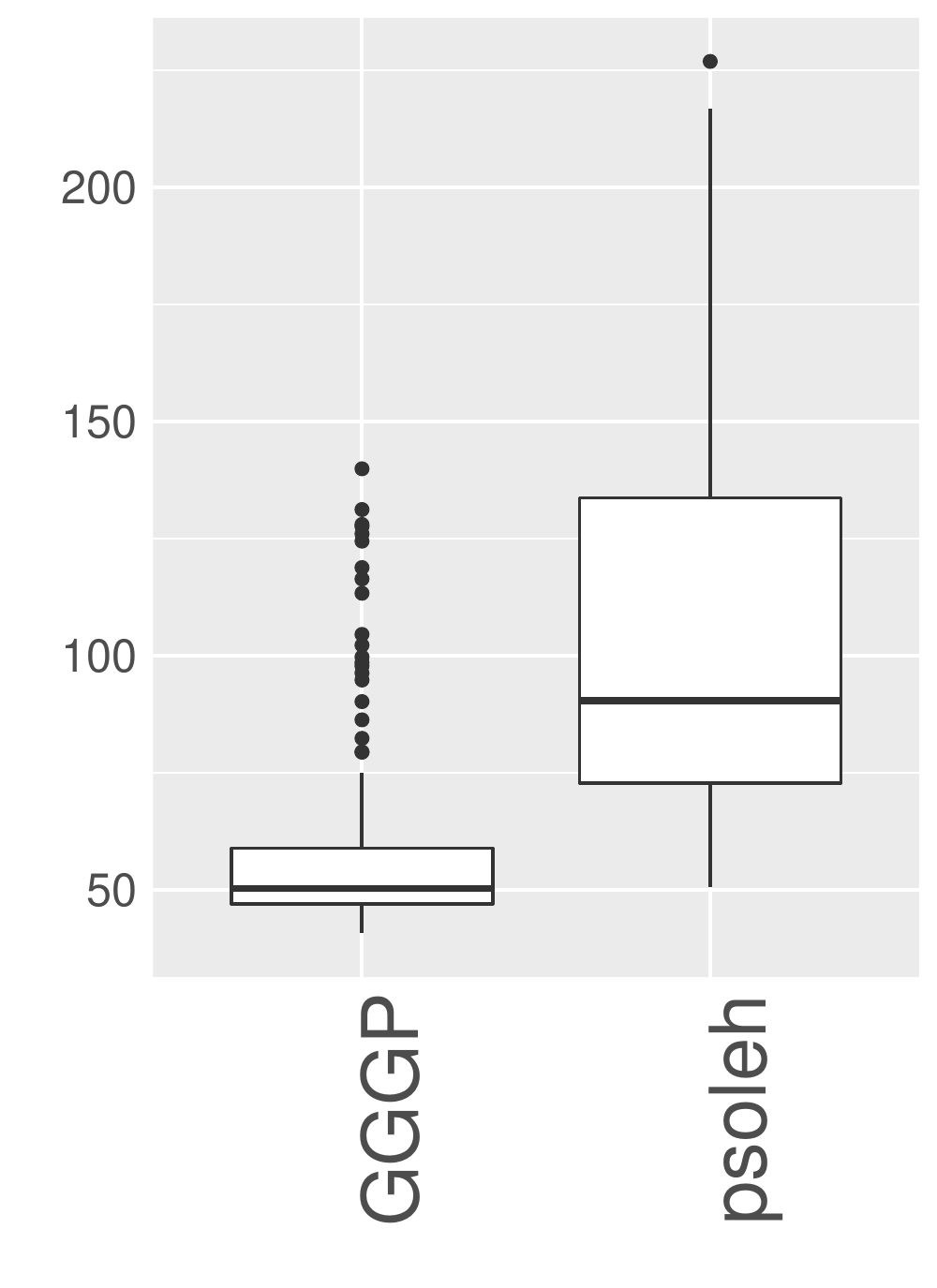}
		\vspace{-4mm} 
  	\caption{\scriptsize{Rosenbrock}} \label{fig:box3053a}
	\end{subfigure}
	\vspace{2mm} 
	\caption{30D best general box plots. 200 sample runs with a budget of 2,000 function evaluations. The comparators are taken from \cite{HughesGoerigkDokka2020a}, where the budget was 5,000 evaluations.}
	\label{fig:box30all}
	
\end{figure}

\begin{figure}[H]
	\centering
	

	\begin{subfigure}[t]{.18\textwidth}
		\includegraphics[width=\textwidth]{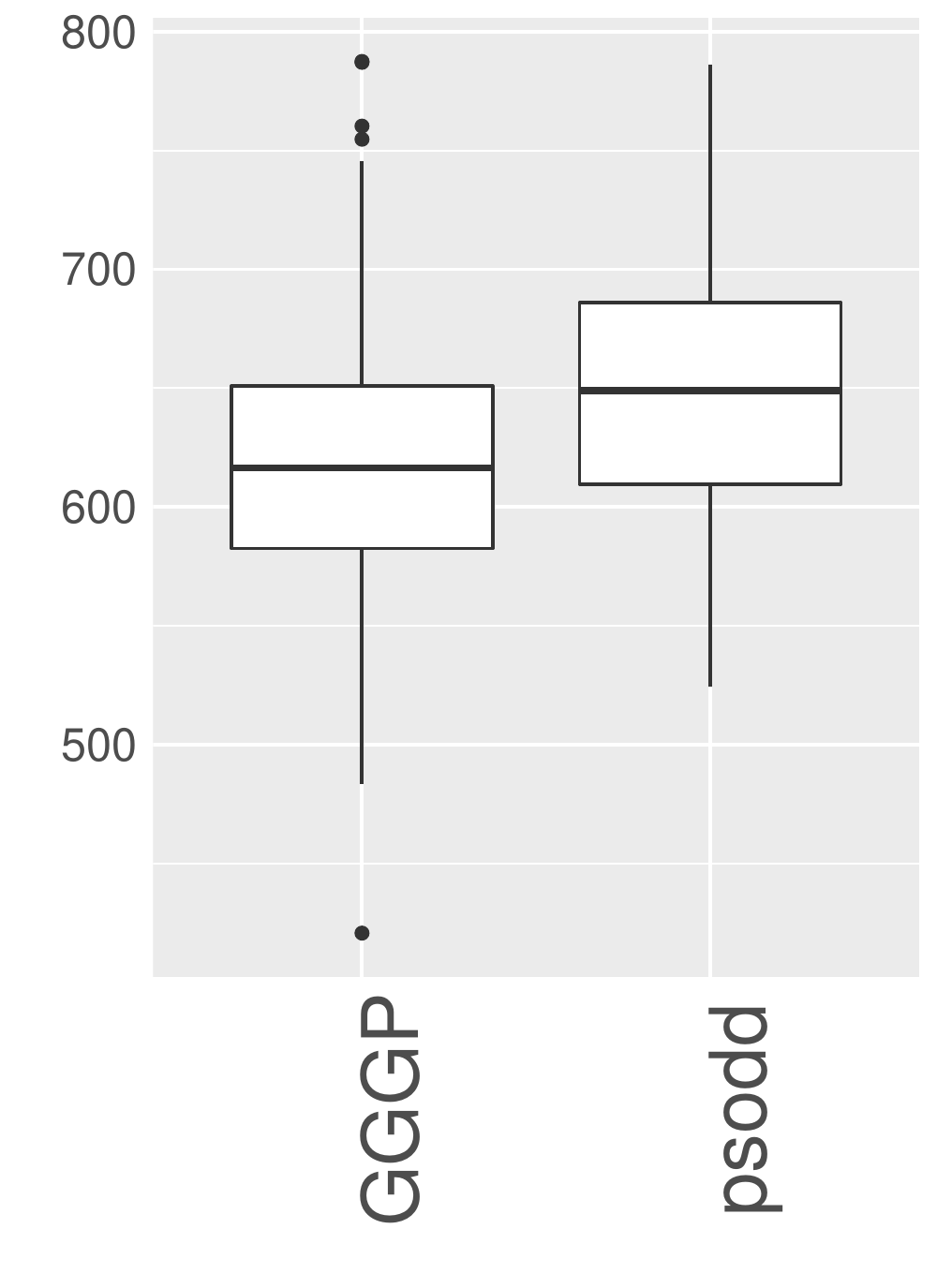}
		\vspace{-4mm} 
  	\caption{\scriptsize{Rastrigin}} \label{fig:box10051a}
	\end{subfigure}%
	\hspace{-1mm} 
	\begin{subfigure}[t]{.18\textwidth}
		\includegraphics[width=\textwidth]{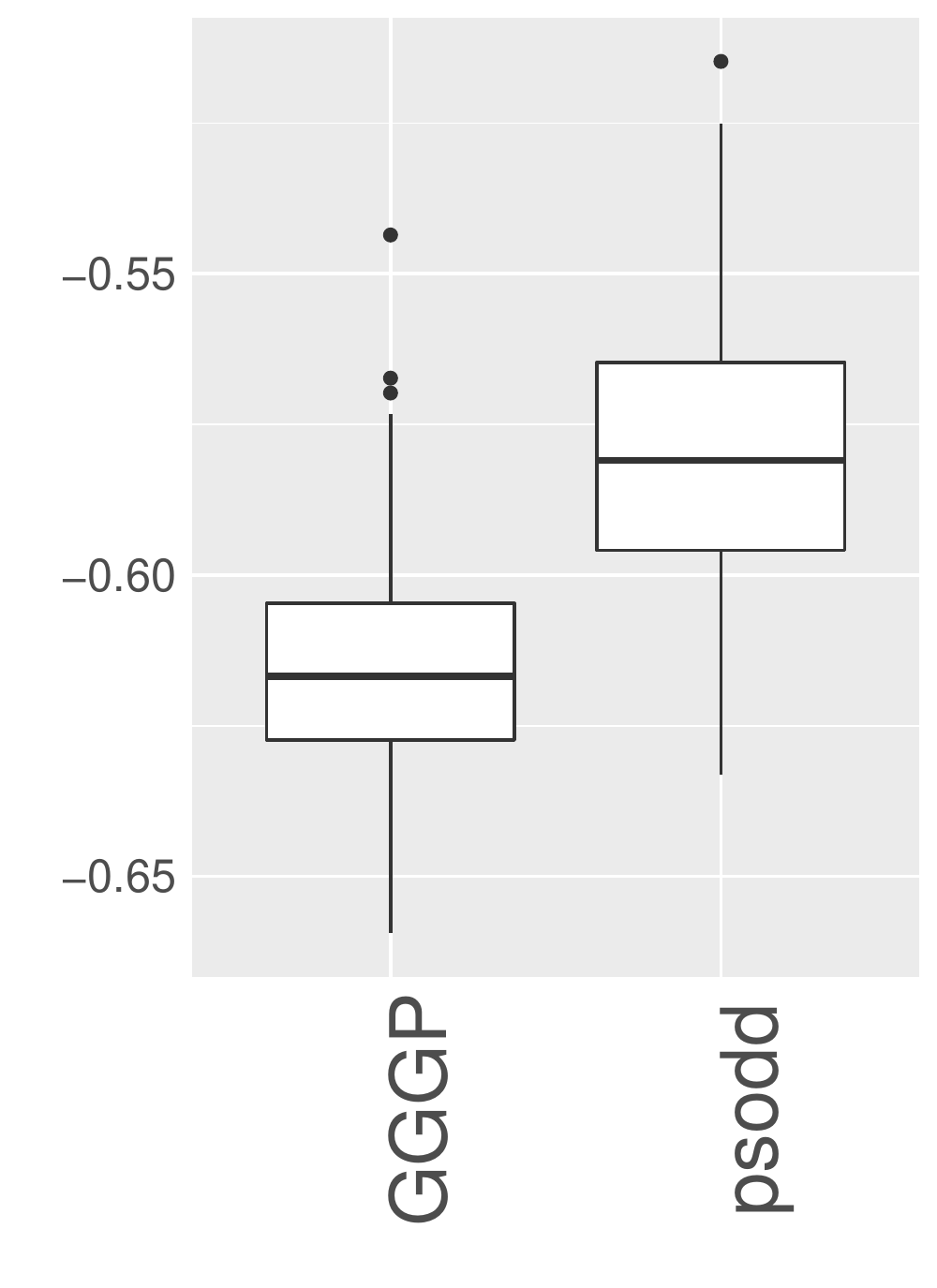}
		\vspace{-4mm} 
  	\caption{\scriptsize{Multipeak F1}} \label{fig:box10057a}
	\end{subfigure}%
	\hspace{-1mm} 
	\begin{subfigure}[t]{.18\textwidth}
		\includegraphics[width=\textwidth]{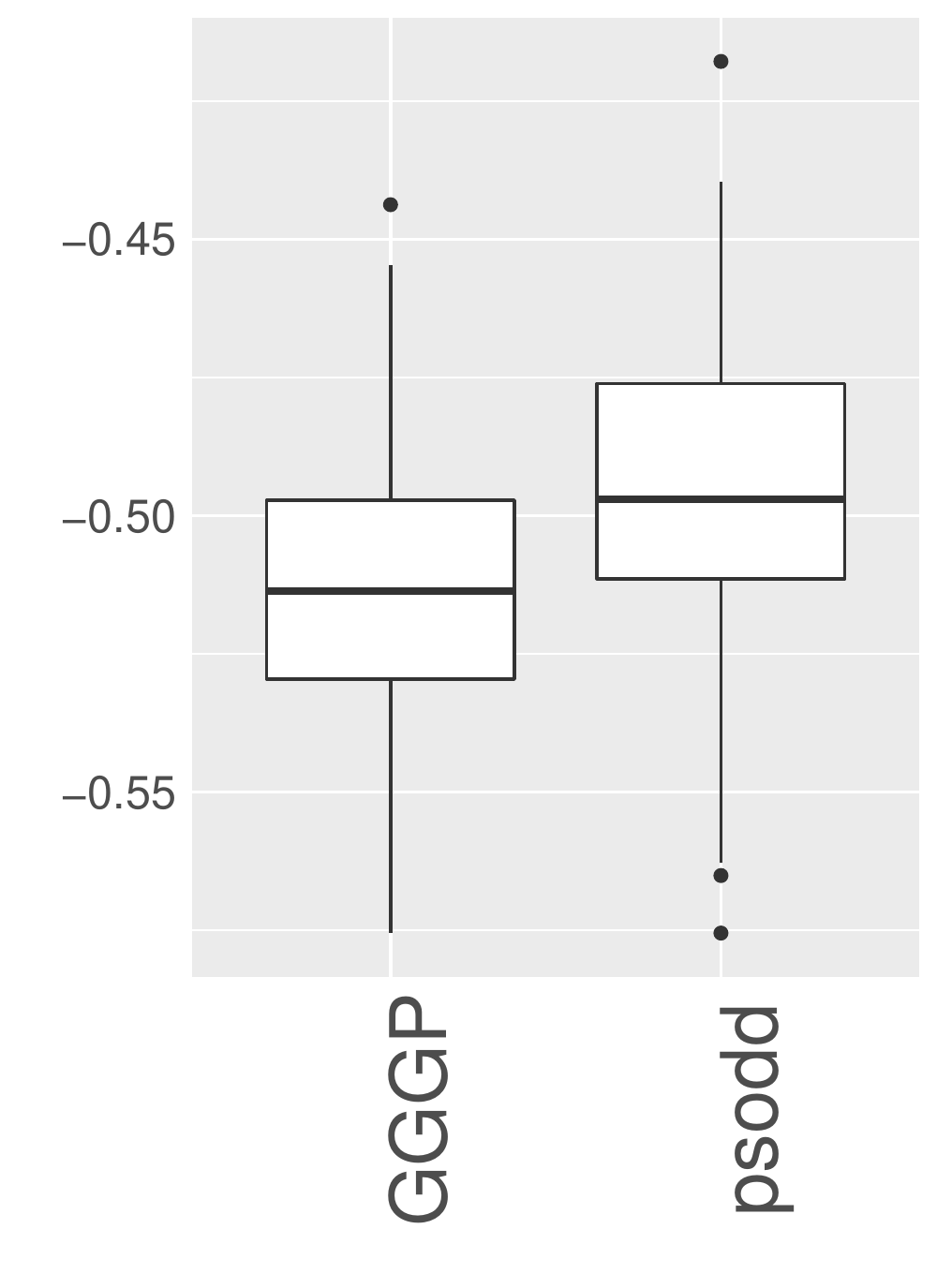}
		\vspace{-4mm} 
  	\caption{\scriptsize{Multipeak F2}} \label{fig:box10058a}
	\end{subfigure}%
	\hspace{-1mm} 
	\begin{subfigure}[t]{.18\textwidth}
		\includegraphics[width=\textwidth]{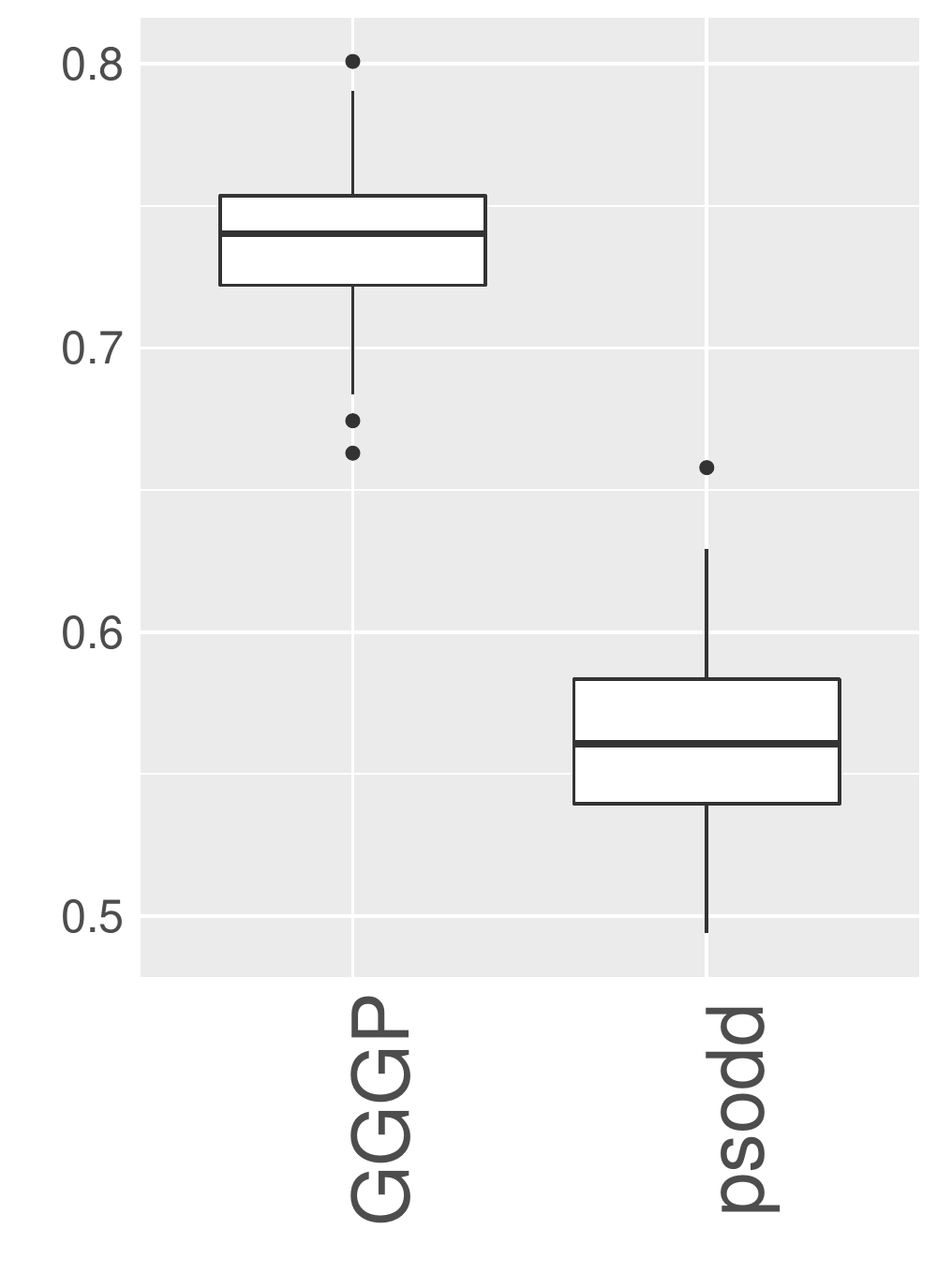}
		\vspace{-4mm} 
  	\caption{\scriptsize{Brankes}} \label{fig:box10067a}
	\end{subfigure}%
	\hspace{-1mm} 
	\begin{subfigure}[t]{.18\textwidth}
		\includegraphics[width=\textwidth]{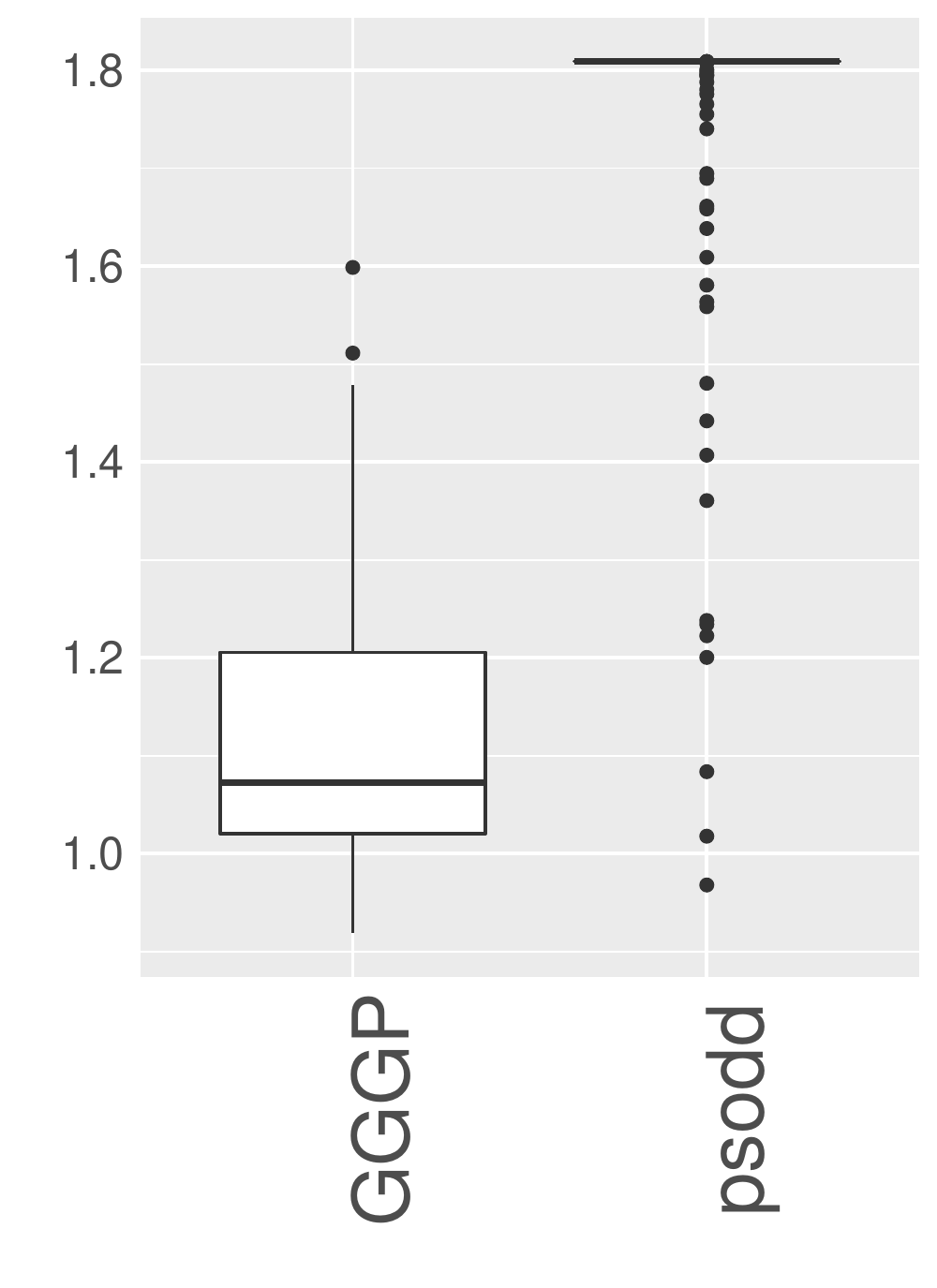}
		\vspace{-4mm} 
  	\caption{\scriptsize{Pickelhaube}} \label{fig:box10065a}
	\end{subfigure}
	
	\vspace{4mm} 
		
	\begin{subfigure}[t]{.18\textwidth}
		\includegraphics[width=\textwidth]{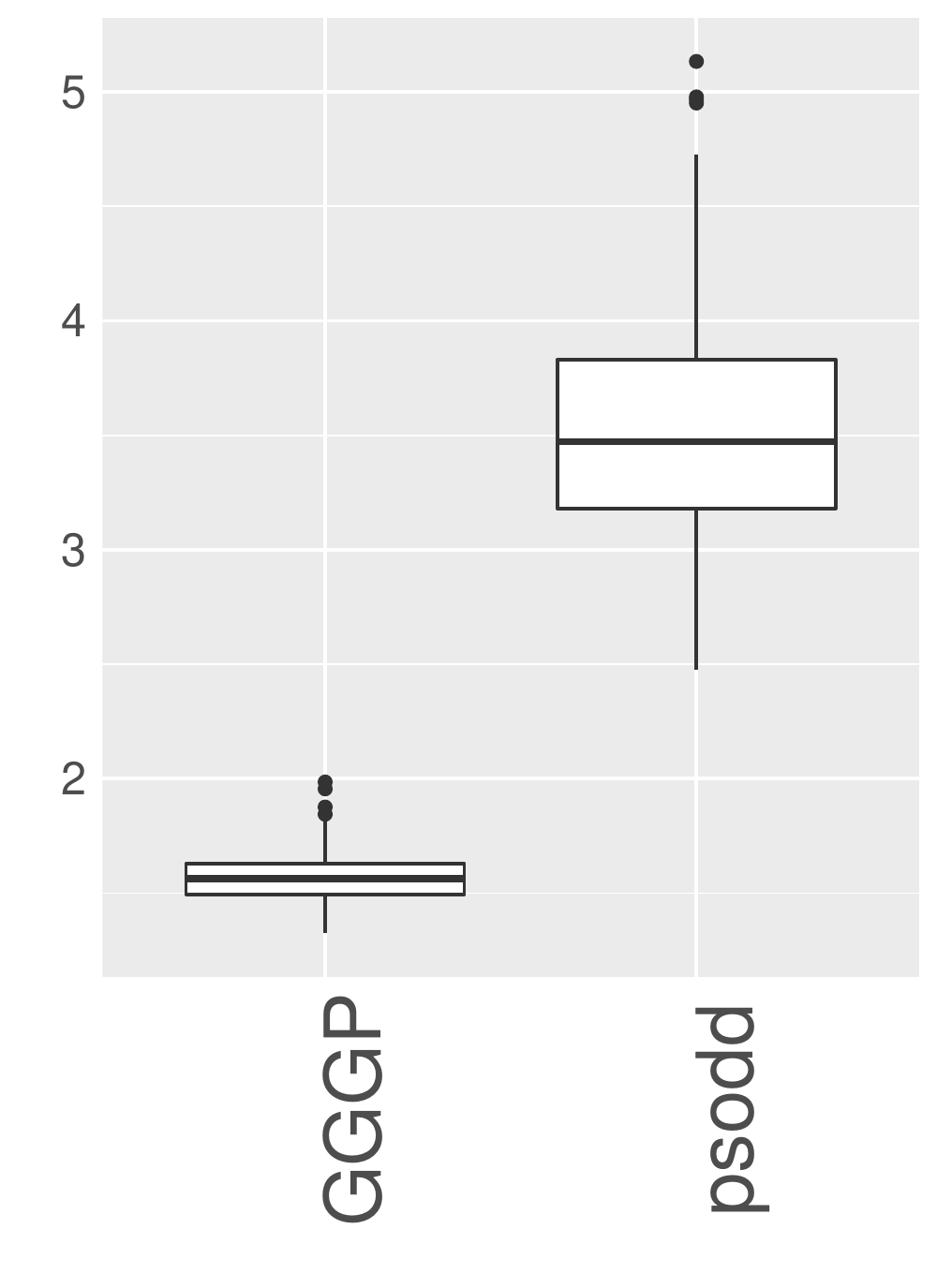}
		\vspace{-4mm} 
  	\caption{\scriptsize{Heaviside}} \label{fig:box10066a}
	\end{subfigure}%
	\hspace{-1mm} 
	\begin{subfigure}[t]{.18\textwidth}
		\includegraphics[width=\textwidth]{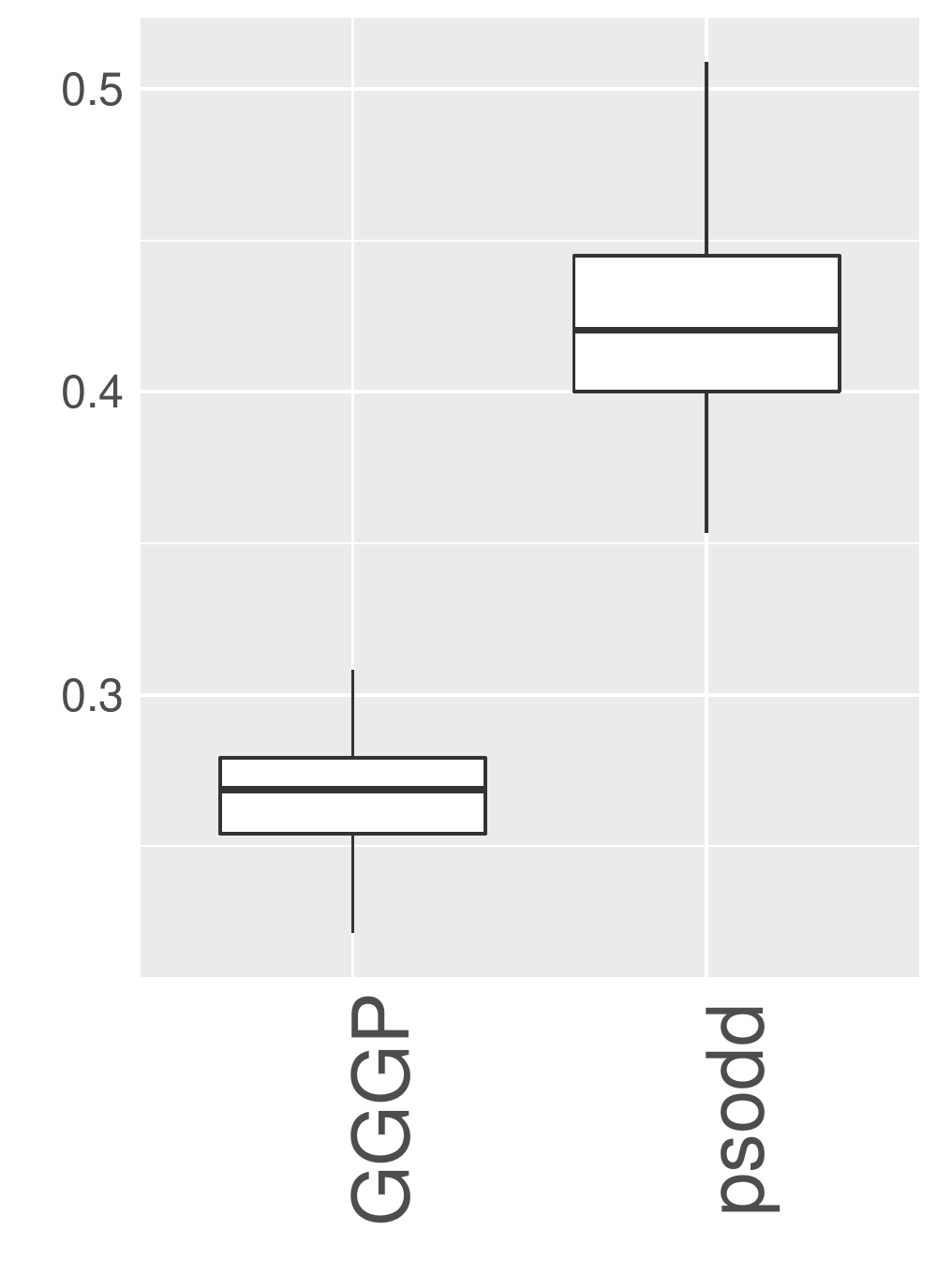}
		\vspace{-4mm} 
  	\caption{\scriptsize{Sawtooth}} \label{fig:box10055a}
	\end{subfigure}%
	\hspace{-1mm} 
	\begin{subfigure}[t]{.18\textwidth}
		\includegraphics[width=\textwidth]{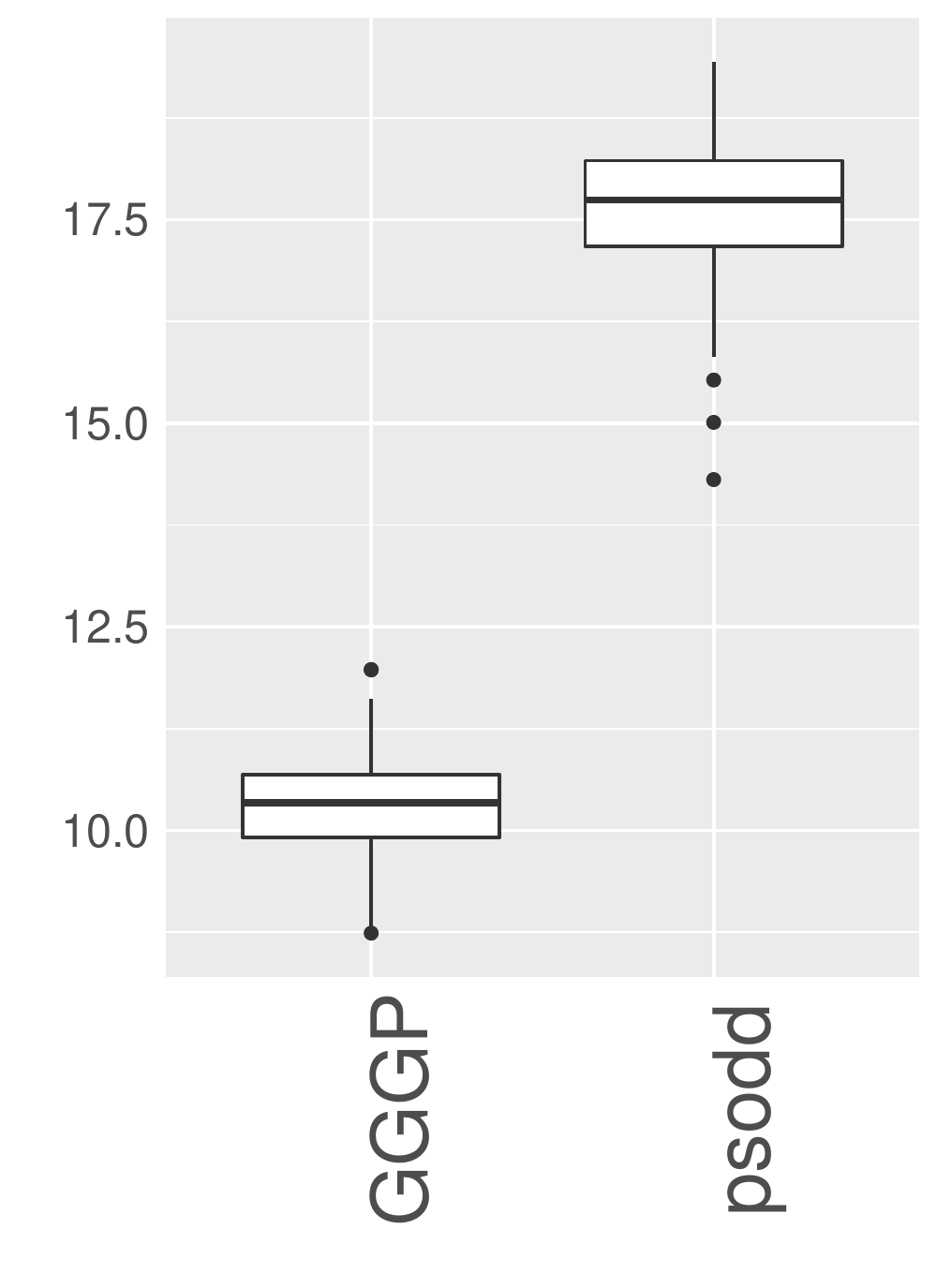}
		\vspace{-4mm} 
  	\caption{\scriptsize{Ackley}} \label{fig:box10052a}
	\end{subfigure}%
	\hspace{-1mm} 
	\begin{subfigure}[t]{.18\textwidth}
		\includegraphics[width=\textwidth]{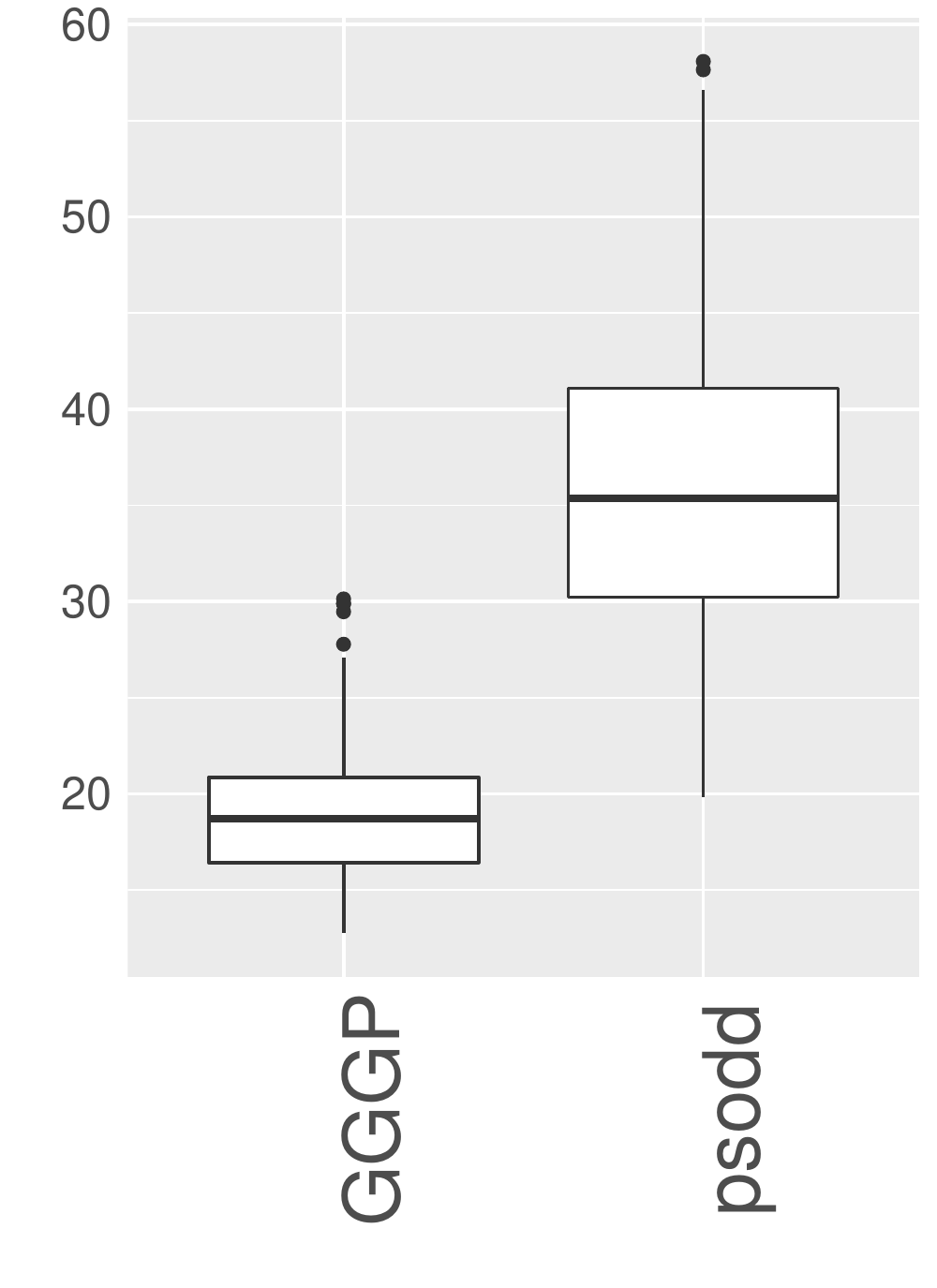}
		\vspace{-4mm} 
  	\caption{\scriptsize{Sphere}} \label{fig:box10054a}
	\end{subfigure}%
	\hspace{-1mm} 
	\begin{subfigure}[t]{.18\textwidth}
		\includegraphics[width=\textwidth]{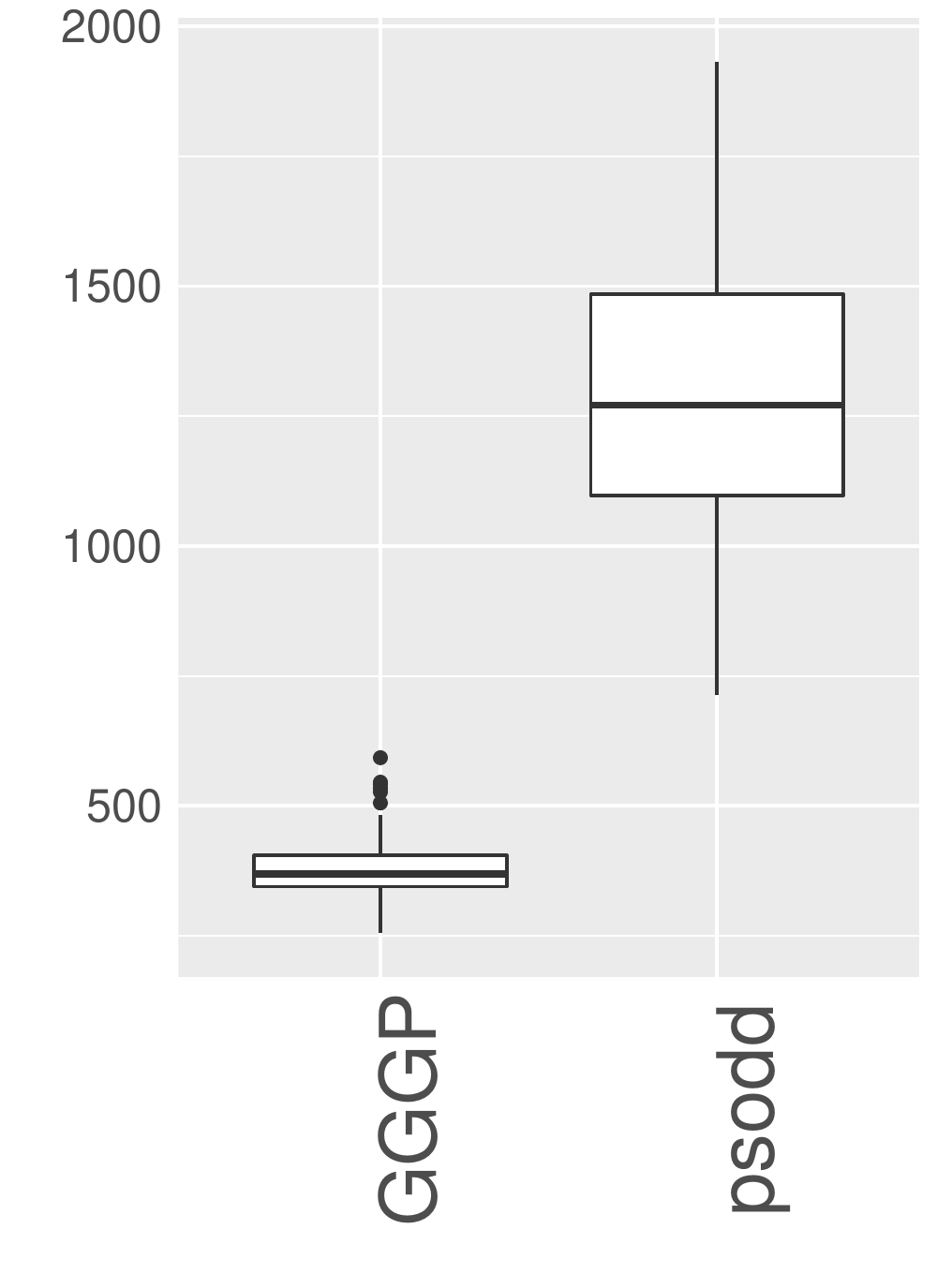}
		\vspace{-4mm} 
  	\caption{\scriptsize{Rosenbrock}} \label{fig:box10053a}
	\end{subfigure}
	\vspace{2mm} 
	\caption{100D best general box plots. 200 sample runs with a budget of 2,000 function evaluations. The comparators are taken from \cite{HughesGoerigkDokka2020a}, where the budget was 5,000 evaluations.}
	\label{fig:box100all}
	
\end{figure}

Mean estimates of the optimum robust values for the best performing heuristics, from the 200 sample runs and following the post-processing stage described in Section~\ref{sec:setUp}, are shown in Table~\ref{fig:means}. Corresponding box plots are shown in Figures~\ref{fig:box30ind} to~\ref{fig:box100all}. Individual case results are for the best heuristic evolved for a given test problem, then applied to that problem. General case results are for the best general case heuristic at 30D as applied to all 10 test problems, and the best at 100D applied to all 10 problems.

\begin{table}[htbp]
\small 
\begin{center}
\vspace{3mm} 
\begin{tabular}{c|cc|cc|cc|cc}
& \multicolumn{2}{c|}{Individual 30D} & \multicolumn{2}{c|}{Individual 100D} & \multicolumn{2}{c|}{General 30D} & \multicolumn{2}{c}{General 100D} \\
& GGGP & Comp & GGGP & Comp & GGGP & Comp & GGGP & Comp \\
\hline

Rastrigin & \textbf{154.93} & 226.57 & \textbf{416.35} & 648.67 & \textbf{216.80} & 226.57 & \textbf{615.40} & 648.67 \\
Multipeak F1 & \textbf{-0.64} & -0.63 & \textbf{-0.64} & -0.58 & -0.58 & \textbf{-0.63} & \textbf{-0.62} & -0.58 \\ 
Multipeak F2 & \textbf{-0.62} & -0.51 & \textbf{-0.54} & -0.49 & \textbf{-0.56} & -0.51 & \textbf{-0.51} & -0.49 \\ 
Branke's & 0.49 & \textbf{0.47} & \textbf{0.55} & 0.56 & 0.57 & \textbf{0.47} & 0.74 & \textbf{0.56} \\ 
Pickelhaube & \textbf{0.43} & \textbf{0.44} & \textbf{1.04} & 1.63 & \textbf{0.42} & 0.44 & \textbf{1.11} & 1.77 \\ 
Heaviside & \textbf{1.03} & 1.03 & \textbf{1.36} & 3.32 & \textbf{1.03} & 1.06 & \textbf{1.57} & 3.55 \\ 
Sawtooth & \textbf{0.30} & 0.35 & \textbf{0.25} & 0.33 & \textbf{0.34} & 0.35 & \textbf{0.27} & 0.42 \\ 
Ackley & \textbf{5.71} & 6.78 & \textbf{10.33} & 14.46 & \textbf{5.91} & 6.78 & \textbf{10.30} & 17.65 \\ 
Sphere & \textbf{1.68} & 2.86 & \textbf{15.60} & 36.11 & \textbf{2.62} & 5.30 & \textbf{19.05} & 36.11 \\ 
Rosenbrock & \textbf{55.23} & 89.24 & \textbf{311.96} & 1288.00 & \textbf{57.38} & 104.00 & \textbf{375.01} & 1288.00

\end{tabular}
\caption{Mean estimates of the optimum robust values for the best performing heuristics, due to 200 sample runs and using a budget of 2,000 functions evaluations. Comparators are taken from \cite{HughesGoerigkDokka2020a} and use a budget of 5,000 functions evaluations. Best results are shown in bold.}
\label{fig:means}

\end{center}
\end{table}

Comparator results taken from \cite{HughesGoerigkDokka2020a} are also shown. There several heuristics were analysed, with each parameter-fitted to 4 of the 10 test problems used here, separately for 30D and 100D. The budget was 5,000 function evaluations. For the individual cases the comparator results shown are due the best performing specific heuristic for an individual problem. For the general case the comparator results are for the best performing heuristic overall in \cite{HughesGoerigkDokka2020a}, as applied to all 10 problems. 

As the individual case comparators were not tuned on specific problems, comparisons with our individual results should be considered indicative. For the general case a direct comparison is reasonable. Comparisons should be interpreted in the context of the use of a budget of 2,000 function evaluations here. Labels on the box plots, Figures~\ref{fig:box30ind} to~\ref{fig:box100all}, give the specific comparator heuristic responsible for each set of results.

In Table~\ref{fig:means} values in bold indicate results which are best or statistically equivalent to the best, based on Wilcoxon rank-sum tests with 95\% confidence. 

At 30D in 8 of the individual cases our GP analysis produces the best heuristic, with one worse than the comparator and one statistically equivalent. For the general case the GP again produces the best heuristic for 8 problems, with the comparator best for 2. In view of our much reduced budget this shows a significantly improved performance for the general case. The individual case comparisons also indicate a good performance.

At 100D the performance of the new heuristics is even better. For all of the individual cases the GP produces the best results. In a number of instances we see substantial improvements, which is encouraging. For the general case the new heuristic is best for 9 problems and worse for one. Again for several problems the new results show significant improvements. In view of the reduced budget this is a strong performance.

\vspace{4mm} 

\vspace{4mm} 

\subsection{Component analysis}
\label{sec:analysis}

\subsubsection{Analysis of best performing heuristics}
\label{sec:analysisBests}

The component breakdown for the best heuristic generated in each of the 22 GP runs are shown in Tables~\ref{fig:best30D} and~\ref{fig:best100D}. The results discussed in Section~\ref{sec:results} are generated by these heuristics. We first consider this snapshot of the components associated with the very best performing heuristics, and then move on to consider the component breakdowns across all heuristics from all GP runs, against heuristic performance.

From Tables~\ref{fig:best30D} and~\ref{fig:best100D} it can be observed that an inner maximisation using random sampling is much preferred, with only 2 heuristics employing an alternative, PSO. Furthermore only a small number of points are typically sampled, with 14 heuristics using 5 points or less in the inner maximisation. In all cases a particle level stopping condition is used.

For the movement formulations, a baseline Inertia velocity formulation is preferred by 20 heuristics. An extended capability is used by all heuristics, with 16 using the full +DD+LEH capability. In all heuristics where dormancy and relocation are used (including +LEH), relocation using the largest empty hypersphere is selected over random relocation. Where a descent direction vector is used (including +DD), a unit vector form of $\pmb{r}_3$ is employed in all but 3 heuristics, rather than a randomised vector. 

The best heuristics typically employ small swarm sizes, with 17 using less than 10 particles. Of the 4 network topologies appearing in Tables~\ref{fig:best30D} and~\ref{fig:best100D}, 14 heuristics use Global, with Hierarchical, von Neumann and Ring also represented. Where dormancy is relevant the use of existing information to inform it is preferred in 14 heuristics. The use of existing information to update a robust value on completion of an inner search is preferred 19 times. Some form of PSO level mutation is employed 12 times, 7 of which are by sampling from a Gaussian distribution. 

Beyond this narrow snapshot of the component breakdowns of the very best performing search algorithms, an assessment of the forms of component included across the large numbers of heuristics generated by our GP runs will give some indication of how each alternative impacts heuristic performance. The alternative forms that a component may take are given by the grammar in Figure~\ref{fig:grammar}. For a given component the levels of representation of each alternative form across all heuristics generated in the GP runs, is driven by evolutionary processes and so will indicate some preference. At a component level Table~\ref{fig:proportions} gives the proportions of each alternative form separately for 30D and 100D, from all heuristics generated here. For a given dimension the heuristics due to all individual cases and the general case are taken in total. 

In Table~\ref{fig:proportions} results for the top third best performing heuristics are also shown. This is a high level indication of the impact of each alternative form on heuristic performance. In Figures~\ref{fig:inType30Comp} to~\ref{fig:mutationComp} we expand on this information for selected components. Each plot relates to all heuristics generated in a single GP run. By arranging all heuristics generated in a run in order of best-to-worst fitness, an assessment of the heuristic component breakdowns with fitness is made, using a decile scale of heuristic performance (x axis). Within each decile of the fitness-sorted heuristics, for each component, the proportion of each alternative choice for that component is calculated. For a given component and a given test problem case, each line on the plot represents one of the alternatives available for that component. A line is generated from the 10 decile point values, representing the proportion of heuristics within that decile which employ that alternative (y-axis). Within any given decile the plotted values add to 100\% as they refer to the proportion of heuristics within that decile.

For example consider Figure~\ref{fig:inType30Comp} comprising 10 individual test problem plots at 30D, for the form of inner maximisation. There are 3 alternatives for the inner maximisation: random sampling (red), particle swarm optimisation (green), or genetic algorithm (blue). For the Rastrigin problem, of all of the heuristics within the first decile, that is the top 10\% performing heuristics, 97\%	employ random sampling whilst 2\% and 1\% employ PSO and GA respectively. These are the values plotted at the 0.1 position on the x-axis. Within the next decile at the 0.2 position on the x-axis (the 10\% -- 20\% range of best performing heuristics), the value are 75\%, 7\% and 18\% for random sampling, PSO and GA respectively.

These plots indicate how component breakdowns relate to heuristic performance. In addition the relative areas under the lines indicate the total proportion of each component category across all heuristics in a single GP run. For example in the Rastrigin plot in Figure~\ref{fig:inType30Comp} the proportions of each alternative for the form of inner maximisation is 48\%, 24\% and 28\% for random sampling, PSO and GA respectively. In Sections~\ref{sec:compInner} to~\ref{sec:compR3} we consider the main components individually.

\newpage
\begin{landscape} %
	
\begin{table}[htbp]
\footnotesize
\begin{center}

\begin{tabular}{c|cccccccccccc}
 & Group & inner & Movement & Network & inExt & stop & nDorm & nPBest & Mutation & Baseline & Relocate & r3 \\
\hline

 Rastrigin & 2 & Random & rPSO+DD+LEH & Global & 48 & Yes & Yes & Yes & None & Inertia & LEH & Unity \\

 Multipeak F1 & 21 & Random & rPSO+DD+LEH & Global & 4 & Yes & No & Yes & Uniform & Inertia & LEH & Unity \\

 Multipeak F2 & 8 & Random & rPSO+DD+LEH & Hierarchical & 5 & Yes & Yes & Yes & Uniform & Inertia & LEH & Unity \\

 Branke's & 9 & Random & rPSO+DD & Hierarchical & 4 & Yes & na & Yes & Uniform & Inertia & na & Random \\

 Pickelhaube & 3 & Random & rPSO+DD+LEH & Global & 4 & Yes & Yes & Yes & None & Constriction & LEH & Unity \\

 Heaviside & 3 & Random & rPSO+DD & Global & 7 & Yes & na & Yes & None & Inertia & na & Unity \\

 Sawtooth & 10 & Random & rPSO+LEH & Global & 6 & Yes & No & Yes & Gaussian & Inertia & LEH & na \\

 Ackley & 3 & Random & rPSO+DD+LEH & Global & 4 & Yes & Yes & Yes & None & Inertia & LEH & Unity \\

 Sphere & 2 & Random & rPSO+DD & Hierarchical & 15 & Yes & na & Yes & None & Inertia & na & Unity \\

 Rosenbrock & 9 & Random & rPSO+DD+LEH & Global & 6 & Yes & Yes & No & None & Constriction & LEH & Unity \\

 General & 4 & Random & rPSO+DD & Ring & 4 & Yes & na & Yes & Gaussian & Inertia & na & Unity 

\end{tabular}
\caption{30D components of best heuristics.}
\label{fig:best30D}
\end{center}
\end{table}

\vspace*{2.0mm}

\begin{table}[htbp]
\footnotesize
\begin{center}

\begin{tabular}{c|cccccccccccc}
 & Group & inner & Movement & Network & inExt & stop & nDorm & nPBest & Mutation & Baseline & Relocate & r3 \\
\hline

 Rastrigin & 4 & Random & rPSO+DD+LEH & Global & 9 & Yes & Yes & Yes & Gaussian & Constriction & LEH & Unity \\

 Multipeak F1 & 8 & Random & rPSO+DD+LEH & Global & 4 & Yes & No & Yes & Uniform & Inertia & LEH & Unity \\

 Multipeak F2 & 8 & Random & rPSO+DD+LEH & von Neumann & 5 & Yes & Yes & Yes & Gaussian & Inertia & LEH & Unity \\

 Branke's & 15 & PSO & rPSO+LEH & Hierarchical & 4 & Yes & No & No & Uniform & Inertia & LEH & na \\

 Pickelhaube & 9 & Random & rPSO+DD+LEH & Hierarchical & 4 & Yes & Yes & Yes & None & Inertia & LEH & Unity \\

 Heaviside & 8 & Random & rPSO+DD+LEH & Global & 4 & Yes & Yes & No & Gaussian & Inertia & LEH & Unity \\

 Sawtooth & 21 & Random & rPSO+DD+LEH & von Neumann & 4 & Yes & Yes & Yes & Gaussian & Inertia & LEH & Unity \\

 Ackley & 11 & PSO & rPSO+DD+LEH & Global & 4 & Yes & Yes & Yes & None & Inertia & LEH & Random \\

 Sphere & 5 & Random & rPSO+DD+LEH & Global & 8 & Yes & Yes & Yes & None & Inertia & LEH & Unity \\

 Rosenbrock & 5 & Random & rPSO+DD+LEH & Global & 8 & Yes & Yes & Yes & Gaussian & Inertia & LEH & Random \\

 General & 8 & Random & rPSO+DD+LEH & Global & 4 & Yes & Yes & Yes & None & Inertia & LEH & Unity 

\end{tabular}
\caption{100D components of best heuristics.}
\label{fig:best100D}
\end{center}
\end{table}

\vspace*{2.0mm}

\end{landscape} %

\newpage

\begin{table}[H]%
\footnotesize 
\begin{center}

\begin{tabular}{cc|cccc}	
										
\textbf{Component}	&		&	\textbf{30D all}	&	\textbf{30D top}	&	\textbf{100D all}	&	\textbf{100D top}	\\
\hline		
\textbf{Form of inner search}	&	\textbf{Random}	&	46.6\%	&	67.5\%	&	42.5\%	&	55.3\%	\\
	&	\textbf{PSO}	&	25.2\%	&	12.8\%	&	29.8\%	&	25.5\%	\\
	&	\textbf{GA}	&	28.1\%	&	19.8\%	&	27.7\%	&	19.2\%	\\
\hline											
\textbf{Extent of inner search}	&	\textbf{[2-10]}	&	49.9\%	&	83.0\%	&	50.3\%	&	84.9\%	\\
	&	\textbf{[11-20]}	&	21.1\%	&	12.8\%	&	20.3\%	&	10.3\%	\\
	&	\textbf{[21-30]}	&	11.2\%	&	2.4\%	&	11.5\%	&	3.5\%	\\
	&	\textbf{[31-40]}	&	5.3\%	&	0.6\%	&	5.7\%	&	0.8\%	\\
	&	\textbf{\textgreater 40}	&	12.5\%	&	1.3\%	&	12.2\%	&	0.6\%	\\
\hline											
\textbf{Form of baseline rPSO formula}	&	\textbf{Constriction}	&	36.2\%	&	17.0\%	&	35.8\%	&	19.5\%	\\
	&	\textbf{Inertia}	&	63.8\%	&	83.0\%	&	64.2\%	&	80.5\%	\\
\hline											
\textbf{Form of movement}	&	\textbf{rPSO}	&	12.1\%	&	3.1\%	&	11.4\%	&	1.6\%	\\
	&	\textbf{+DD}	&	25.0\%	&	25.2\%	&	15.0\%	&	7.8\%	\\
	&	\textbf{+LEH}	&	25.8\%	&	25.0\%	&	32.2\%	&	36.3\%	\\
	&	\textbf{+DD+LEH}	&	37.1\%	&	46.8\%	&	41.4\%	&	54.3\%	\\
\hline											
\textbf{Form of network}	&	\textbf{Global}	&	24.5\%	&	41.2\%	&	25.6\%	&	43.0\%	\\
	&	\textbf{Focal}	&	9.8\%	&	3.1\%	&	9.9\%	&	3.2\%	\\
	&	\textbf{Ring (size=2)}	&	11.3\%	&	6.9\%	&	11.1\%	&	6.4\%	\\
	&	\textbf{von Neumann}	&	13.6\%	&	12.8\%	&	14.3\%	&	12.9\%	\\
	&	\textbf{Clan}	&	13.2\%	&	10.9\%	&	13.6\%	&	13.5\%	\\
	&	\textbf{Cluster}	&	12.1\%	&	8.8\%	&	11.5\%	&	7.9\%	\\
	&	\textbf{Hierarchy}	&	15.5\%	&	16.5\%	&	14.1\%	&	13.1\%	\\
\hline											
\textbf{Group (swarm) size}	&	\textbf{[2-10]}	&	25.6\%	&	43.2\%	&	22.4\%	&	33.9\%	\\
	&	\textbf{[11-20]}	&	19.4\%	&	17.0\%	&	21.3\%	&	22.1\%	\\
	&	\textbf{[21-30]}	&	20.2\%	&	18.0\%	&	20.3\%	&	18.1\%	\\
	&	\textbf{[31-40]}	&	17.6\%	&	11.9\%	&	18.2\%	&	13.8\%	\\
	&	\textbf{\textgreater 40}	&	17.2\%	&	10.0\%	&	17.8\%	&	12.1\%	\\
\hline											
\textbf{Inclusion of stopping condition}	&	\textbf{No}	&	36.8\%	&	16.2\%	&	40.6\%	&	26.9\%	\\
	&	\textbf{Yes}	&	63.2\%	&	83.8\%	&	59.4\%	&	73.1\%	\\
\hline											
\textbf{Use of existing info. for dormancy}	&	\textbf{No}	&	30.4\%	&	32.3\%	&	32.3\%	&	33.4\%	\\
	&	\textbf{Yes}	&	32.6\%	&	39.4\%	&	41.3\%	&	57.3\%	\\
	&	\textbf{Not applicable}	&	37.0\%	&	28.3\%	&	26.4\%	&	9.4\%	\\	
\hline											
\textbf{Use of existing info. for personal best}	&	\textbf{No}	&	46.3\%	&	39.9\%	&	44.6\%	&	36.5\%	\\
	&	\textbf{Yes}	&	53.7\%	&	60.1\%	&	55.4\%	&	63.5\%	\\
\hline											
\textbf{Form of mutation}	&	\textbf{None}	&	37.2\%	&	42.4\%	&	39.1\%	&	48.0\%	\\
	&	\textbf{Random}	&	31.0\%	&	27.2\%	&	31.9\%	&	29.9\%	\\
	&	\textbf{Gaussian}	&	31.8\%	&	30.4\%	&	29.0\%	&	22.1\%	\\
\hline											
\textbf{Form of relocation due to dormancy}	&	\textbf{LEH}	&	55.5\%	&	71.7\%	&	65.7\%	&	90.6\%	\\
	&	\textbf{Random}	&	7.5\%	&	0.0\%	&	7.9\%	&	0.0\%	\\
	&	\textbf{Not applicable}	&	37.0\%	&	28.3\%	&	26.4\%	&	9.4\%	\\
\hline											
\textbf{Form of r3 vector}	&	\textbf{Random}	&	27.3\%	&	27.9\%	&	29.6\%	&	32.3\%	\\
	&	\textbf{Unity}	&	34.8\%	&	44.0\%	&	26.8\%	&	29.8\%	\\
	&	\textbf{Not applicable}	&	37.9\%	&	28.1\%	&	43.6\%	&	37.9\%	\\
\hline													

\end{tabular}

\caption{Proportions of component make ups over all heuristics. Here 'top' refers to the top one third of heuristics when sorted best to worst.}
\label{fig:proportions}
\end{center}
\end{table}

\subsubsection{Form and extent of inner maximisation}
\label{sec:compInner}

\begin{figure}[H]
	\centering
	
	\vspace{-5mm} 

	\begin{subfigure}[t]{.2\textwidth}
		\includegraphics[width=\textwidth]{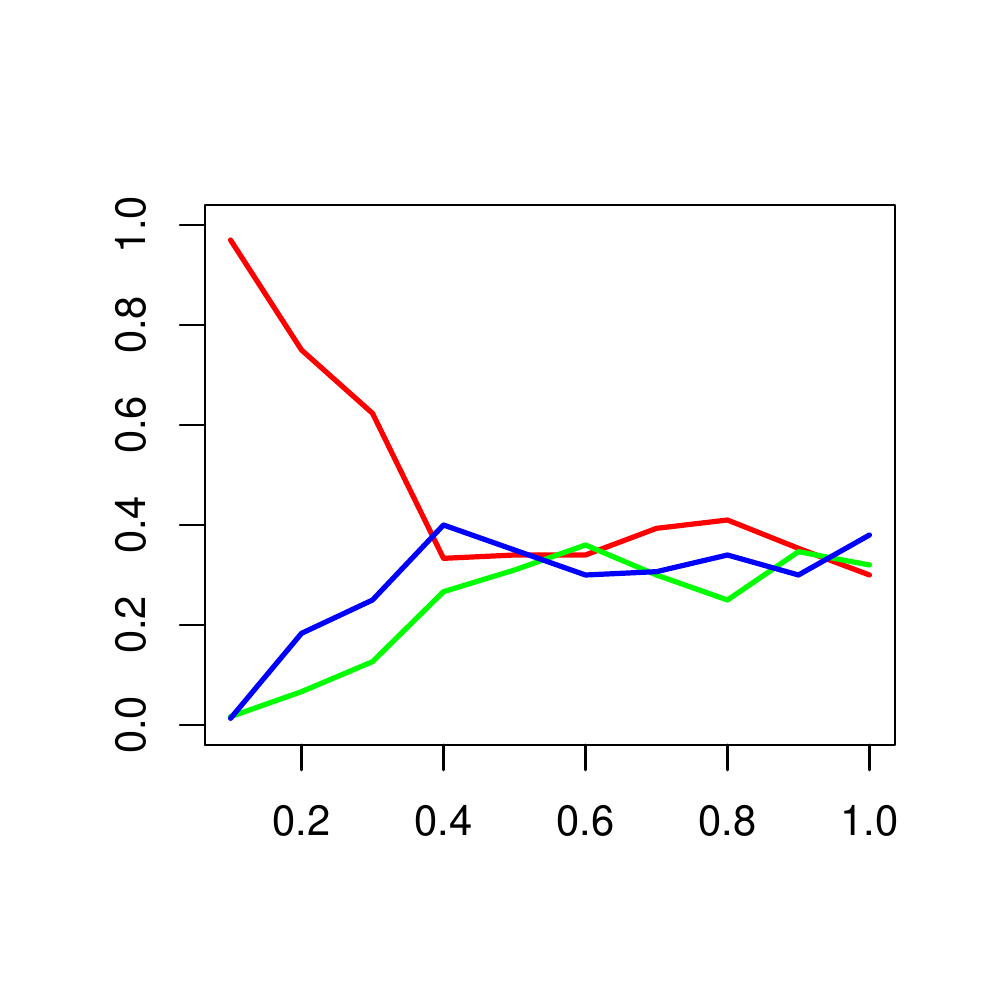}
		\vspace{-10mm} 
  	\caption{\scriptsize{Rastrigin}} \label{fig:iT30Ra}
	\end{subfigure}%
	\hspace{-6mm} 
	\begin{subfigure}[t]{.2\textwidth}
		\includegraphics[width=\textwidth]{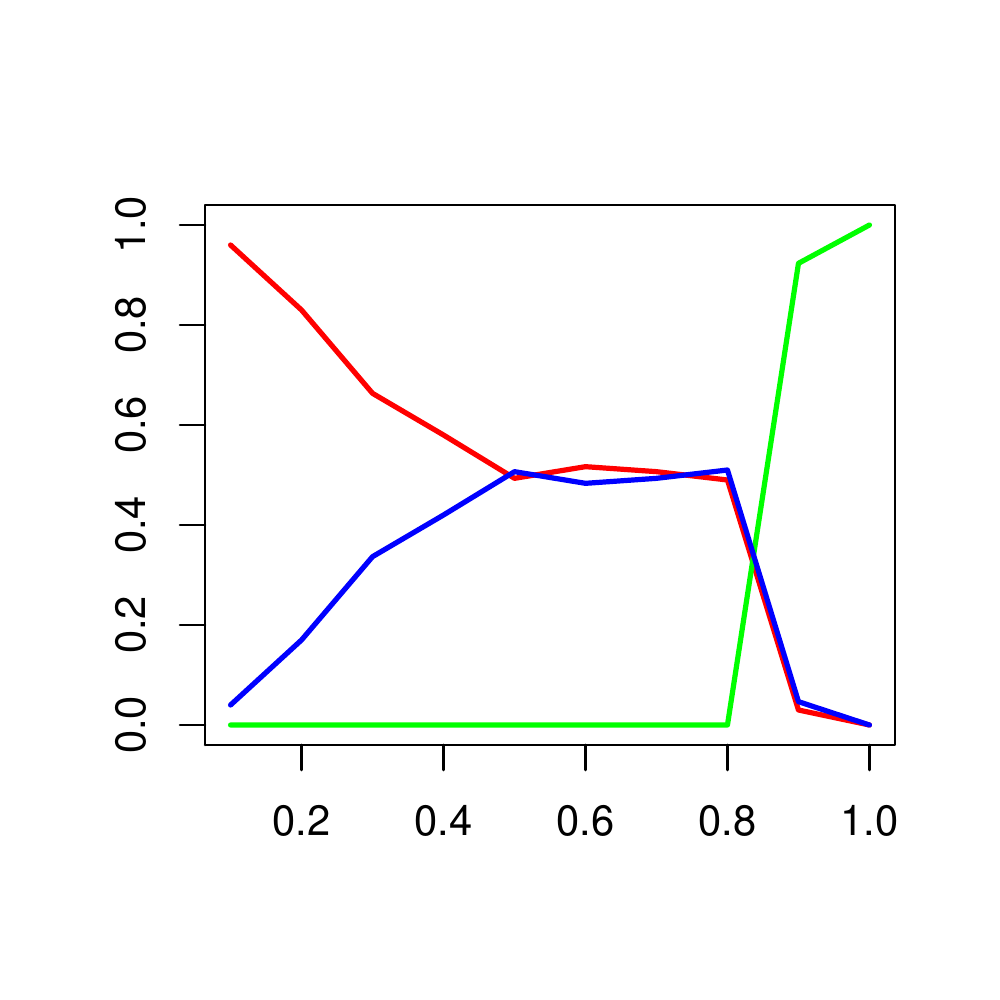}
		\vspace{-10mm} 
  	\caption{\scriptsize{Multipeak F1}} \label{fig:iT30M1}
	\end{subfigure}%
	\hspace{-6mm} 
	\begin{subfigure}[t]{.2\textwidth}
		\includegraphics[width=\textwidth]{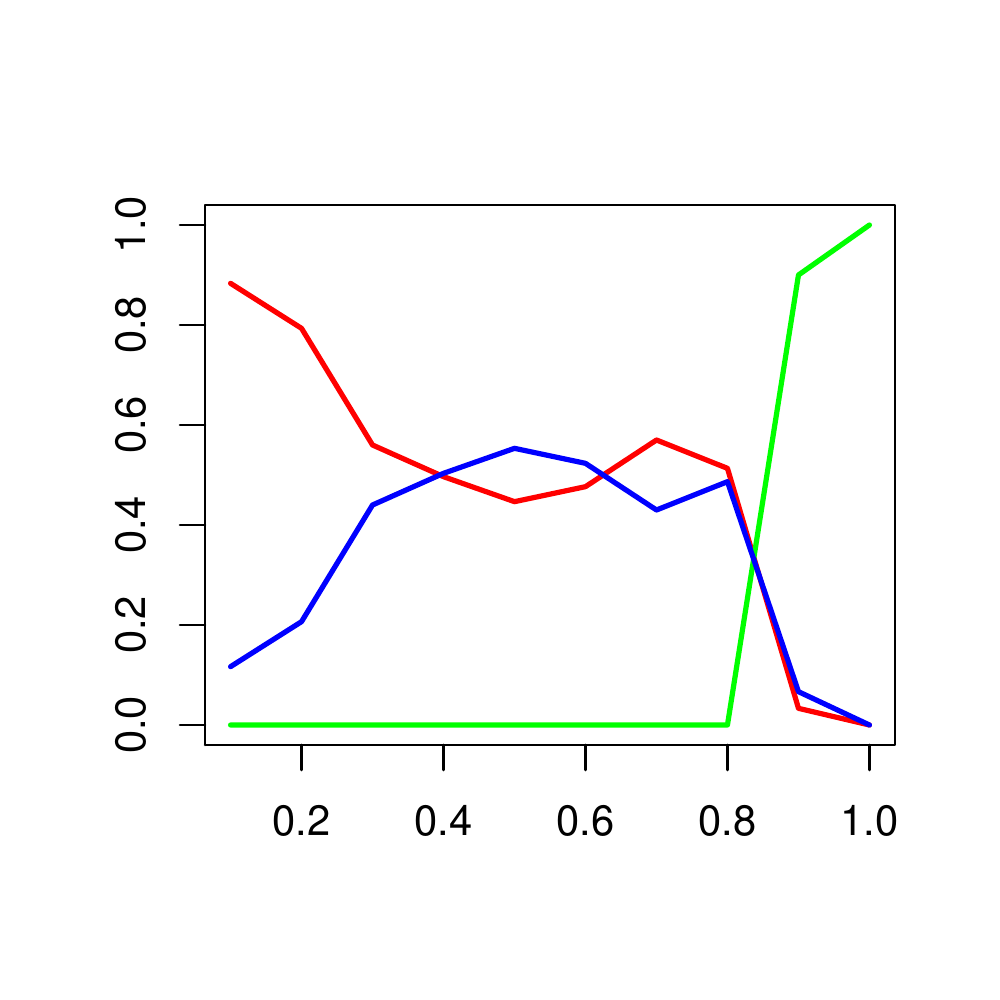}
		\vspace{-10mm} 
  	\caption{\scriptsize{Multipeak F2}} \label{fig:iT30M2}
	\end{subfigure}%
	\hspace{-6mm} 
	\begin{subfigure}[t]{.2\textwidth}
		\includegraphics[width=\textwidth]{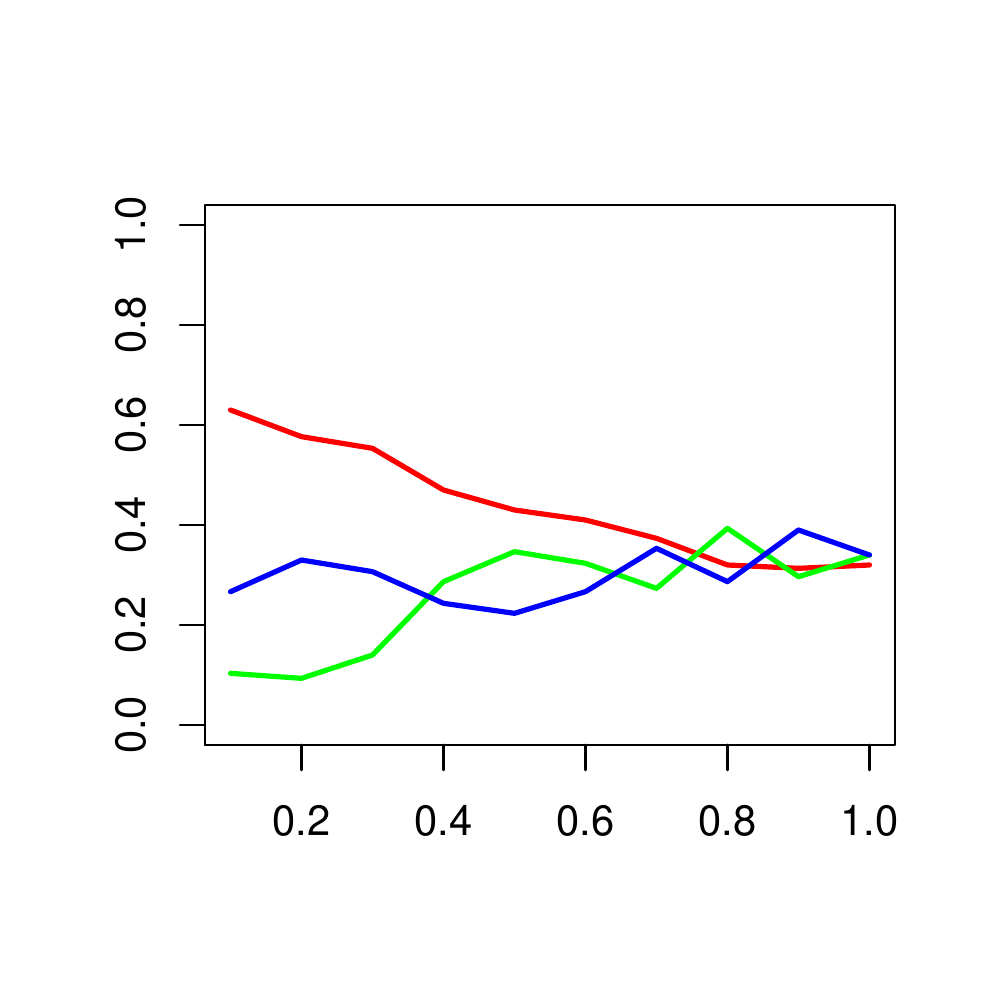}
		\vspace{-10mm} 
  	\caption{\scriptsize{Brankes}} \label{fig:iT30Br}
	\end{subfigure}%
	\hspace{-6mm} 
	\begin{subfigure}[t]{.2\textwidth}
		\includegraphics[width=\textwidth]{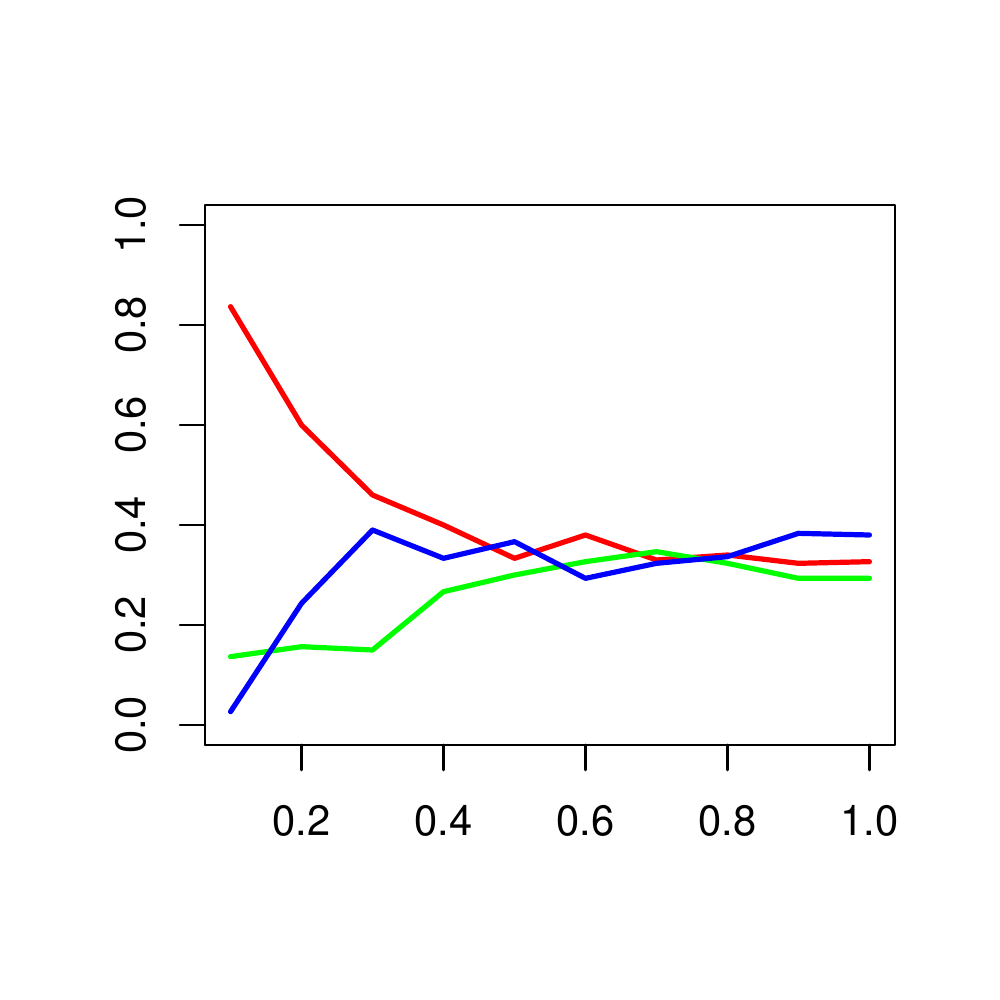}
		\vspace{-10mm} 
  	\caption{\scriptsize{Pickelhaube}} \label{fig:iT30Pi}
	\end{subfigure}
	
	\vspace{-2mm} 
		
	\begin{subfigure}[t]{.2\textwidth}
		\includegraphics[width=\textwidth]{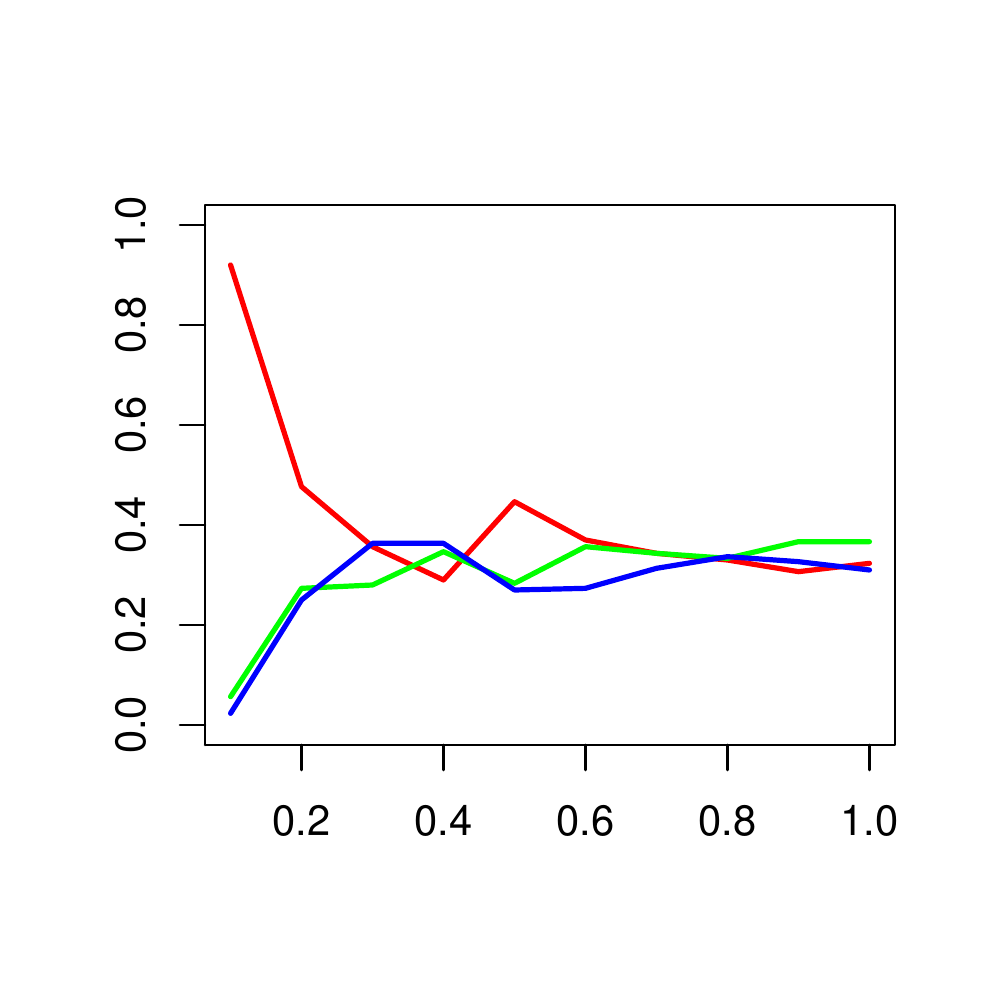}
		\vspace{-10mm} 
  	\caption{\scriptsize{Heaviside}} \label{fig:iT30Hv}
	\end{subfigure}%
	\hspace{-6mm} 
	\begin{subfigure}[t]{.2\textwidth}
		\includegraphics[width=\textwidth]{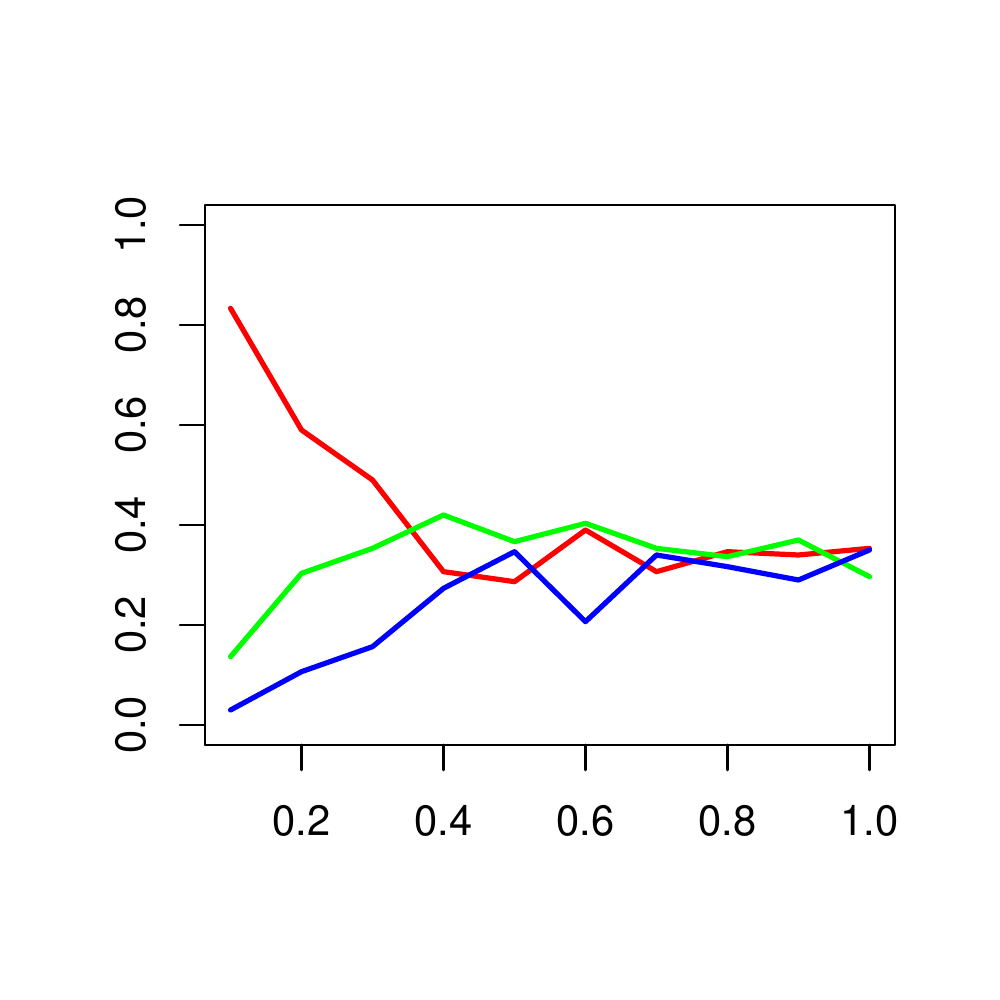}
		\vspace{-10mm} 
  	\caption{\scriptsize{Sawtooth}} \label{fig:iT30Sa}
	\end{subfigure}%
	\hspace{-6mm} 
	\begin{subfigure}[t]{.2\textwidth}
		\includegraphics[width=\textwidth]{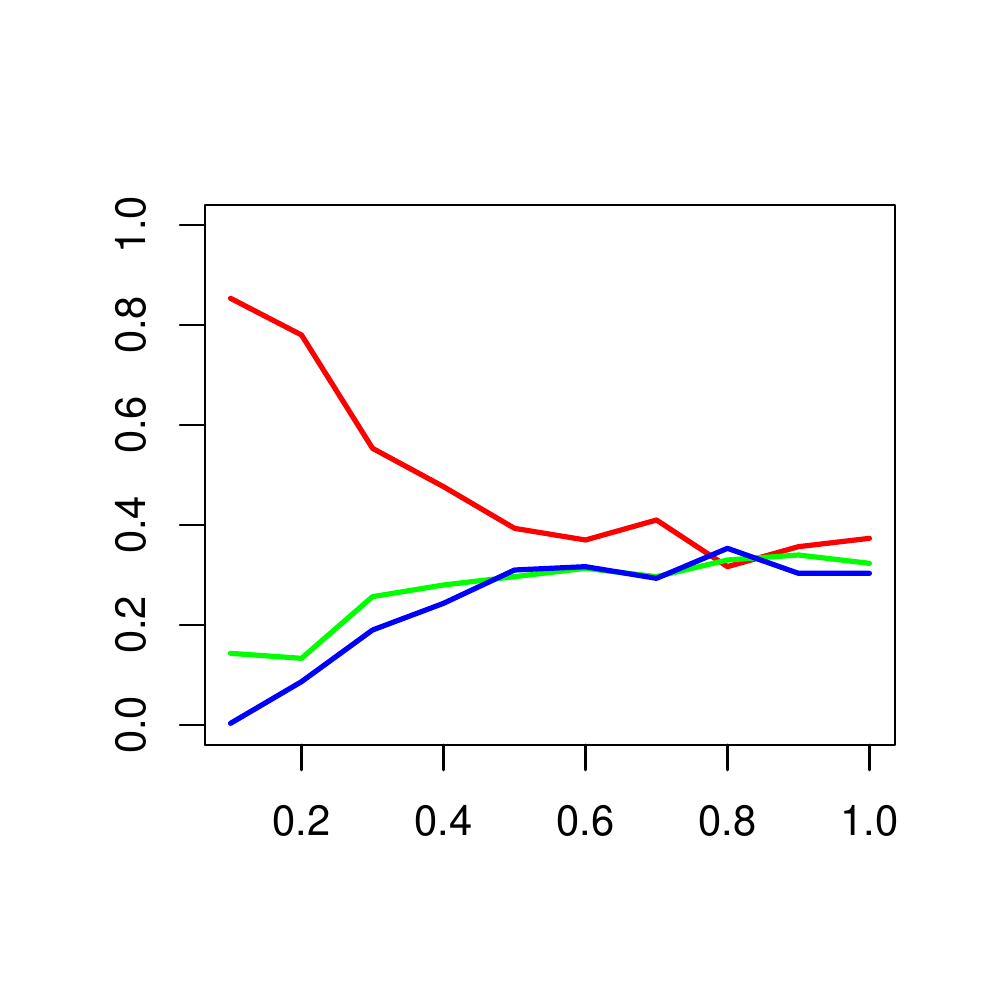}
		\vspace{-10mm} 
  	\caption{\scriptsize{Ackley}} \label{fig:iT30Ac}
	\end{subfigure}%
	\hspace{-6mm} 
	\begin{subfigure}[t]{.2\textwidth}
		\includegraphics[width=\textwidth]{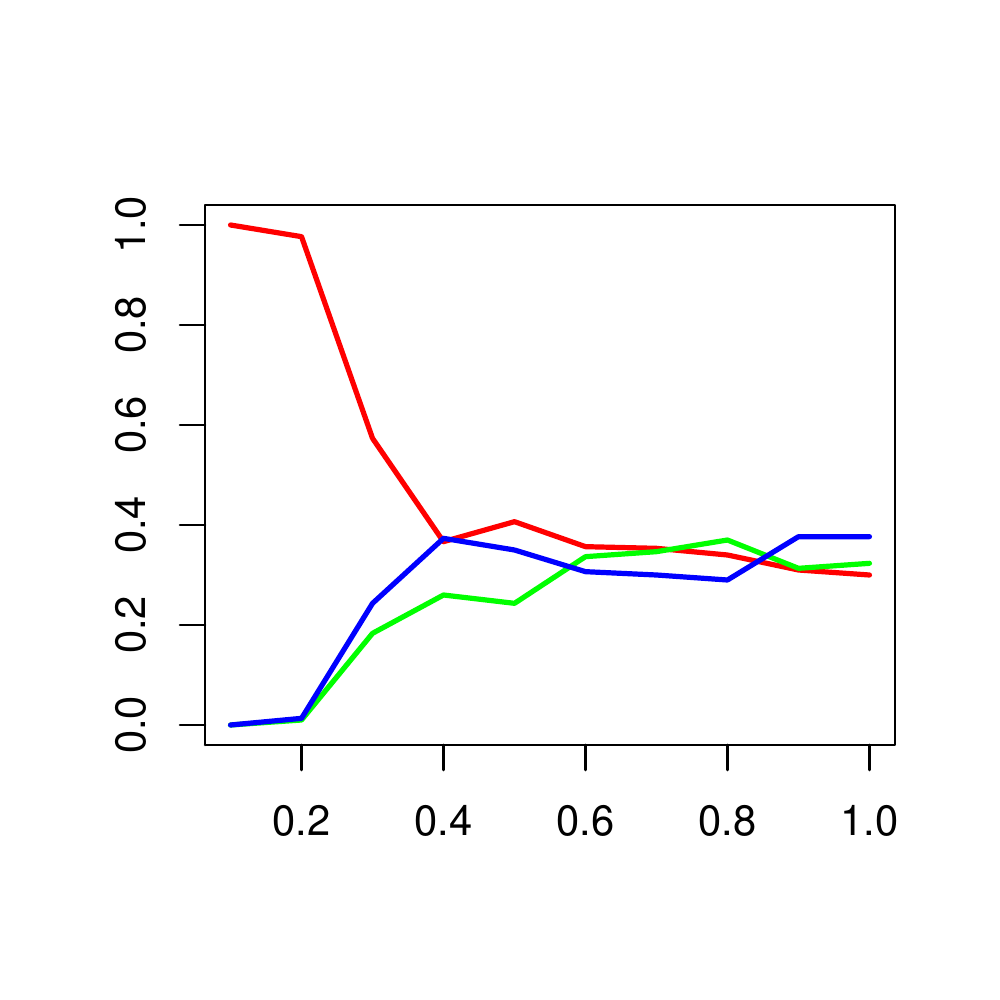}
		\vspace{-10mm} 
  	\caption{\scriptsize{Sphere}} \label{fig:iT30Sp}
	\end{subfigure}%
	\hspace{-6mm} 
	\begin{subfigure}[t]{.2\textwidth}
		\includegraphics[width=\textwidth]{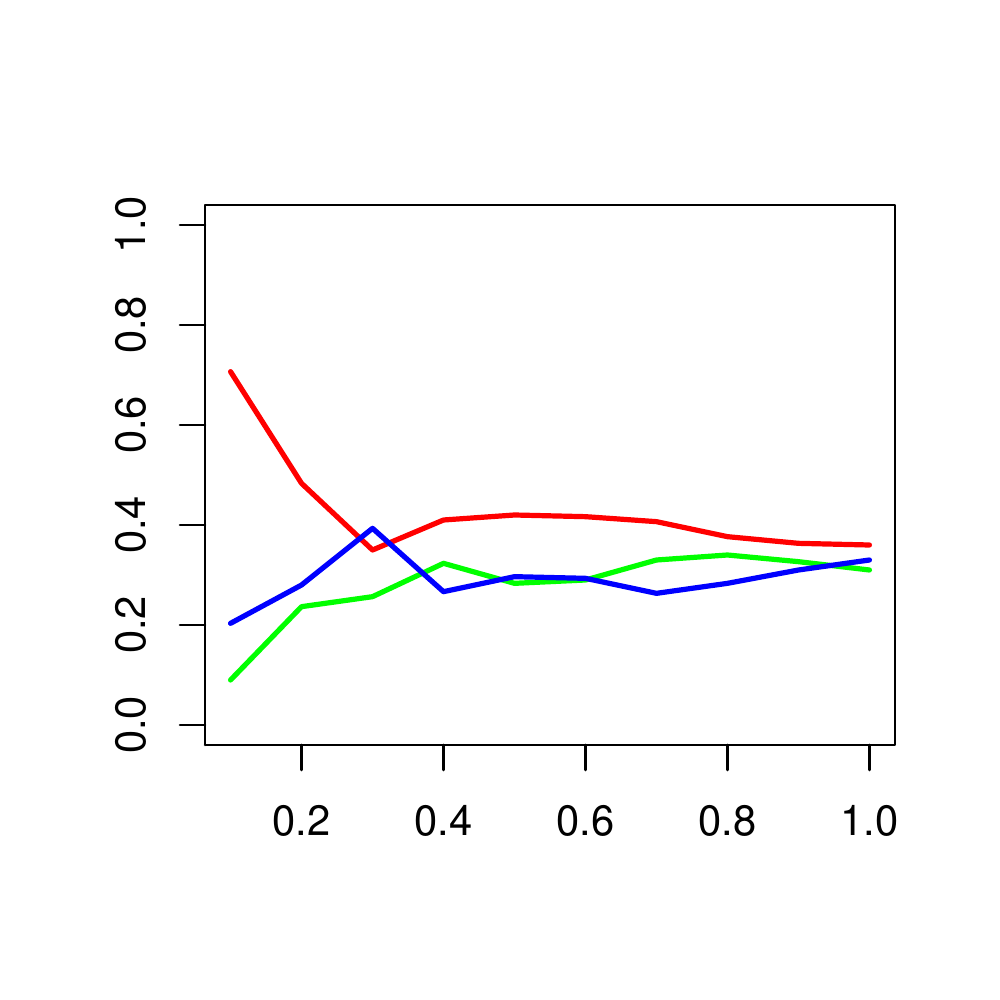}
		\vspace{-10mm} 
  	\caption{\scriptsize{Rosenbrock}} \label{fig:iT30Ro}
	\end{subfigure}
		
	\caption{Component -- decile  breakdowns for the form of inner search, across all GGGP heuristics at 30D. Components: Random (red), PSO (green), GA (blue).}
	\label{fig:inType30Comp}
	
\end{figure}

\begin{figure}[H]
	\centering
	
	\vspace{-5mm} 

	\begin{subfigure}[t]{.2\textwidth}
		\includegraphics[width=\textwidth]{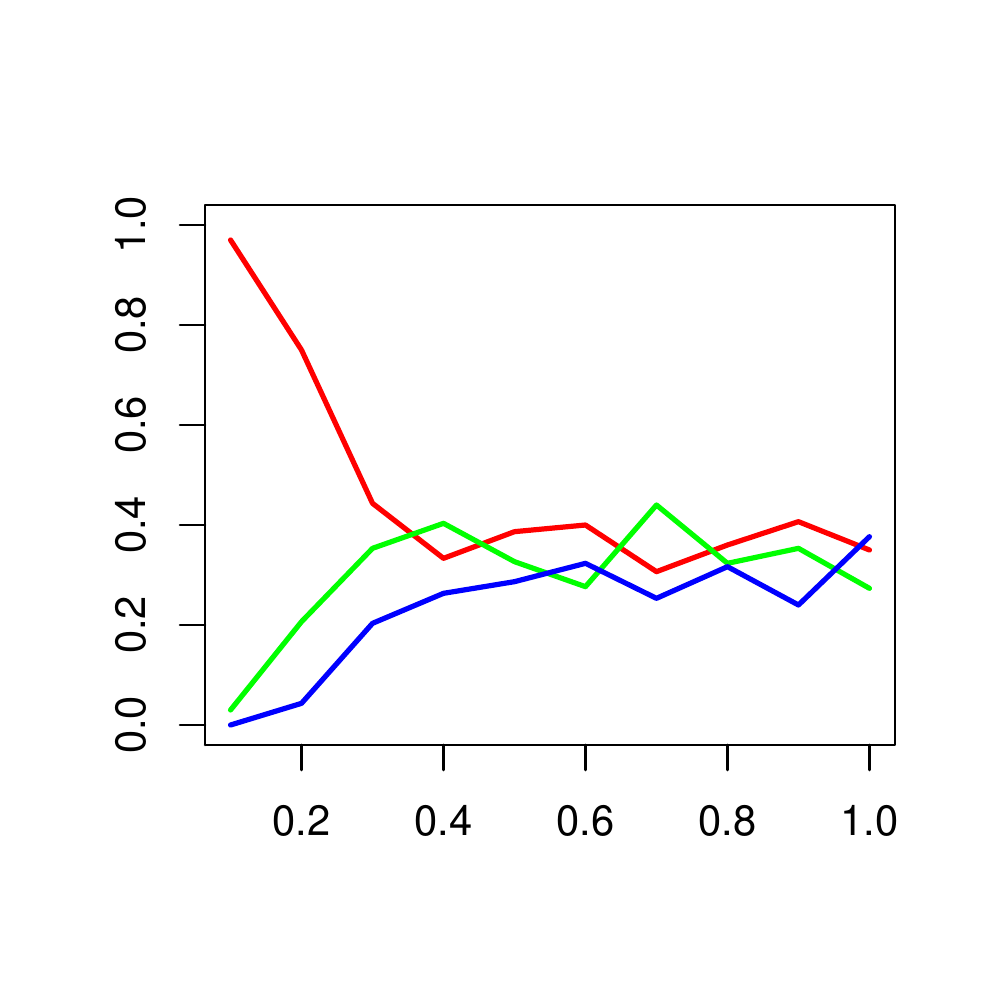}
		\vspace{-10mm} 
  	\caption{\scriptsize{Rastrigin}} \label{fig:iT100Ra}
	\end{subfigure}%
	\hspace{-6mm} 
	\begin{subfigure}[t]{.2\textwidth}
		\includegraphics[width=\textwidth]{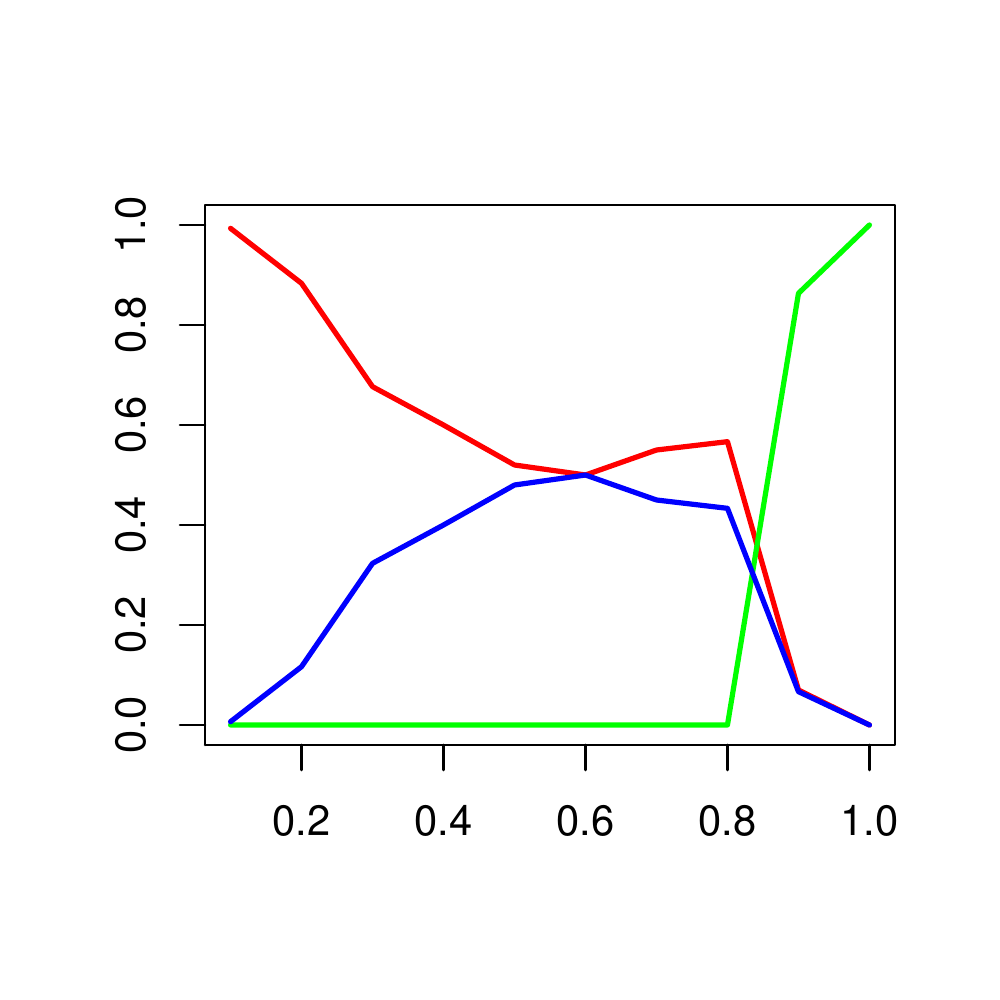}
		\vspace{-10mm} 
  	\caption{\scriptsize{Multipeak F1}} \label{fig:iT100M1}
	\end{subfigure}%
	\hspace{-6mm} 
	\begin{subfigure}[t]{.2\textwidth}
		\includegraphics[width=\textwidth]{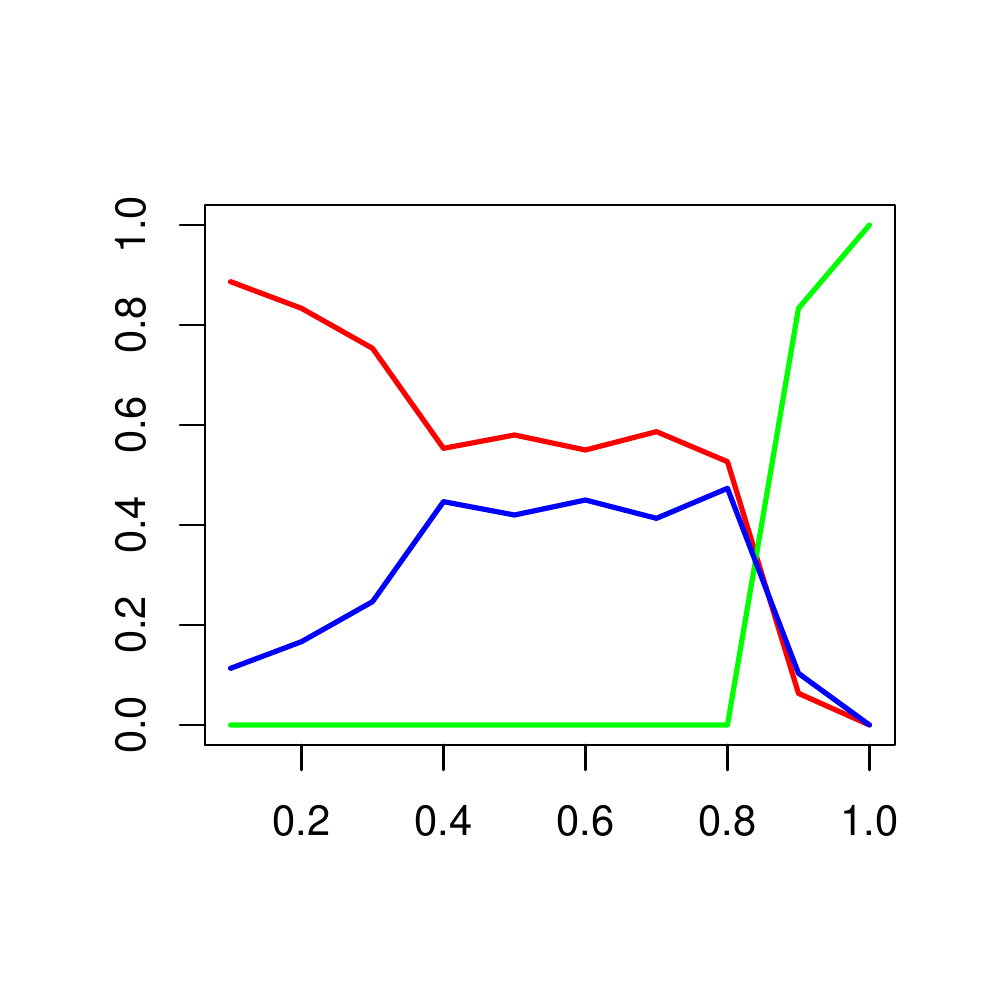}
		\vspace{-10mm} 
  	\caption{\scriptsize{Multipeak F2}} \label{fig:iT100M2}
	\end{subfigure}%
	\hspace{-6mm} 
	\begin{subfigure}[t]{.2\textwidth}
		\includegraphics[width=\textwidth]{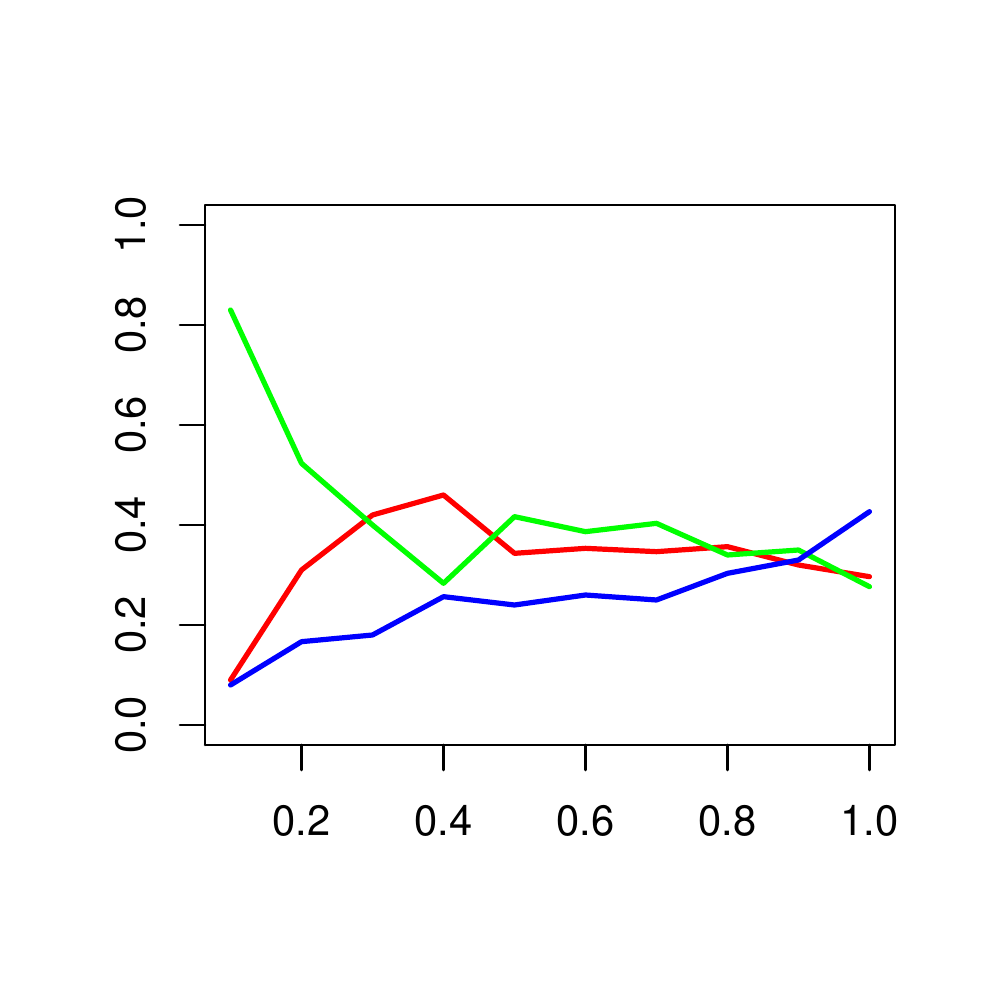}
		\vspace{-10mm} 
  	\caption{\scriptsize{Brankes}} \label{fig:iT100Br}
	\end{subfigure}%
	\hspace{-6mm} 
	\begin{subfigure}[t]{.2\textwidth}
		\includegraphics[width=\textwidth]{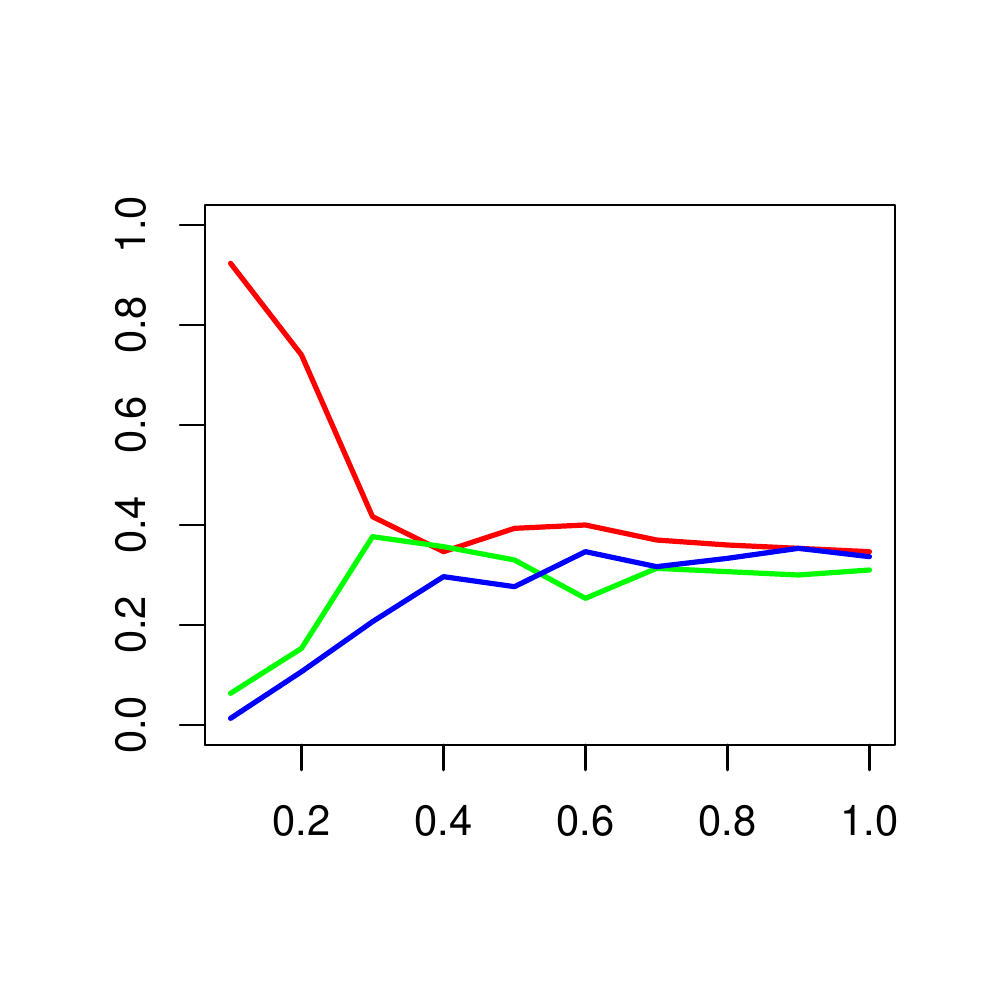}
		\vspace{-10mm} 
  	\caption{\scriptsize{Pickelhaube}} \label{fig:iT100Pi}
	\end{subfigure}
	
	\vspace{-2mm} 
		
	\begin{subfigure}[t]{.2\textwidth}
		\includegraphics[width=\textwidth]{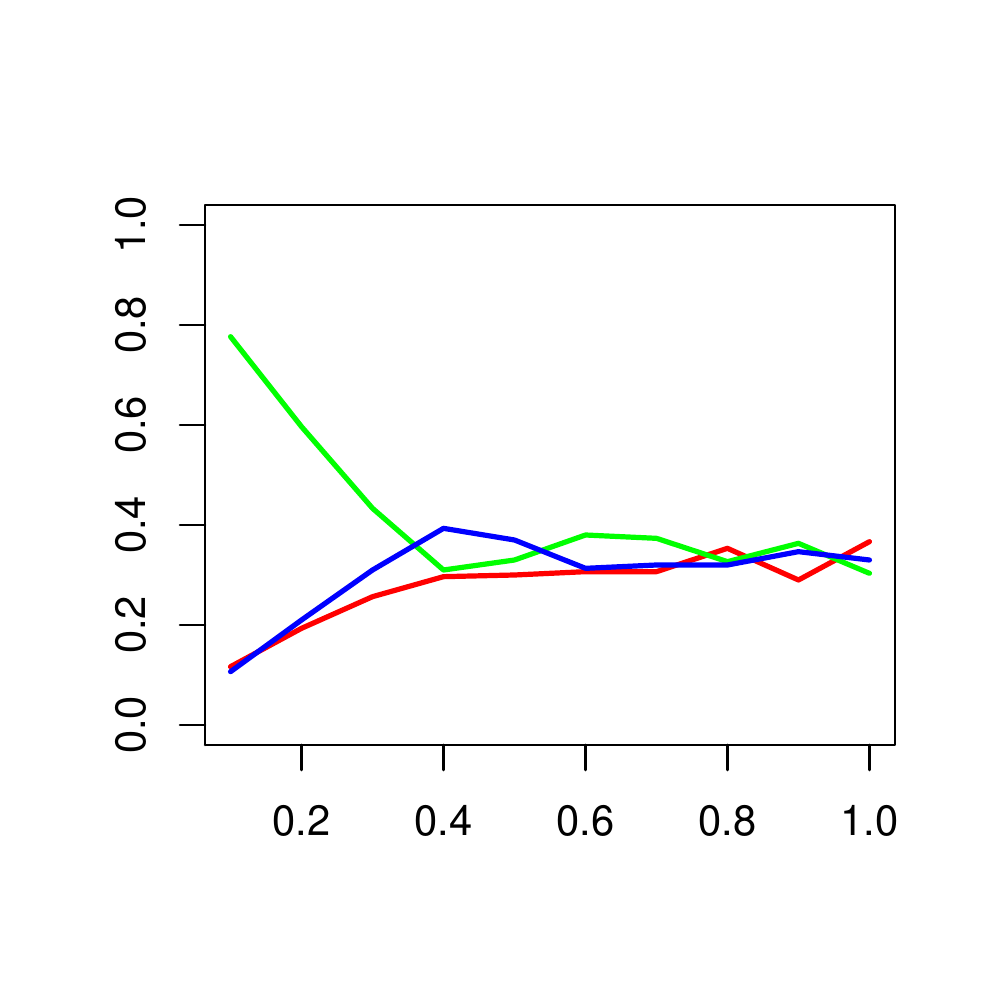}
		\vspace{-10mm} 
  	\caption{\scriptsize{Heaviside}} \label{fig:iT100Hv}
	\end{subfigure}%
	\hspace{-6mm} 
	\begin{subfigure}[t]{.2\textwidth}
		\includegraphics[width=\textwidth]{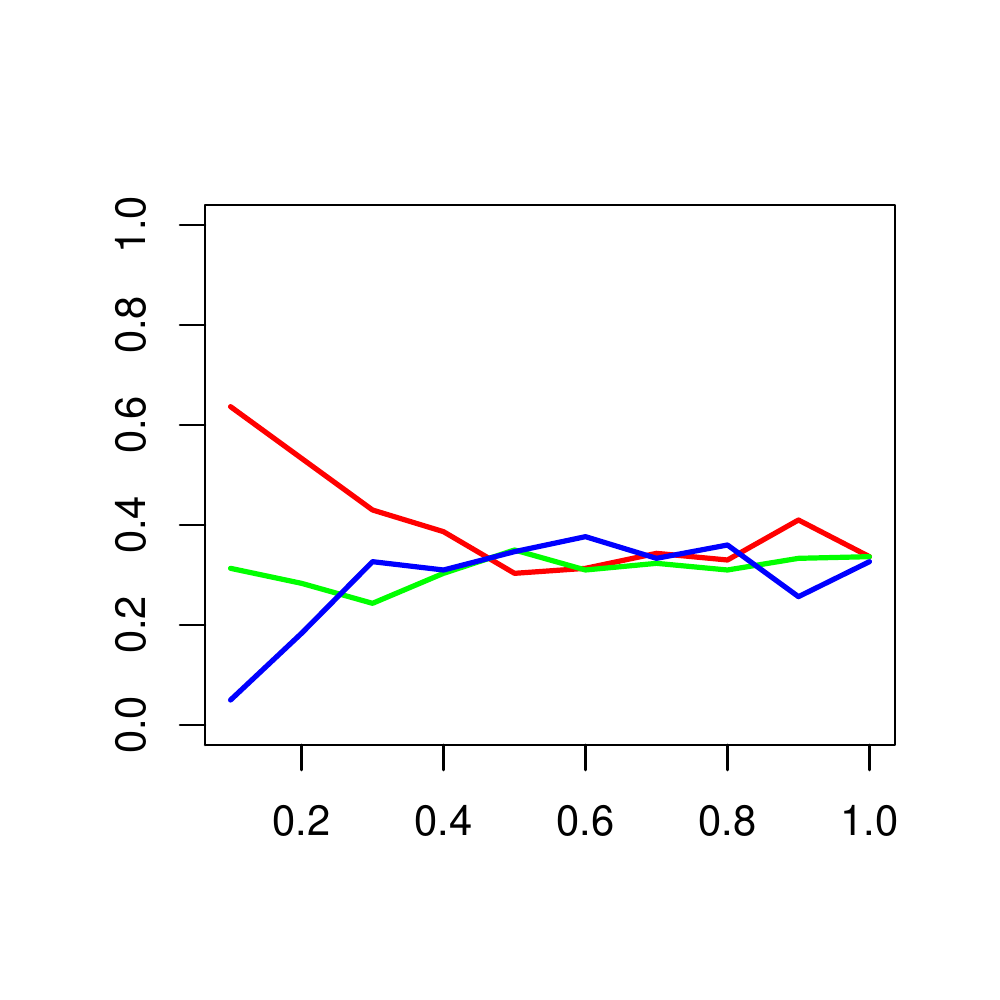}
		\vspace{-10mm} 
  	\caption{\scriptsize{Sawtooth}} \label{fig:iT100Sa}
	\end{subfigure}%
	\hspace{-6mm} 
	\begin{subfigure}[t]{.2\textwidth}
		\includegraphics[width=\textwidth]{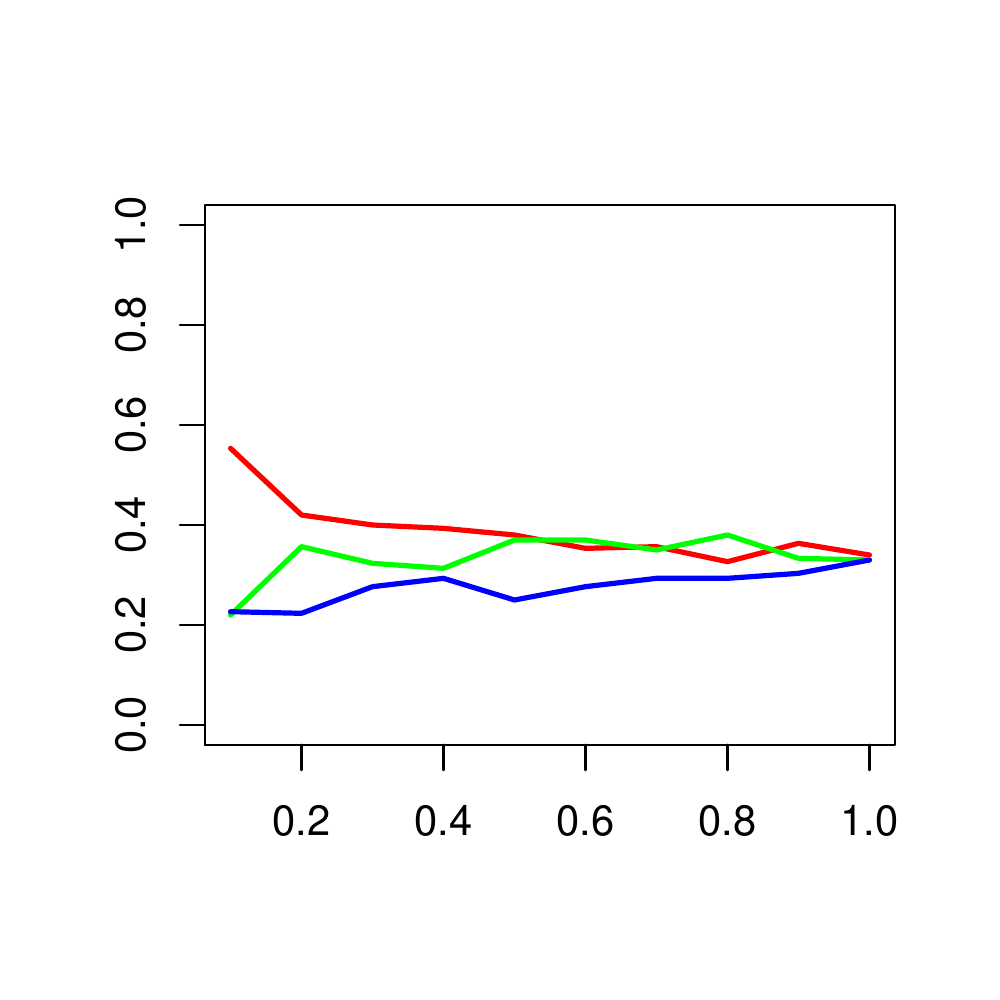}
		\vspace{-10mm} 
  	\caption{\scriptsize{Ackley}} \label{fig:iT100Ac}
	\end{subfigure}%
	\hspace{-6mm} 
	\begin{subfigure}[t]{.2\textwidth}
		\includegraphics[width=\textwidth]{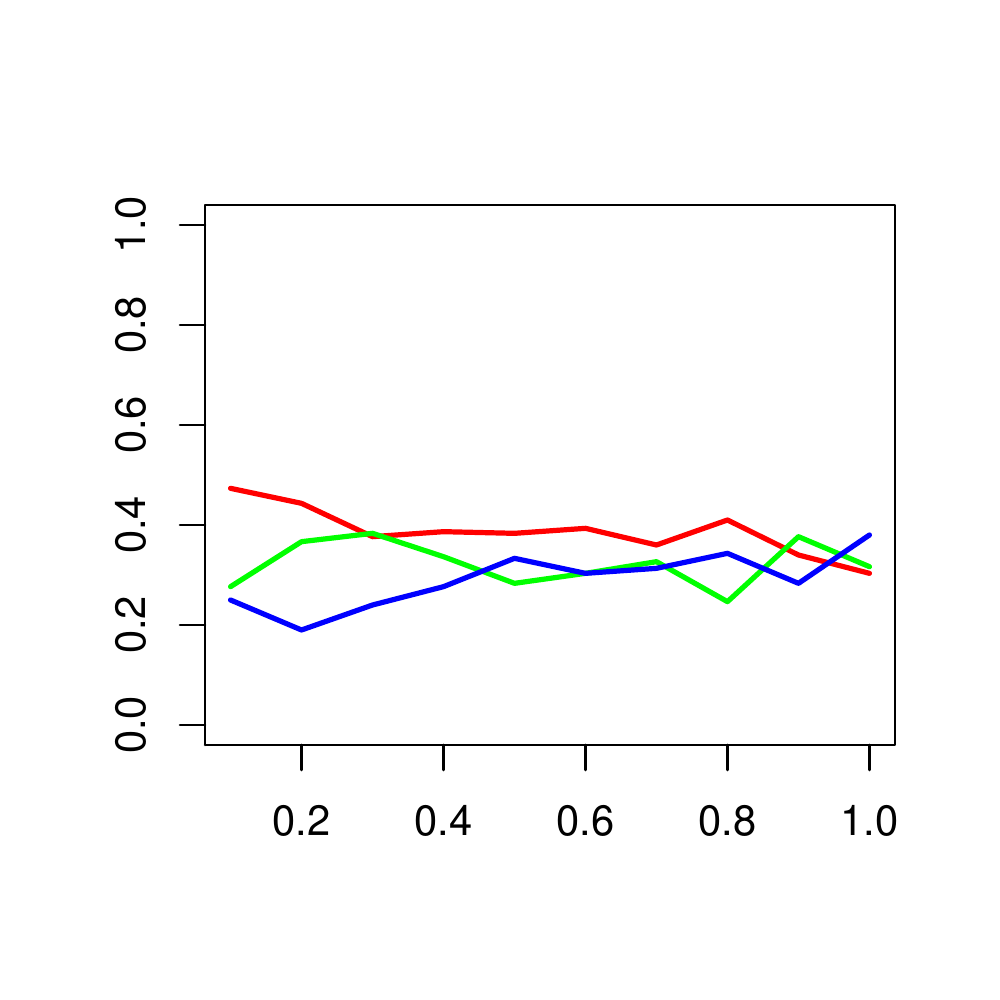}
		\vspace{-10mm} 
  	\caption{\scriptsize{Sphere}} \label{fig:iT100Sp}
	\end{subfigure}%
	\hspace{-6mm} 
	\begin{subfigure}[t]{.2\textwidth}
		\includegraphics[width=\textwidth]{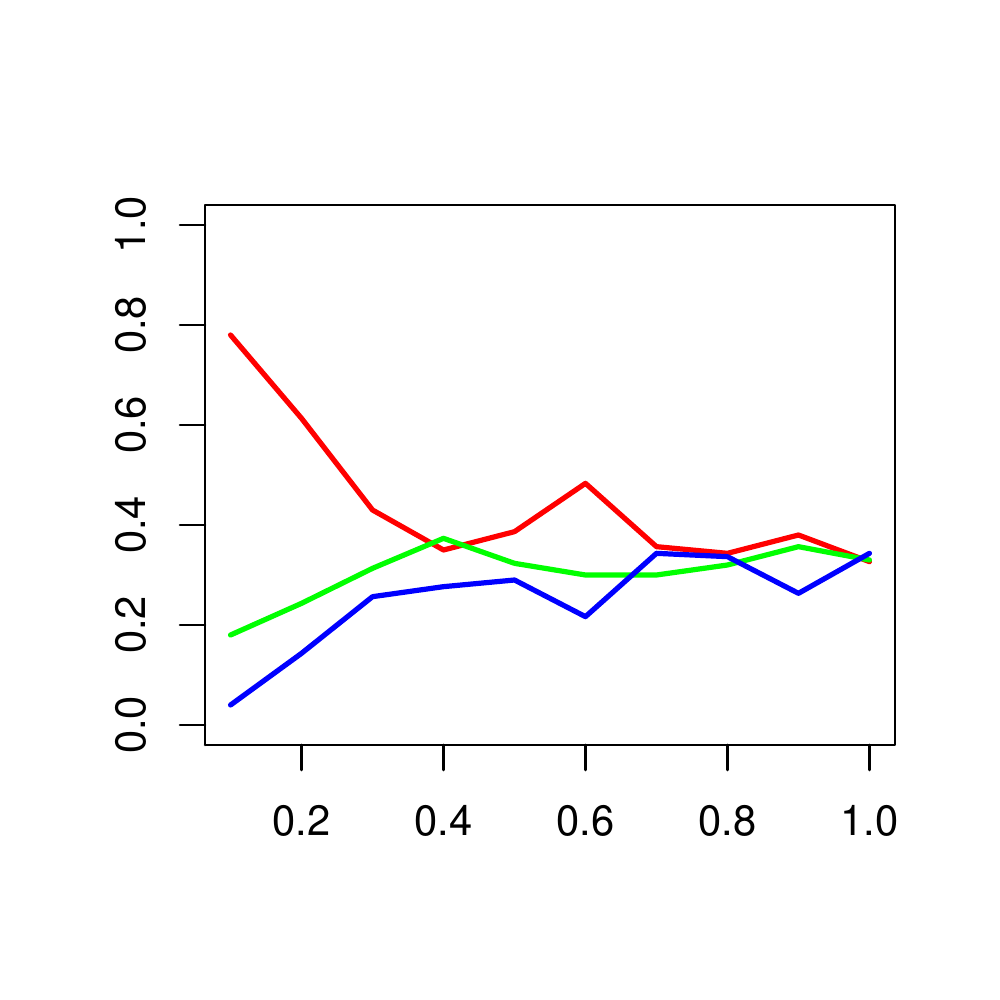}
		\vspace{-10mm} 
  	\caption{\scriptsize{Rosenbrock}} \label{fig:iT100Ro}
	\end{subfigure}
		
	\caption{Component -- decile  breakdowns for the form of inner search, across all GGGP heuristics at 100D. Components: Random (red), PSO (green), GA (blue).}
	\label{fig:inType100Comp}
	
\end{figure}

\begin{figure}[H]
	\centering
	
	\vspace{-5mm} 

	\begin{subfigure}[t]{.24\textwidth}
		\includegraphics[width=\textwidth]{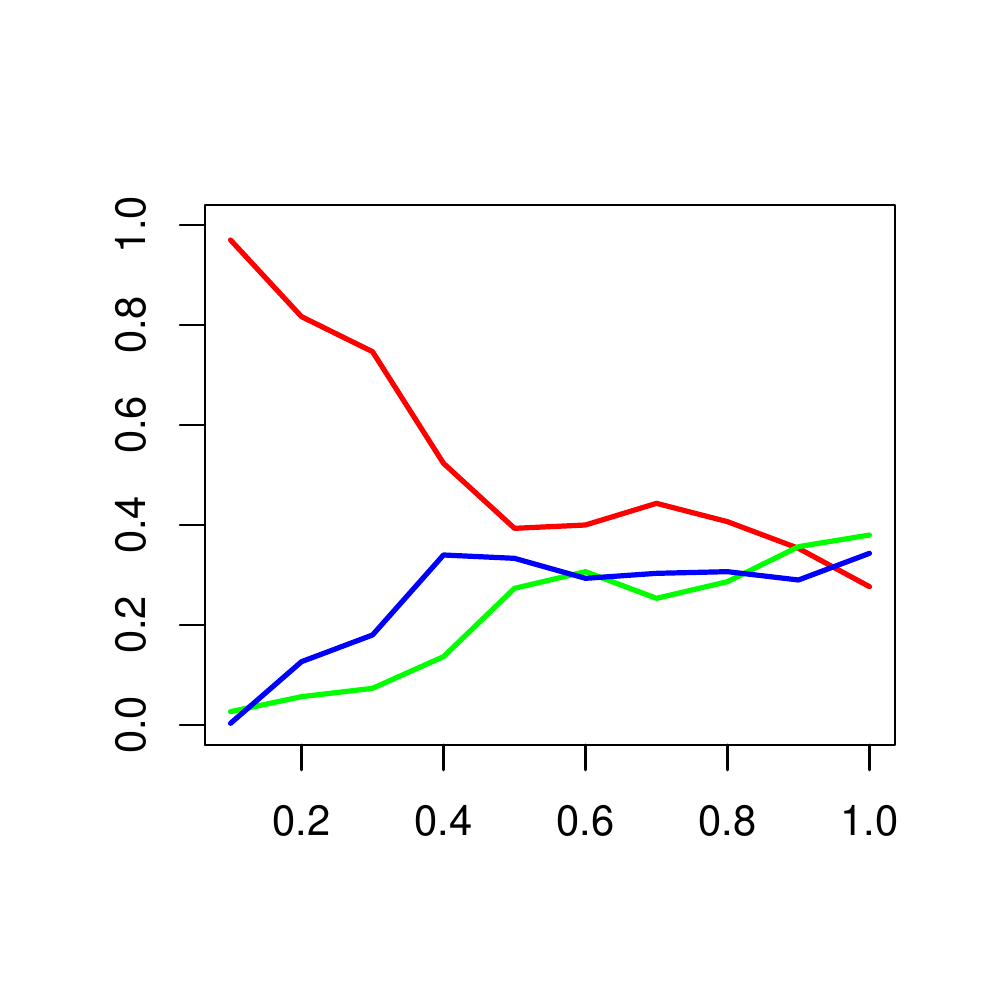}
		\vspace{-10mm} 
  	\caption{\scriptsize{30D}} \label{fig:iT30}
	\end{subfigure}%
	\hspace{-6mm} 
	\begin{subfigure}[t]{.24\textwidth}
		\includegraphics[width=\textwidth]{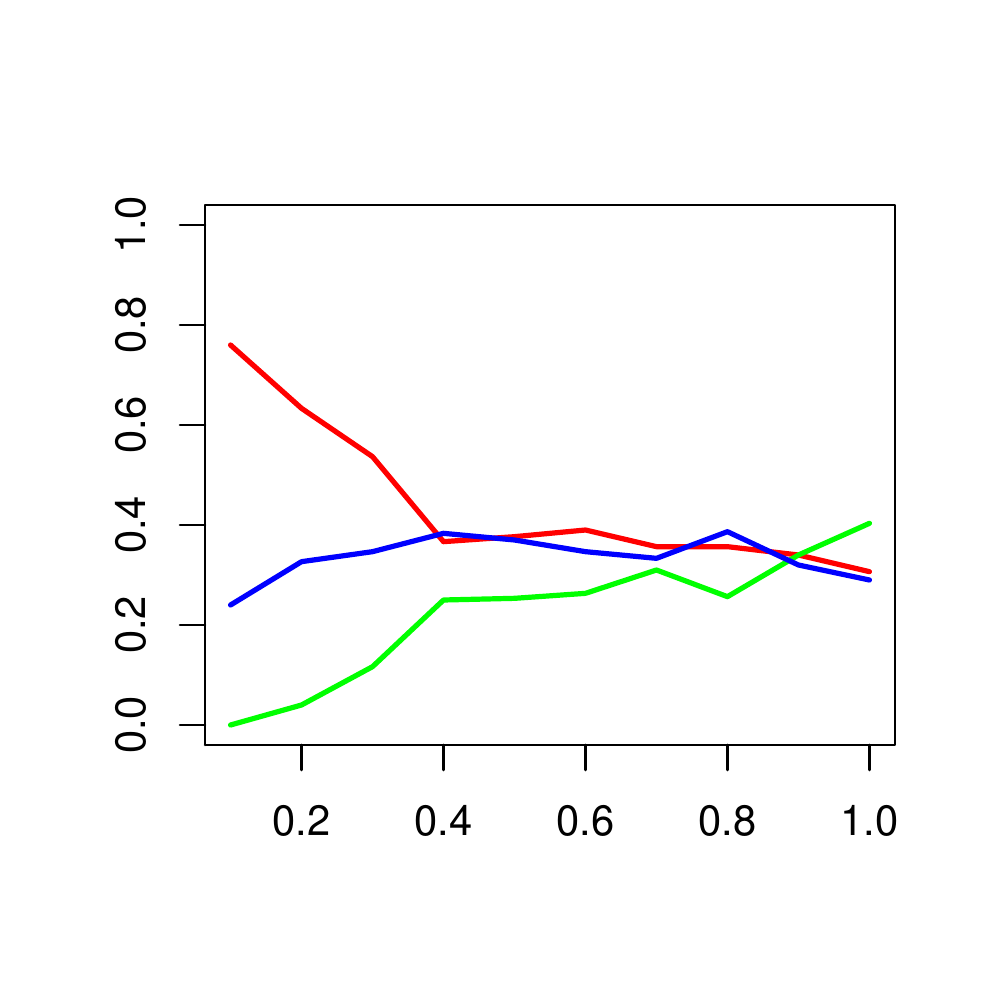}
		\vspace{-10mm} 
  	\caption{\scriptsize{100D}} \label{fig:iT100}
	\end{subfigure}

	\caption{Component -- decile breakdowns for the form of inner search, across all GGGP heuristics at 30D and 100D for the general heuristics. Components: Random (red), PSO (green), GA (blue).}
	\label{fig:inTyComp}
	
\end{figure}

From Table~\ref{fig:proportions} and Figures~\ref{fig:inType30Comp} to~\ref{fig:inTyComp} it can be seen that the most used form of inner maximisation search is random sampling. For the 30D individual cases random sampling is the most commonly associated choice with the best performing heuristics, dominating for most problems. At 100D random sampling is again most typically associated with the best heuristics, although it is much less dominant, with PSO the most used form of inner search in the best performing heuristics for the Branke and Heaviside problems. For the general cases random sampling dominates the best performing heuristics.

\begin{figure}[H]
	\centering
	

	\begin{subfigure}[t]{.2\textwidth}
		\includegraphics[width=\textwidth]{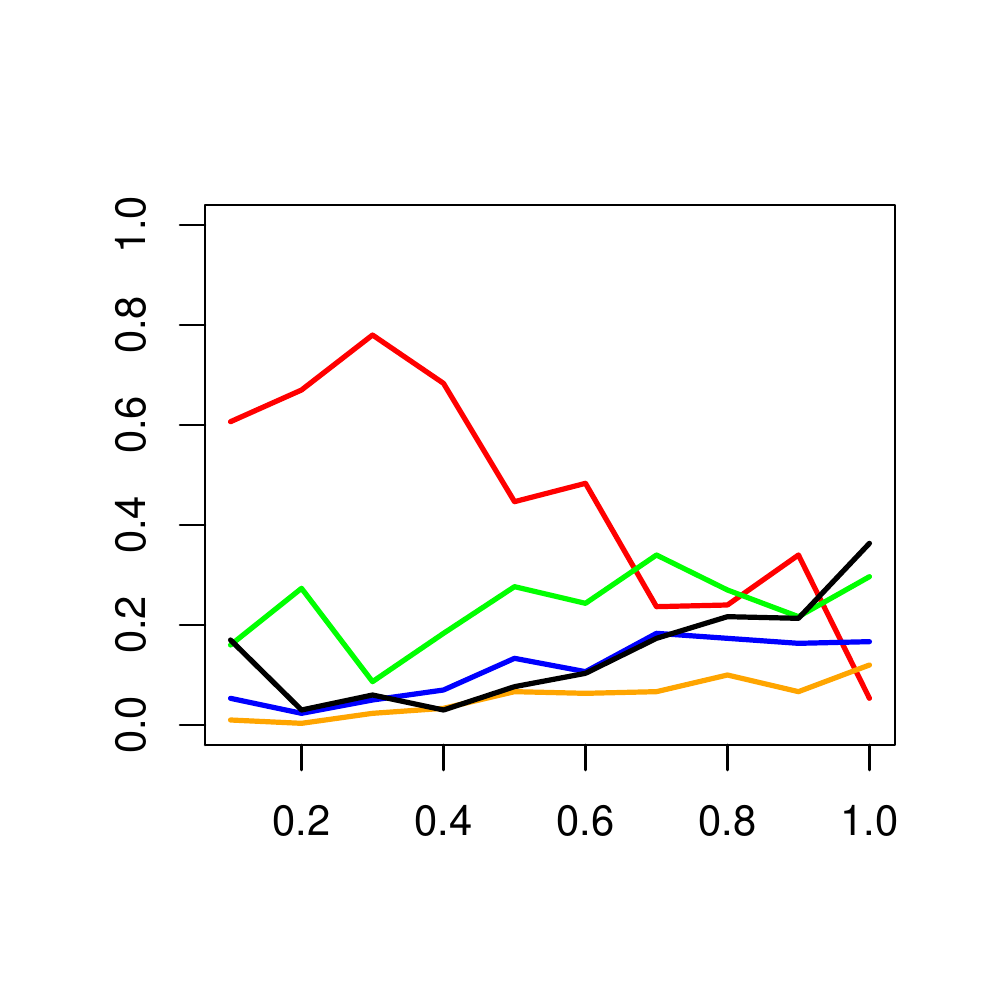}
		\vspace{-10mm} 
  	\caption{\scriptsize{Rastrigin}} \label{fig:iX30Ra}
	\end{subfigure}%
	\hspace{-6mm} 
	\begin{subfigure}[t]{.2\textwidth}
		\includegraphics[width=\textwidth]{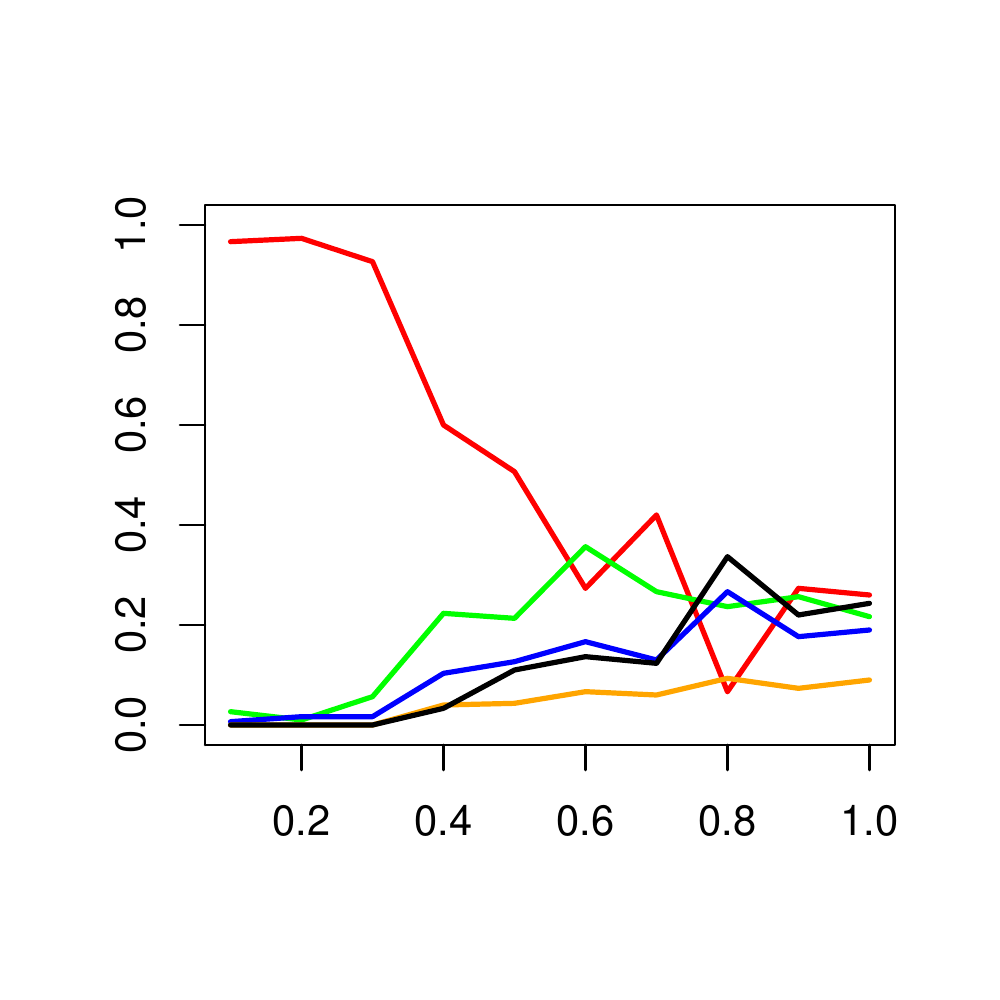}
		\vspace{-10mm} 
  	\caption{\scriptsize{Multipeak F1}} \label{fig:iX30M1}
	\end{subfigure}%
	\hspace{-6mm} 
	\begin{subfigure}[t]{.2\textwidth}
		\includegraphics[width=\textwidth]{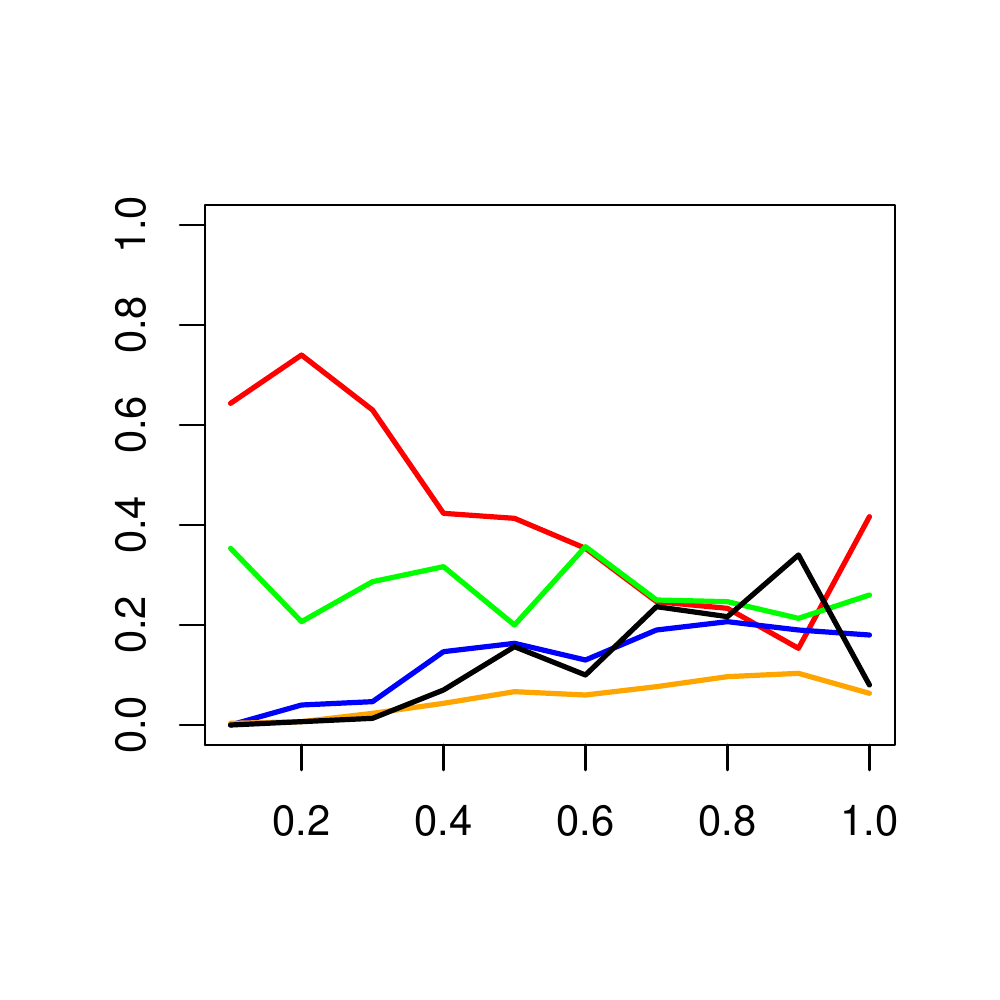}
		\vspace{-10mm} 
  	\caption{\scriptsize{Multipeak F2}} \label{fig:iX30M2}
	\end{subfigure}%
	\hspace{-6mm} 
	\begin{subfigure}[t]{.2\textwidth}
		\includegraphics[width=\textwidth]{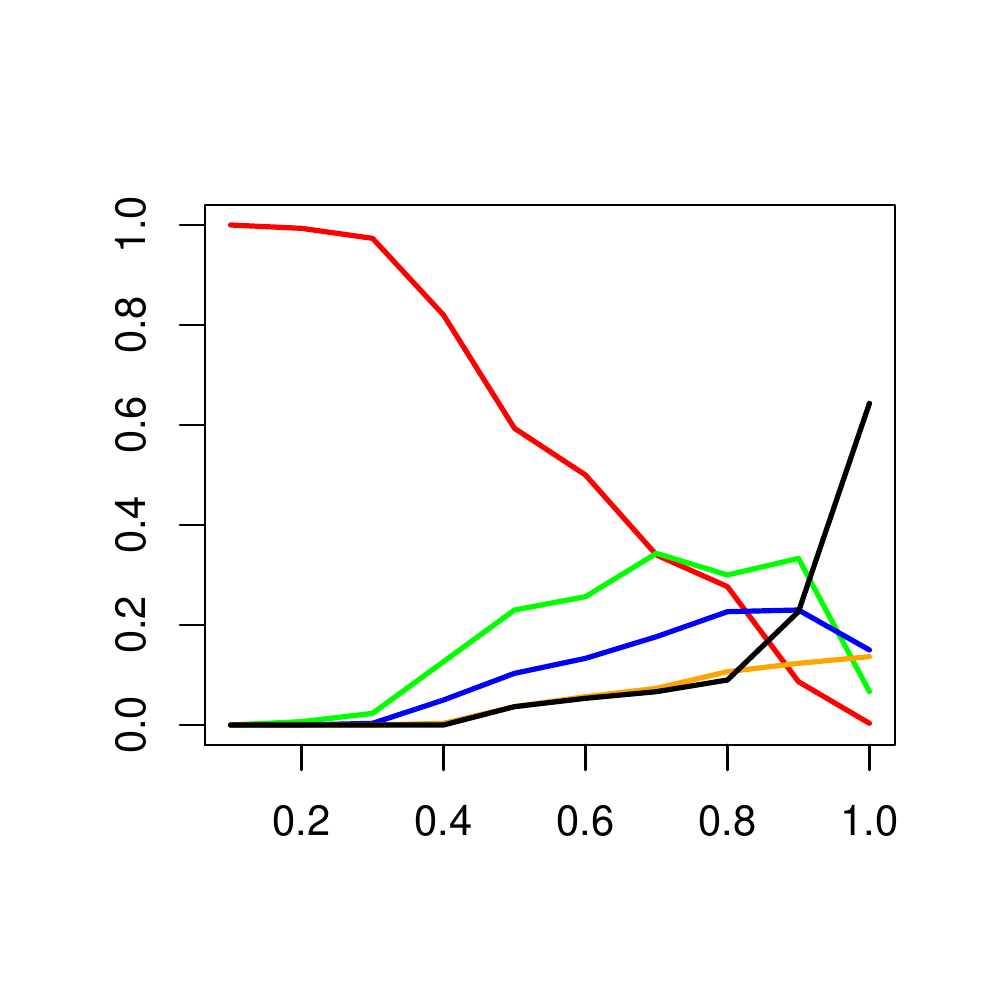}
		\vspace{-10mm} 
  	\caption{\scriptsize{Brankes}} \label{fig:iX30Br}
	\end{subfigure}%
	\hspace{-6mm} 
	\begin{subfigure}[t]{.2\textwidth}
		\includegraphics[width=\textwidth]{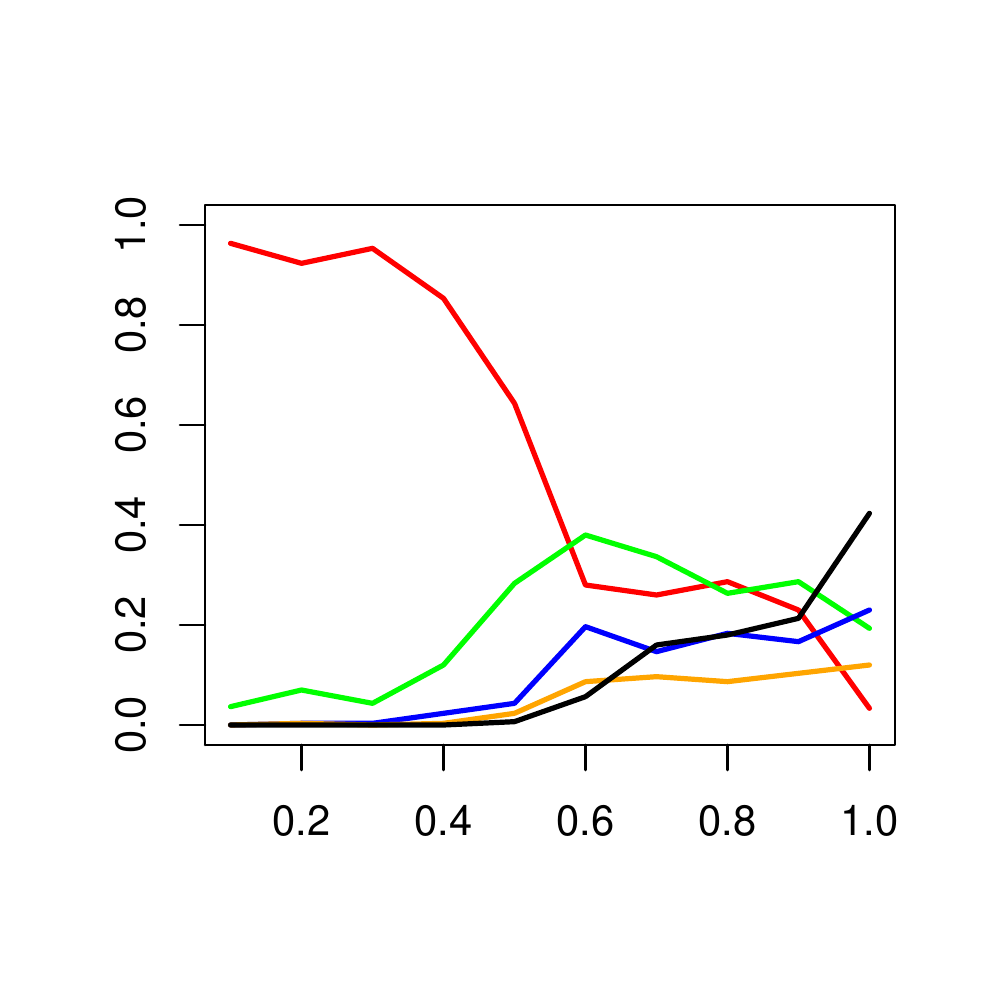}
		\vspace{-10mm} 
  	\caption{\scriptsize{Pickelhaube}} \label{fig:iX30Pi}
	\end{subfigure}
	
	\vspace{-2mm} 
		
	\begin{subfigure}[t]{.2\textwidth}
		\includegraphics[width=\textwidth]{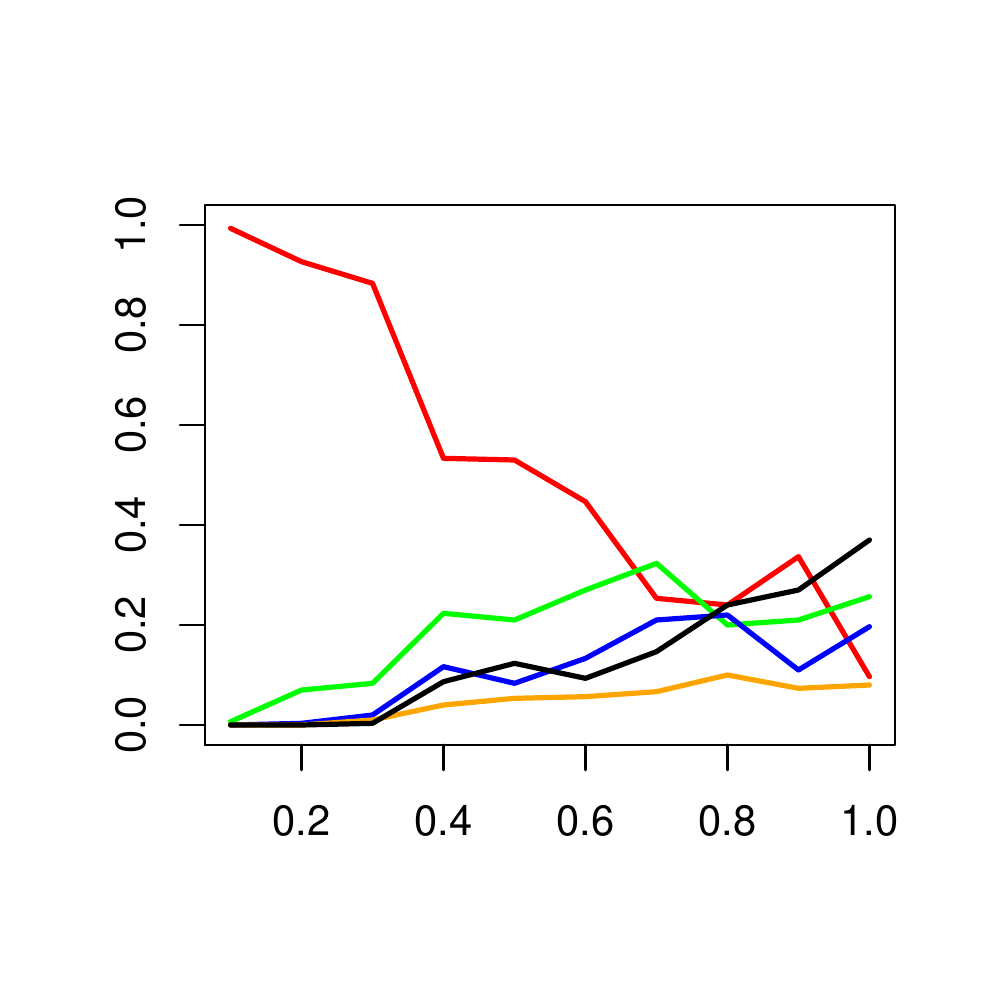}
		\vspace{-10mm} 
  	\caption{\scriptsize{Heaviside}} \label{fig:iX30Hv}
	\end{subfigure}%
	\hspace{-6mm} 
	\begin{subfigure}[t]{.2\textwidth}
		\includegraphics[width=\textwidth]{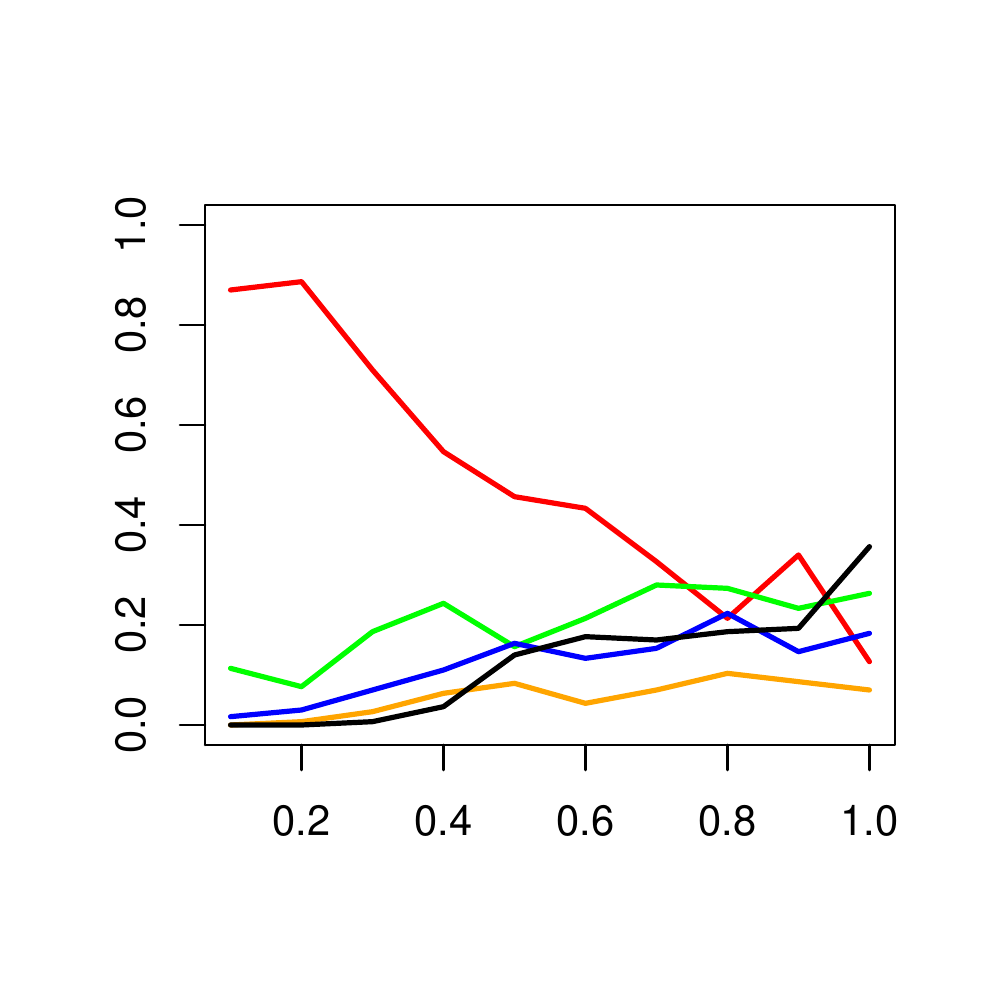}
		\vspace{-10mm} 
  	\caption{\scriptsize{Sawtooth}} \label{fig:iX30Sa}
	\end{subfigure}%
	\hspace{-6mm} 
	\begin{subfigure}[t]{.2\textwidth}
		\includegraphics[width=\textwidth]{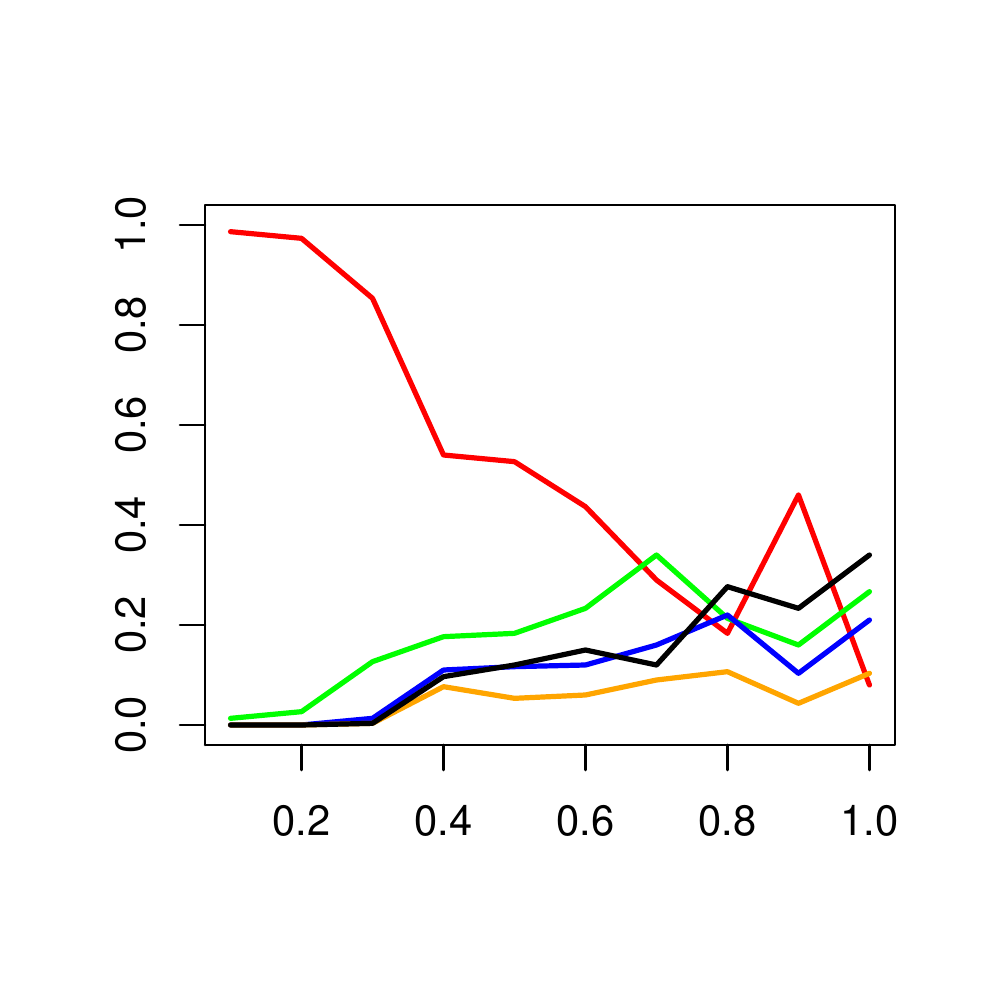}
		\vspace{-10mm} 
  	\caption{\scriptsize{Ackley}} \label{fig:iX30Ac}
	\end{subfigure}%
	\hspace{-6mm} 
	\begin{subfigure}[t]{.2\textwidth}
		\includegraphics[width=\textwidth]{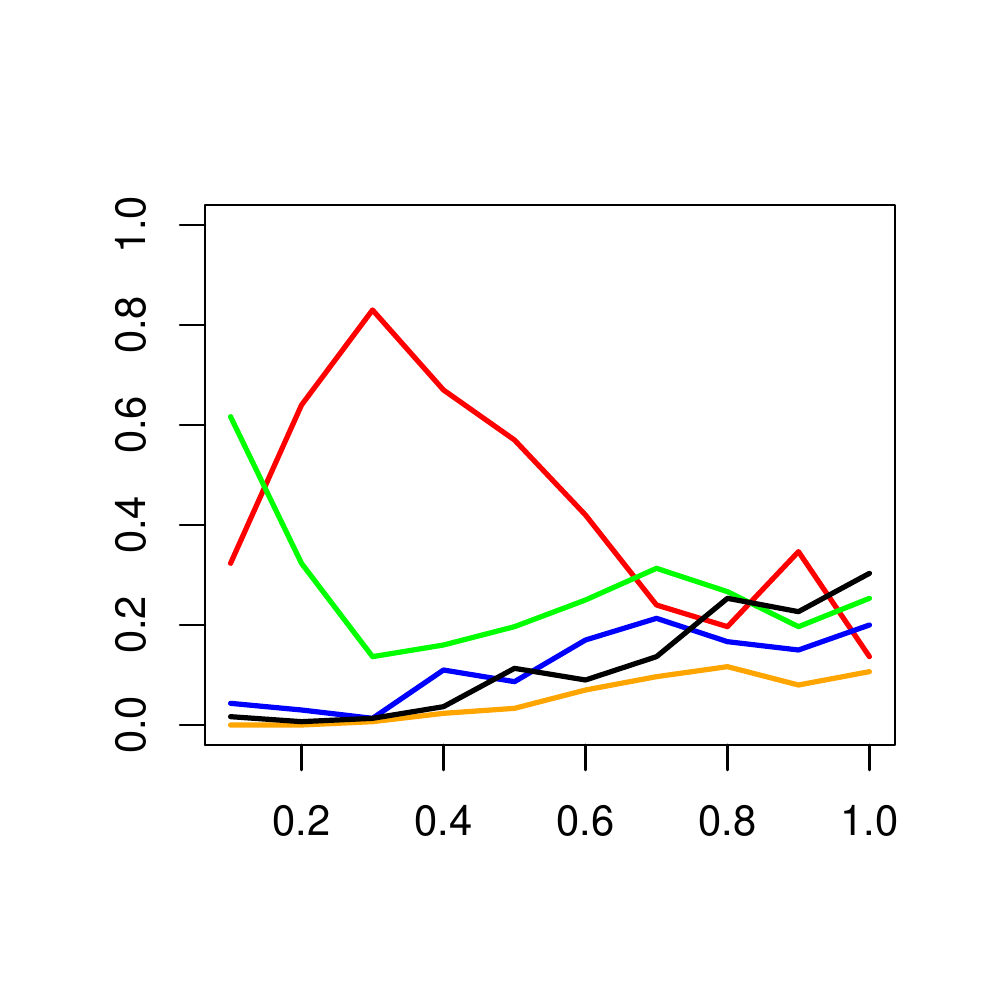}
		\vspace{-10mm} 
  	\caption{\scriptsize{Sphere}} \label{fig:iX30Sp}
	\end{subfigure}%
	\hspace{-6mm} 
	\begin{subfigure}[t]{.2\textwidth}
		\includegraphics[width=\textwidth]{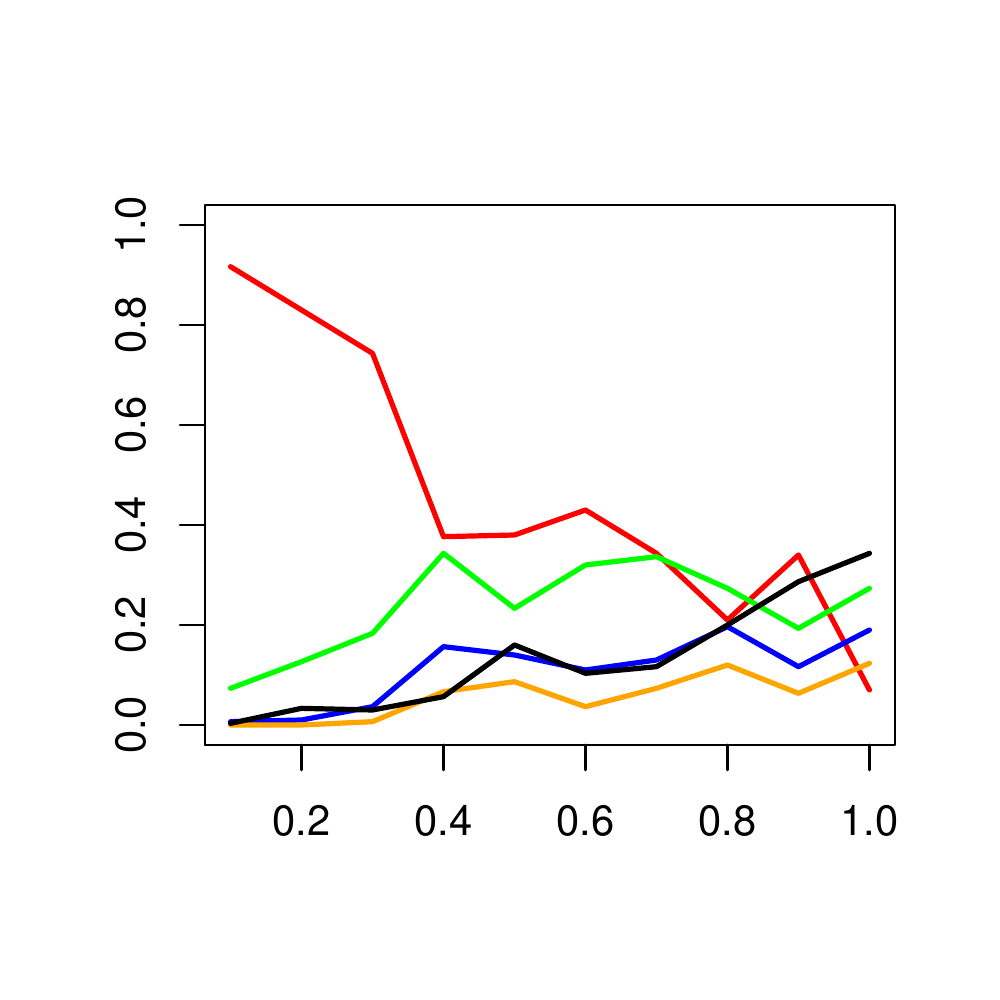}
		\vspace{-10mm} 
  	\caption{\scriptsize{Rosenbrock}} \label{fig:iX30Ro}
	\end{subfigure}
		
	\caption{Component -- decile breakdowns for the extent (size) of the inner maximisation search, across all GGGP heuristics at 30D. Components: [2 -- 10] (red), [11 -- 20] (green), [21 -- 30] (blue), [31 -- 40] (orange), \textgreater \ 40 (black).}
	\label{fig:inExt30Comp}
	
\end{figure}

\begin{figure}[H]
	\centering
	
	\vspace{-5mm} 

	\begin{subfigure}[t]{.2\textwidth}
		\includegraphics[width=\textwidth]{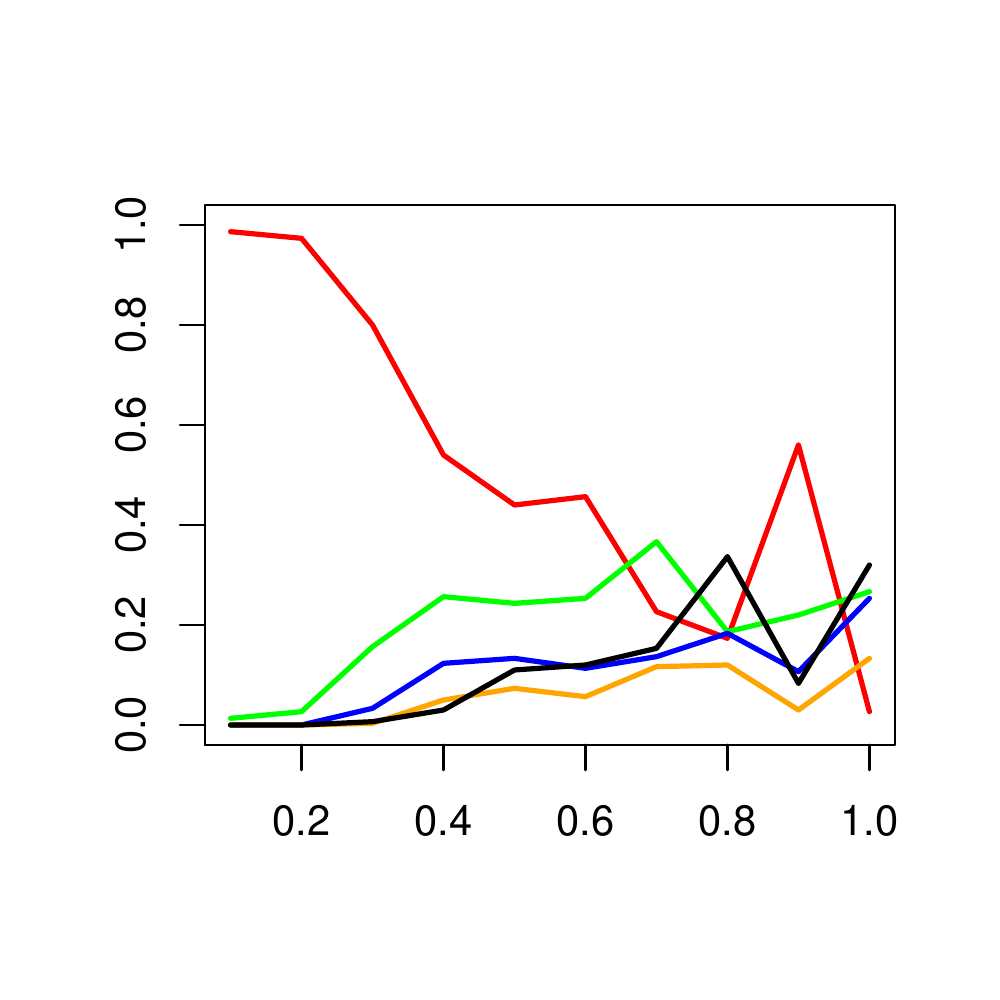}
		\vspace{-10mm} 
  	\caption{\scriptsize{Rastrigin}} \label{fig:iX100Ra}
	\end{subfigure}%
	\hspace{-6mm} 
	\begin{subfigure}[t]{.2\textwidth}
		\includegraphics[width=\textwidth]{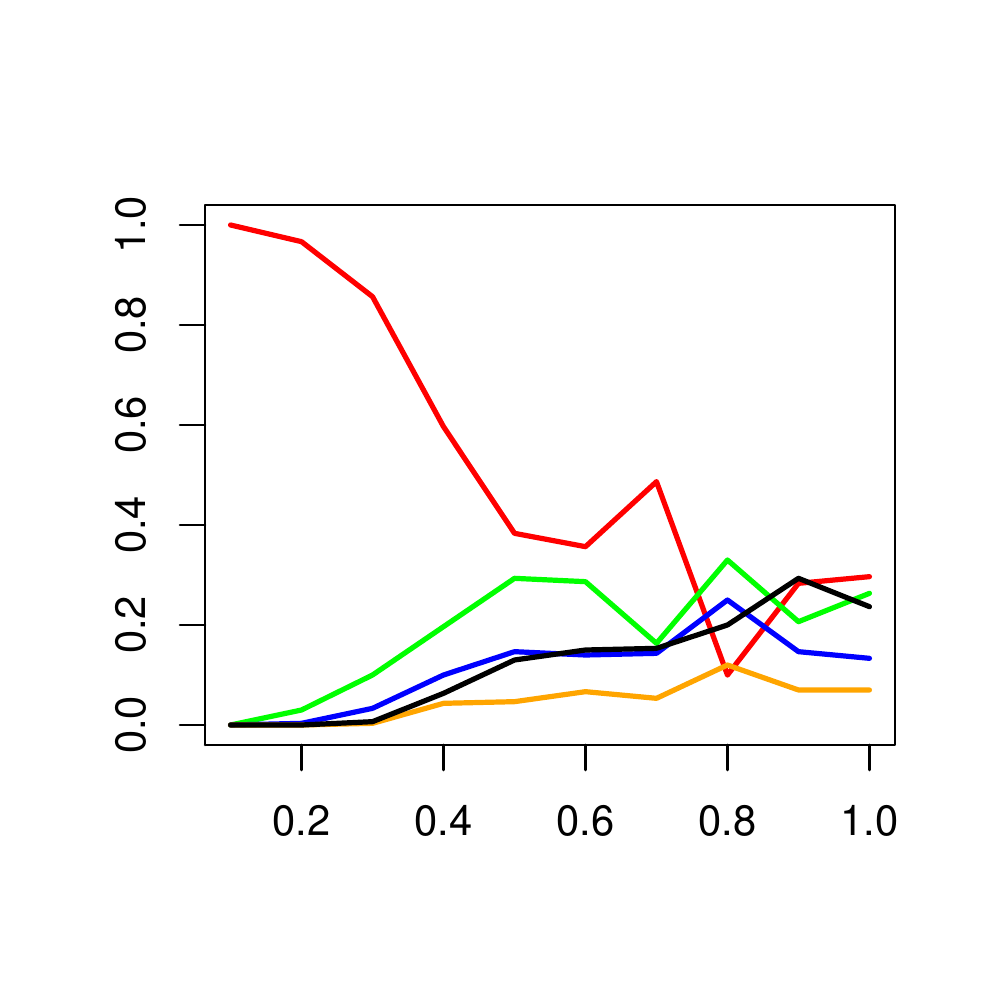}
		\vspace{-10mm} 
  	\caption{\scriptsize{Multipeak F1}} \label{fig:iX100M1}
	\end{subfigure}%
	\hspace{-6mm} 
	\begin{subfigure}[t]{.2\textwidth}
		\includegraphics[width=\textwidth]{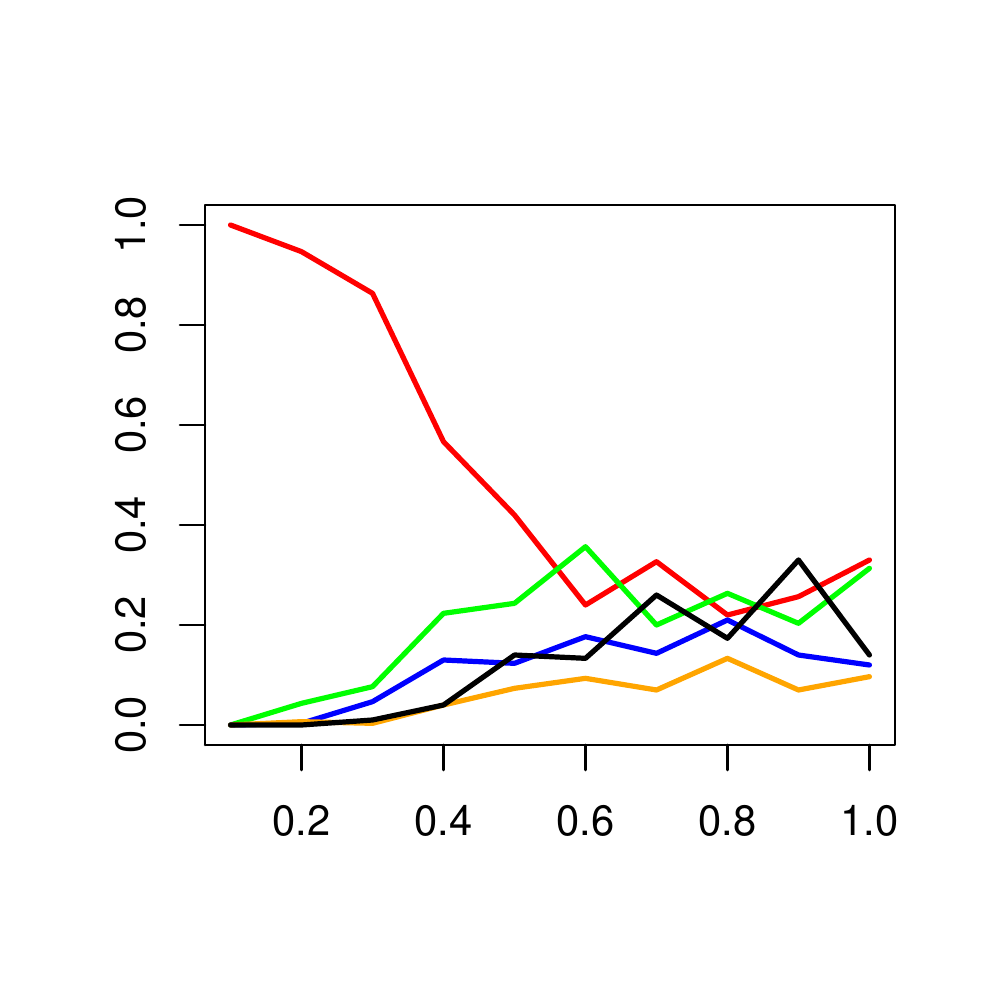}
		\vspace{-10mm} 
  	\caption{\scriptsize{Multipeak F2}} \label{fig:iX100M2}
	\end{subfigure}%
	\hspace{-6mm} 
	\begin{subfigure}[t]{.2\textwidth}
		\includegraphics[width=\textwidth]{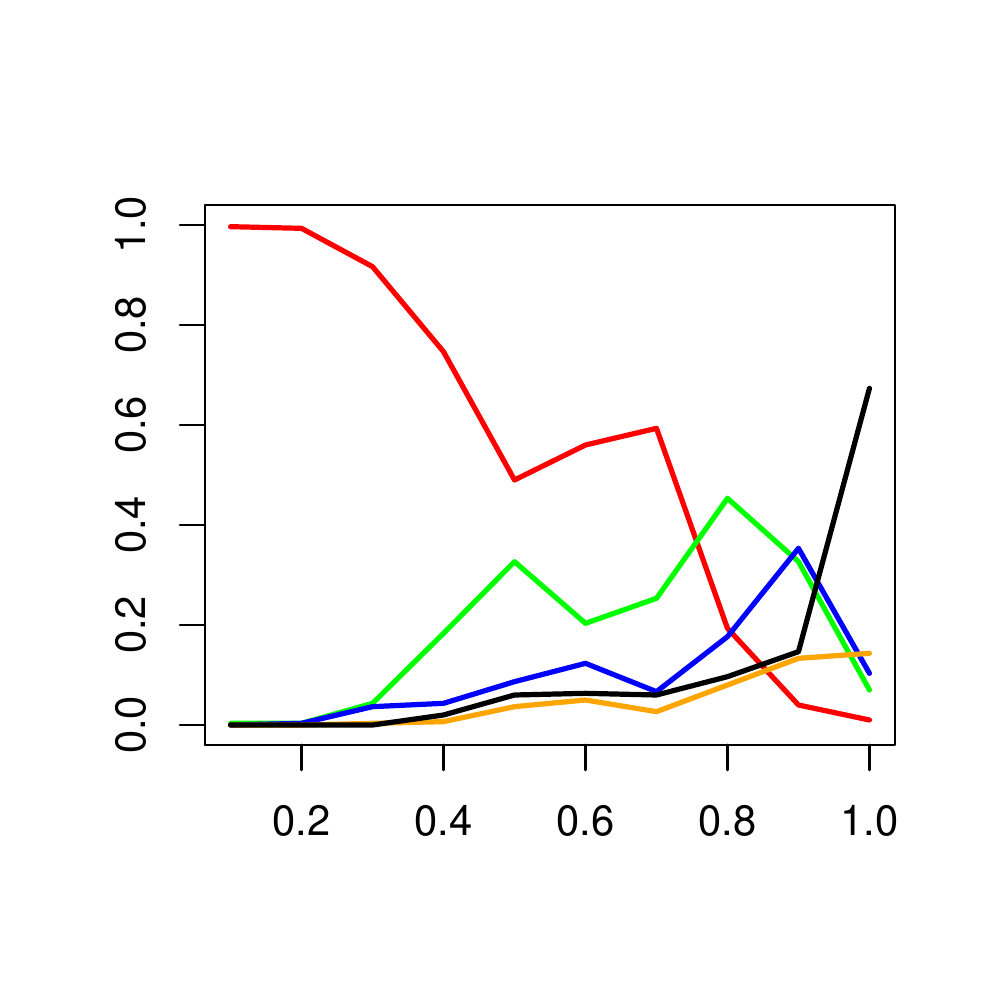}
		\vspace{-10mm} 
  	\caption{\scriptsize{Brankes}} \label{fig:iX100Br}
	\end{subfigure}%
	\hspace{-6mm} 
	\begin{subfigure}[t]{.2\textwidth}
		\includegraphics[width=\textwidth]{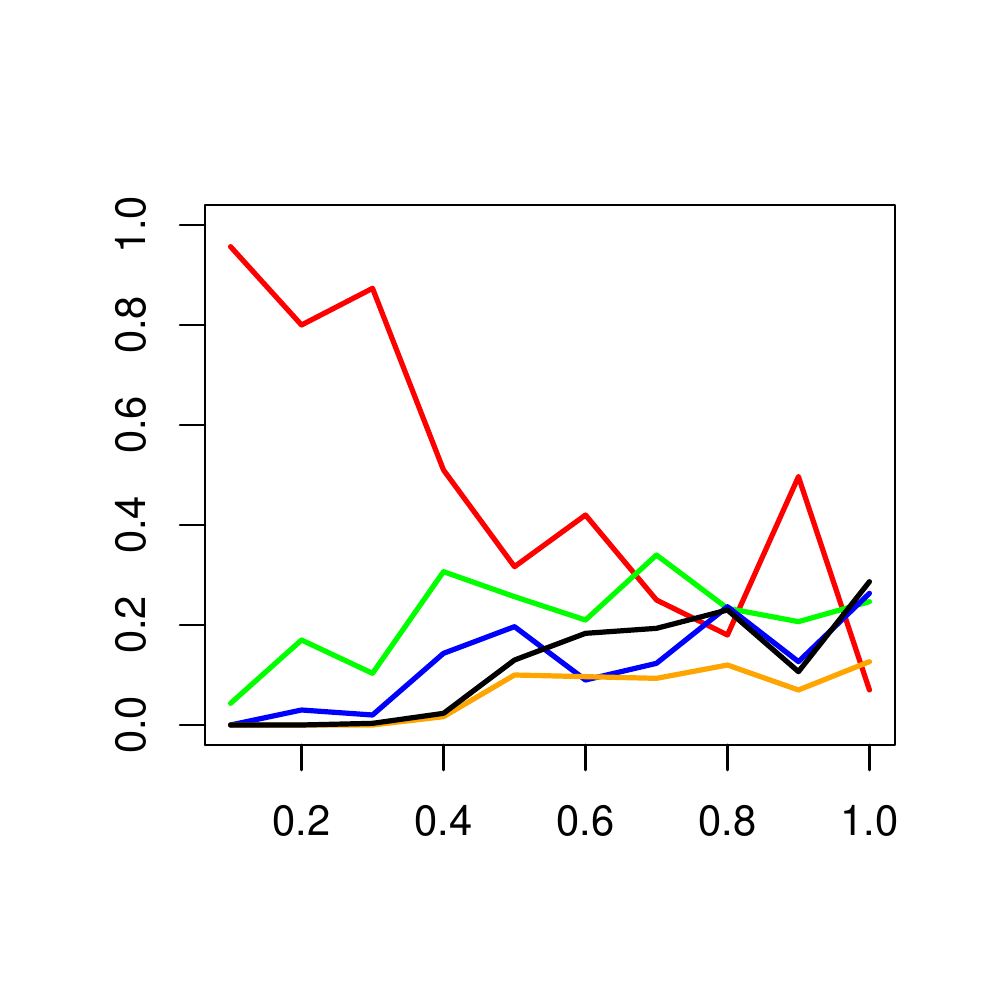}
		\vspace{-10mm} 
  	\caption{\scriptsize{Pickelhaube}} \label{fig:iX100Pi}
	\end{subfigure}
	
	\vspace{-2mm} 
		
	\begin{subfigure}[t]{.2\textwidth}
		\includegraphics[width=\textwidth]{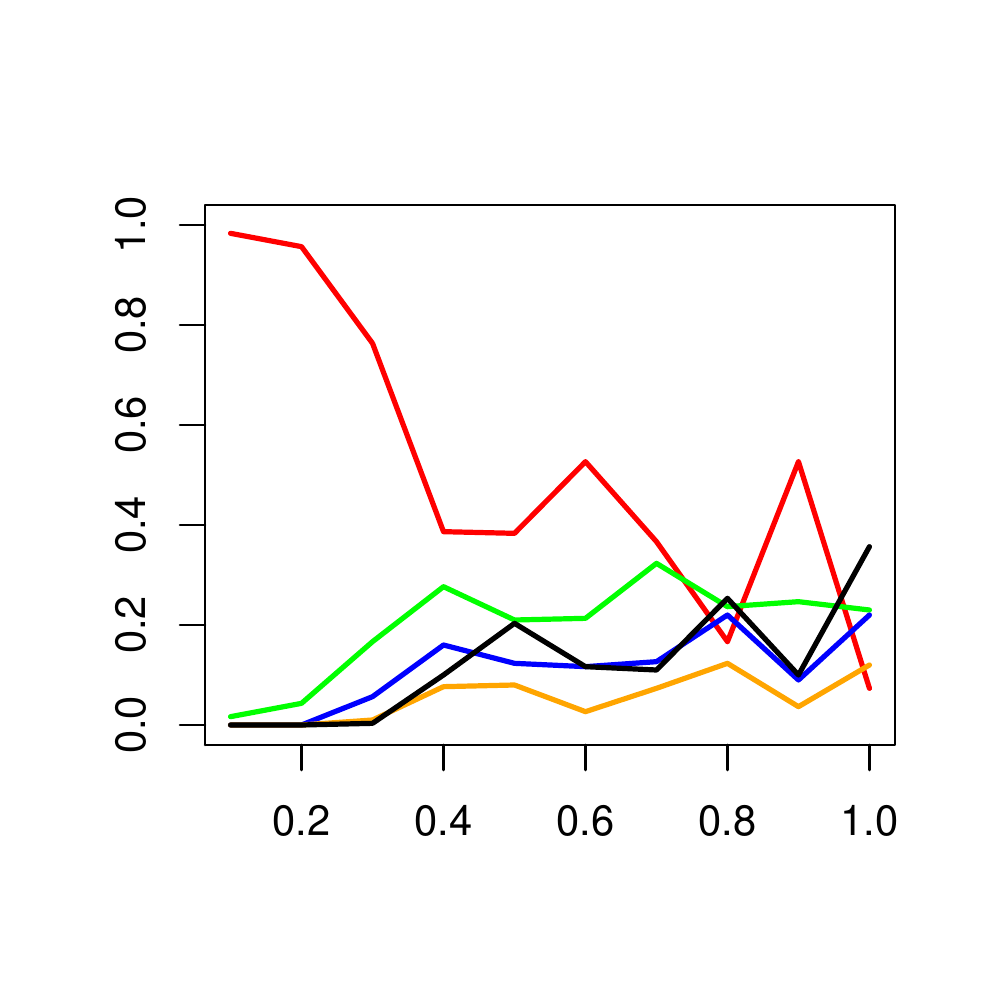}
		\vspace{-10mm} 
  	\caption{\scriptsize{Heaviside}} \label{fig:iX100Hv}
	\end{subfigure}%
	\hspace{-6mm} 
	\begin{subfigure}[t]{.2\textwidth}
		\includegraphics[width=\textwidth]{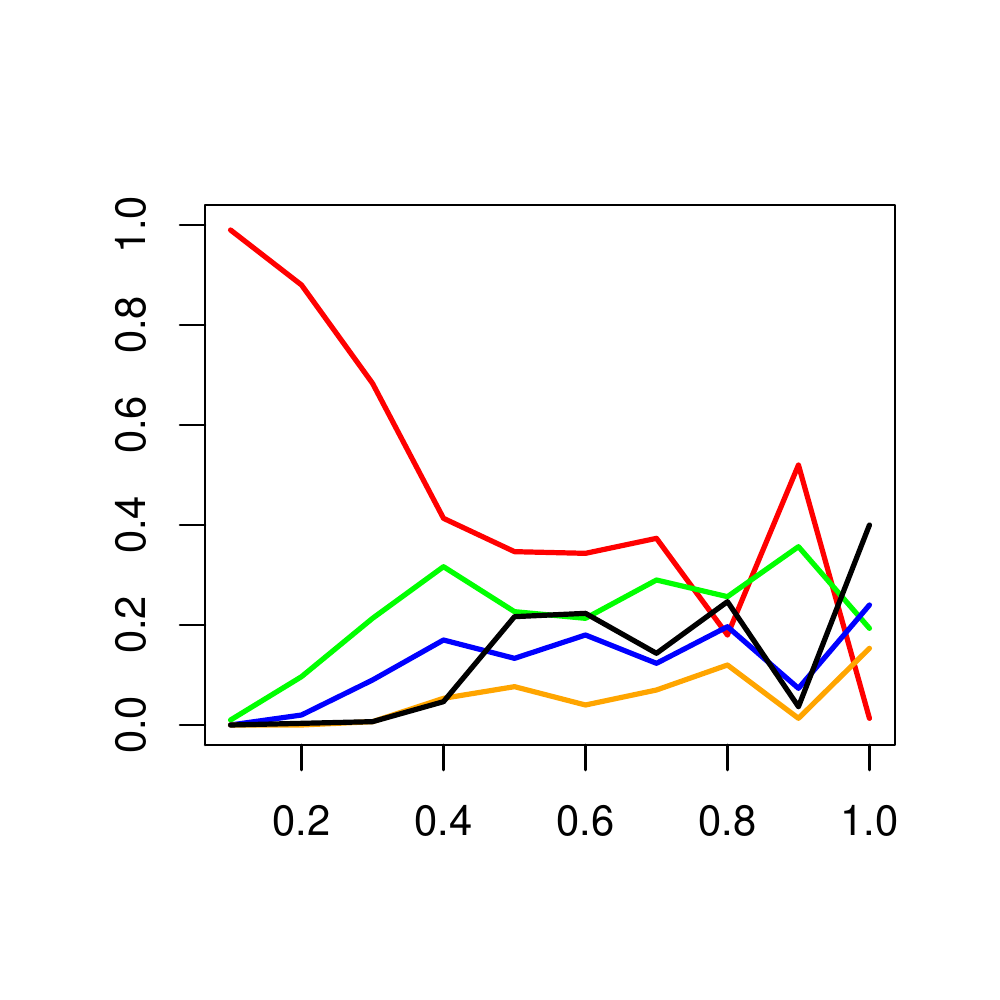}
		\vspace{-10mm} 
  	\caption{\scriptsize{Sawtooth}} \label{fig:iX100Sa}
	\end{subfigure}%
	\hspace{-6mm} 
	\begin{subfigure}[t]{.2\textwidth}
		\includegraphics[width=\textwidth]{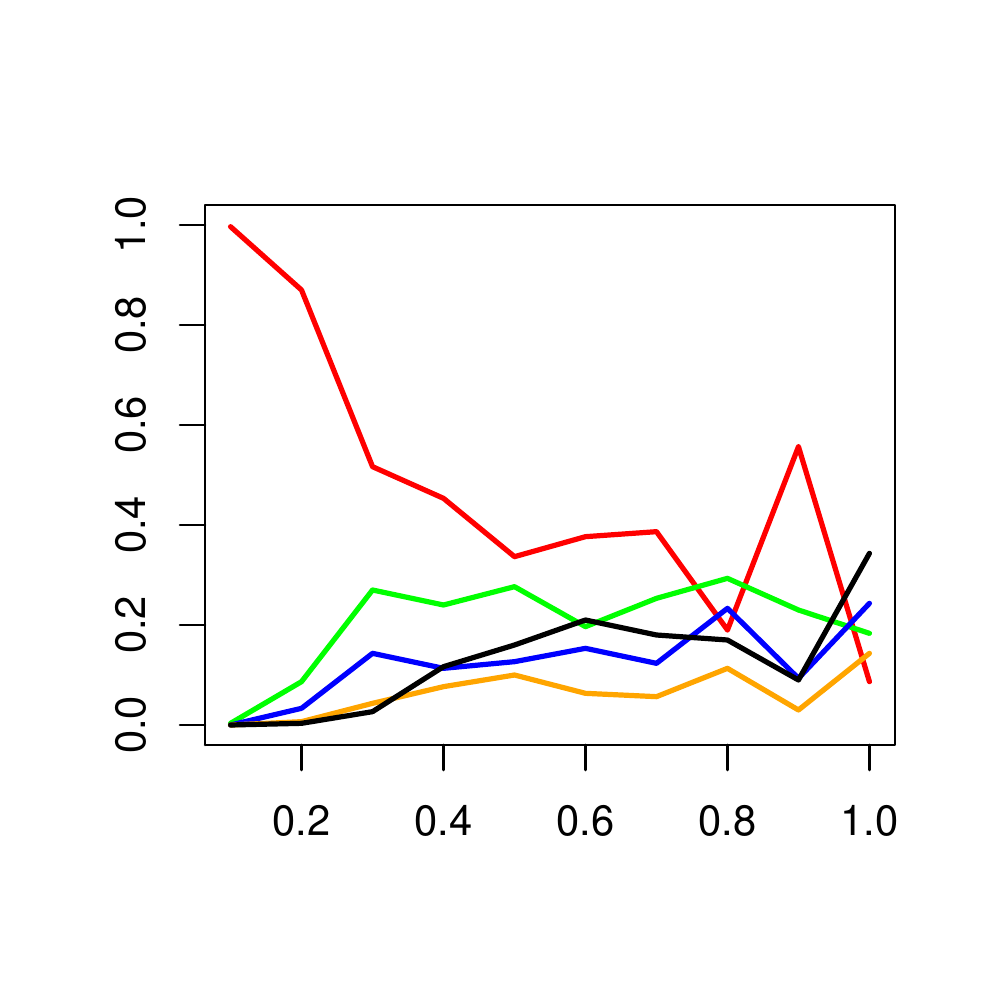}
		\vspace{-10mm} 
  	\caption{\scriptsize{Ackley}} \label{fig:iX100Ac}
	\end{subfigure}%
	\hspace{-6mm} 
	\begin{subfigure}[t]{.2\textwidth}
		\includegraphics[width=\textwidth]{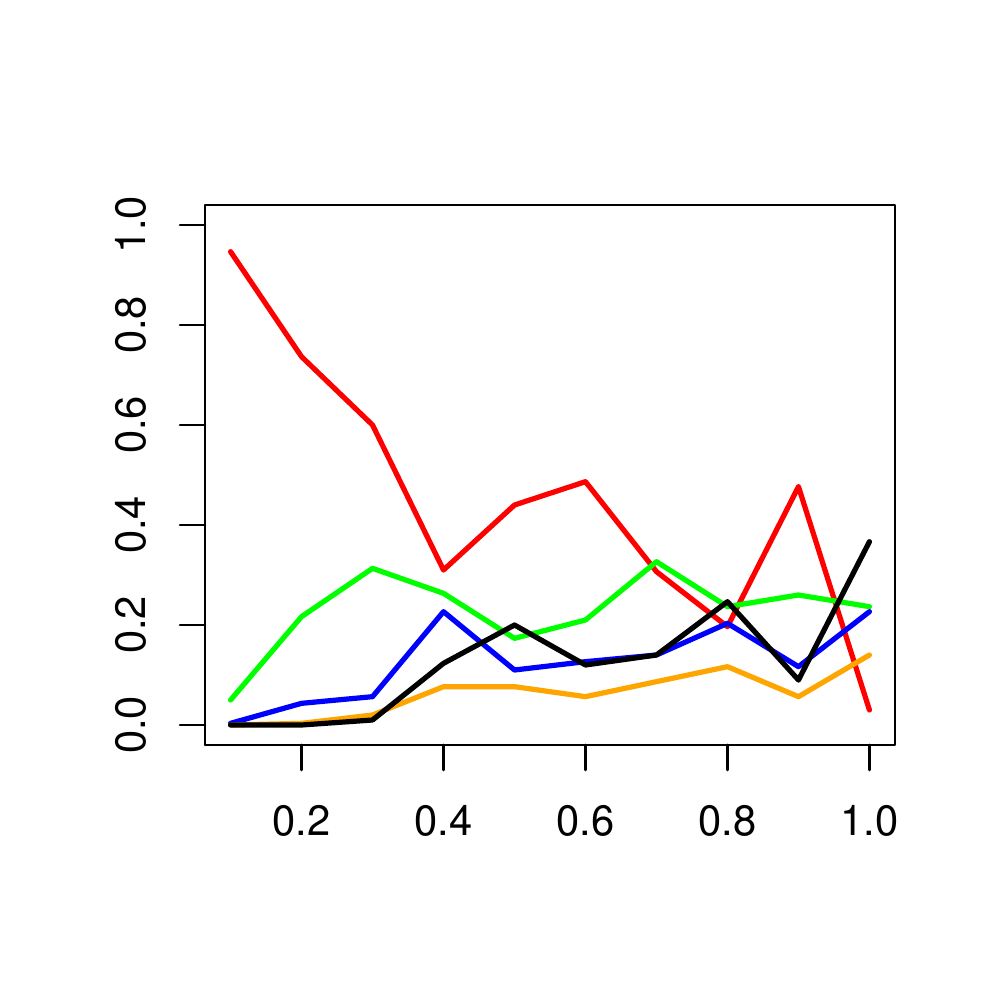}
		\vspace{-10mm} 
  	\caption{\scriptsize{Sphere}} \label{fig:iX100Sp}
	\end{subfigure}%
	\hspace{-6mm} 
	\begin{subfigure}[t]{.2\textwidth}
		\includegraphics[width=\textwidth]{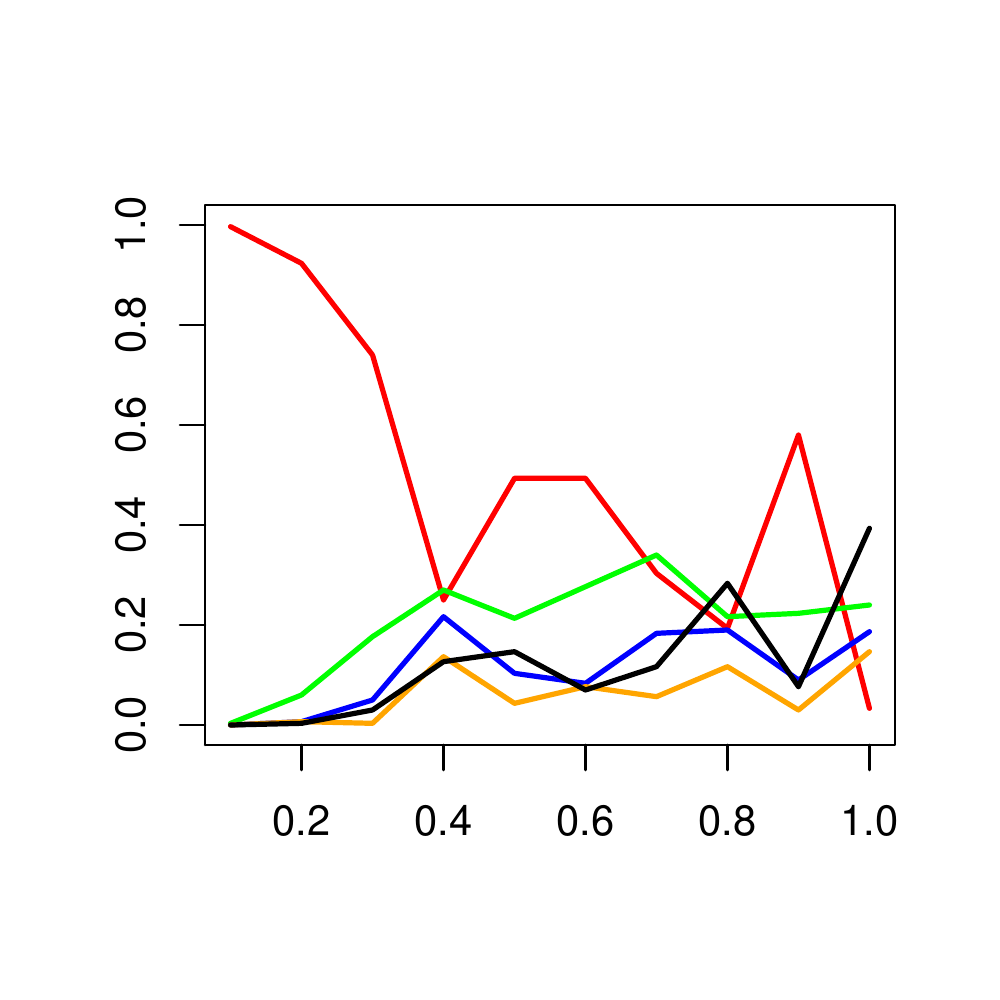}
		\vspace{-10mm} 
  	\caption{\scriptsize{Rosenbrock}} \label{fig:iX100Ro}
	\end{subfigure}
		
	\caption{Component -- decile breakdowns for the extent (size) of the inner maximisation search, across all GGGP heuristics at 100D. Components: [2 -- 10] (red), [11 -- 20] (green), [21 -- 30] (blue), [31 -- 40] (orange), \textgreater \ 40 (black).}
	\label{fig:inExt100Comp}
	
\end{figure}

\begin{figure}[H]
	\centering
	
	\vspace{-5mm} 

	\begin{subfigure}[t]{.25\textwidth}
		\includegraphics[width=\textwidth]{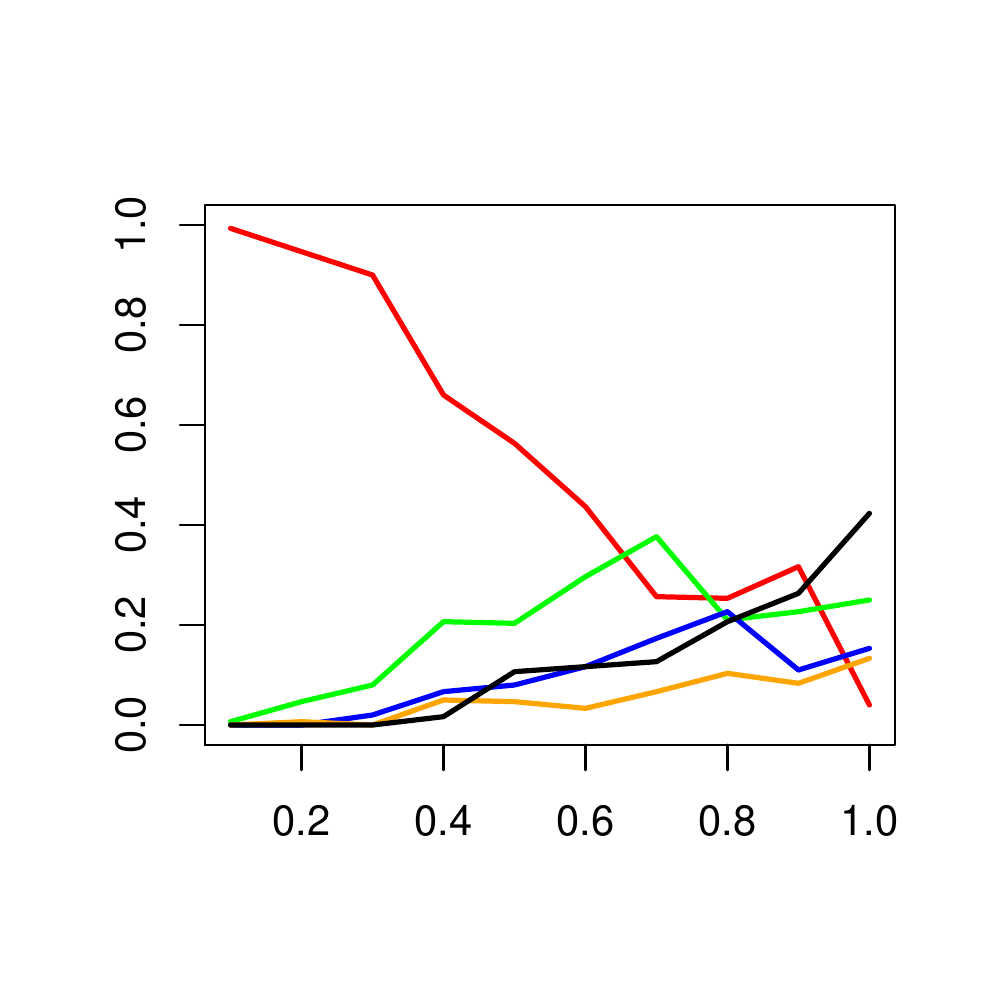}
		\vspace{-10mm} 
  	\caption{\scriptsize{30D}} \label{fig:iX30}
	\end{subfigure}%
	\hspace{-6mm} 
	\begin{subfigure}[t]{.25\textwidth}
		\includegraphics[width=\textwidth]{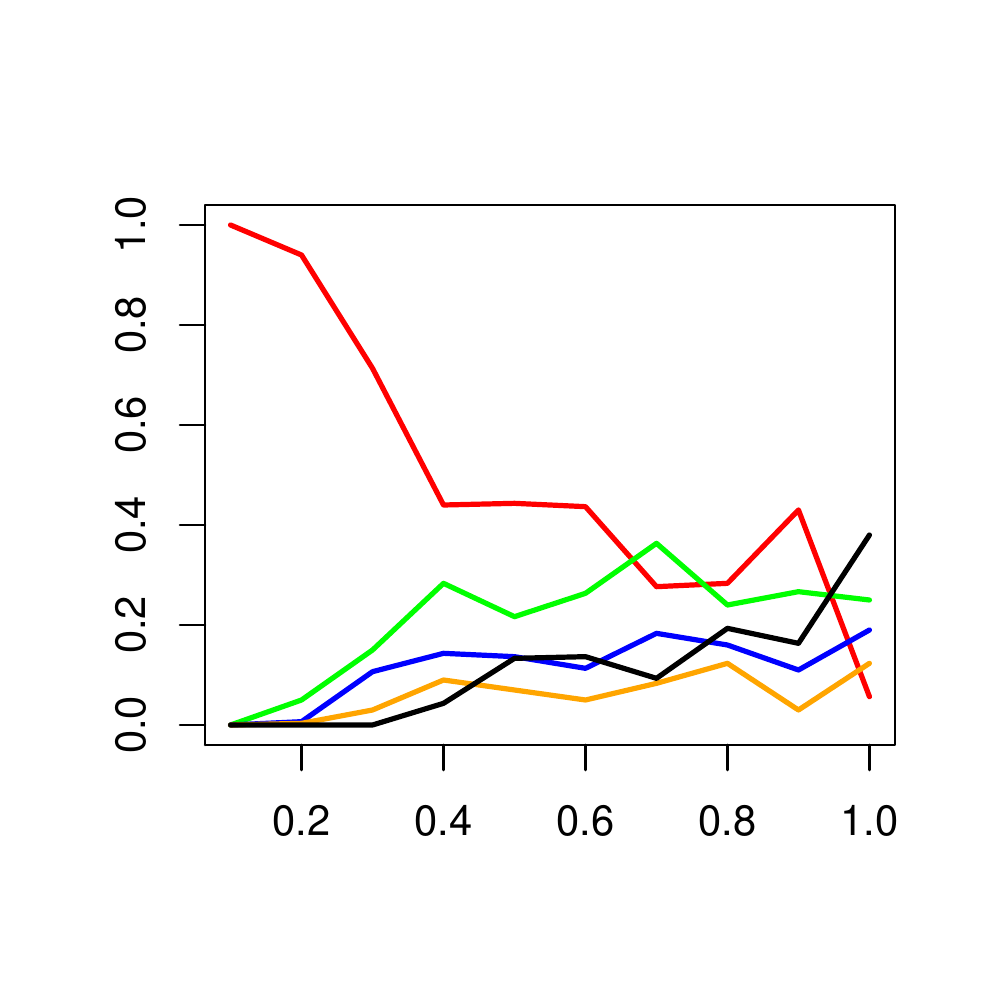}
		\vspace{-10mm} 
  	\caption{\scriptsize{100D}} \label{fig:iX100}
	\end{subfigure}

	\caption{Component -- decile breakdowns for the extent (size) of the inner maximisation search, across all GGGP heuristics at 30 D and 100D for the general heuristics. Components: [2 -- 10] (red), [11 -- 20] (green), [21 -- 30] (blue), [31 -- 40] (orange), \textgreater \ 40 (black).}
	\label{fig:inExtComp}
	
\end{figure}

In terms of the number of points evaluated in an inner maximisation, Table~\ref{fig:proportions} and Figures~\ref{fig:inExt30Comp} to~\ref{fig:inExtComp} show a clear dominance for low numbers, primarily in the range 2-10 (red). This is both for all heuristics generated and those performing best. Additional analysis shows that for the best performing third of heuristics, the proportions that employ random sampling using a number of points in the 2-10 range is 56\% and 49\% for 30D and 100D respectively.

\subsubsection{Baseline and extended movement capabilities}
\label{sec:compMove}

Table~\ref{fig:proportions} shows that the inertia formulation of the baseline particle velocity equation appears in 64\% of all heuristics for both 30D and 100D, whilst for the best performing third it dominates, appearing in over 80\% of heuristics. A decile level analysis confirms this dominance in the best performing heuristics. 

\begin{figure}[H]
	\centering
	

	\begin{subfigure}[t]{.2\textwidth}
		\includegraphics[width=\textwidth]{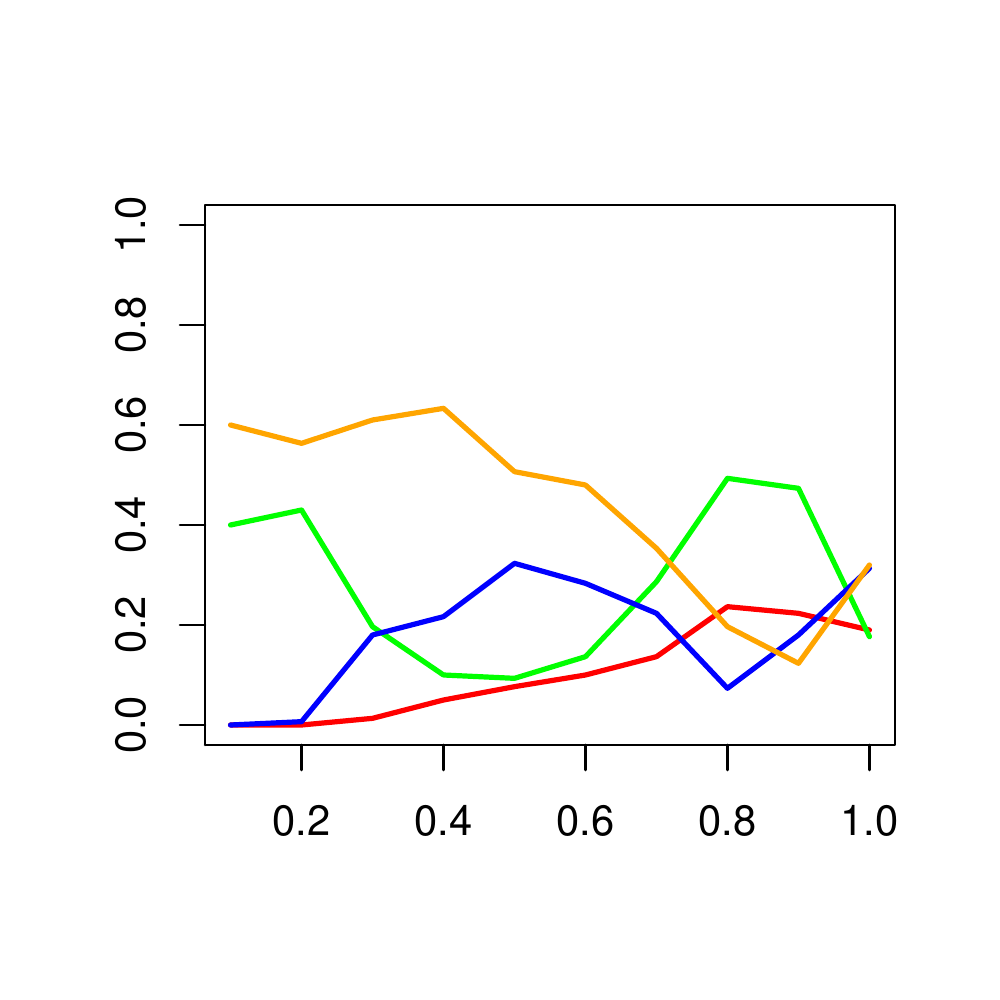}
		\vspace{-10mm} 
  	\caption{\scriptsize{Rastrigin}} \label{fig:mv30Ra}
	\end{subfigure}%
	\hspace{-6mm} 
	\begin{subfigure}[t]{.2\textwidth}
		\includegraphics[width=\textwidth]{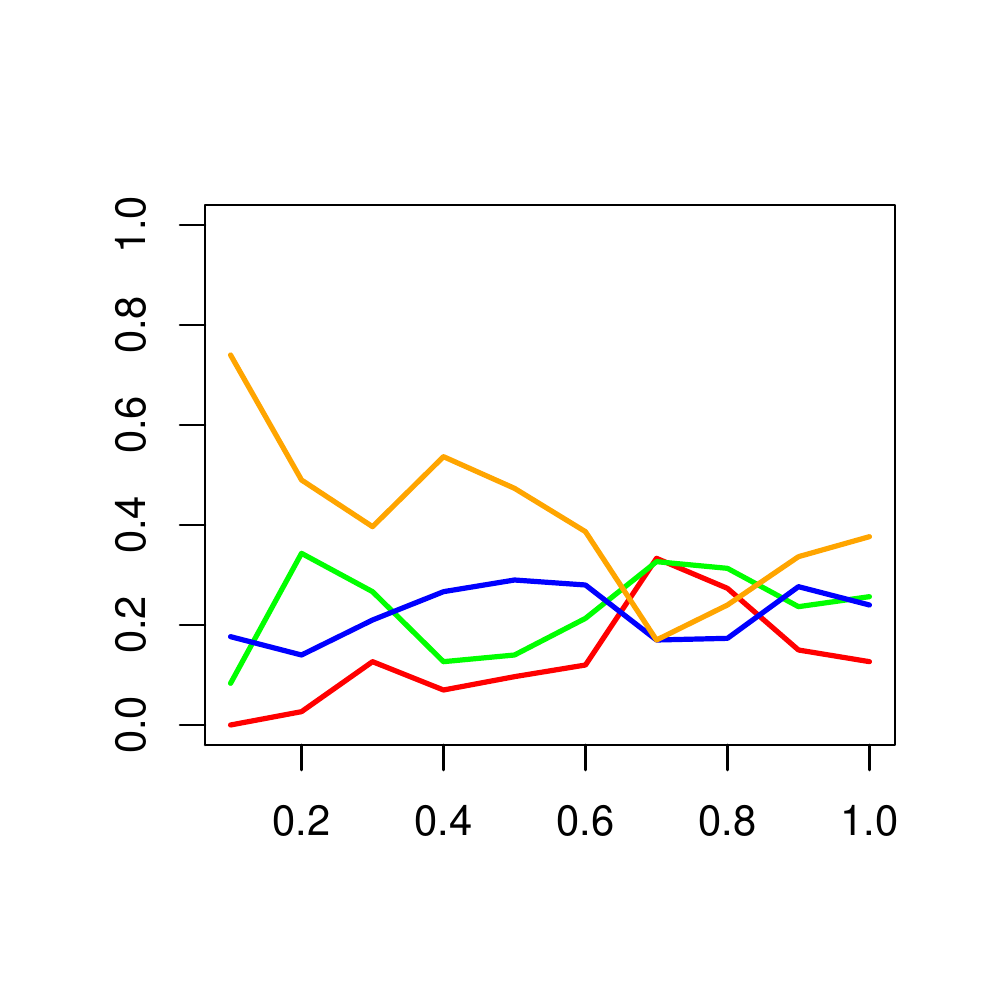}
		\vspace{-10mm} 
  	\caption{\scriptsize{Multipeak F1}} \label{fig:mv30M1}
	\end{subfigure}%
	\hspace{-6mm} 
	\begin{subfigure}[t]{.2\textwidth}
		\includegraphics[width=\textwidth]{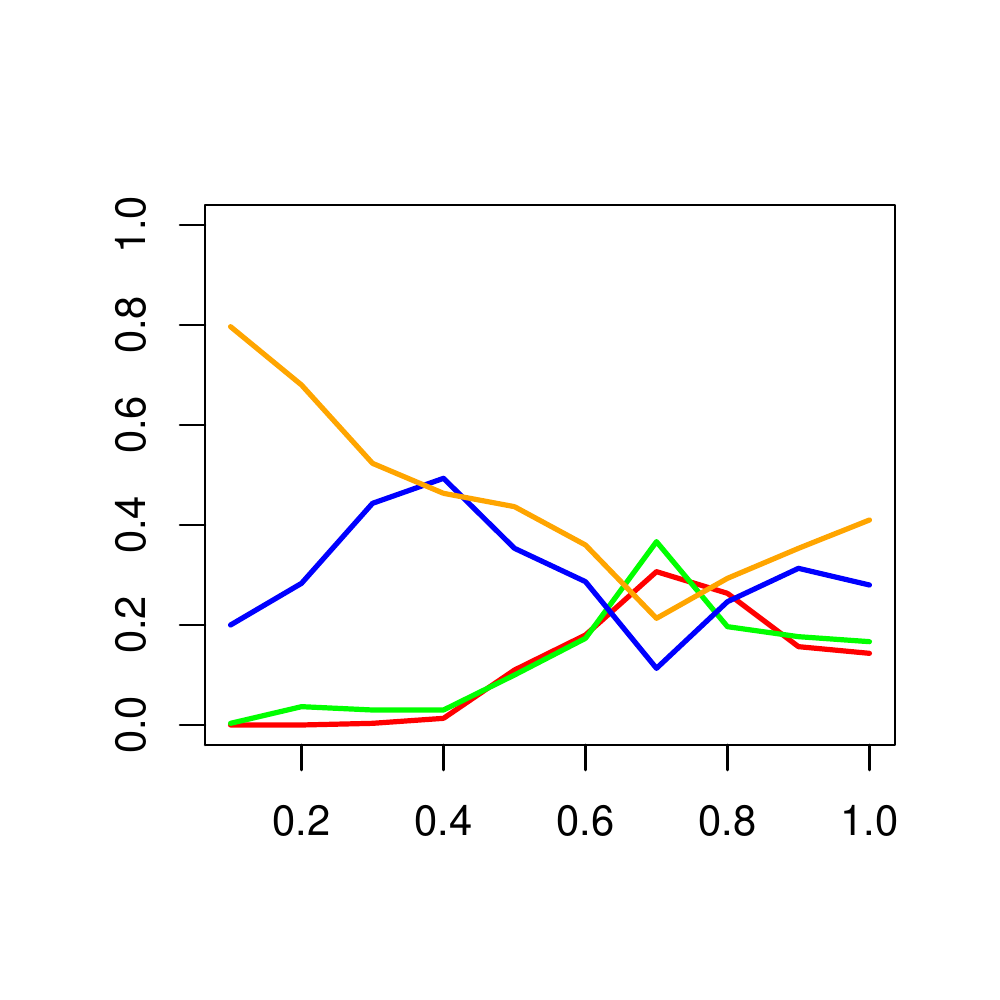}
		\vspace{-10mm} 
  	\caption{\scriptsize{Multipeak F2}} \label{fig:mv30M2}
	\end{subfigure}%
	\hspace{-6mm} 
	\begin{subfigure}[t]{.2\textwidth}
		\includegraphics[width=\textwidth]{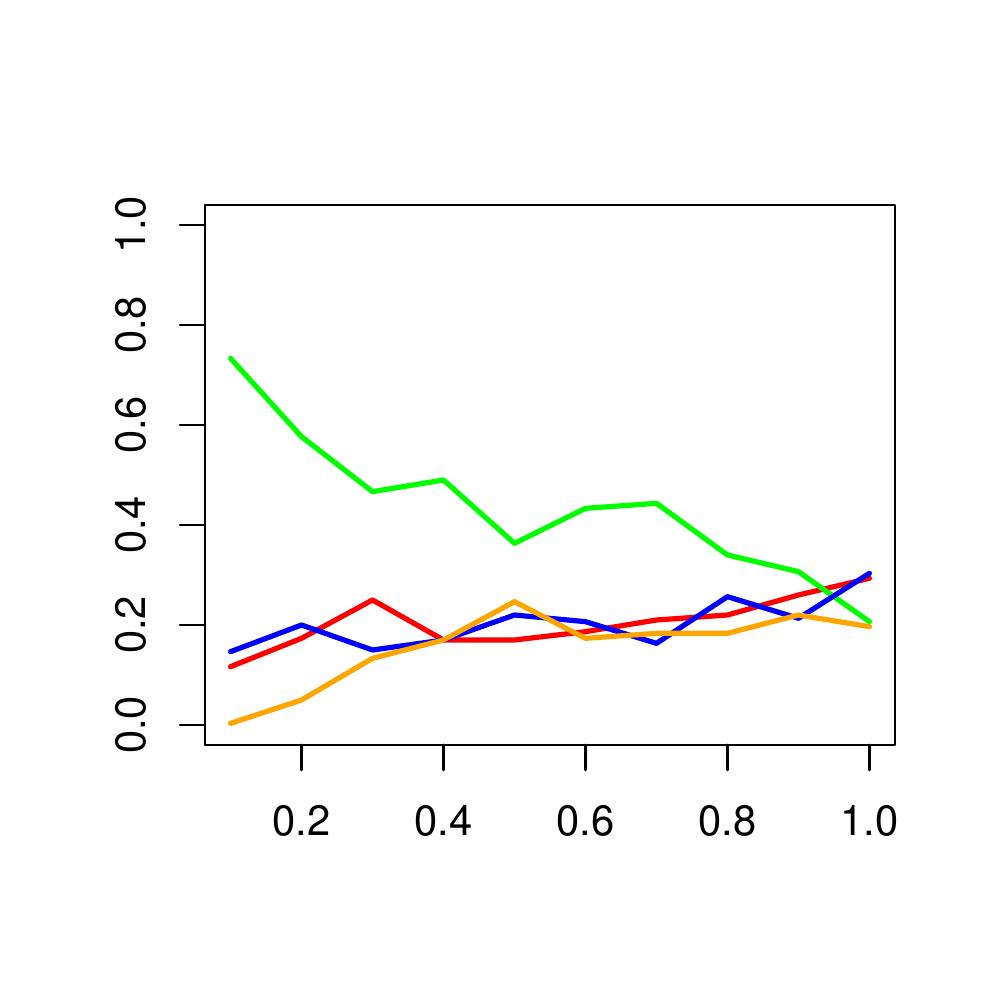}
		\vspace{-10mm} 
  	\caption{\scriptsize{Brankes}} \label{fig:mv30Br}
	\end{subfigure}%
	\hspace{-6mm} 
	\begin{subfigure}[t]{.2\textwidth}
		\includegraphics[width=\textwidth]{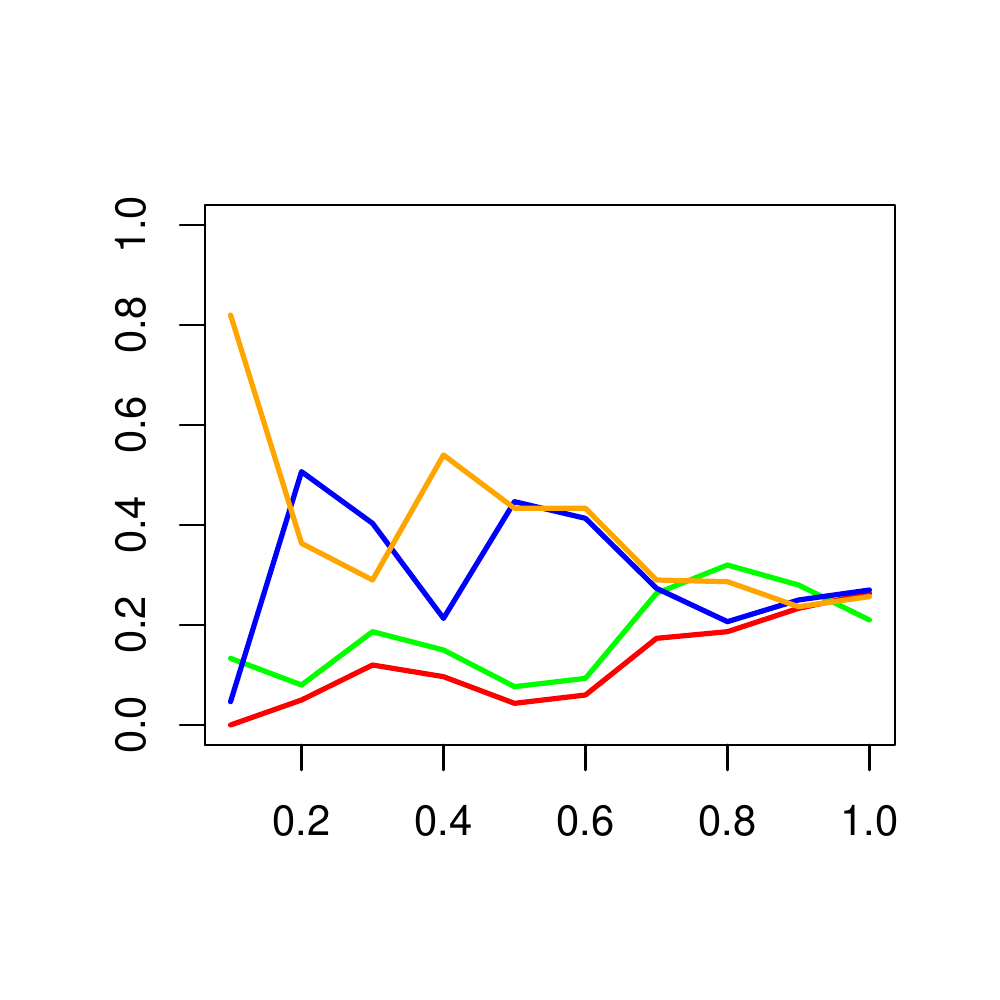}
		\vspace{-10mm} 
  	\caption{\scriptsize{Pickelhaube}} \label{fig:mv30Pi}
	\end{subfigure}
	
	\vspace{-2mm} 
		
	\begin{subfigure}[t]{.2\textwidth}
		\includegraphics[width=\textwidth]{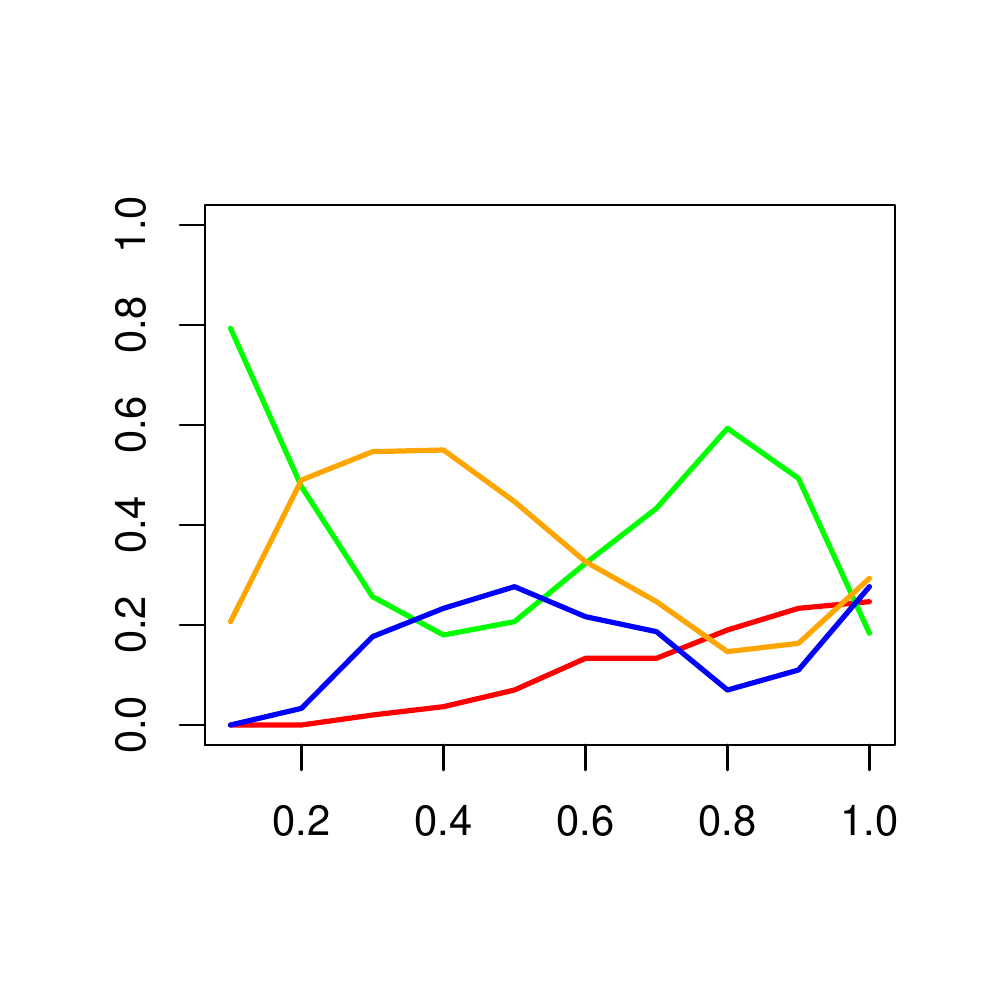}
		\vspace{-10mm} 
  	\caption{\scriptsize{Heaviside}} \label{fig:mv30Hv}
	\end{subfigure}%
	\hspace{-6mm} 
	\begin{subfigure}[t]{.2\textwidth}
		\includegraphics[width=\textwidth]{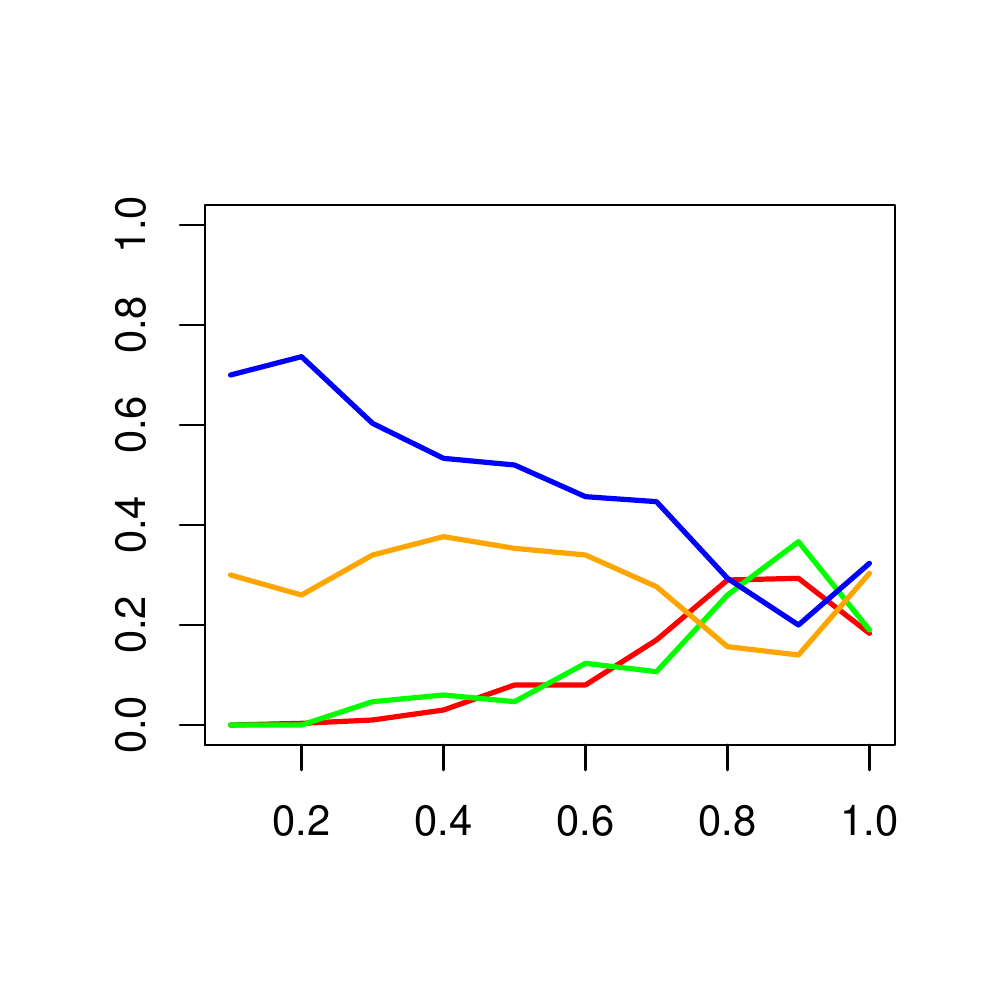}
		\vspace{-10mm} 
  	\caption{\scriptsize{Sawtooth}} \label{fig:mv30Sa}
	\end{subfigure}%
	\hspace{-6mm} 
	\begin{subfigure}[t]{.2\textwidth}
		\includegraphics[width=\textwidth]{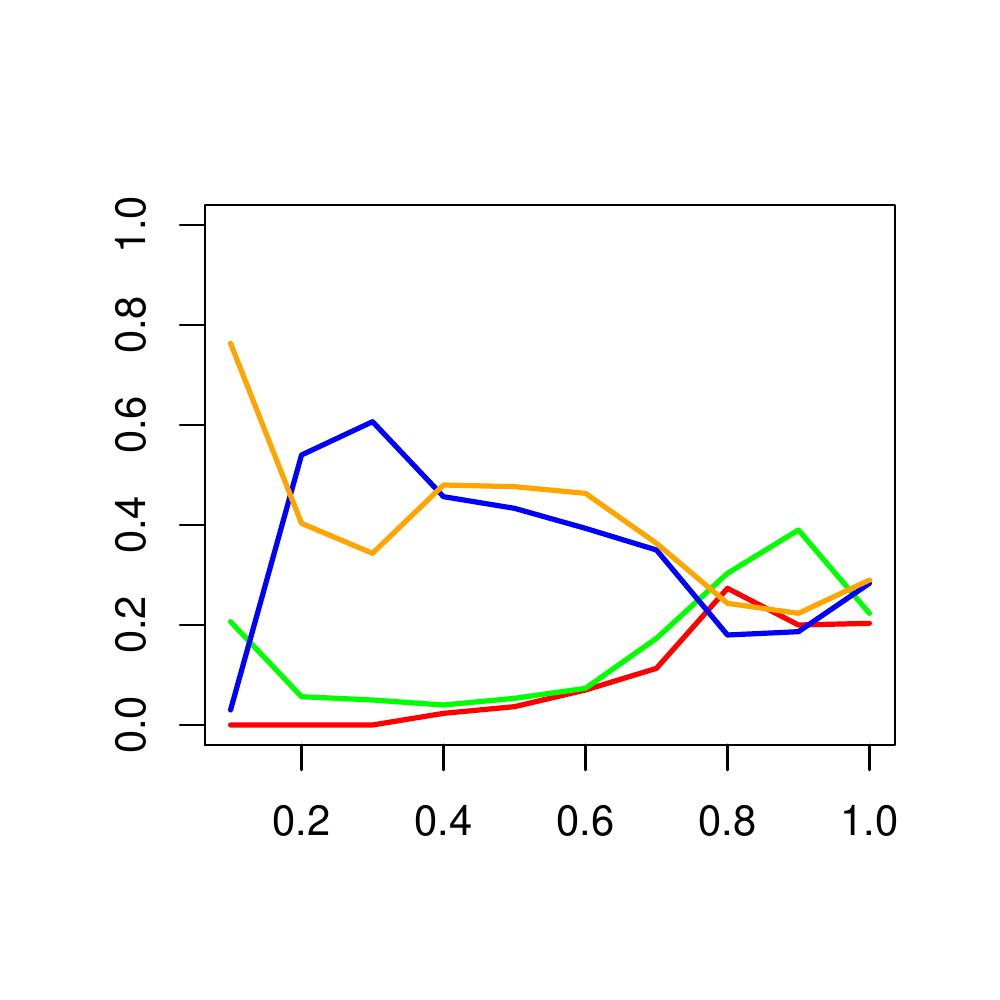}
		\vspace{-10mm} 
  	\caption{\scriptsize{Ackley}} \label{fig:mv30Ac}
	\end{subfigure}%
	\hspace{-6mm} 
	\begin{subfigure}[t]{.2\textwidth}
		\includegraphics[width=\textwidth]{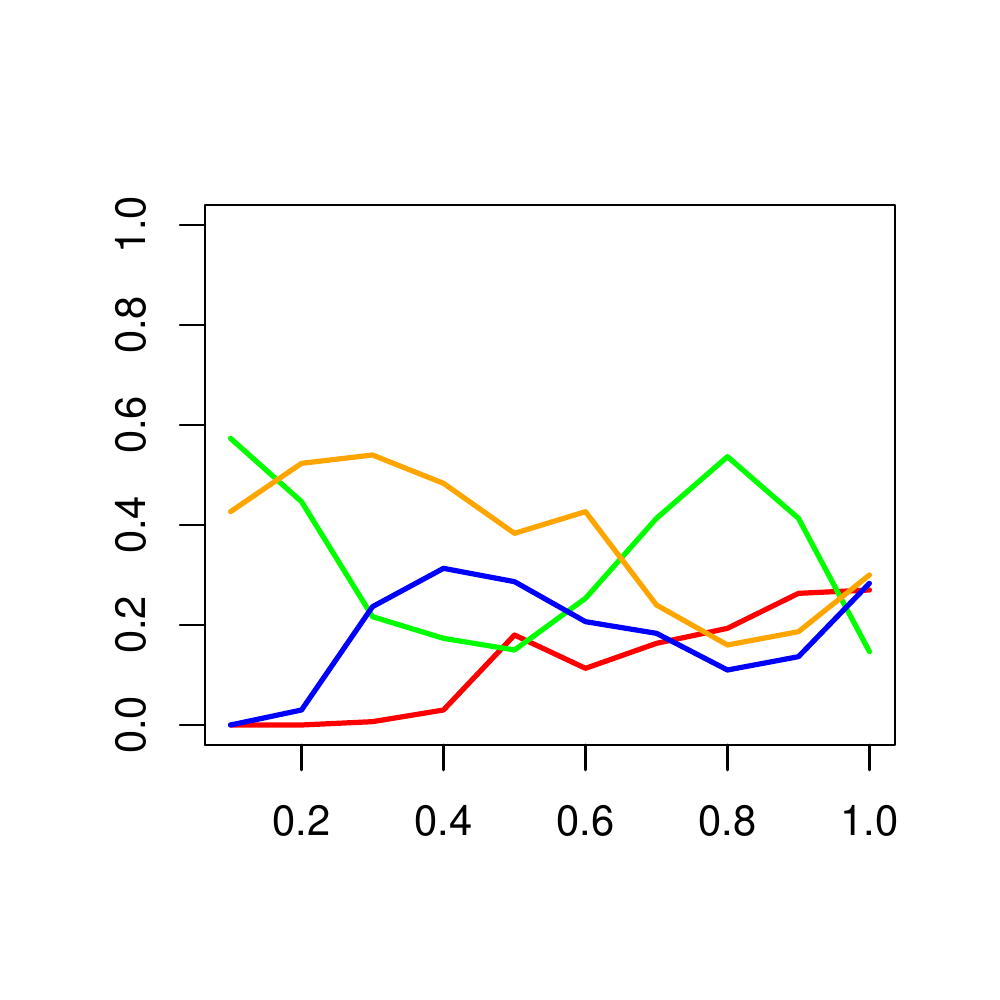}
		\vspace{-10mm} 
  	\caption{\scriptsize{Sphere}} \label{fig:mv30Sp}
	\end{subfigure}%
	\hspace{-6mm} 
	\begin{subfigure}[t]{.2\textwidth}
		\includegraphics[width=\textwidth]{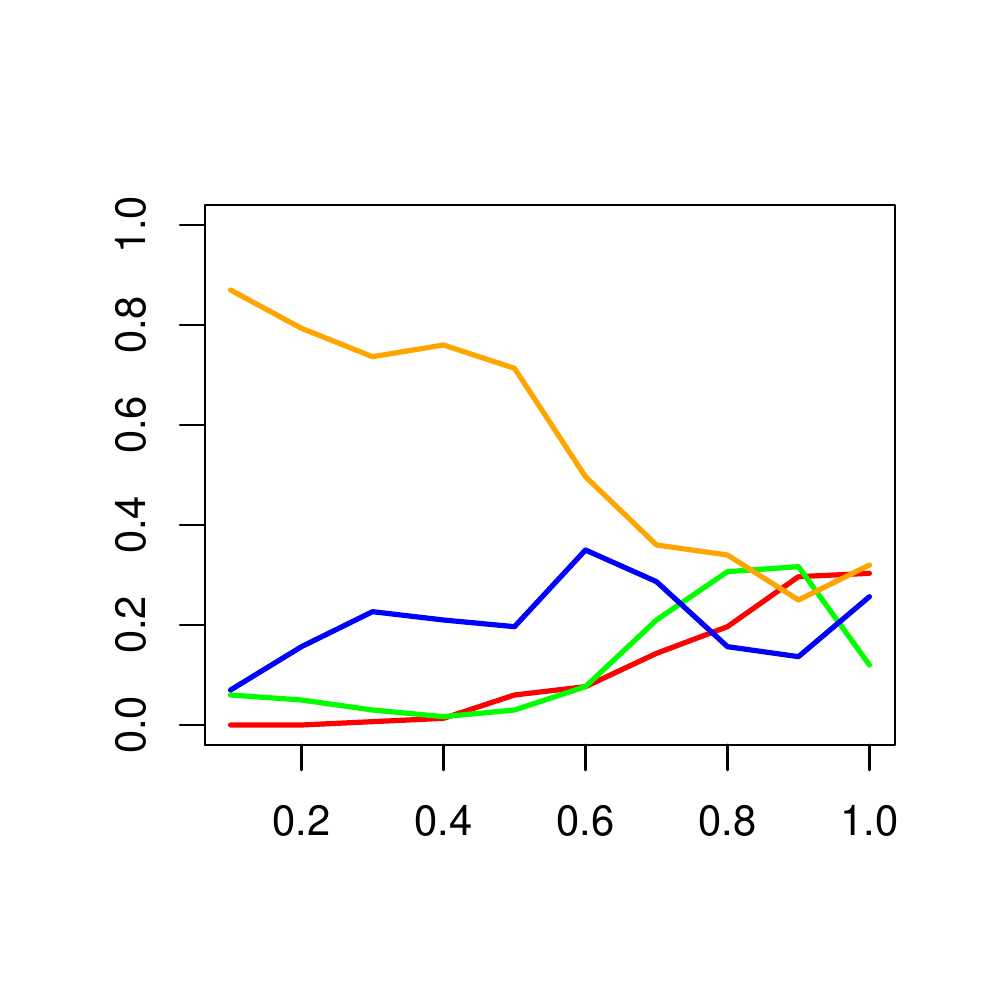}
		\vspace{-10mm} 
  	\caption{\scriptsize{Rosenbrock}} \label{fig:mv30Ro}
	\end{subfigure}
		
	\caption{Component -- decile breakdowns for the form of movement capability, across all GGGP heuristics at 30D. Components: Baseline (red), DD (green), LEH (blue), LEH+DD (orange).}
	\label{fig:move30Comp}

\end{figure}

\begin{figure}[H]
	\centering
	
	\vspace{-5mm} 

	\begin{subfigure}[t]{.2\textwidth}
		\includegraphics[width=\textwidth]{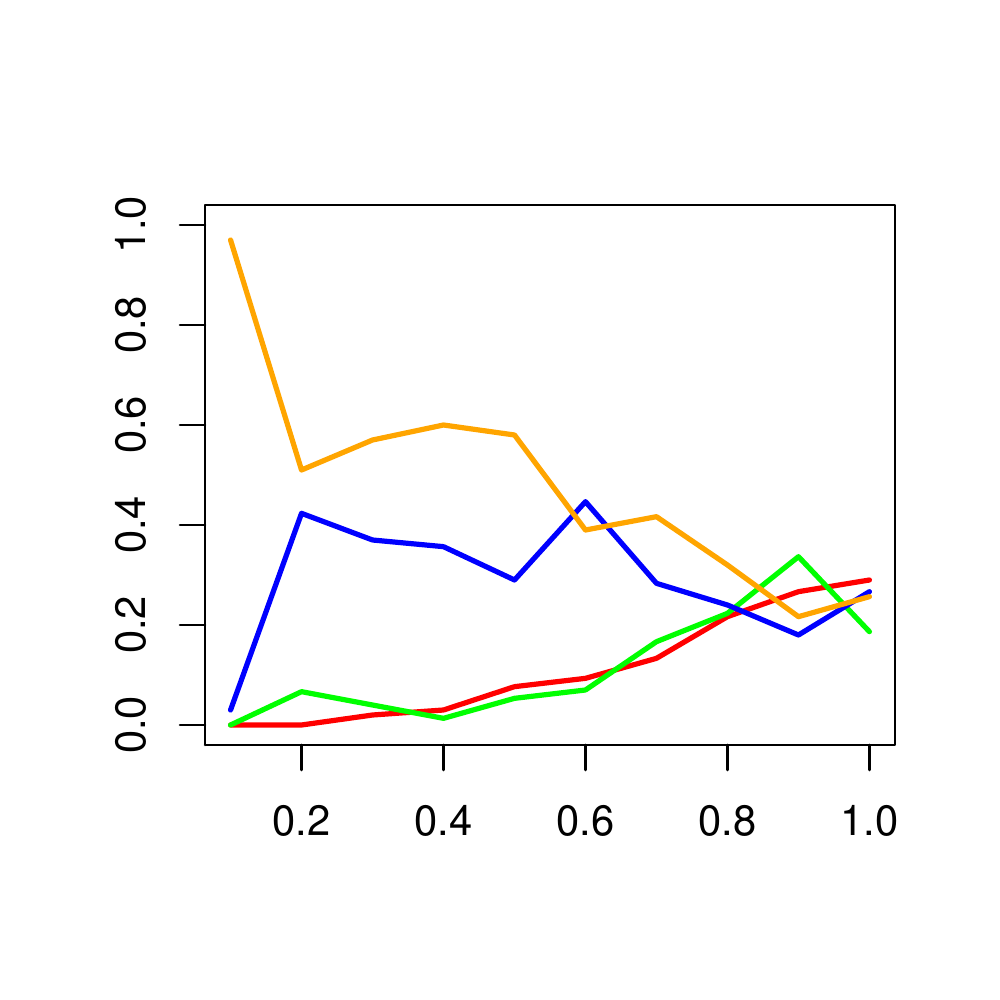}
		\vspace{-10mm} 
  	\caption{\scriptsize{Rastrigin}} \label{fig:mv100Ra}
	\end{subfigure}%
	\hspace{-6mm} 
	\begin{subfigure}[t]{.2\textwidth}
		\includegraphics[width=\textwidth]{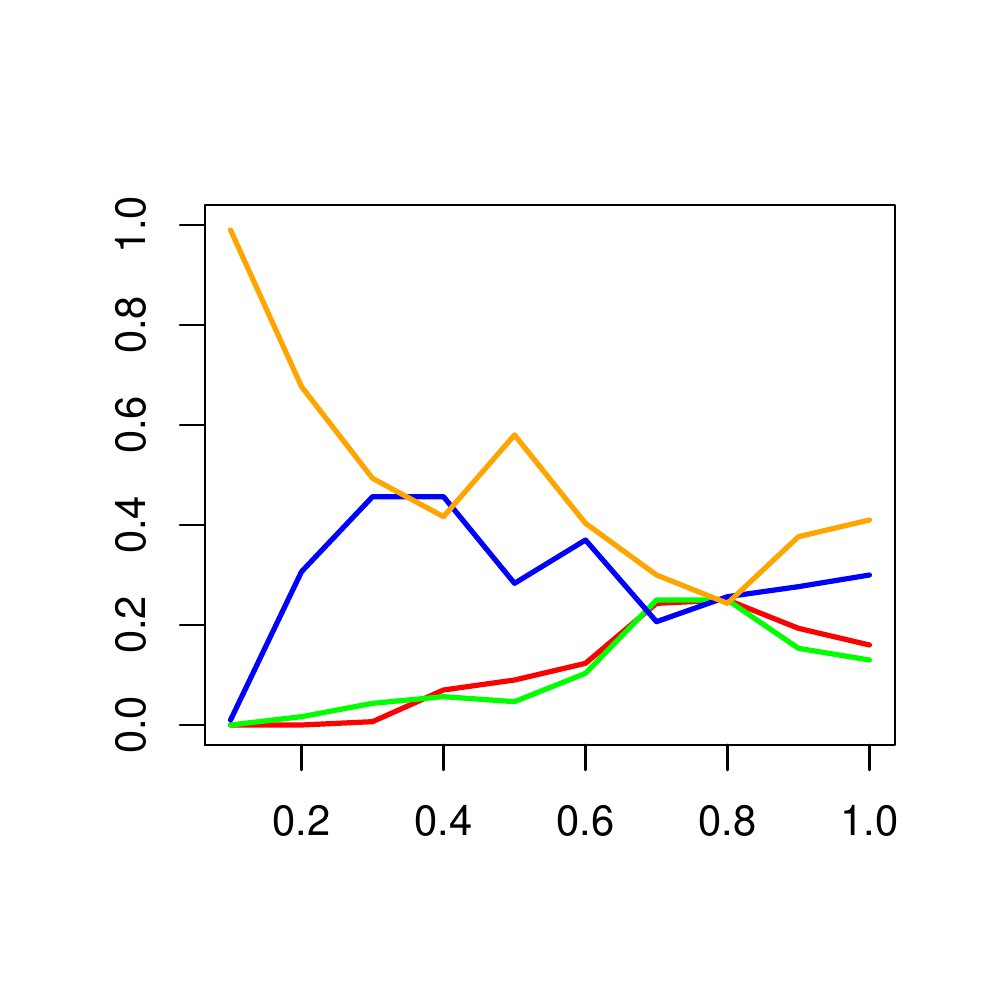}
		\vspace{-10mm} 
  	\caption{\scriptsize{Multipeak F1}} \label{fig:mv100M1}
	\end{subfigure}%
	\hspace{-6mm} 
	\begin{subfigure}[t]{.2\textwidth}
		\includegraphics[width=\textwidth]{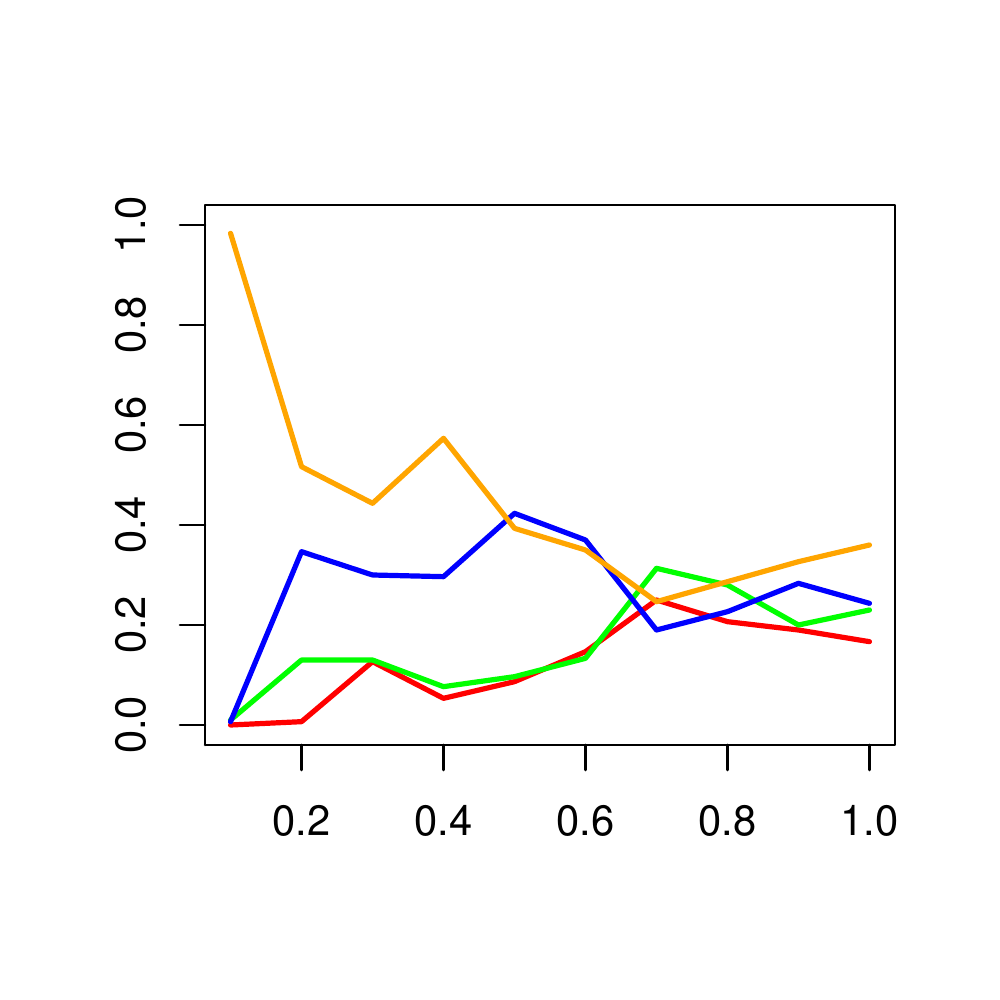}
		\vspace{-10mm} 
  	\caption{\scriptsize{Multipeak F2}} \label{fig:mv100M2}
	\end{subfigure}%
	\hspace{-6mm} 
	\begin{subfigure}[t]{.2\textwidth}
		\includegraphics[width=\textwidth]{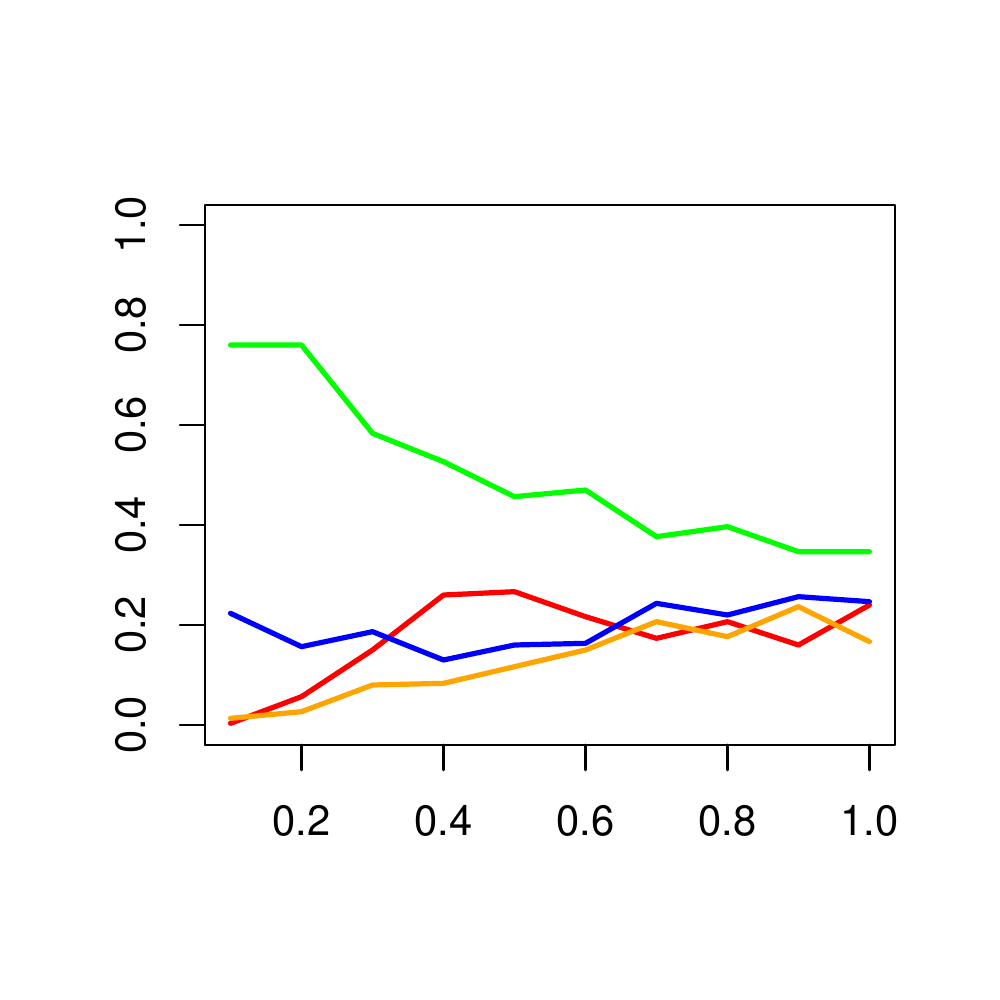}
		\vspace{-10mm} 
  	\caption{\scriptsize{Brankes}} \label{fig:mv100Br}
	\end{subfigure}%
	\hspace{-6mm} 
	\begin{subfigure}[t]{.2\textwidth}
		\includegraphics[width=\textwidth]{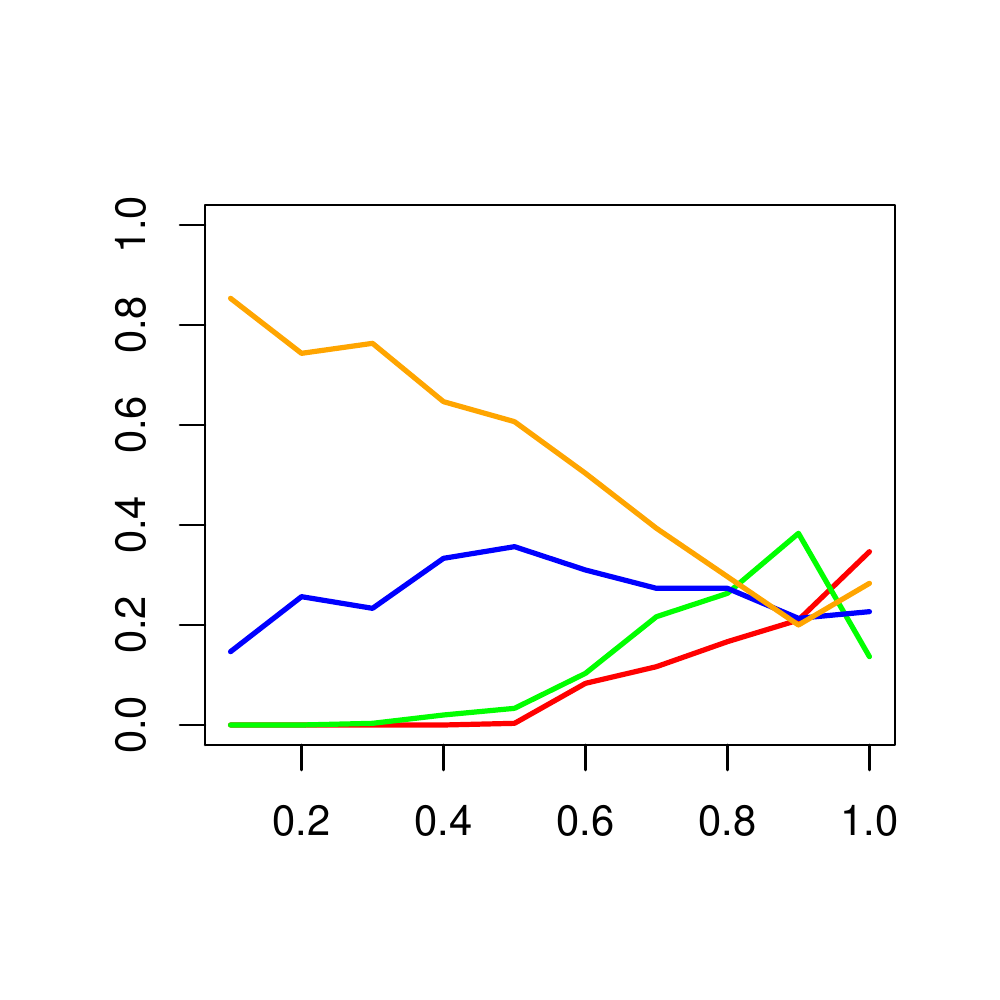}
		\vspace{-10mm} 
  	\caption{\scriptsize{Pickelhaube}} \label{fig:mv100Pi}
	\end{subfigure}
	
	\vspace{-2mm} 
		
	\begin{subfigure}[t]{.2\textwidth}
		\includegraphics[width=\textwidth]{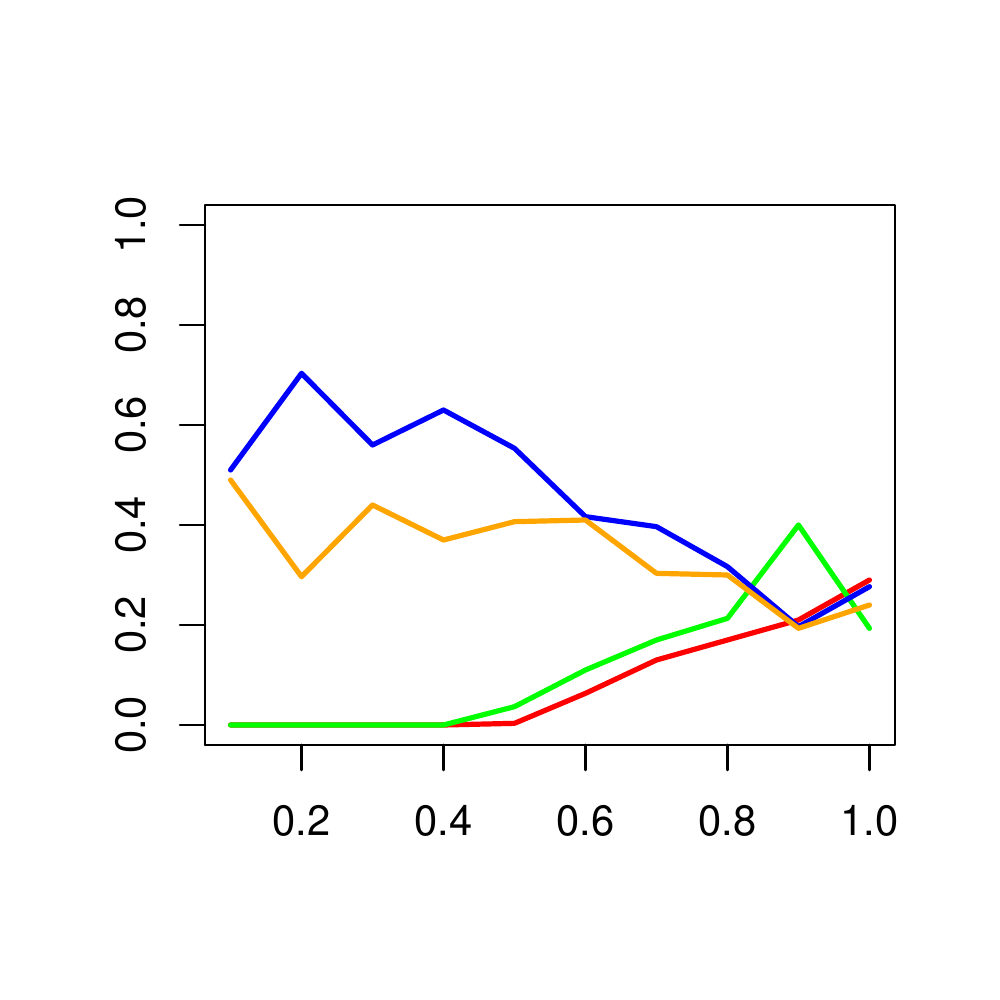}
		\vspace{-10mm} 
  	\caption{\scriptsize{Heaviside}} \label{fig:mv100Hv}
	\end{subfigure}%
	\hspace{-6mm} 
	\begin{subfigure}[t]{.2\textwidth}
		\includegraphics[width=\textwidth]{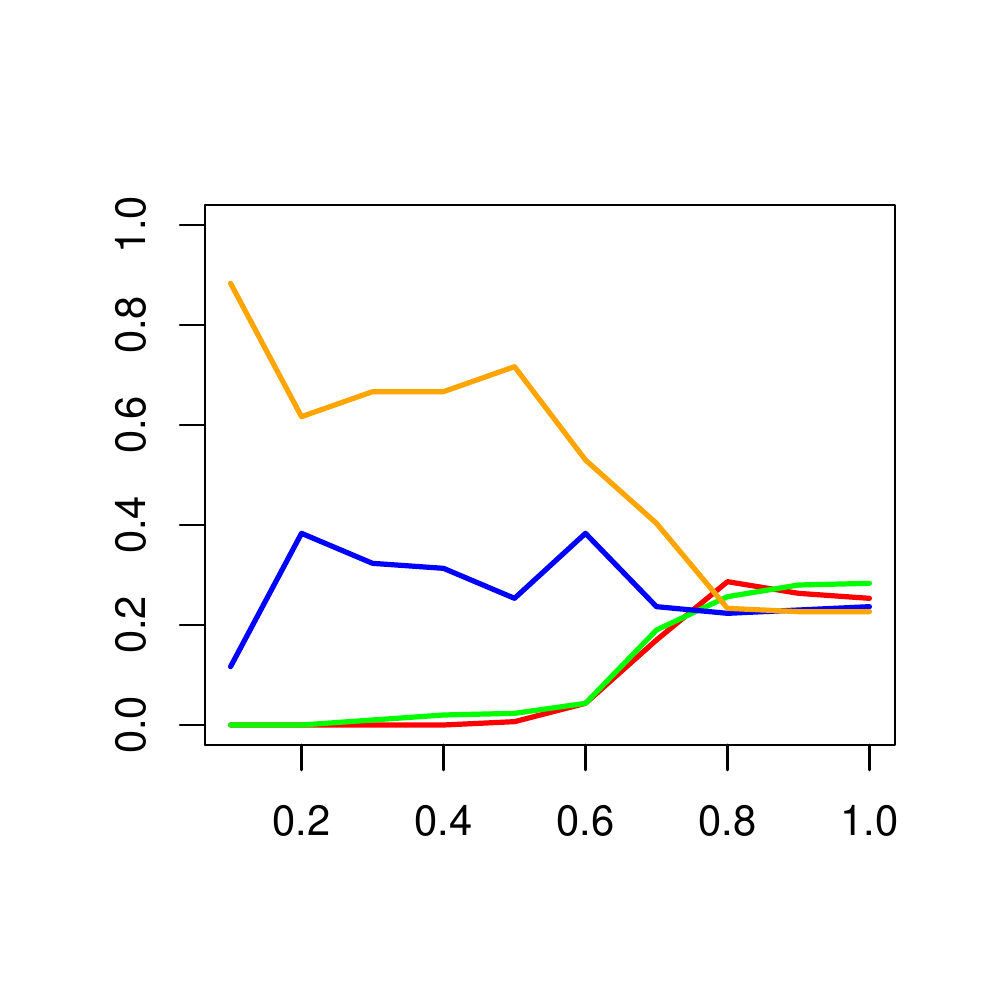}
		\vspace{-10mm} 
  	\caption{\scriptsize{Sawtooth}} \label{fig:mv100Sa}
	\end{subfigure}%
	\hspace{-6mm} 
	\begin{subfigure}[t]{.2\textwidth}
		\includegraphics[width=\textwidth]{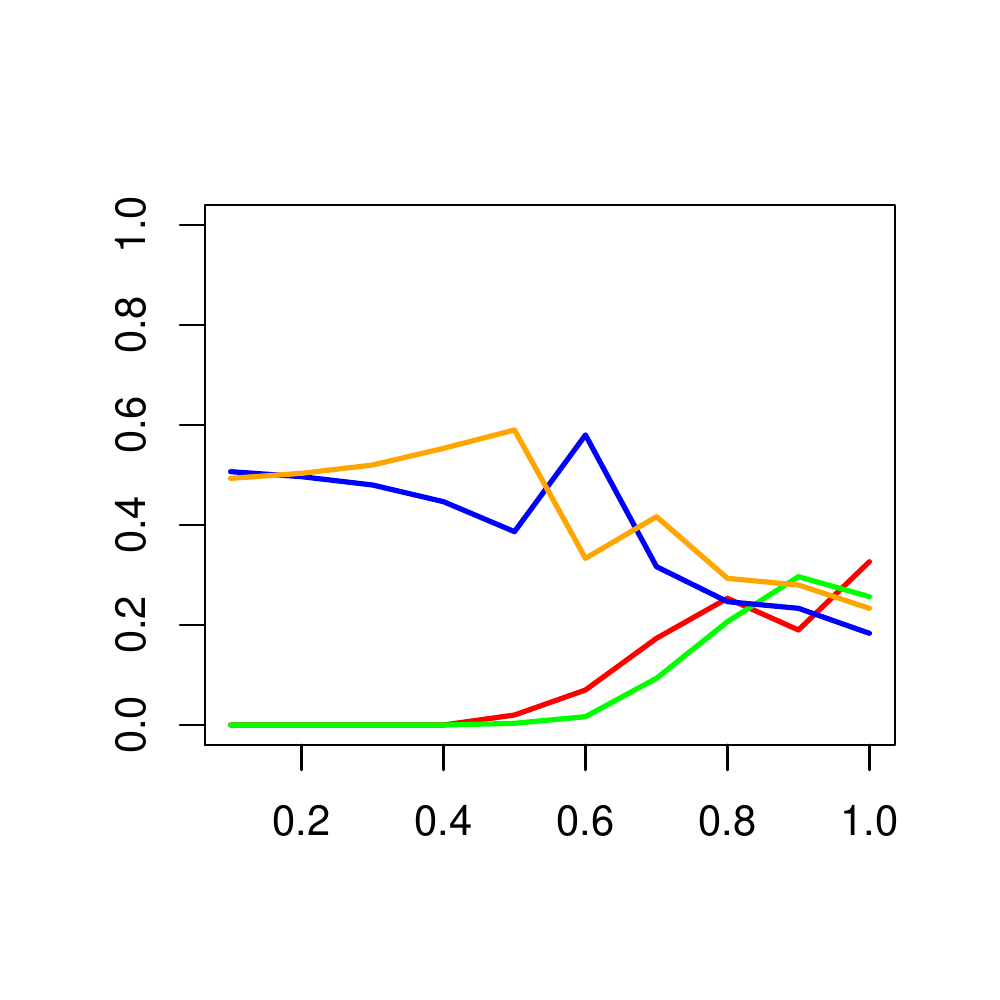}
		\vspace{-10mm} 
  	\caption{\scriptsize{Ackley}} \label{fig:mv100Ac}
	\end{subfigure}%
	\hspace{-6mm} 
	\begin{subfigure}[t]{.2\textwidth}
		\includegraphics[width=\textwidth]{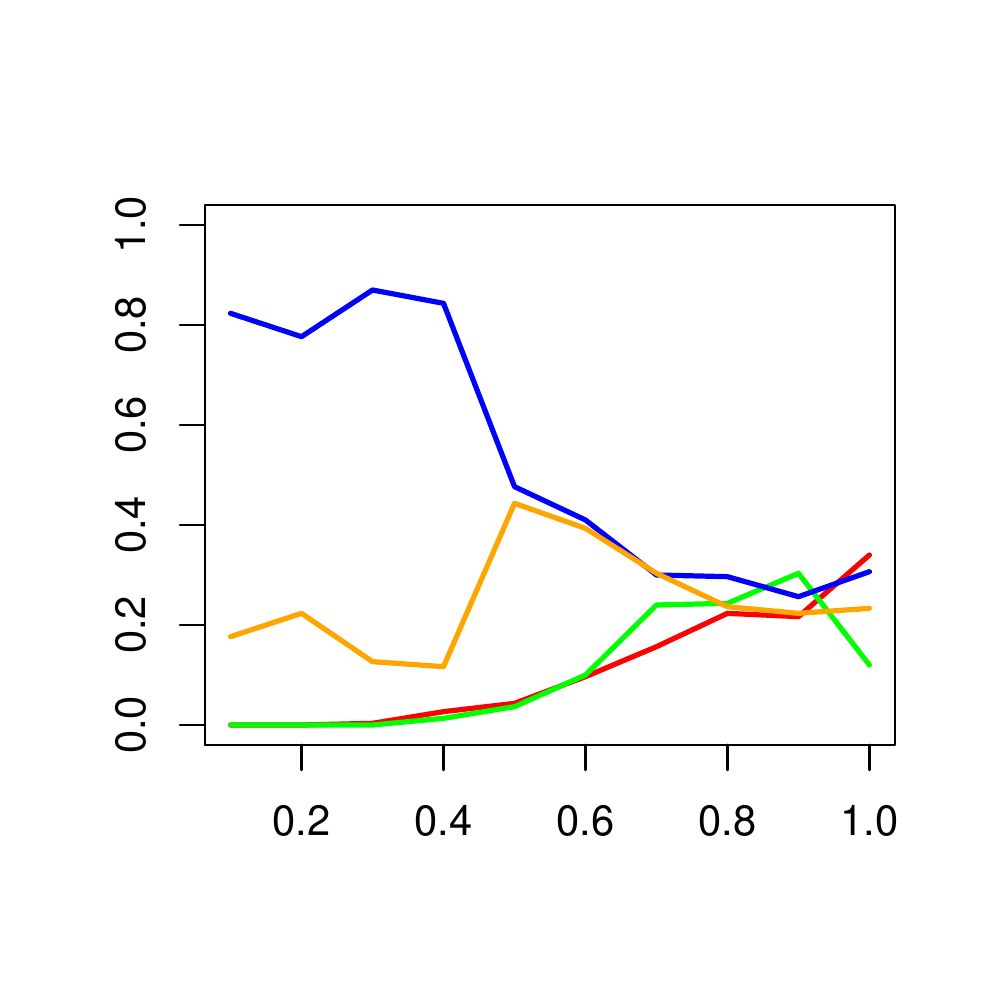}
		\vspace{-10mm} 
  	\caption{\scriptsize{Sphere}} \label{fig:mv100Sp}
	\end{subfigure}%
	\hspace{-6mm} 
	\begin{subfigure}[t]{.2\textwidth}
		\includegraphics[width=\textwidth]{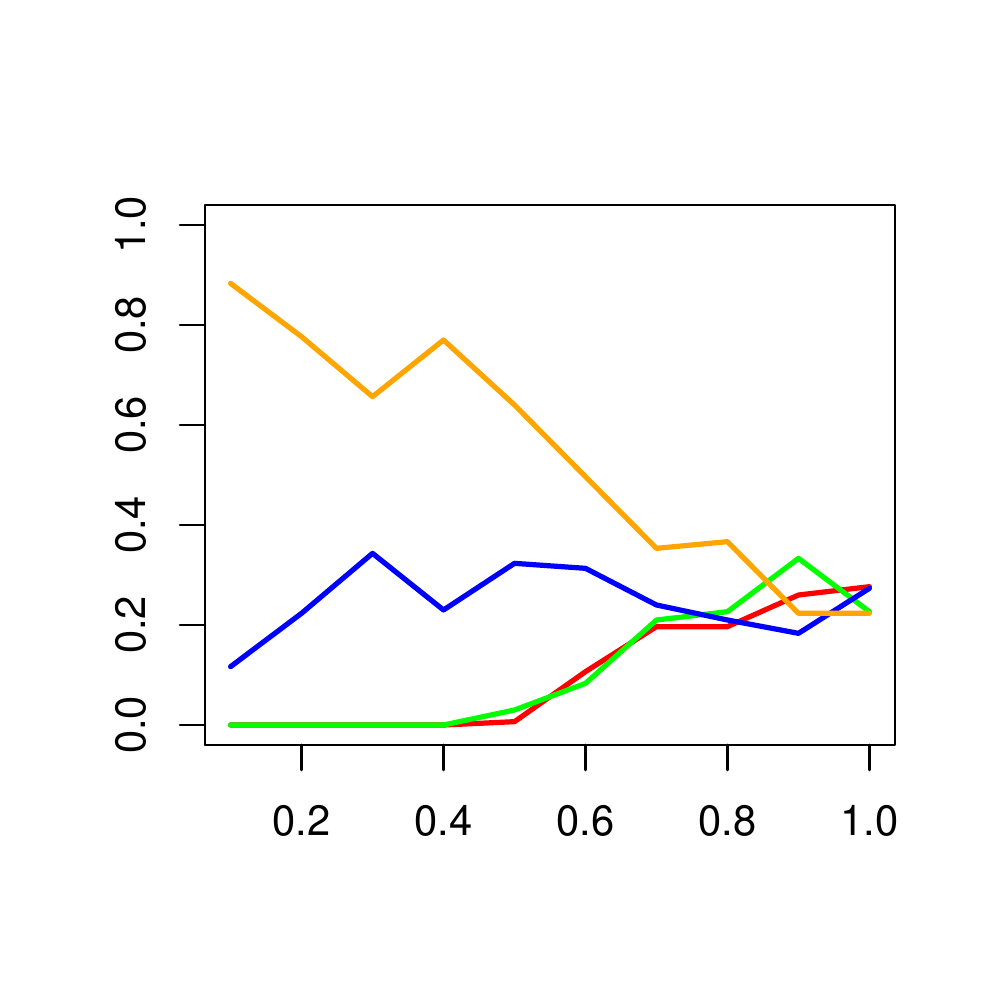}
		\vspace{-10mm} 
  	\caption{\scriptsize{Rosenbrock}} \label{fig:mv100Ro}
	\end{subfigure}
		
	\caption{Component -- decile breakdowns for the form of movement capability, across all GGGP heuristics at 100D. Components: Baseline (red), DD (green), LEH (blue), LEH+DD (orange).}
	\label{fig:move100Comp}
	
\end{figure}

\begin{figure}[H]
	\centering
	
	\vspace{-5mm} 

	\begin{subfigure}[t]{.24\textwidth}
		\includegraphics[width=\textwidth]{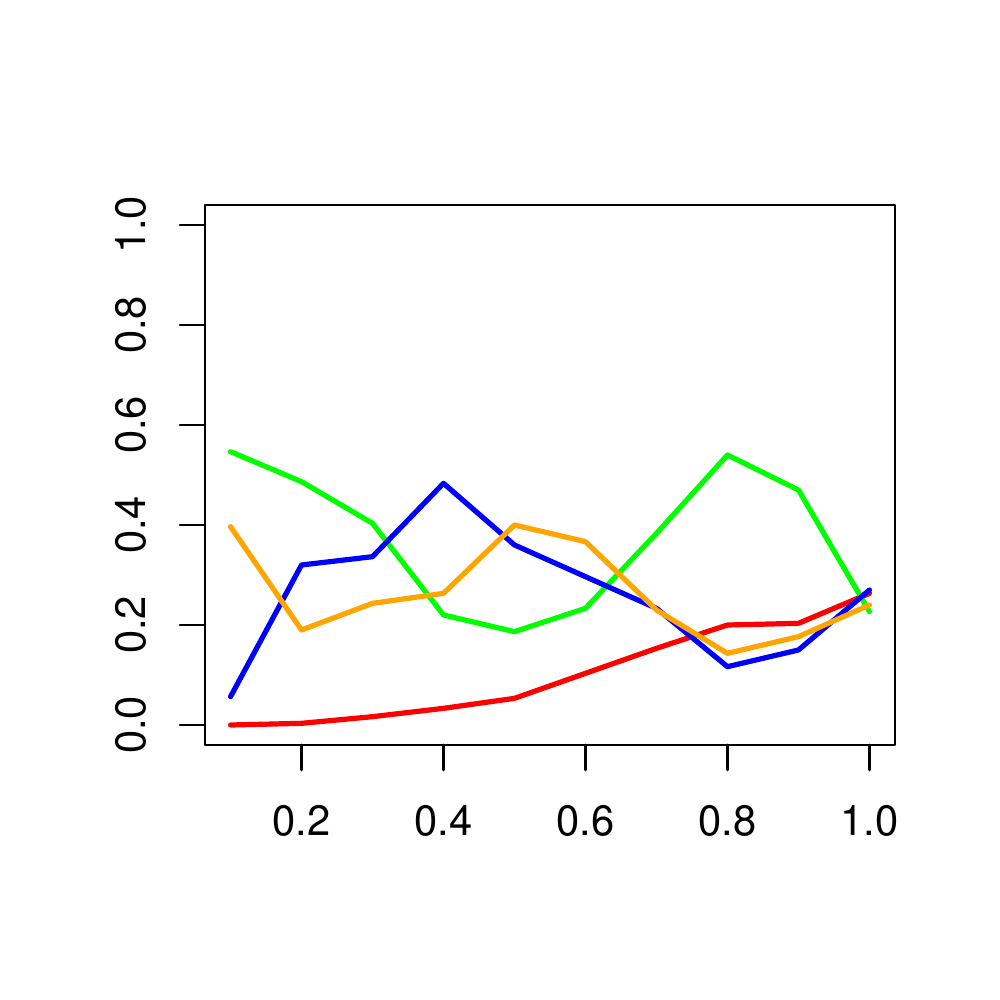}
		\vspace{-10mm} 
  	\caption{\scriptsize{30D}} \label{fig:mv30}
	\end{subfigure}%
	\hspace{-6mm} 
	\begin{subfigure}[t]{.24\textwidth}
		\includegraphics[width=\textwidth]{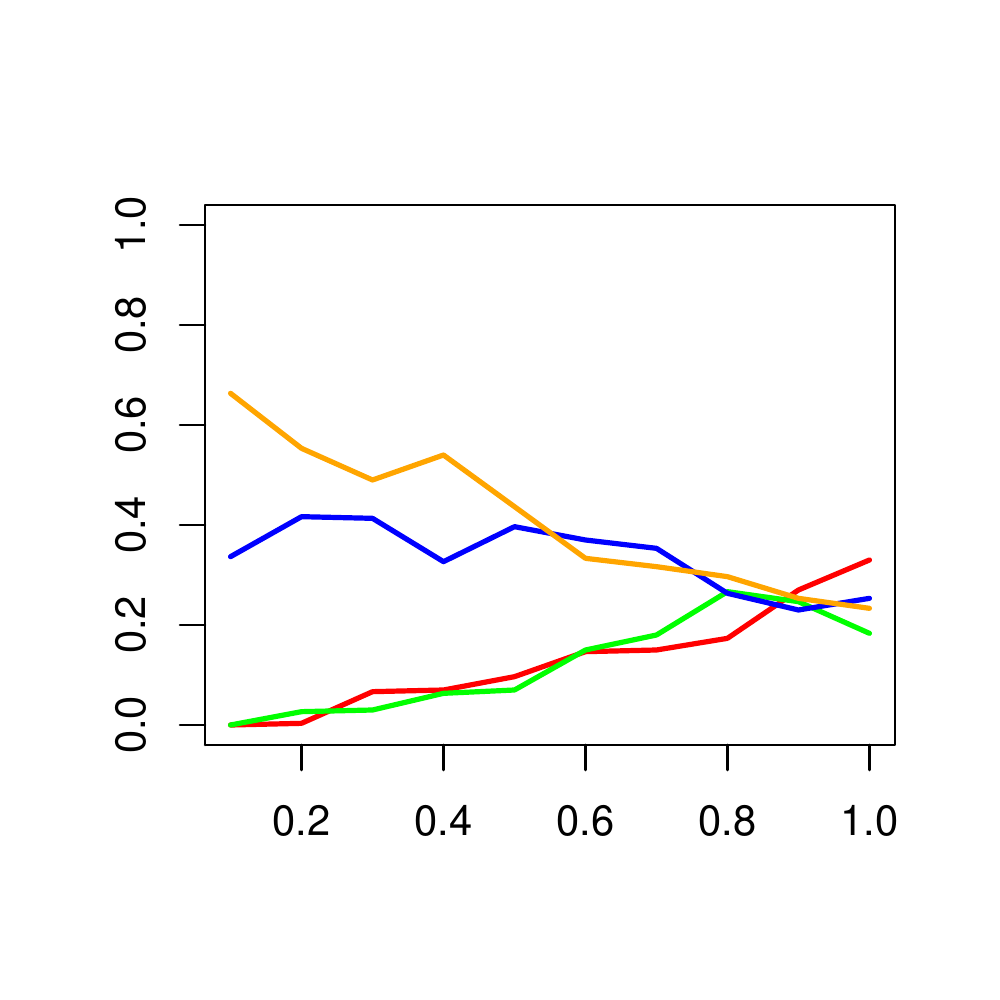}
		\vspace{-10mm} 
  	\caption{\scriptsize{100D}} \label{fig:mv100}
	\end{subfigure}

	\caption{Component -- decile breakdowns for the form of movement capability, across all GGGP heuristics at 30D and 100D for the general heuristics. Components: Baseline (red), DD (green), LEH (blue), LEH+DD (orange).}
	\label{fig:moveComp}
	
\end{figure}

However for the extended form of particle level movement, things are less clear. From Table~\ref{fig:proportions} it can be seen that the most used form of extended capability includes both descent direction and LEH dormancy-relocation (+DD+LEH) for both 30D and 100D. For both dimensions this increases in the top third performing heuristics. Nevertheless both the +LEH and +DD individual capabilities are also well represented. Figures~\ref{fig:move30Comp} to~\ref{fig:moveComp} show the decile analysis for the extended movement capability: no additional capability (red: rPSO), +DD (green), +LEH (blue), or +DD+LEH (orange). The plots are in accord with the high level results, with +DD+LEH most associated with the best performing heuristics but both +DD and +LEH also performing well.

\vspace{2mm} 

\vspace{2mm} 

\subsubsection{Network topology}
\label{sec:compNet}

\begin{figure}[H]
	\centering
	
	\vspace{-10mm} 

	\begin{subfigure}[t]{.2\textwidth}
		\includegraphics[width=\textwidth]{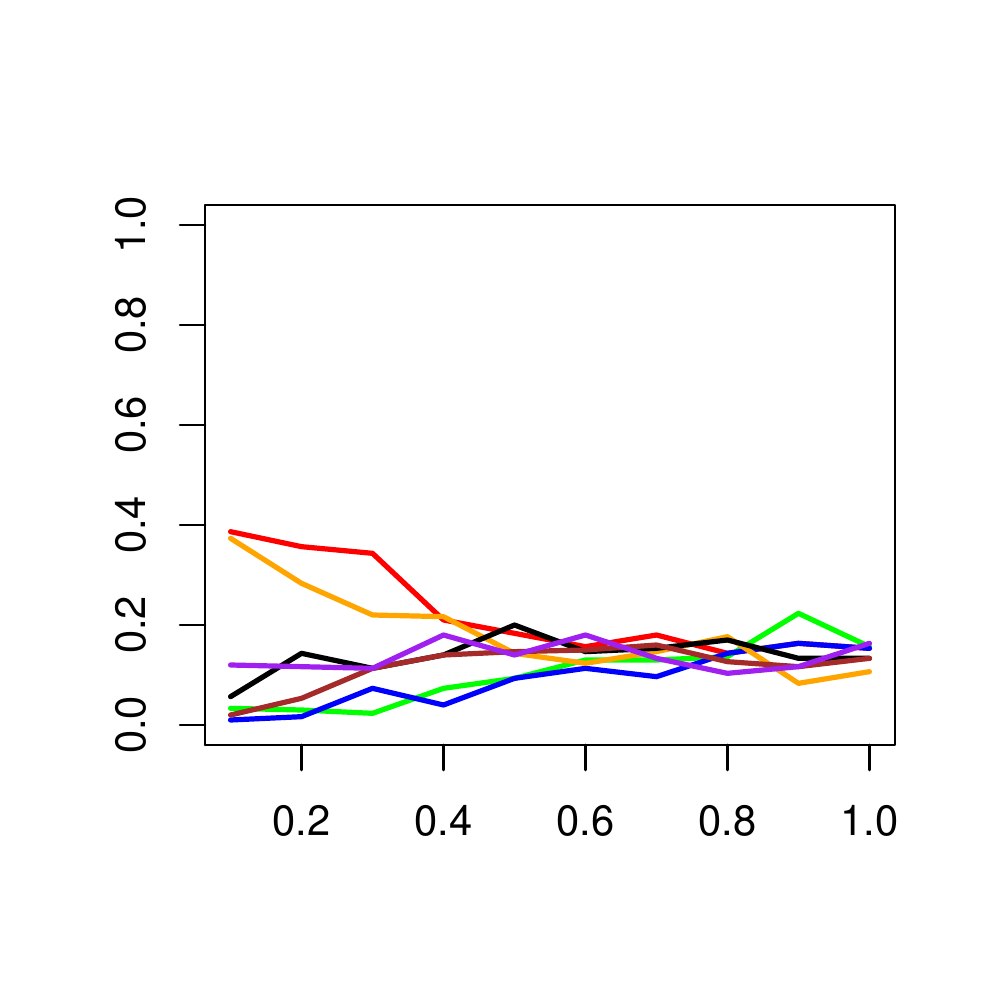}
		\vspace{-10mm} 
  	\caption{\scriptsize{Rastrigin}} \label{fig:nw30Ra}
	\end{subfigure}%
	\hspace{-6mm} 
	\begin{subfigure}[t]{.2\textwidth}
		\includegraphics[width=\textwidth]{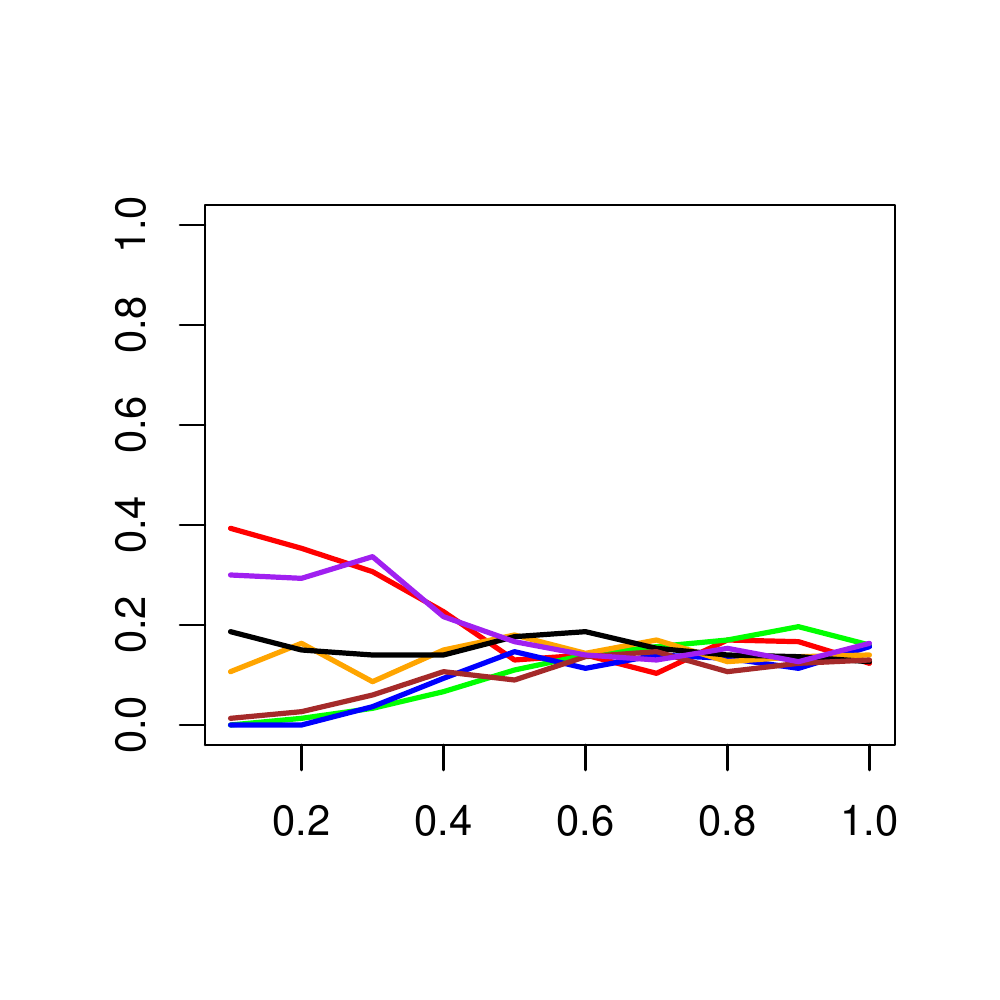}
		\vspace{-10mm} 
  	\caption{\scriptsize{Multipeak F1}} \label{fig:nw30M1}
	\end{subfigure}%
	\hspace{-6mm} 
	\begin{subfigure}[t]{.2\textwidth}
		\includegraphics[width=\textwidth]{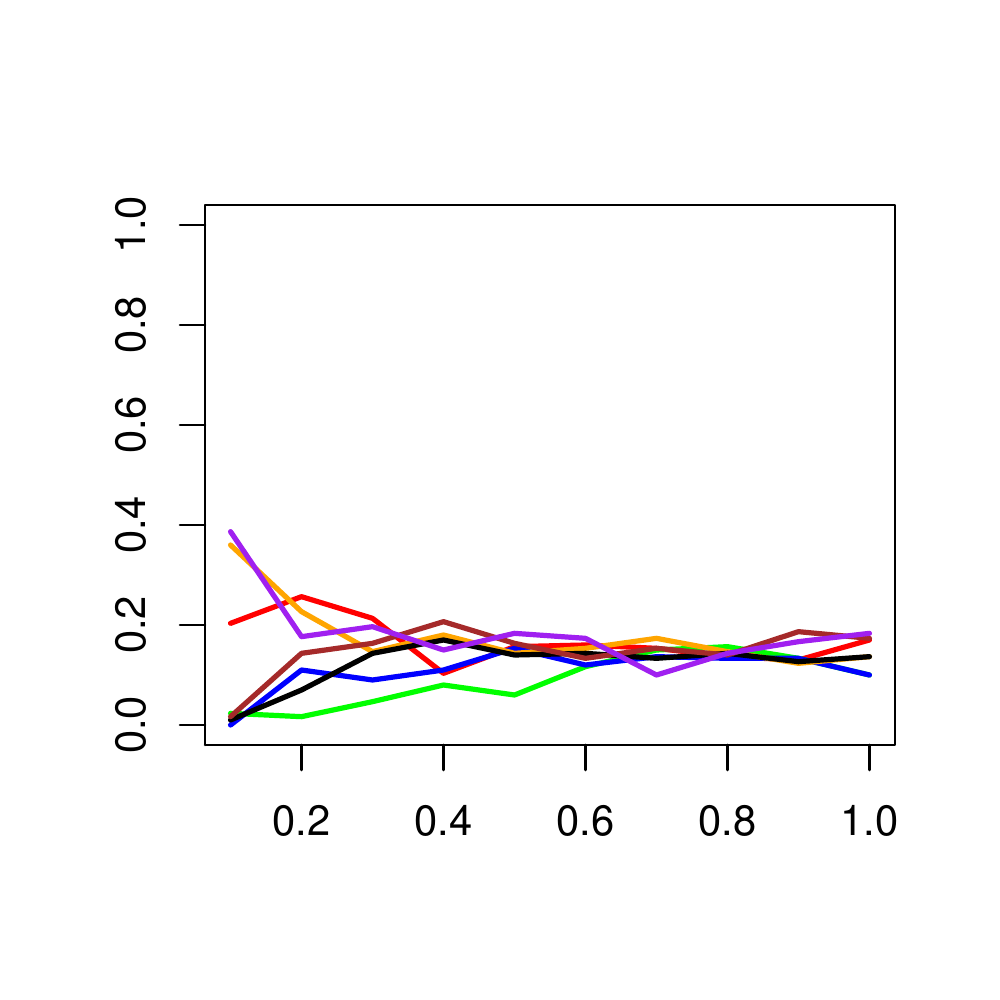}
		\vspace{-10mm} 
  	\caption{\scriptsize{Multipeak F2}} \label{fig:nw30M2}
	\end{subfigure}%
	\hspace{-6mm} 
	\begin{subfigure}[t]{.2\textwidth}
		\includegraphics[width=\textwidth]{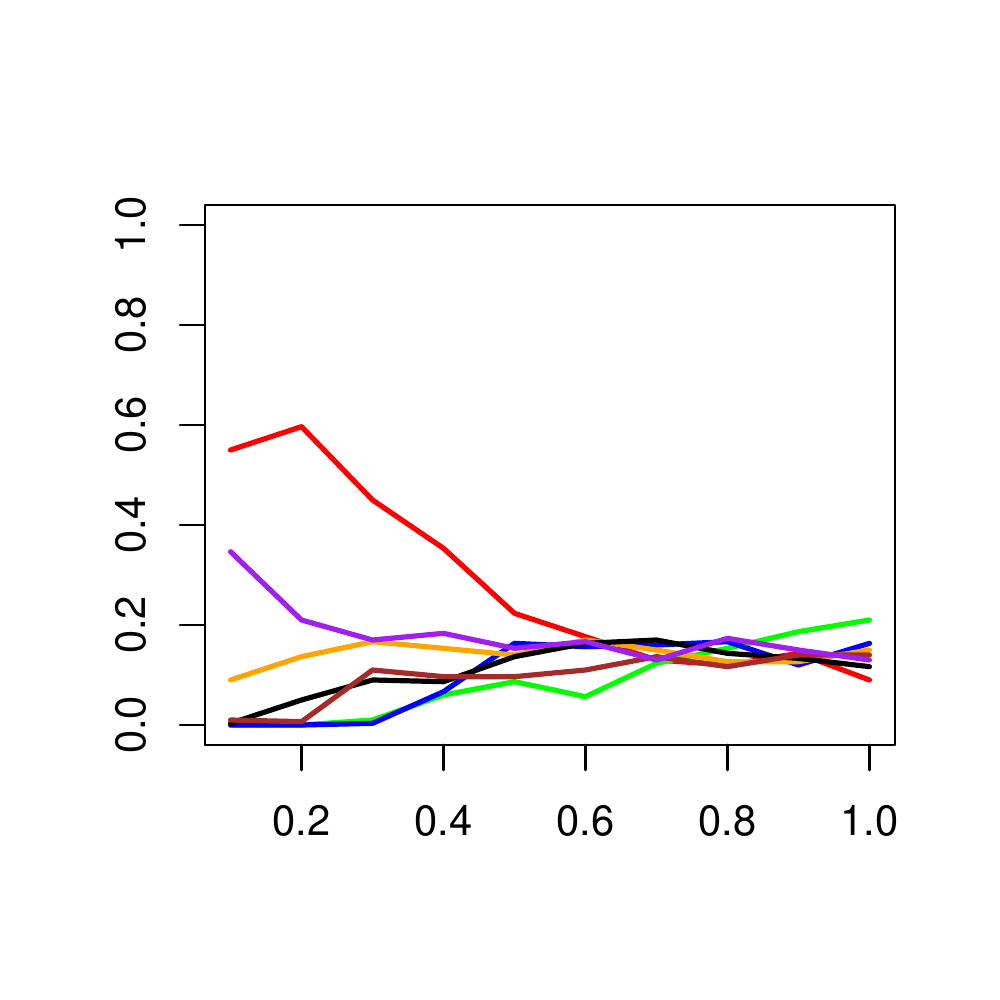}
		\vspace{-10mm} 
  	\caption{\scriptsize{Brankes}} \label{fig:nw30Br}
	\end{subfigure}%
	\hspace{-6mm} 
	\begin{subfigure}[t]{.2\textwidth}
		\includegraphics[width=\textwidth]{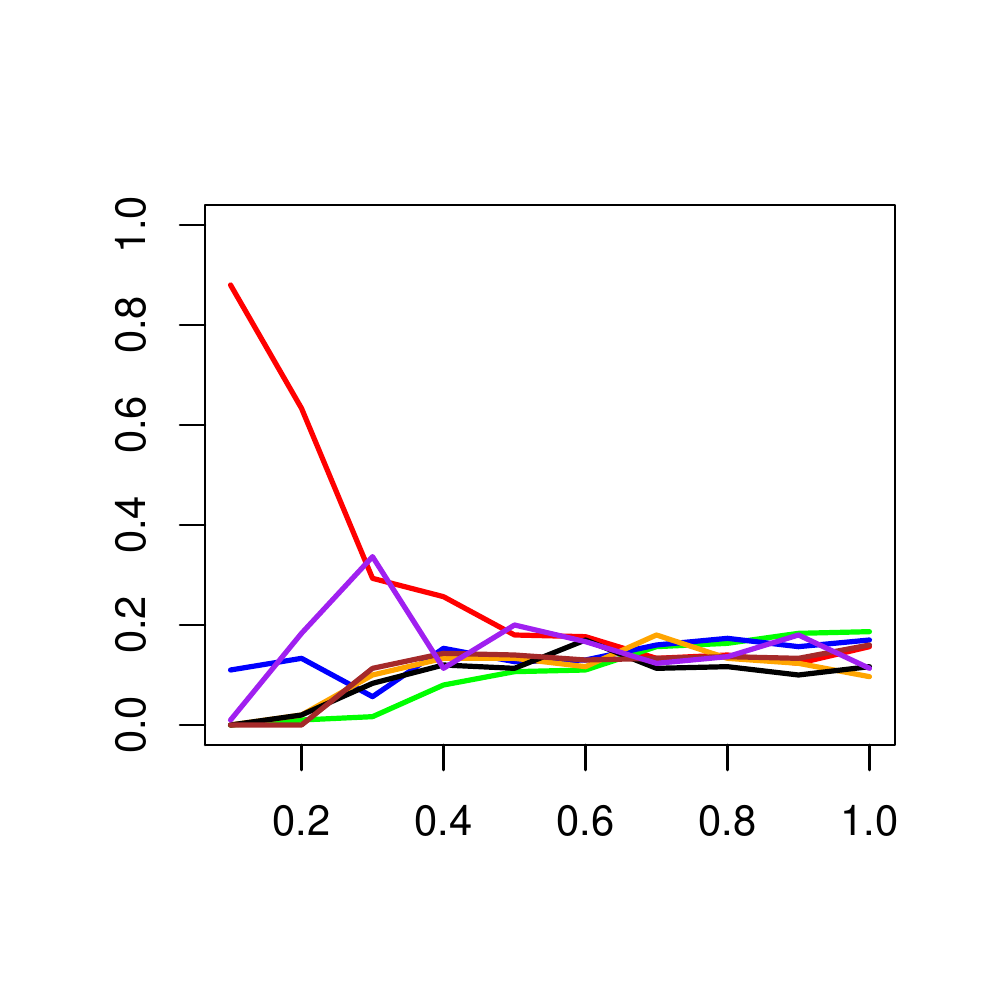}
		\vspace{-10mm} 
  	\caption{\scriptsize{Pickelhaube}} \label{fig:nw30Pi}
	\end{subfigure}
	
	\vspace{-2mm} 
		
	\begin{subfigure}[t]{.2\textwidth}
		\includegraphics[width=\textwidth]{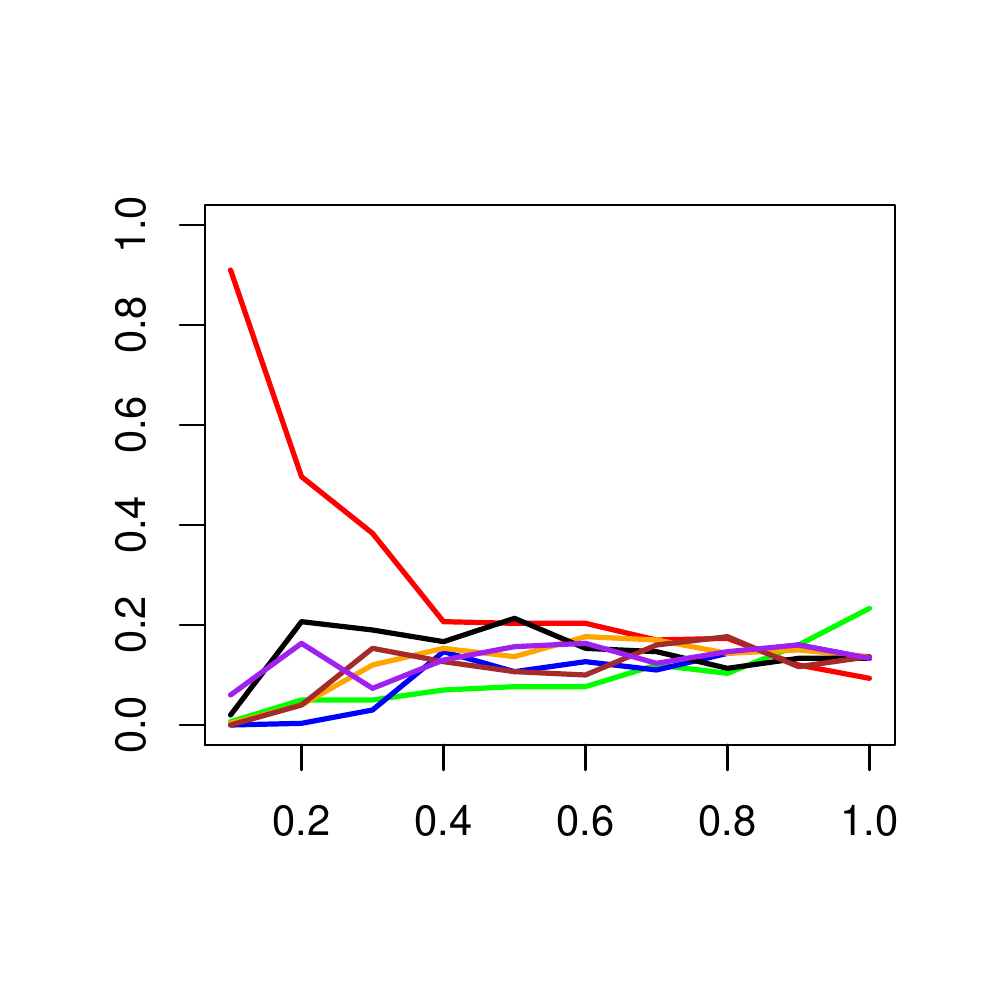}
		\vspace{-10mm} 
  	\caption{\scriptsize{Heaviside}} \label{fig:nw30Hv}
	\end{subfigure}%
	\hspace{-6mm} 
	\begin{subfigure}[t]{.2\textwidth}
		\includegraphics[width=\textwidth]{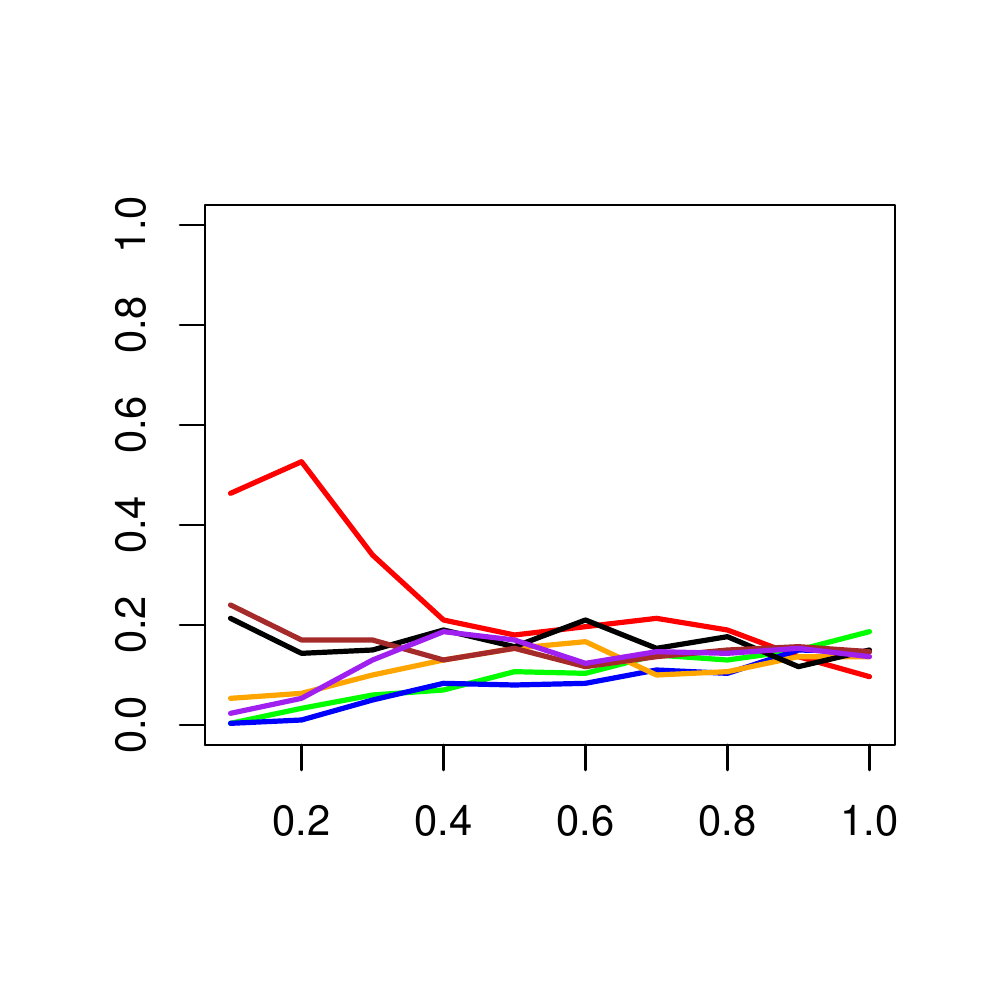}
		\vspace{-10mm} 
  	\caption{\scriptsize{Sawtooth}} \label{fig:nw30Sa}
	\end{subfigure}%
	\hspace{-6mm} 
	\begin{subfigure}[t]{.2\textwidth}
		\includegraphics[width=\textwidth]{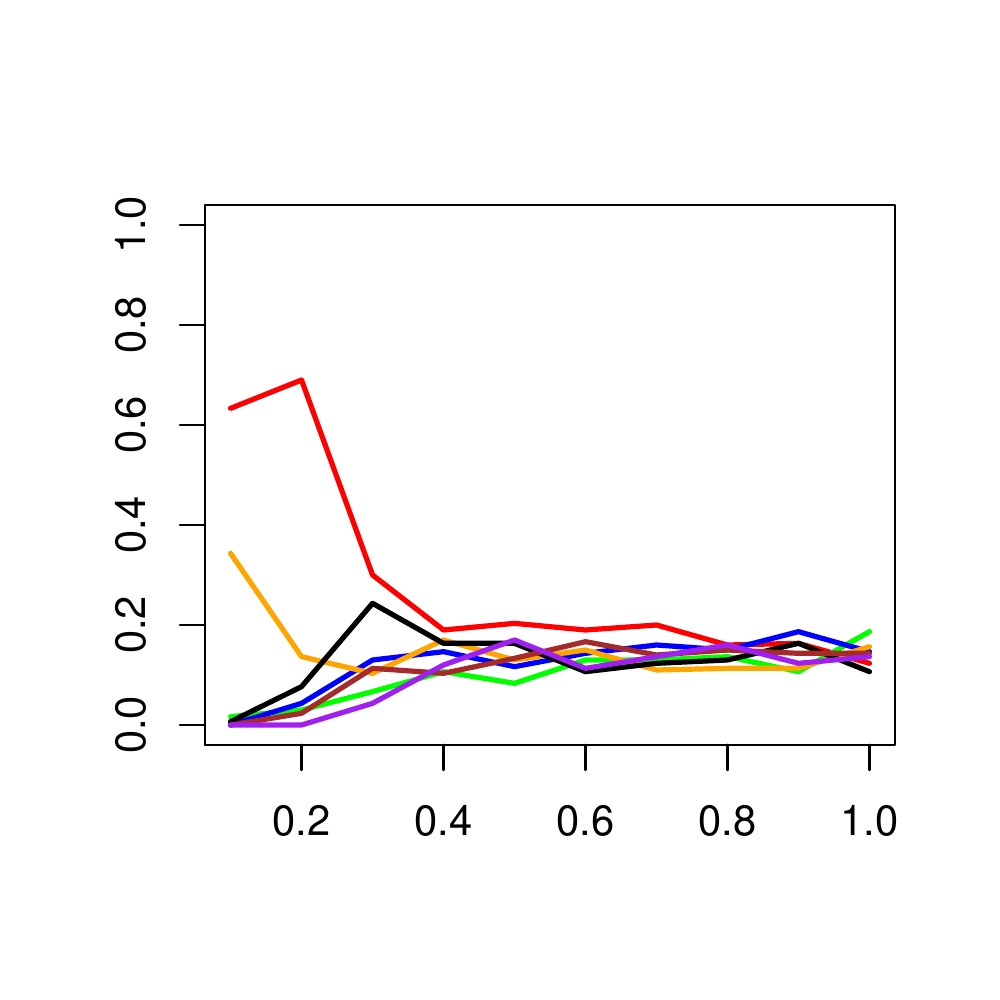}
		\vspace{-10mm} 
  	\caption{\scriptsize{Ackley}} \label{fig:nw30Ac}
	\end{subfigure}%
	\hspace{-6mm} 
	\begin{subfigure}[t]{.2\textwidth}
		\includegraphics[width=\textwidth]{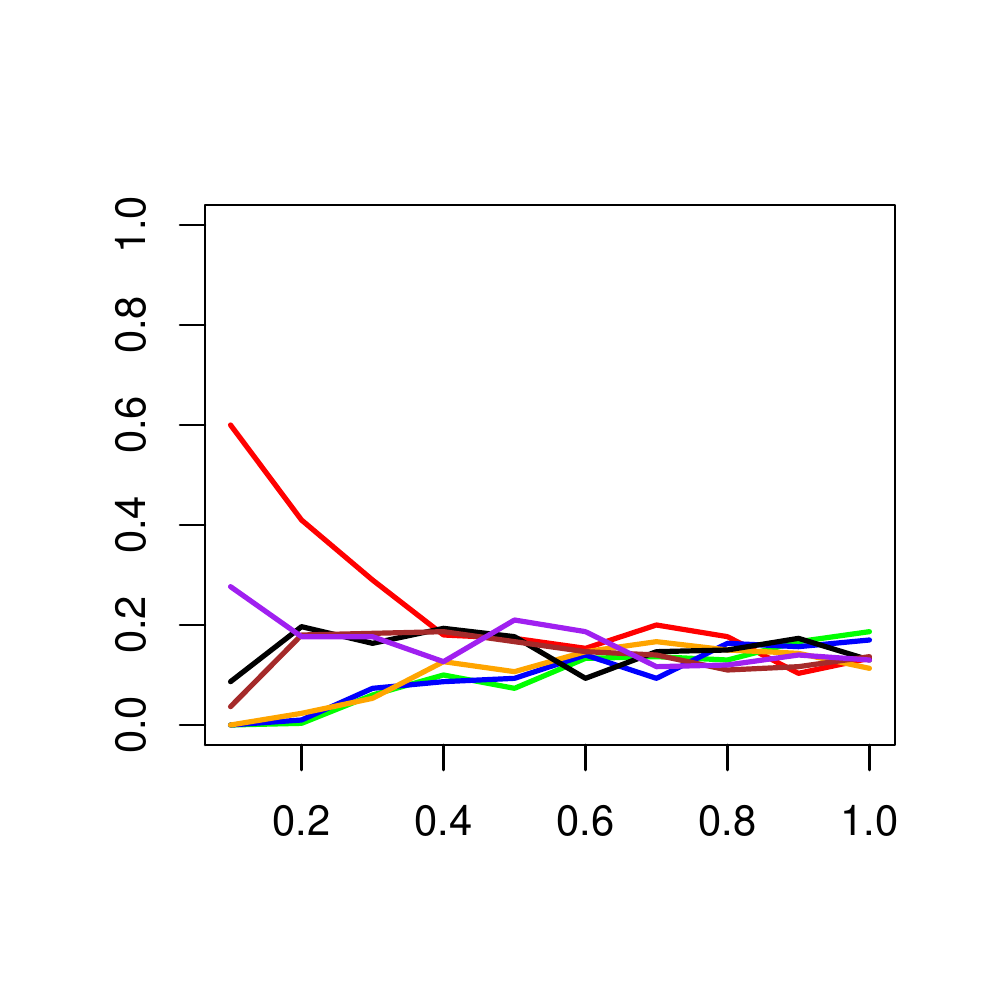}
		\vspace{-10mm} 
  	\caption{\scriptsize{Sphere}} \label{fig:nw30Sp}
	\end{subfigure}%
	\hspace{-6mm} 
	\begin{subfigure}[t]{.2\textwidth}
		\includegraphics[width=\textwidth]{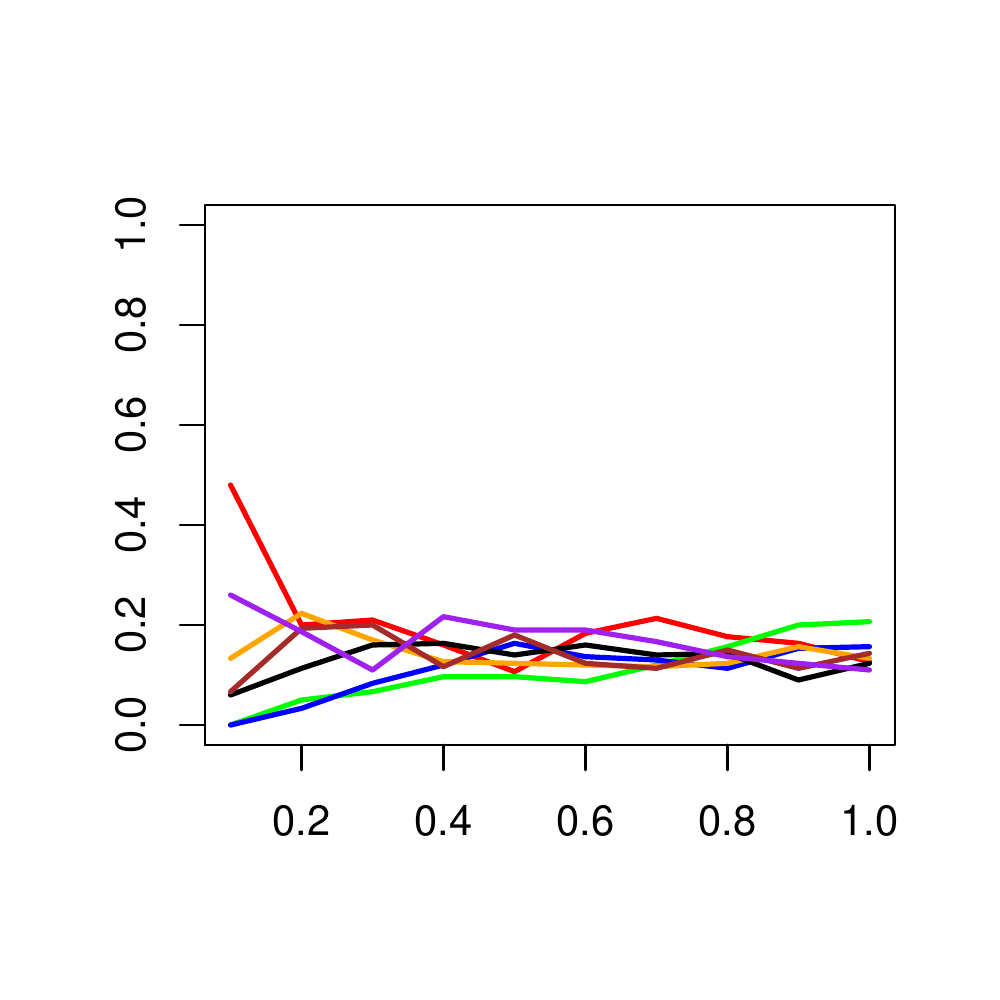}
		\vspace{-10mm} 
  	\caption{\scriptsize{Rosenbrock}} \label{fig:nw30Ro}
	\end{subfigure}
		
	\caption{Component -- decile breakdowns for the form of network for particle information sharing, across all GGGP heuristics at 30D. Components: Global (red), Focal (green), Ring (size=2) (blue), von Neumann (orange), Clan (black), Cluster (brown), Hierarchy (purple).}
	\label{fig:network30Comp}
	
\end{figure}

\begin{figure}[H]
	\centering
	
	\vspace{-5mm} 

	\begin{subfigure}[t]{.2\textwidth}
		\includegraphics[width=\textwidth]{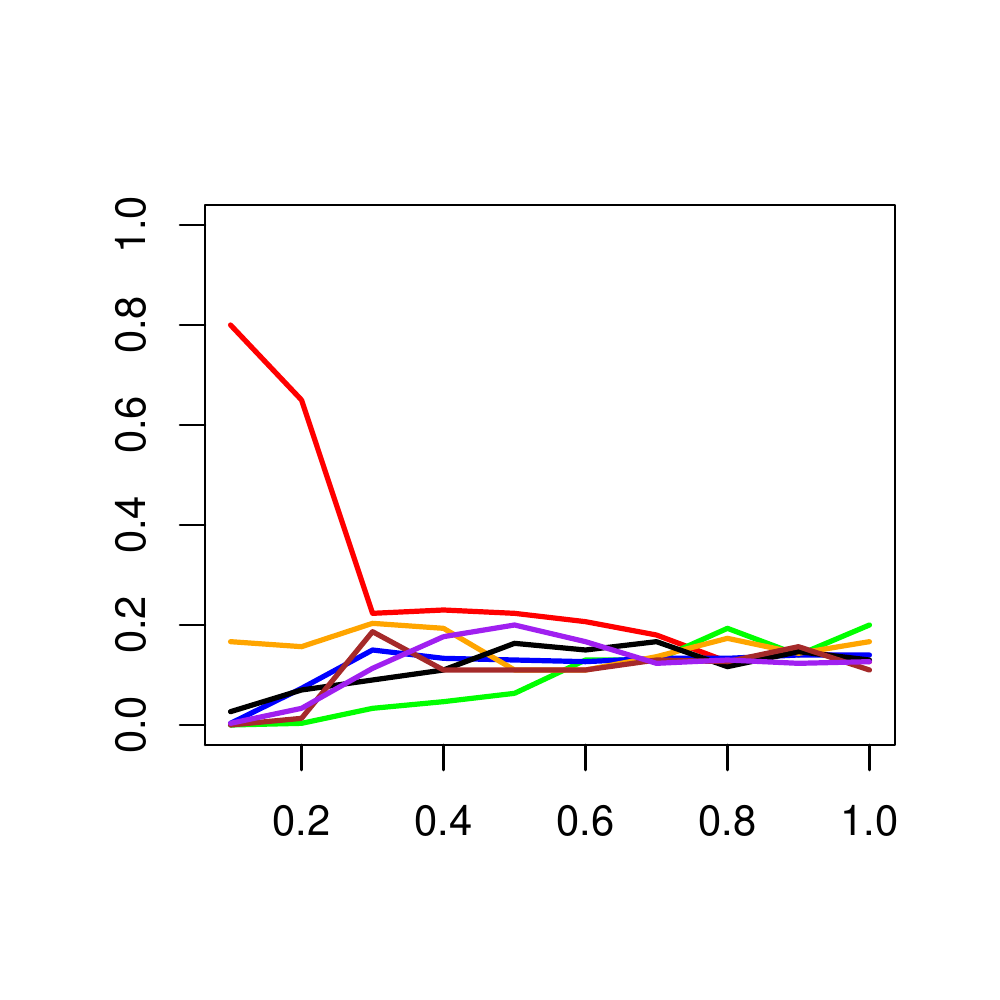}
		\vspace{-10mm} 
  	\caption{\scriptsize{Rastrigin}} \label{fig:nw100Ra}
	\end{subfigure}%
	\hspace{-6mm} 
	\begin{subfigure}[t]{.2\textwidth}
		\includegraphics[width=\textwidth]{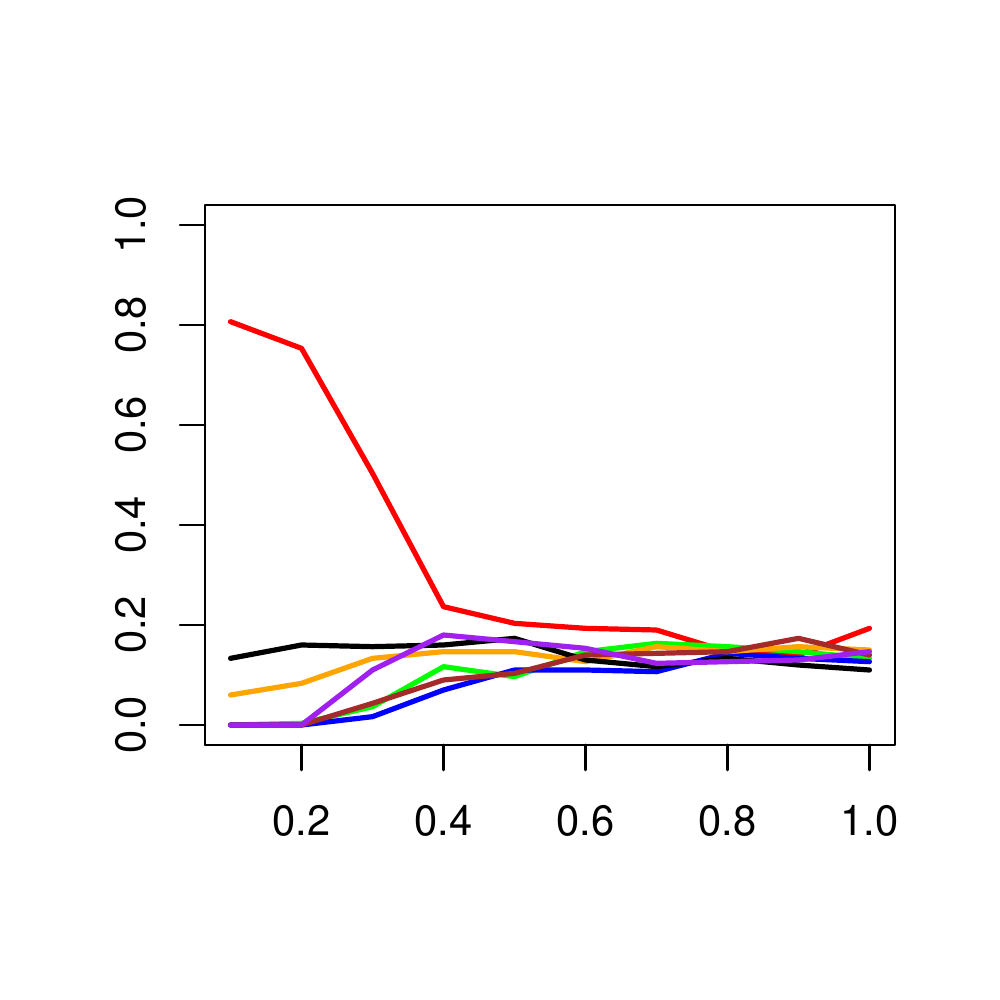}
		\vspace{-10mm} 
  	\caption{\scriptsize{Multipeak F1}} \label{fig:nw100M1}
	\end{subfigure}%
	\hspace{-6mm} 
	\begin{subfigure}[t]{.2\textwidth}
		\includegraphics[width=\textwidth]{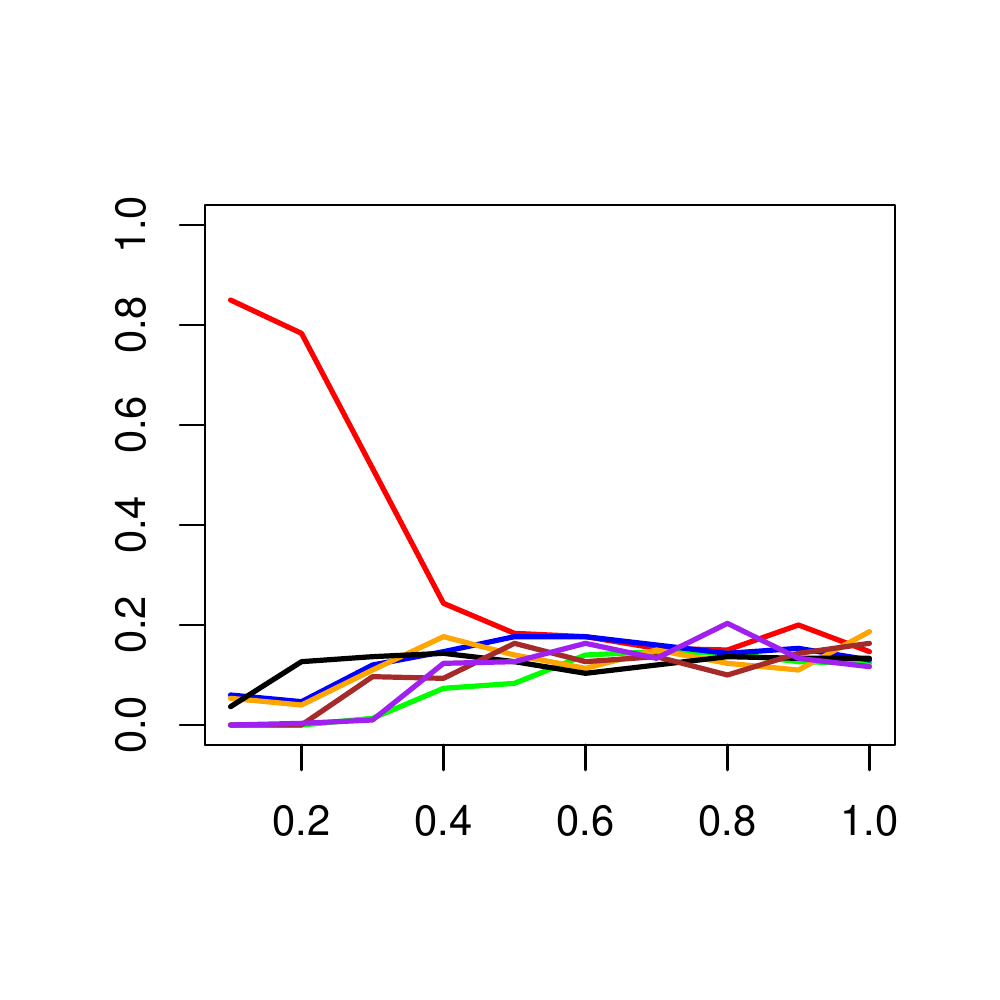}
		\vspace{-10mm} 
  	\caption{\scriptsize{Multipeak F2}} \label{fig:nw100M2}
	\end{subfigure}%
	\hspace{-6mm} 
	\begin{subfigure}[t]{.2\textwidth}
		\includegraphics[width=\textwidth]{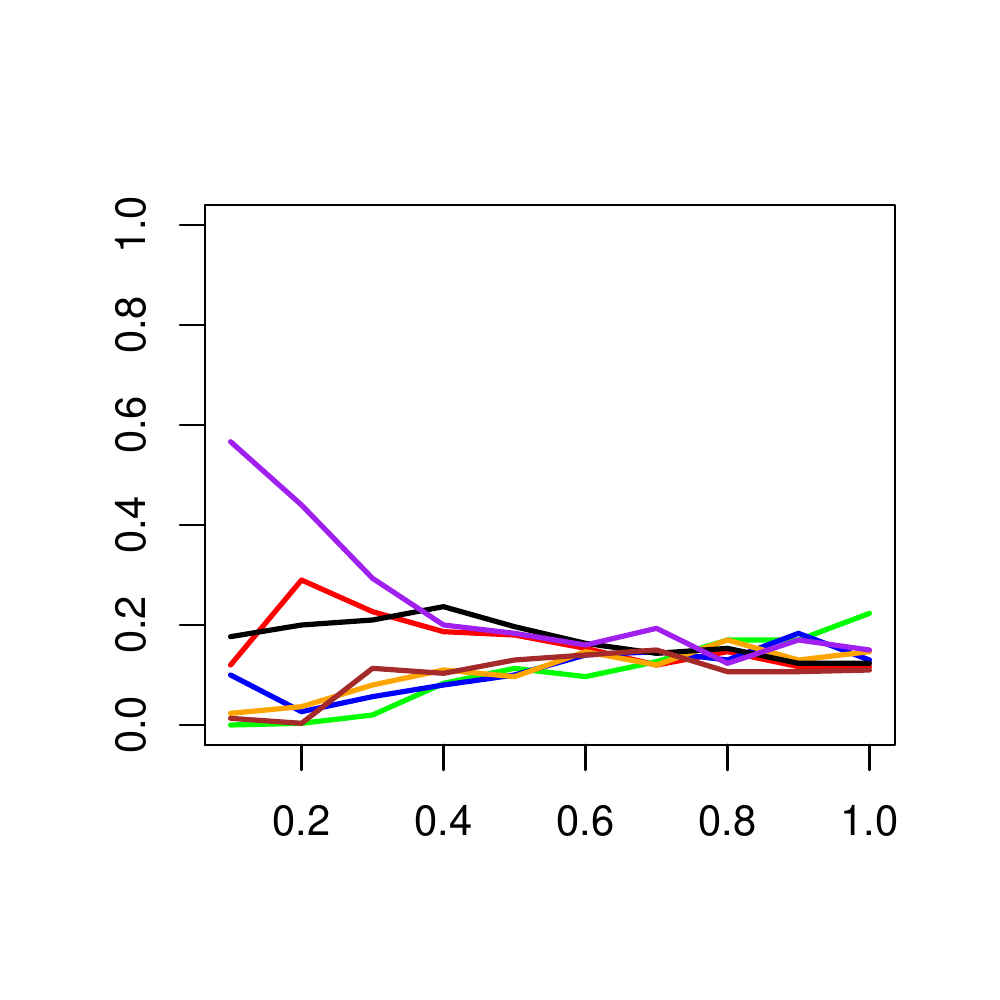}
		\vspace{-10mm} 
  	\caption{\scriptsize{Brankes}} \label{fig:nw100Br}
	\end{subfigure}%
	\hspace{-6mm} 
	\begin{subfigure}[t]{.2\textwidth}
		\includegraphics[width=\textwidth]{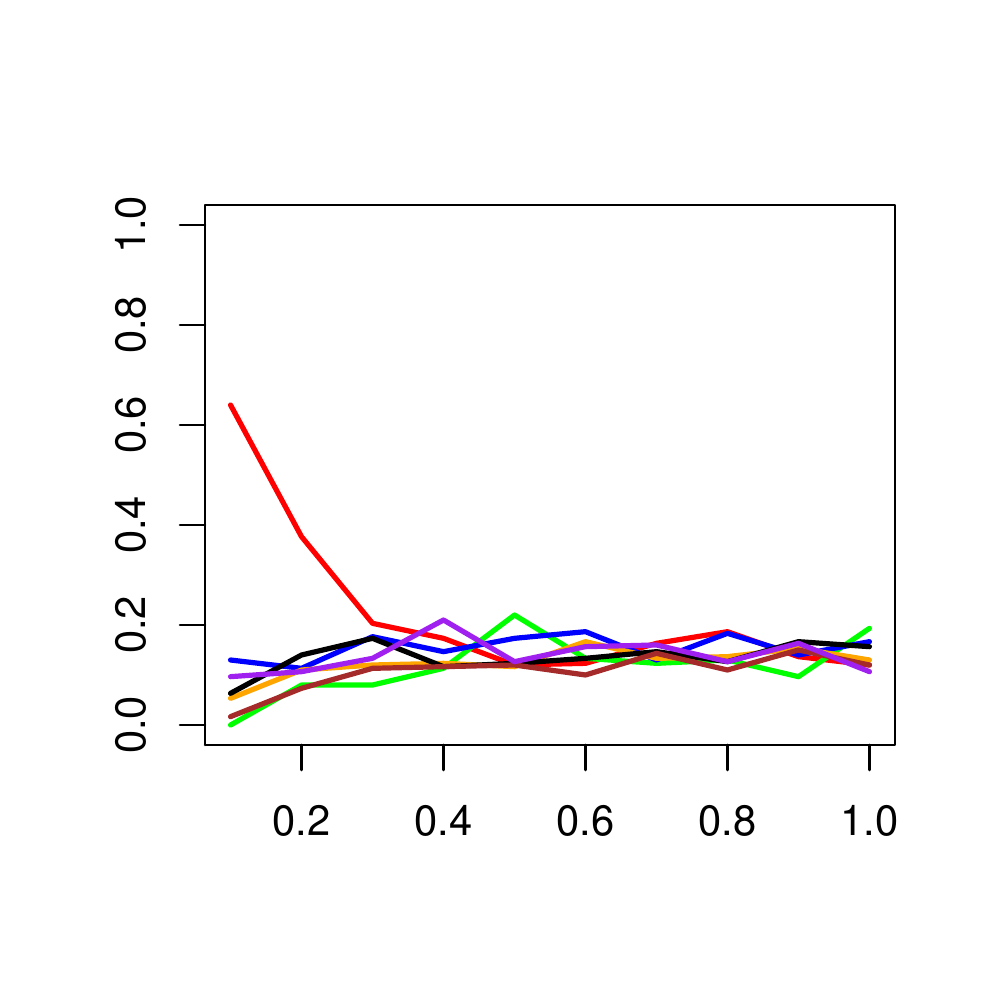}
		\vspace{-10mm} 
  	\caption{\scriptsize{Pickelhaube}} \label{fig:nw100Pi}
	\end{subfigure}
	
	\vspace{-2mm} 
		
	\begin{subfigure}[t]{.2\textwidth}
		\includegraphics[width=\textwidth]{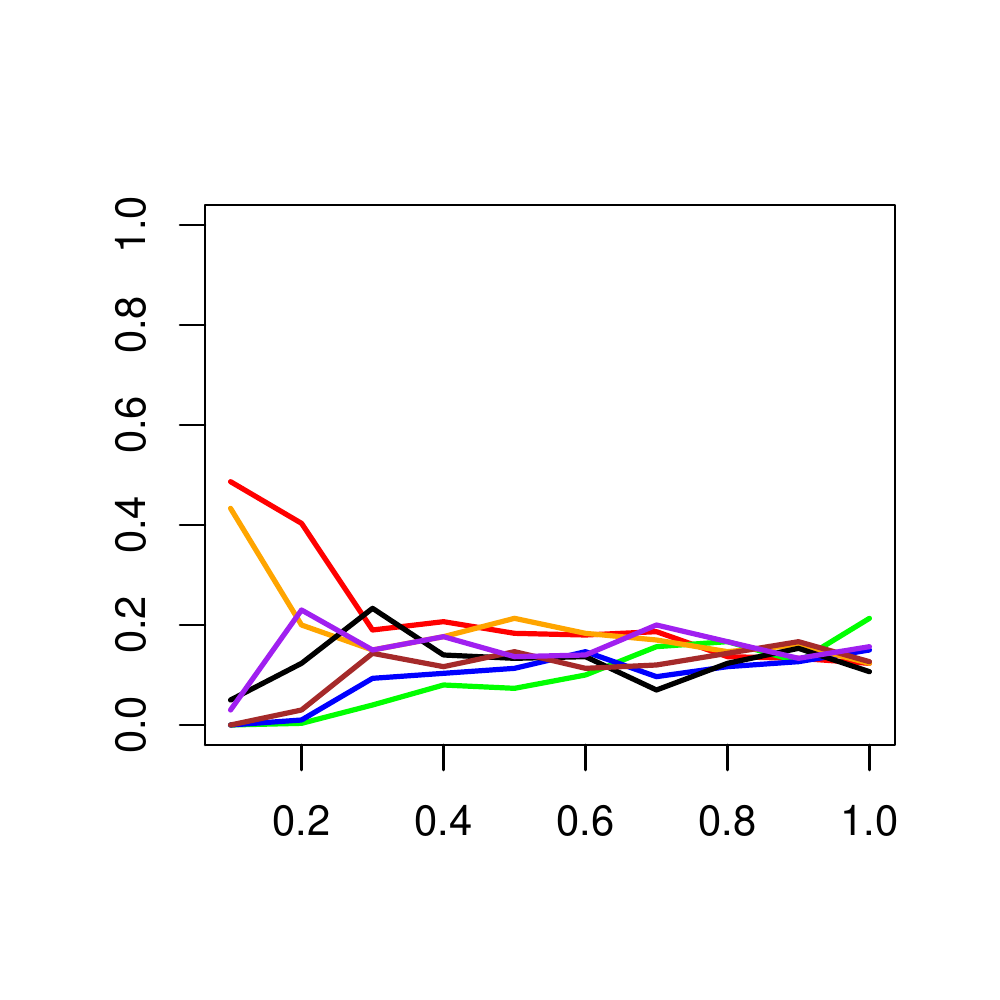}
		\vspace{-10mm} 
  	\caption{\scriptsize{Heaviside}} \label{fig:nw100Hv}
	\end{subfigure}%
	\hspace{-6mm} 
	\begin{subfigure}[t]{.2\textwidth}
		\includegraphics[width=\textwidth]{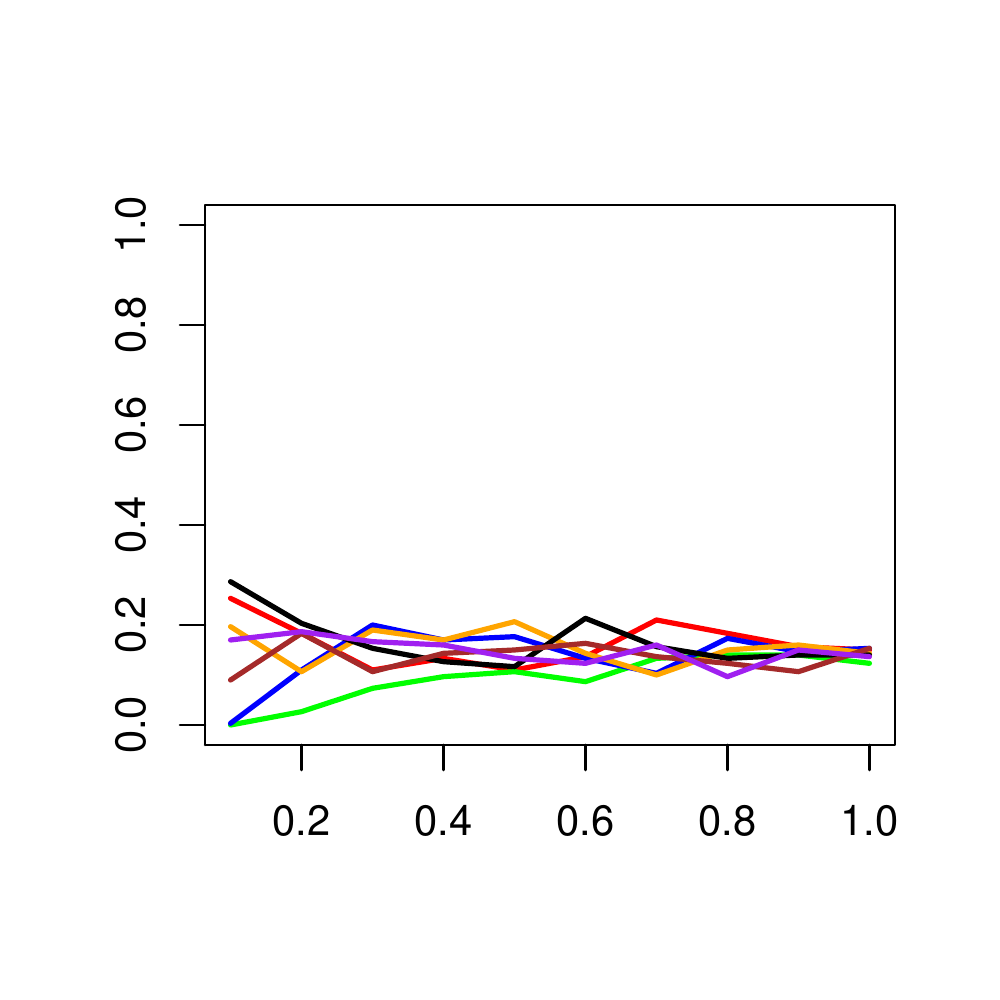}
		\vspace{-10mm} 
  	\caption{\scriptsize{Sawtooth}} \label{fig:nw100Sa}
	\end{subfigure}%
	\hspace{-6mm} 
	\begin{subfigure}[t]{.2\textwidth}
		\includegraphics[width=\textwidth]{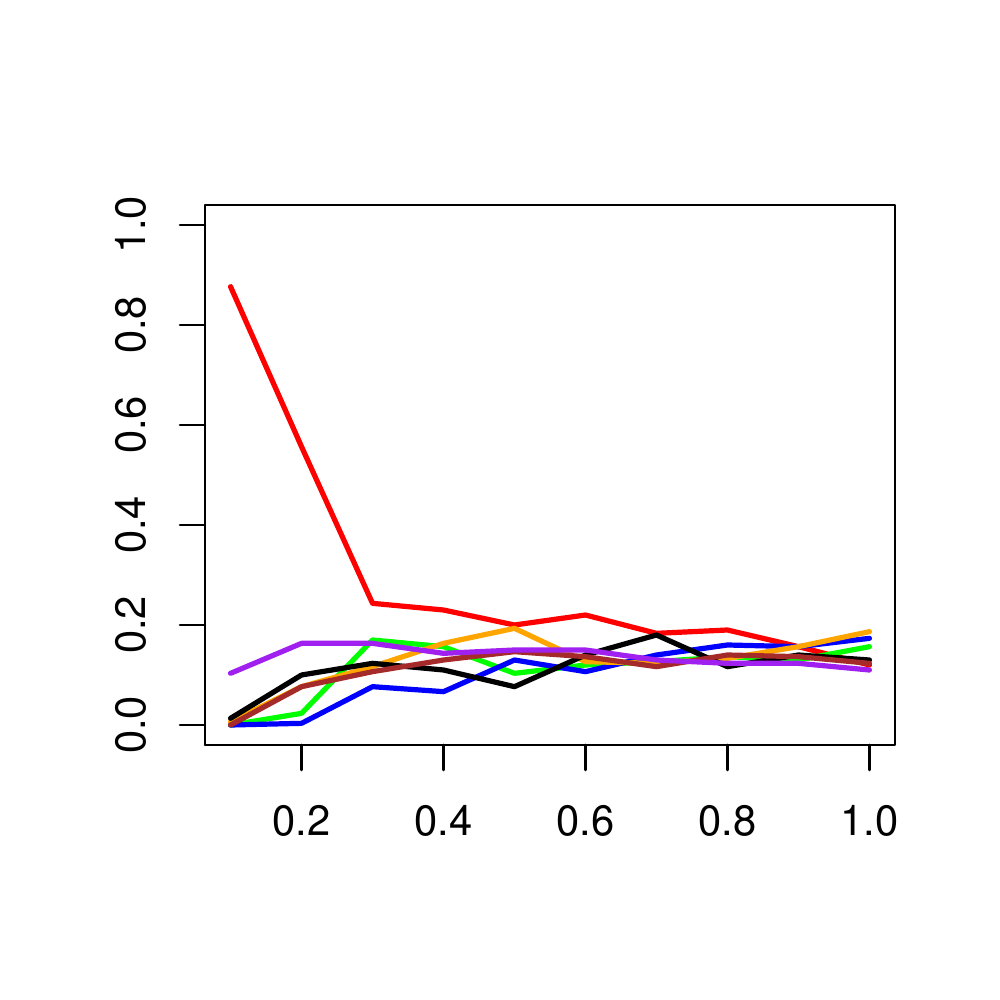}
		\vspace{-10mm} 
  	\caption{\scriptsize{Ackley}} \label{fig:nw100Ac}
	\end{subfigure}%
	\hspace{-6mm} 
	\begin{subfigure}[t]{.2\textwidth}
		\includegraphics[width=\textwidth]{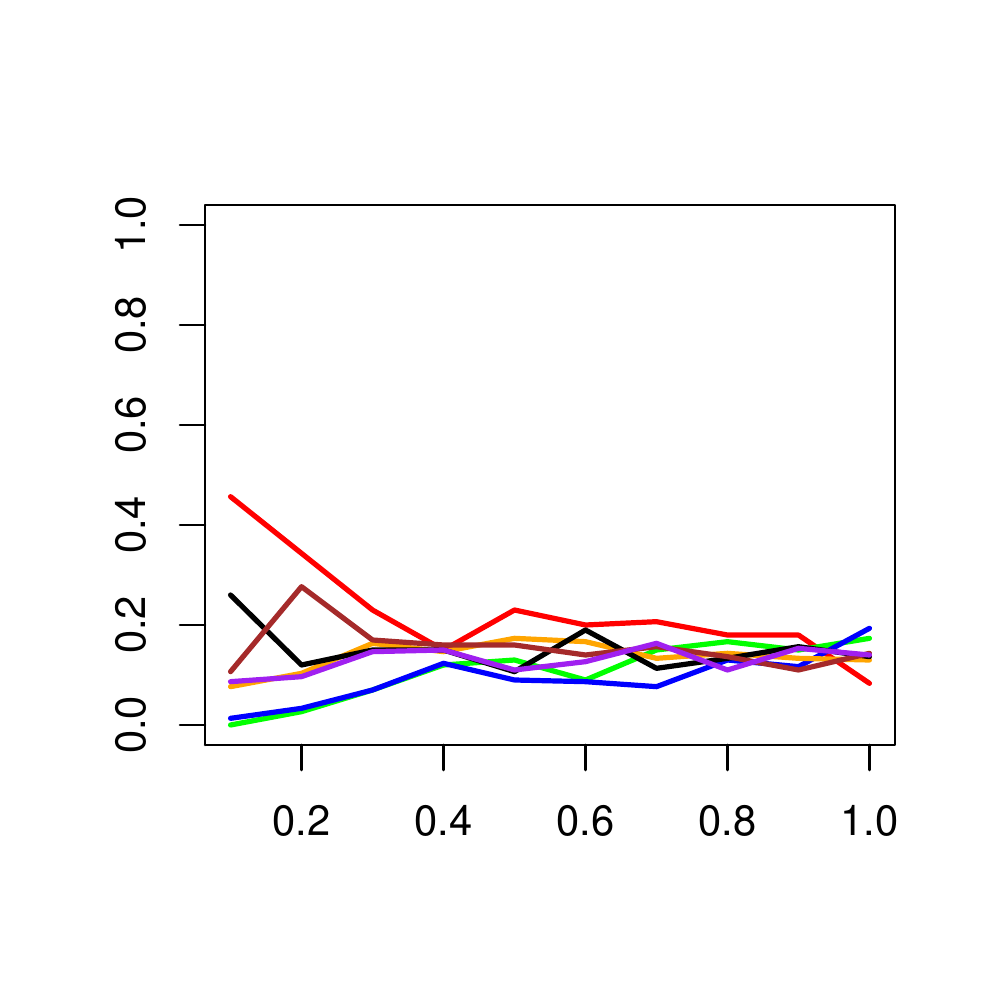}
		\vspace{-10mm} 
  	\caption{\scriptsize{Sphere}} \label{fig:nw100Sp}
	\end{subfigure}%
	\hspace{-6mm} 
	\begin{subfigure}[t]{.2\textwidth}
		\includegraphics[width=\textwidth]{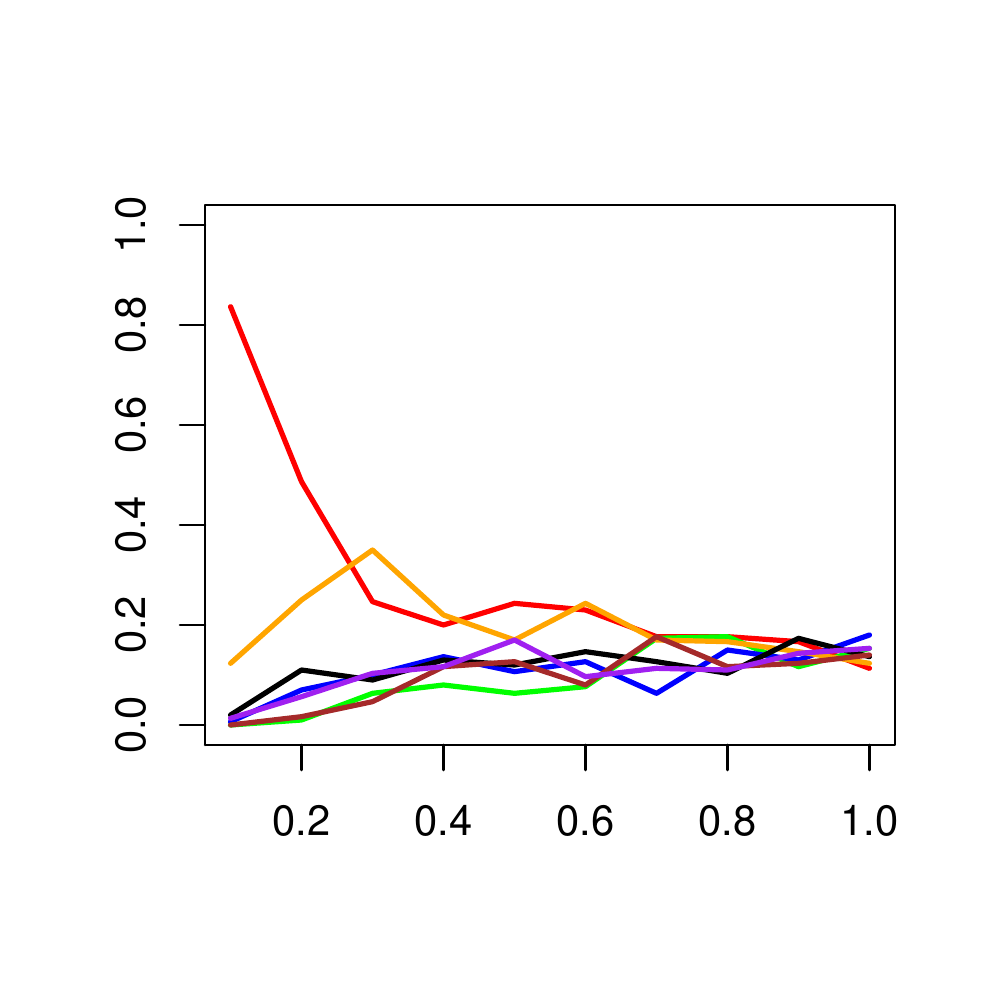}
		\vspace{-10mm} 
  	\caption{\scriptsize{Rosenbrock}} \label{fig:nw100Ro}
	\end{subfigure}
		
	\caption{Component -- decile breakdowns for the form of network for particle information sharing, across all GGGP heuristics at 100D. Components: Global (red), Focal (green), Ring (size=2) (blue), von Neumann (orange), Clan (black), Cluster (brown), Hierarchy (purple).}
	\label{fig:network100Comp}
	
\end{figure}

\begin{figure}[H]
	\centering
	
	\vspace{-5mm} 

	\begin{subfigure}[t]{.24\textwidth}
		\includegraphics[width=\textwidth]{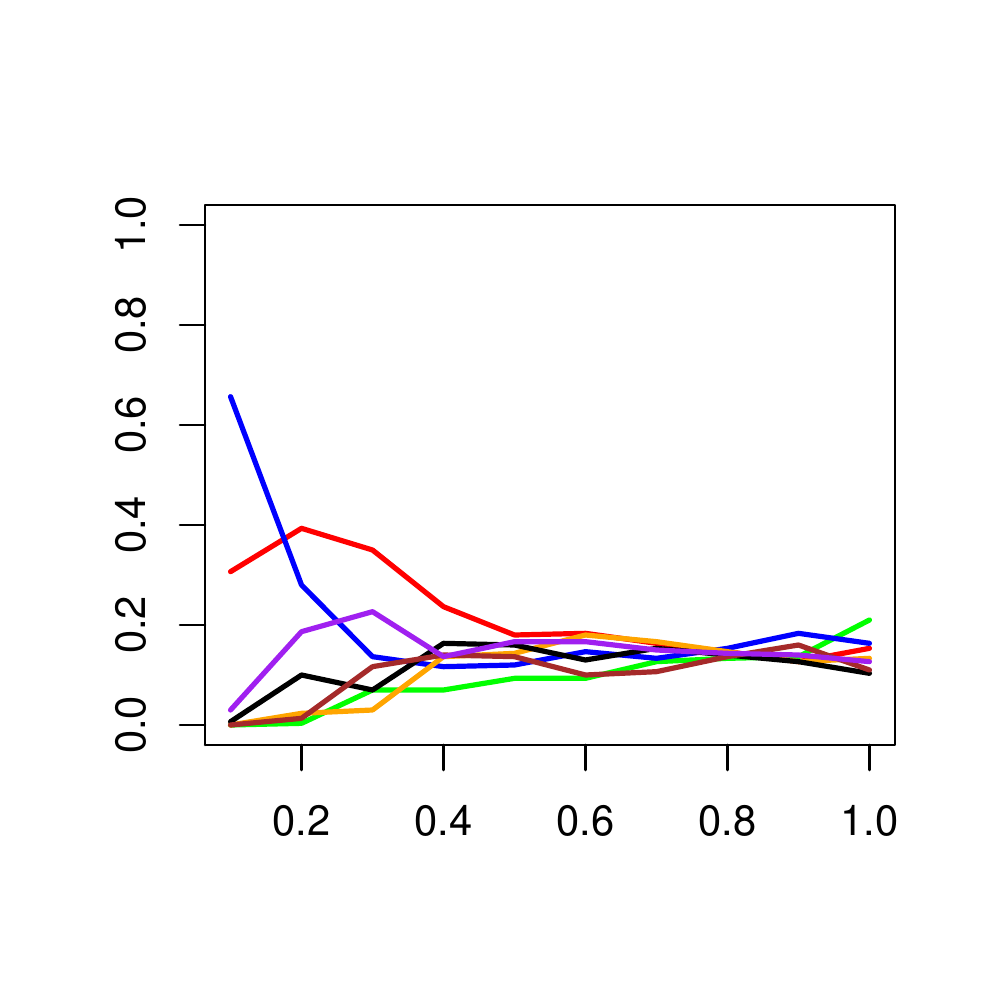}
		\vspace{-10mm} 
  	\caption{\scriptsize{30D}} \label{fig:nw30}
	\end{subfigure}%
	\hspace{-6mm} 
	\begin{subfigure}[t]{.24\textwidth}
		\includegraphics[width=\textwidth]{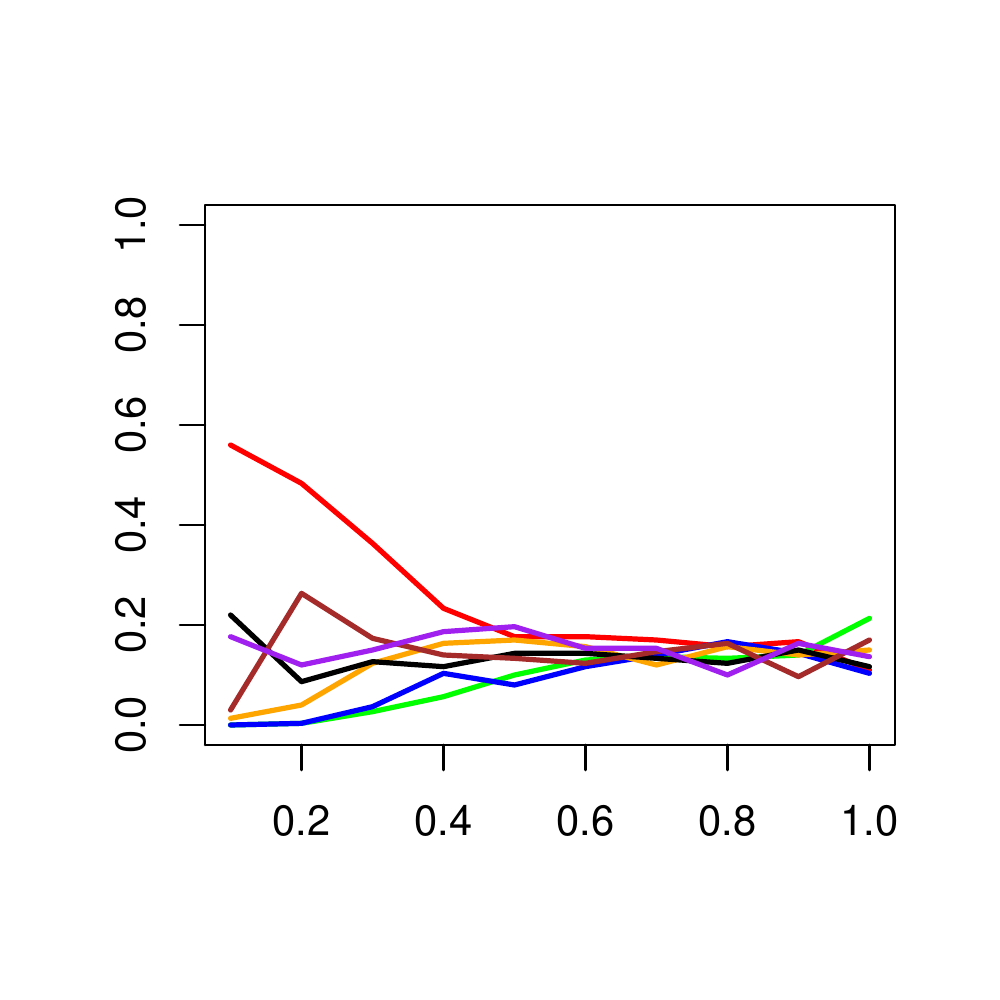}
		\vspace{-10mm} 
  	\caption{\scriptsize{100D}} \label{fig:nw100}
	\end{subfigure}

	\caption{Component -- decile breakdowns for the form of network for particle information sharing, across all GGGP heuristics at 30D and 100D for the general heuristics. Components: Global (red), Focal (green), Ring (size=2) (blue), von Neumann (orange), Clan (black), Cluster (brown), Hierarchy (purple).}
	\label{fig:networkComp}
	
\end{figure}

For the form of network for particle information sharing, Table~\ref{fig:proportions} shows that a Global (red) network is most preferred for both 30D and 100D, with this increasing in the top third performing heuristics. No others forms of network particularly stand out. This is reflected in the decile plots for the individual case heuristics, Figures~\ref{fig:network30Comp} and~\ref{fig:network100Comp}, although both Hierarchical (purple) and von Neumann (orange) networks are also represented in the lowest deciles for a few problems. In the general cases, Figure~\ref{fig:networkComp}, at 30D the Ring (blue) network outperforms Global for the lowest deciles.

\subsubsection{Group size}
\label{sec:compGroup}

\begin{figure}[H]
	\centering
	

	\begin{subfigure}[t]{.2\textwidth}
		\includegraphics[width=\textwidth]{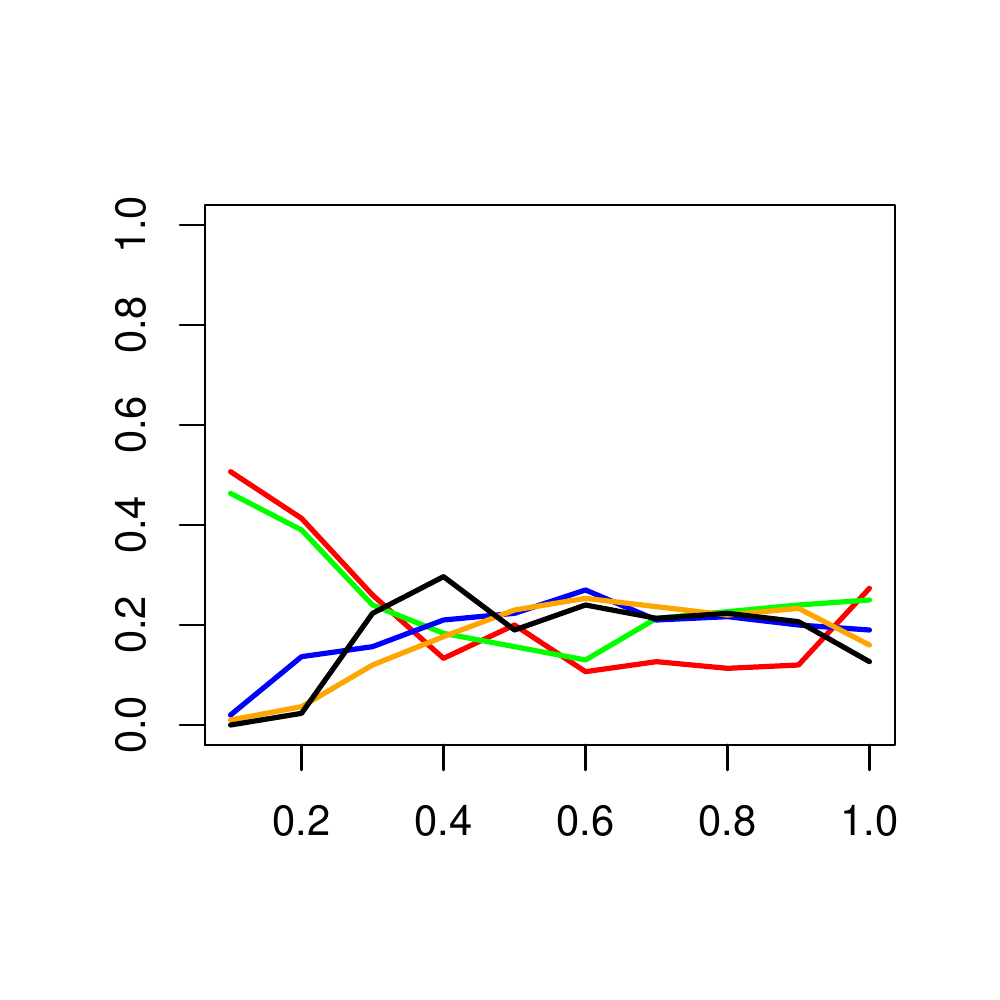}
		\vspace{-10mm} 
  	\caption{\scriptsize{Rastrigin}} \label{fig:gr30Ra}
	\end{subfigure}%
	\hspace{-6mm} 
	\begin{subfigure}[t]{.2\textwidth}
		\includegraphics[width=\textwidth]{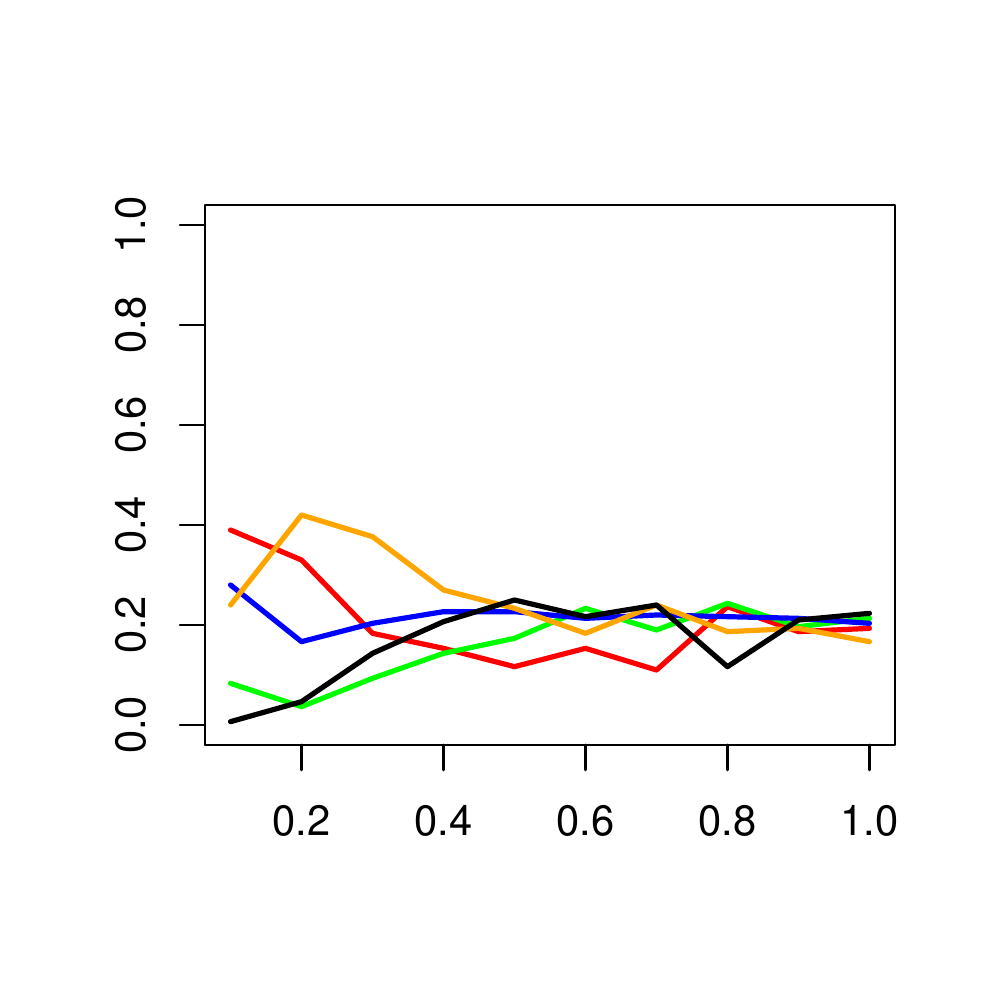}
		\vspace{-10mm} 
  	\caption{\scriptsize{Multipeak F1}} \label{fig:gr30M1}
	\end{subfigure}%
	\hspace{-6mm} 
	\begin{subfigure}[t]{.2\textwidth}
		\includegraphics[width=\textwidth]{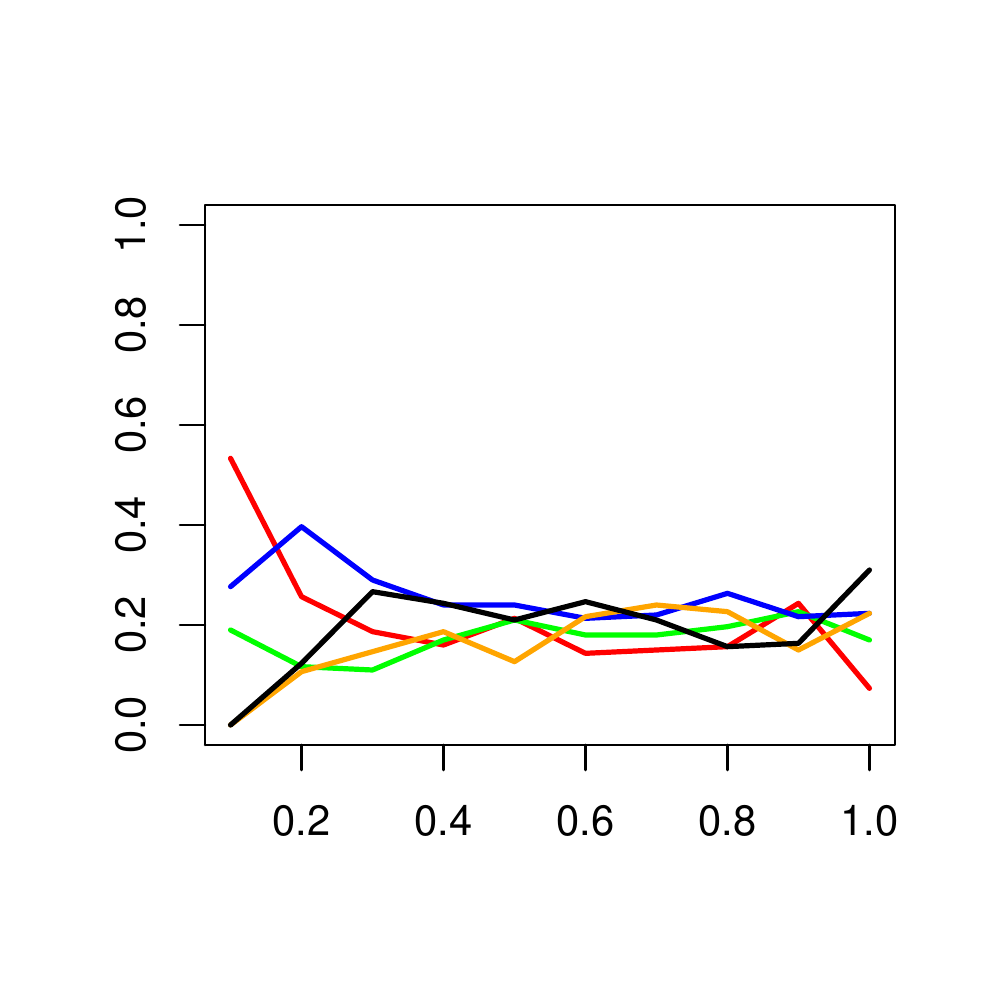}
		\vspace{-10mm} 
  	\caption{\scriptsize{Multipeak F2}} \label{fig:gr30M2}
	\end{subfigure}%
	\hspace{-6mm} 
	\begin{subfigure}[t]{.2\textwidth}
		\includegraphics[width=\textwidth]{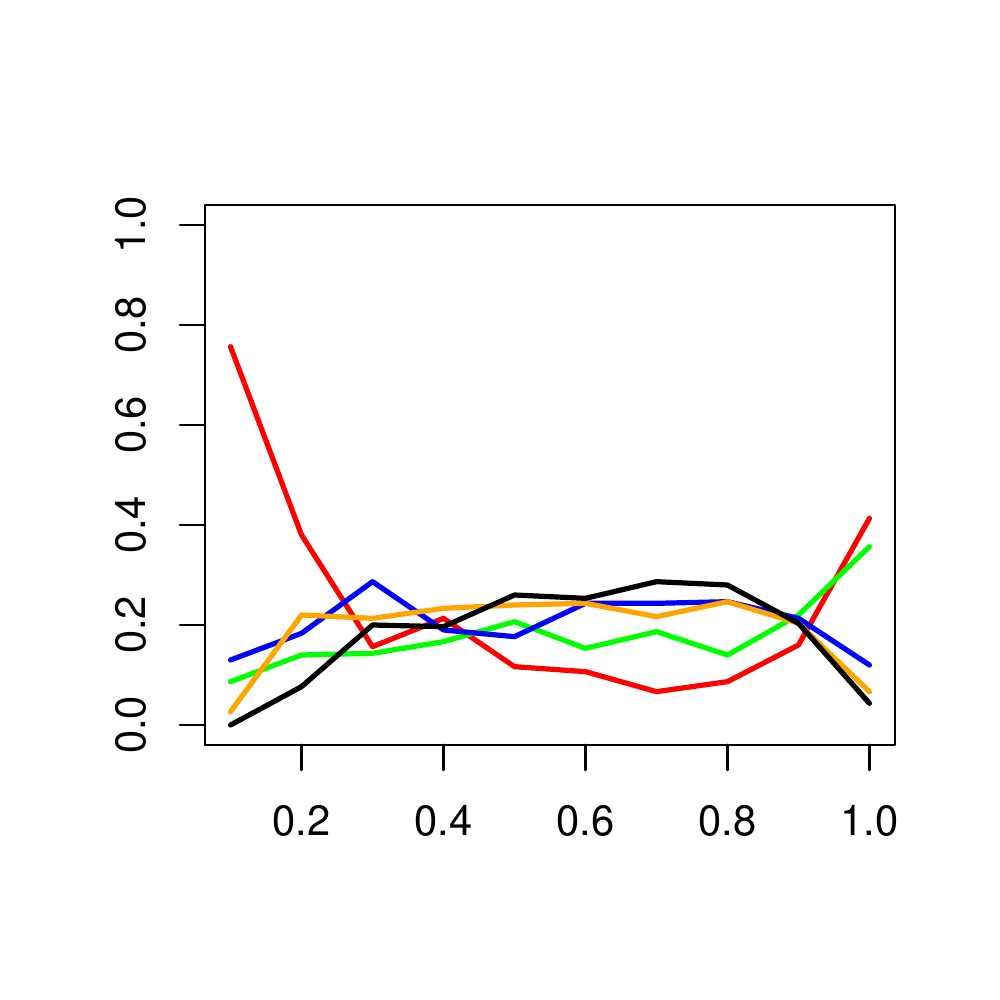}
		\vspace{-10mm} 
  	\caption{\scriptsize{Brankes}} \label{fig:gr30Br}
	\end{subfigure}%
	\hspace{-6mm} 
	\begin{subfigure}[t]{.2\textwidth}
		\includegraphics[width=\textwidth]{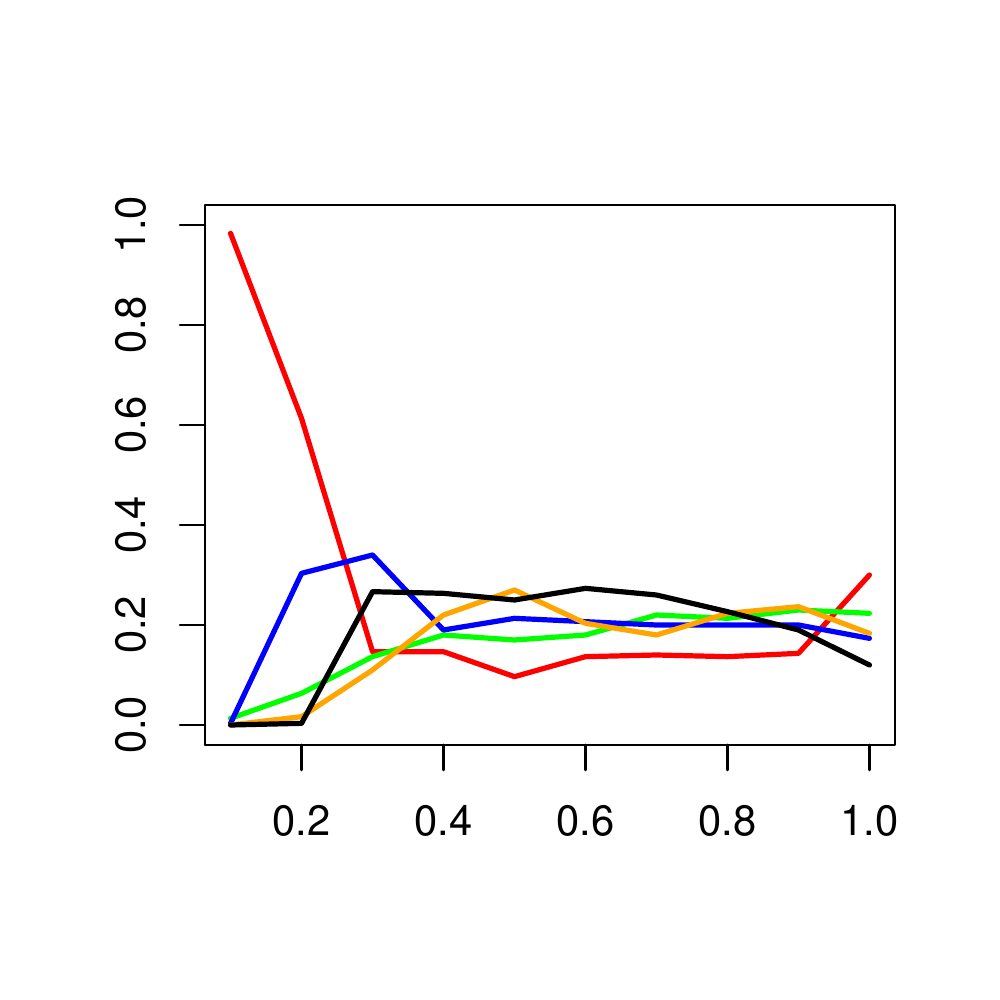}
		\vspace{-10mm} 
  	\caption{\scriptsize{Pickelhaube}} \label{fig:gr30Pi}
	\end{subfigure}
	
	\vspace{-2mm} 
		
	\begin{subfigure}[t]{.2\textwidth}
		\includegraphics[width=\textwidth]{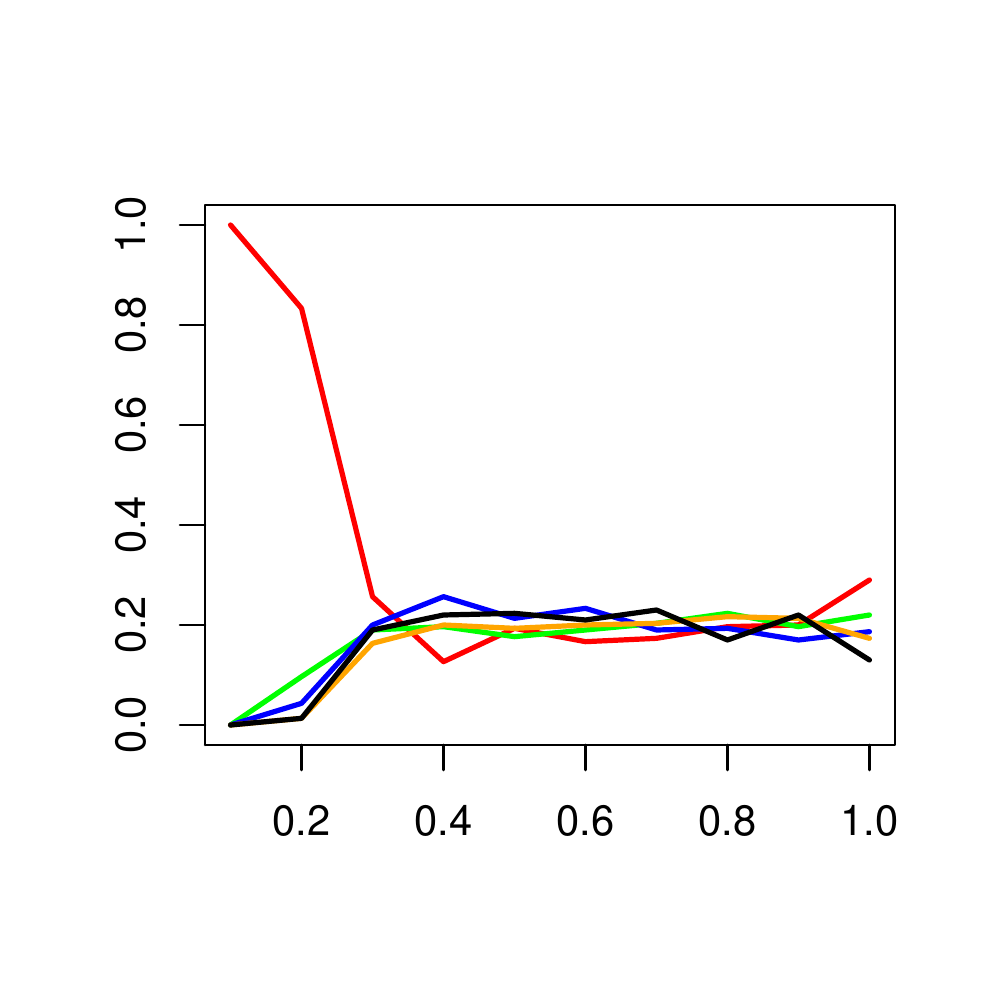}
		\vspace{-10mm} 
  	\caption{\scriptsize{Heaviside}} \label{fig:gr30Hv}
	\end{subfigure}%
	\hspace{-6mm} 
	\begin{subfigure}[t]{.2\textwidth}
		\includegraphics[width=\textwidth]{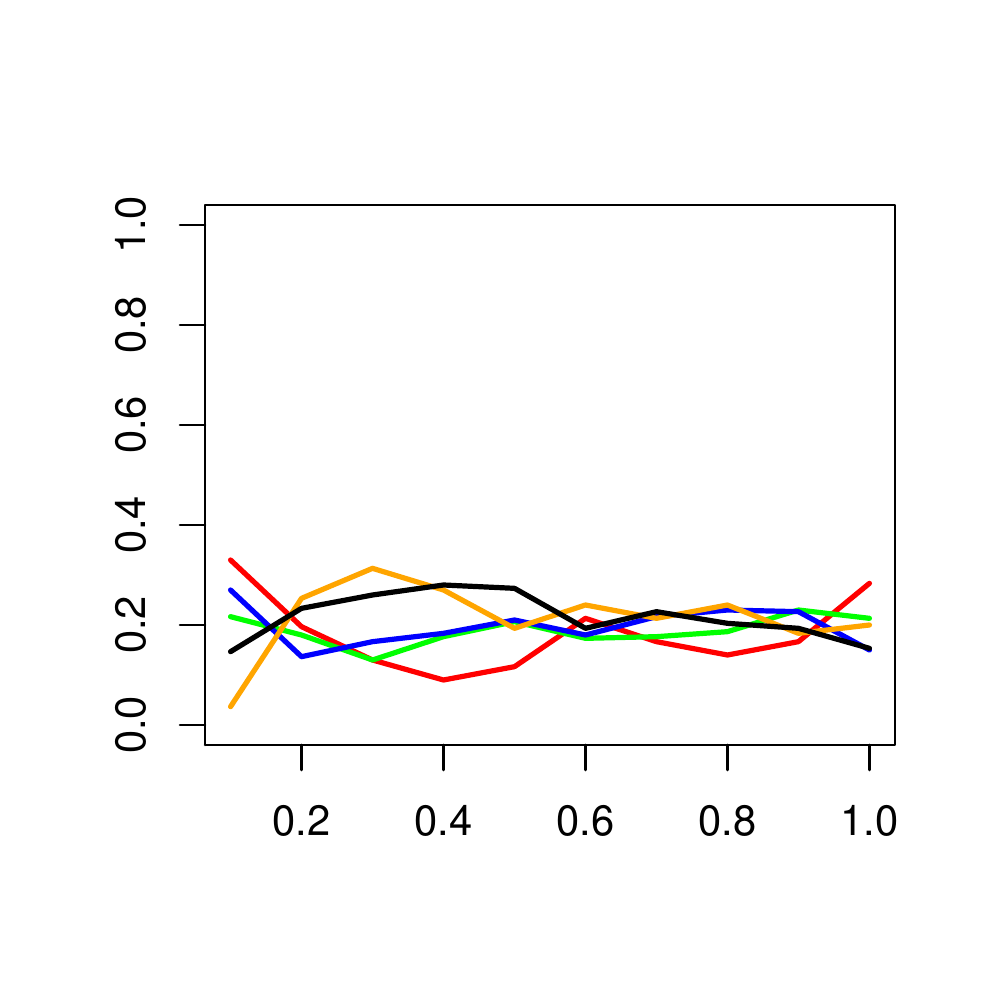}
		\vspace{-10mm} 
  	\caption{\scriptsize{Sawtooth}} \label{fig:gr30Sa}
	\end{subfigure}%
	\hspace{-6mm} 
	\begin{subfigure}[t]{.2\textwidth}
		\includegraphics[width=\textwidth]{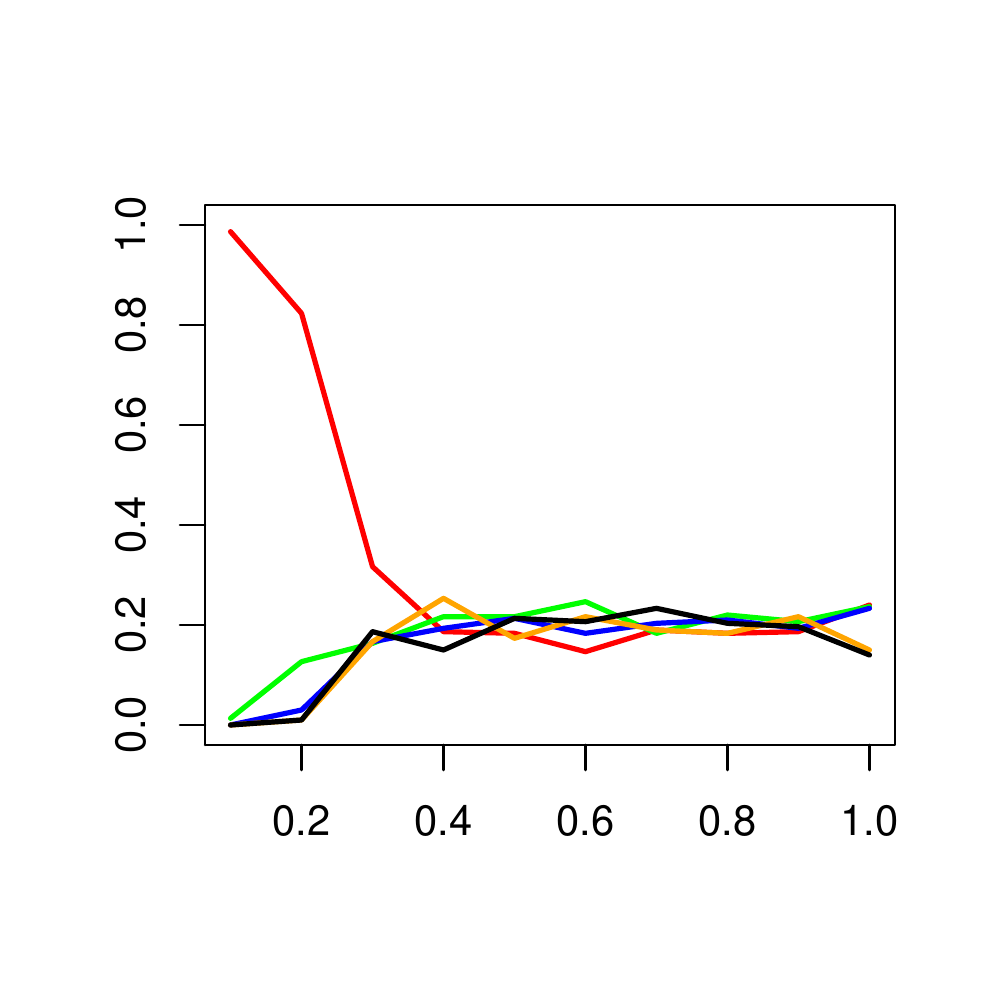}
		\vspace{-10mm} 
  	\caption{\scriptsize{Ackley}} \label{fig:gr30Ac}
	\end{subfigure}%
	\hspace{-6mm} 
	\begin{subfigure}[t]{.2\textwidth}
		\includegraphics[width=\textwidth]{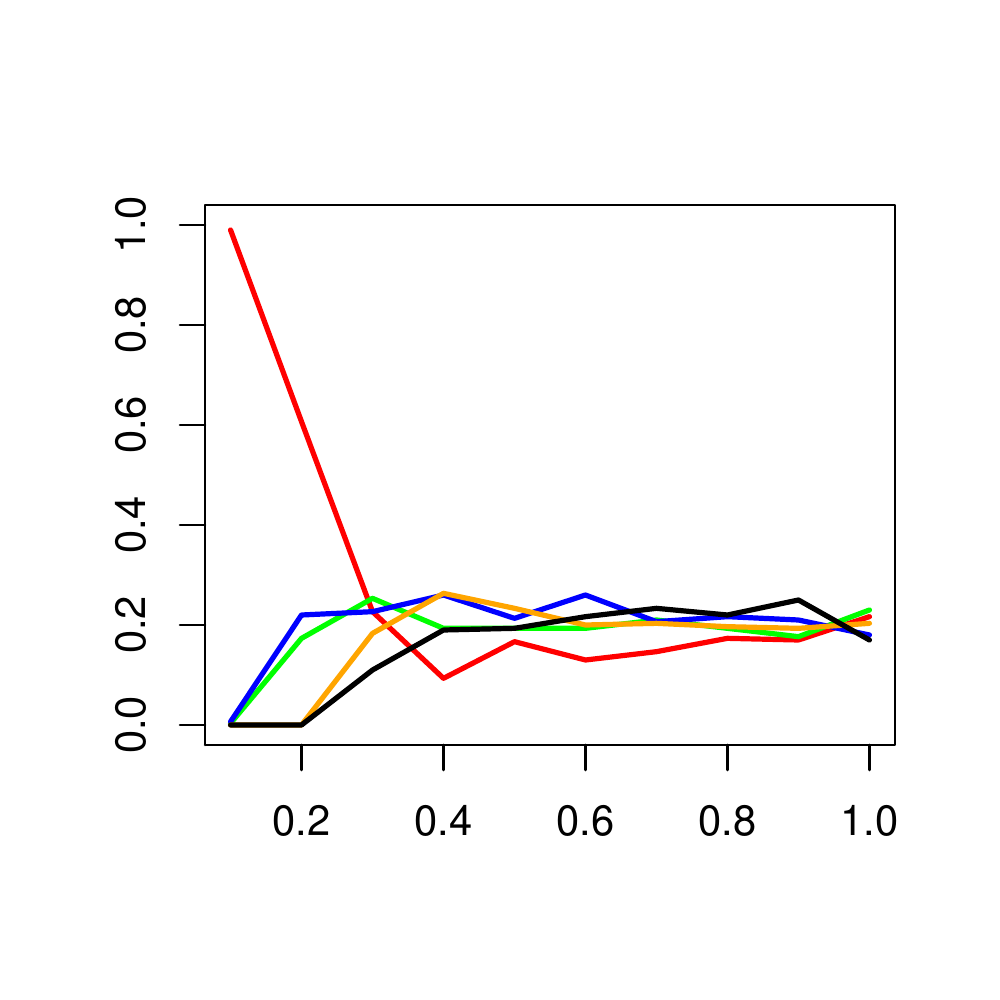}
		\vspace{-10mm} 
  	\caption{\scriptsize{Sphere}} \label{fig:gr30Sp}
	\end{subfigure}%
	\hspace{-6mm} 
	\begin{subfigure}[t]{.2\textwidth}
		\includegraphics[width=\textwidth]{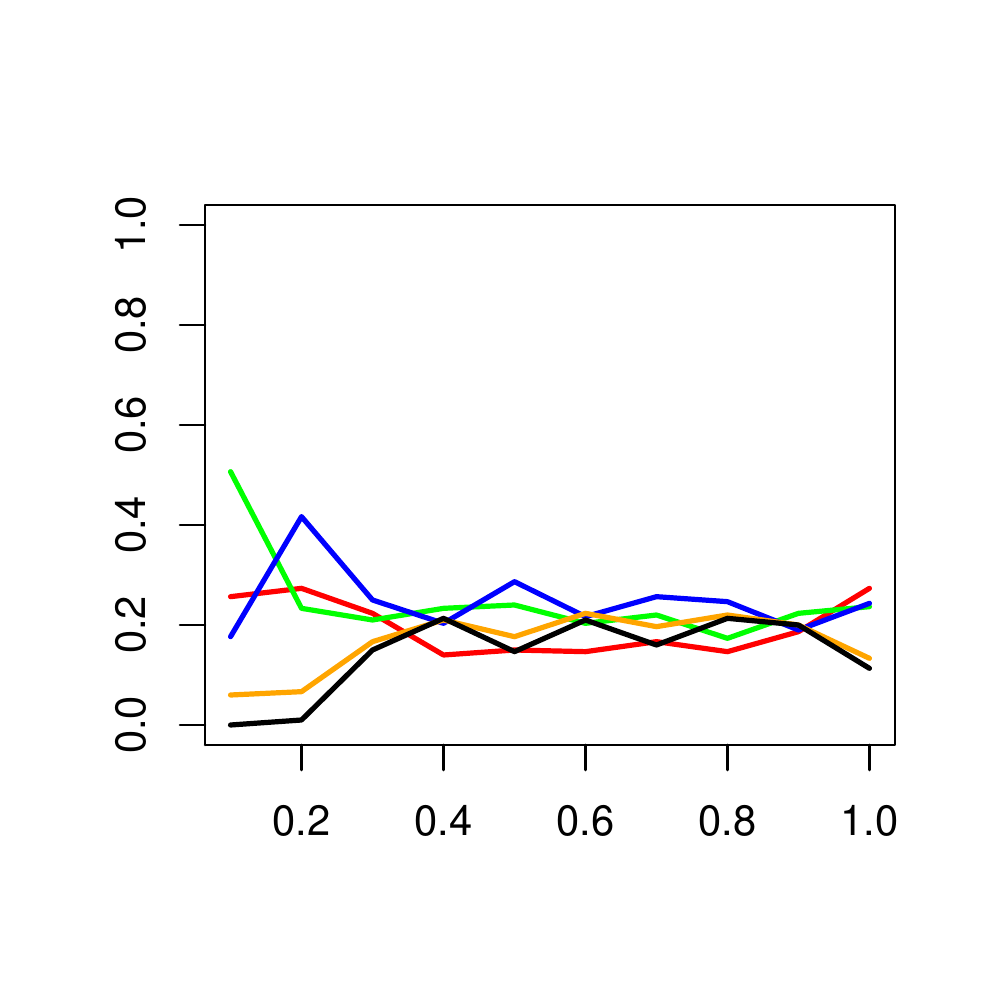}
		\vspace{-10mm} 
  	\caption{\scriptsize{Rosenbrock}} \label{fig:gr30Ro}
	\end{subfigure}
		
	\caption{Component -- decile breakdowns for the group (swarm) size, across all GGGP heuristics at 30D. Components: [2 -- 10] (red), [11 -- 20] (green), [21 -- 30] (blue), [31 -- 40] (orange), \textgreater \ 40 (black).}
	\label{fig:group30Comp}
	
\end{figure}

\begin{figure}[H]
	\centering
	
	\vspace{-5mm} 

	\begin{subfigure}[t]{.2\textwidth}
		\includegraphics[width=\textwidth]{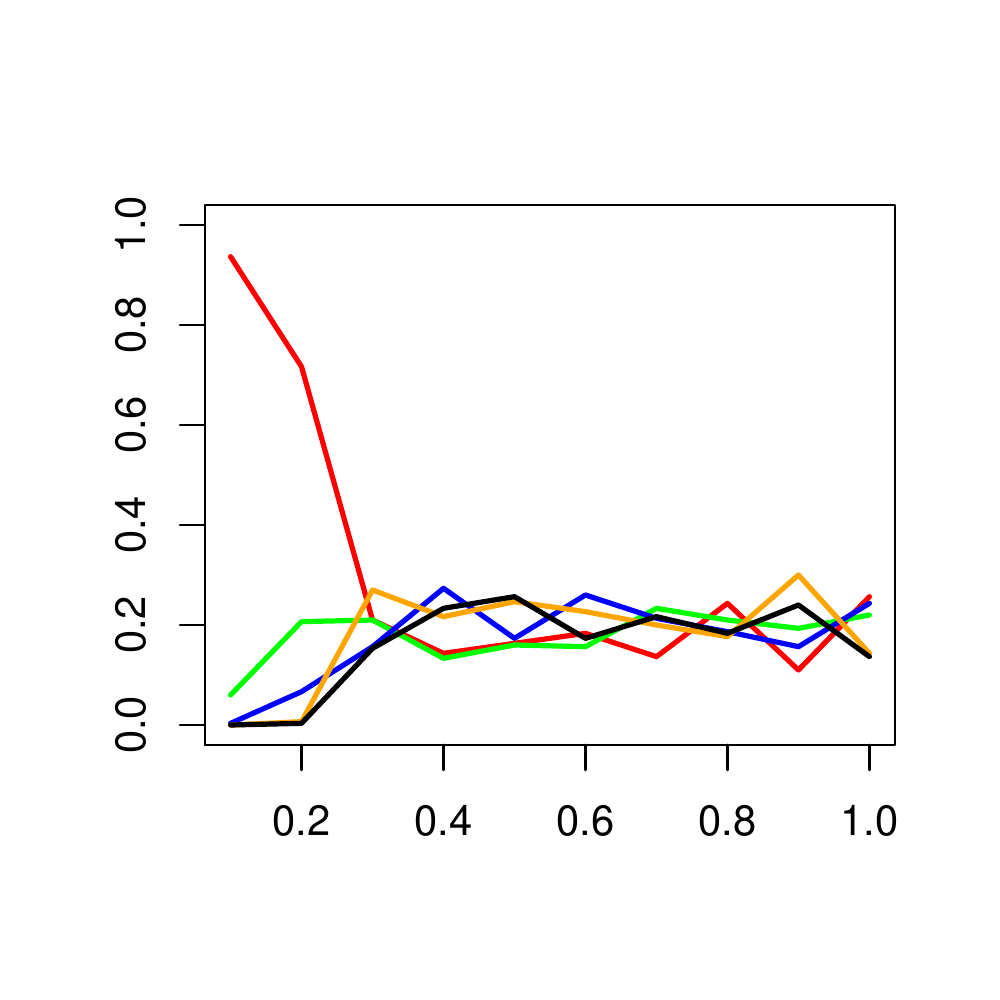}
		\vspace{-10mm} 
  	\caption{\scriptsize{Rastrigin}} \label{fig:gr100Ra}
	\end{subfigure}%
	\hspace{-6mm} 
	\begin{subfigure}[t]{.2\textwidth}
		\includegraphics[width=\textwidth]{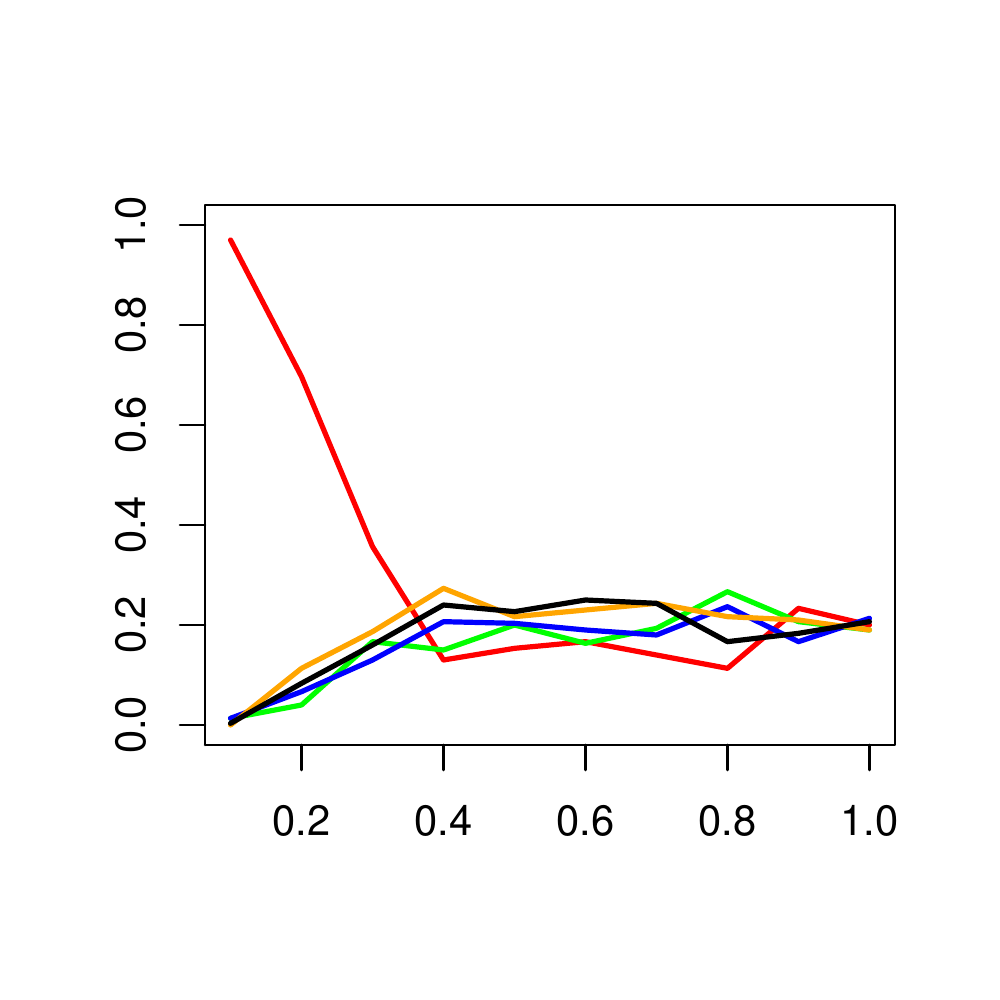}
		\vspace{-10mm} 
  	\caption{\scriptsize{Multipeak F1}} \label{fig:gr100M1}
	\end{subfigure}%
	\hspace{-6mm} 
	\begin{subfigure}[t]{.2\textwidth}
		\includegraphics[width=\textwidth]{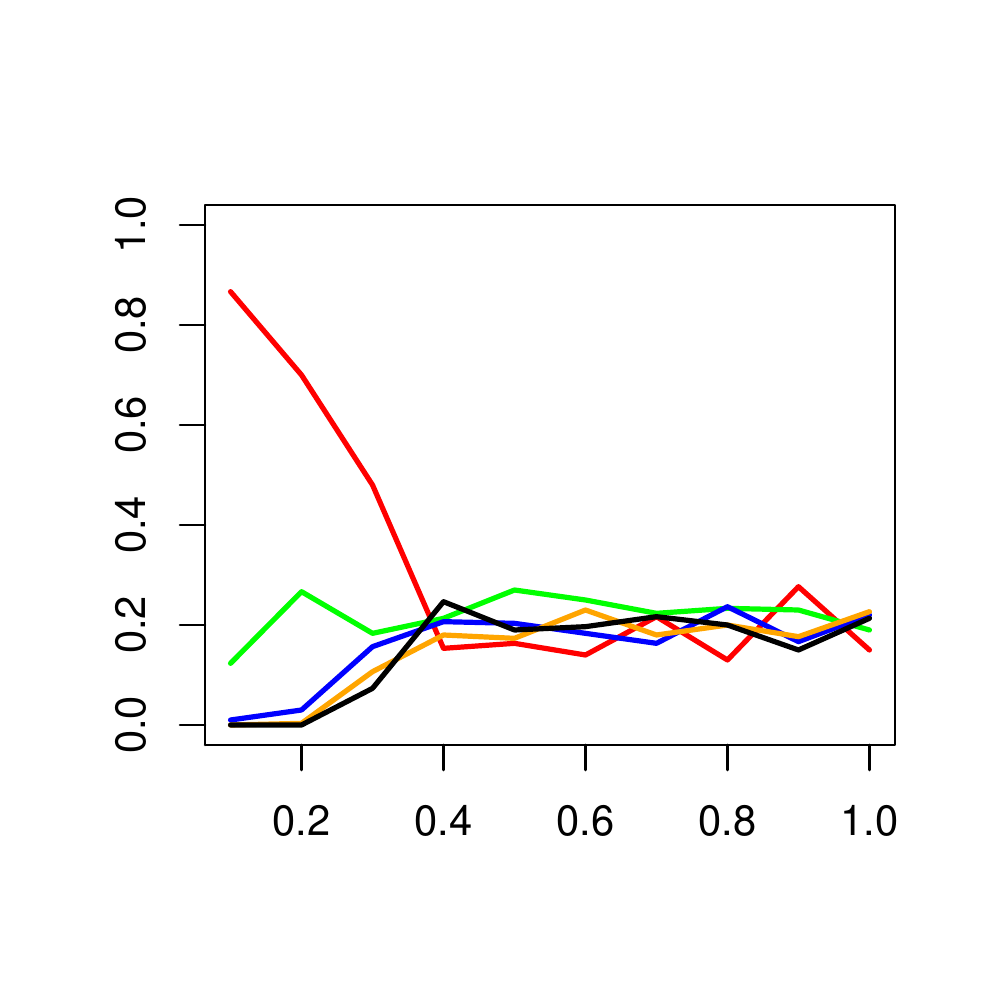}
		\vspace{-10mm} 
  	\caption{\scriptsize{Multipeak F2}} \label{fig:gr100M2}
	\end{subfigure}%
	\hspace{-6mm} 
	\begin{subfigure}[t]{.2\textwidth}
		\includegraphics[width=\textwidth]{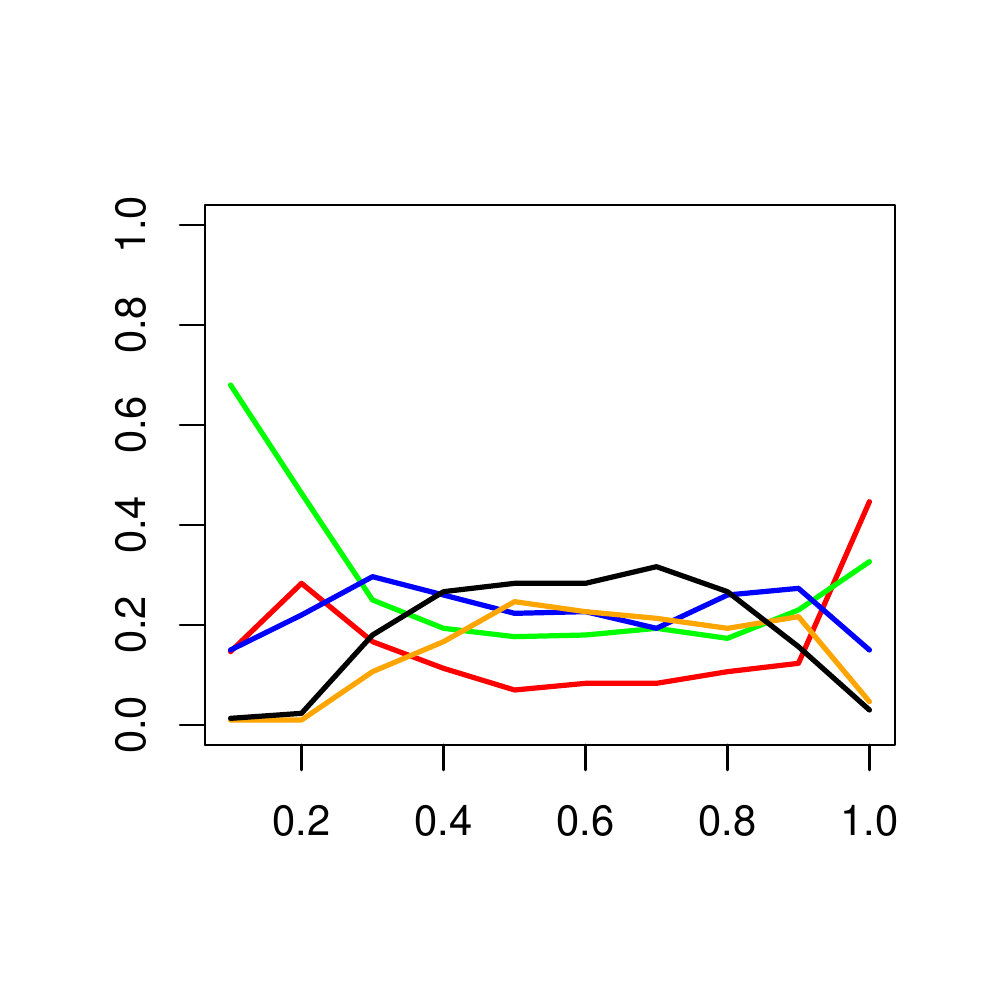}
		\vspace{-10mm} 
  	\caption{\scriptsize{Brankes}} \label{fig:gr100Br}
	\end{subfigure}%
	\hspace{-6mm} 
	\begin{subfigure}[t]{.2\textwidth}
		\includegraphics[width=\textwidth]{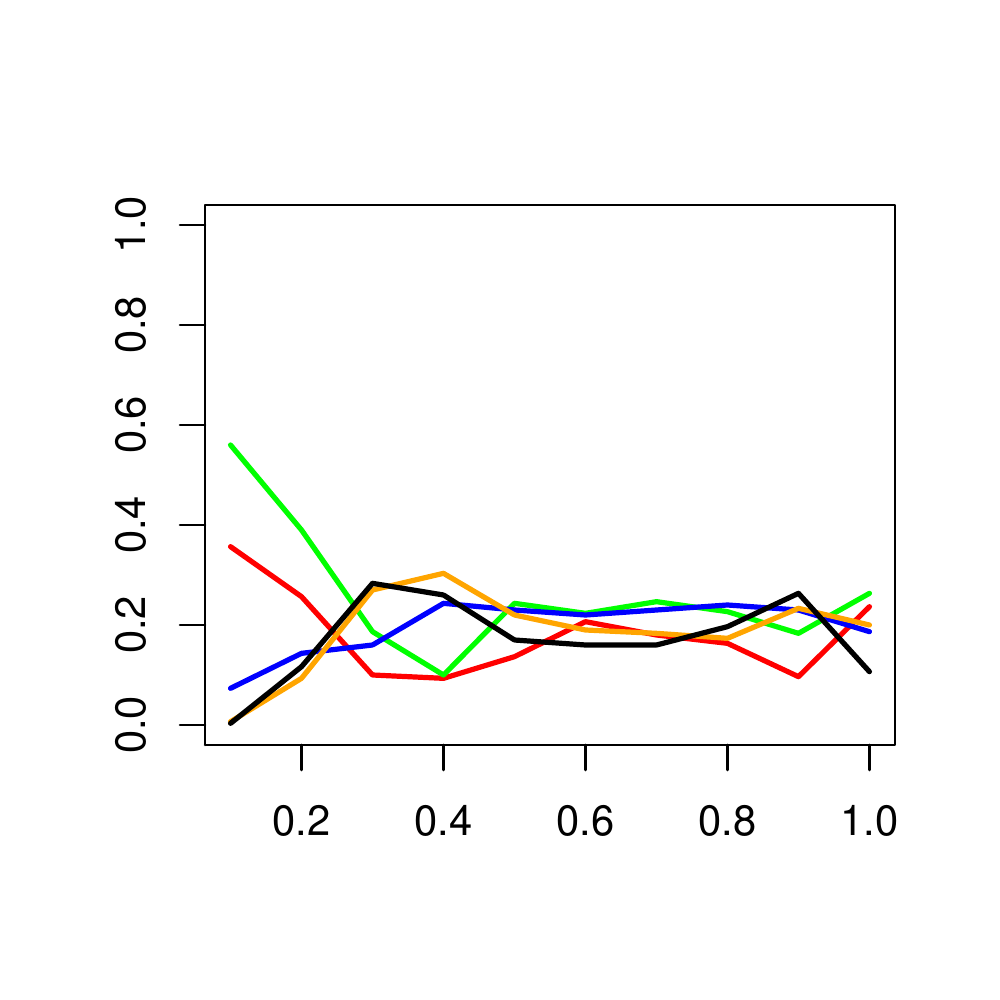}
		\vspace{-10mm} 
  	\caption{\scriptsize{Pickelhaube}} \label{fig:gr100Pi}
	\end{subfigure}
	
	\vspace{-2mm} 
		
	\begin{subfigure}[t]{.2\textwidth}
		\includegraphics[width=\textwidth]{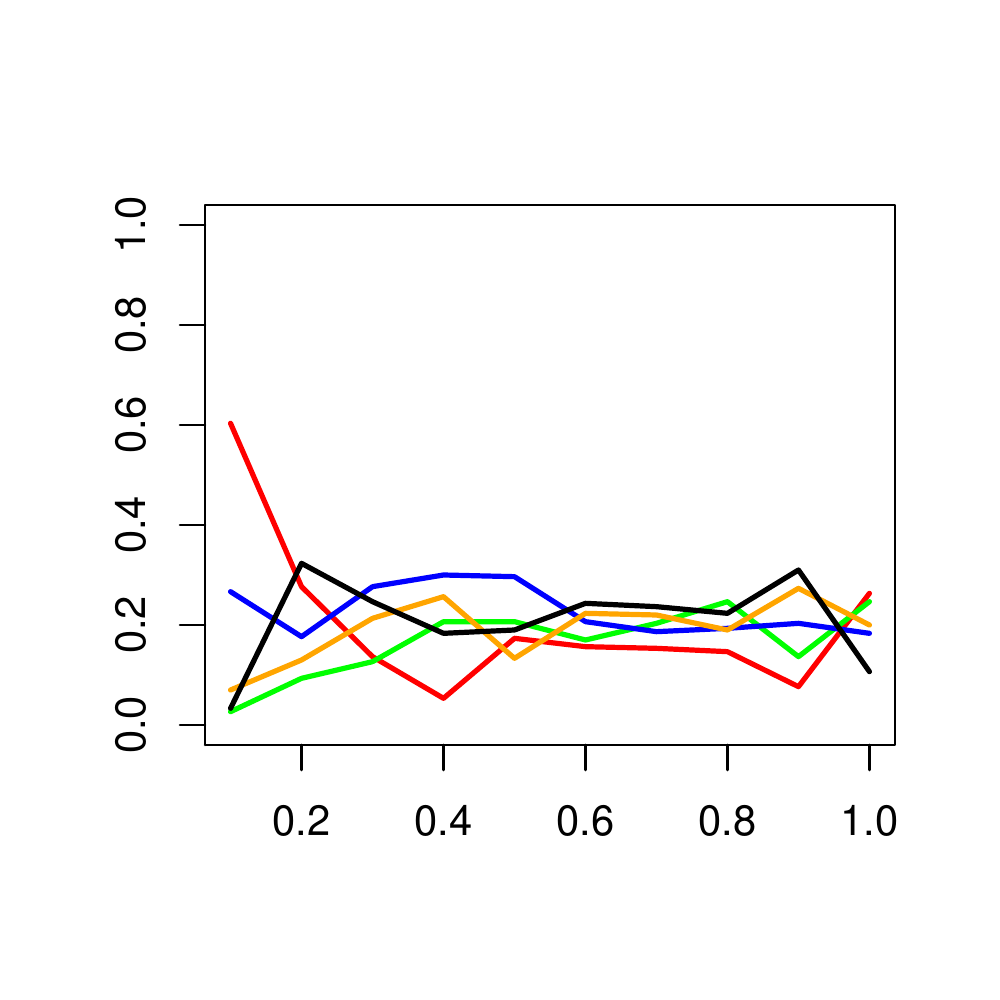}
		\vspace{-10mm} 
  	\caption{\scriptsize{Heaviside}} \label{fig:gr100Hv}
	\end{subfigure}%
	\hspace{-6mm} 
	\begin{subfigure}[t]{.2\textwidth}
		\includegraphics[width=\textwidth]{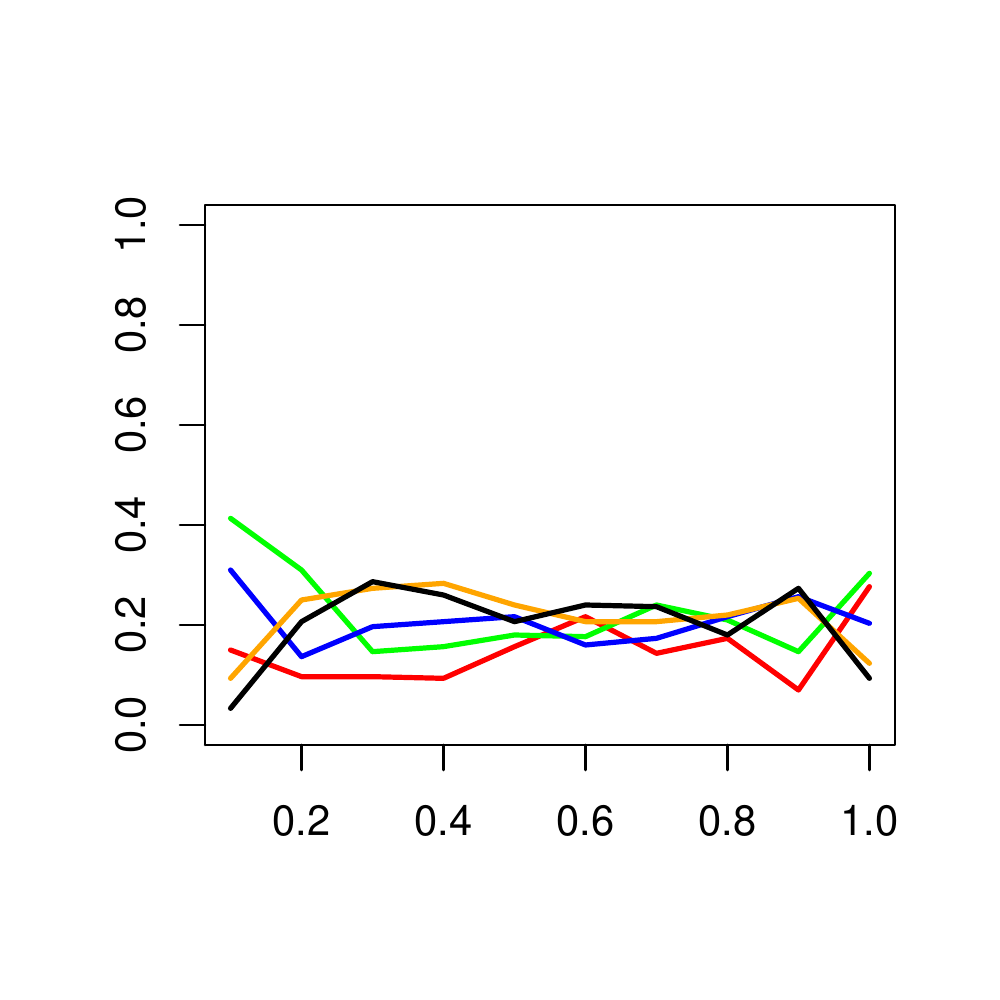}
		\vspace{-10mm} 
  	\caption{\scriptsize{Sawtooth}} \label{fig:gr100Sa}
	\end{subfigure}%
	\hspace{-6mm} 
	\begin{subfigure}[t]{.2\textwidth}
		\includegraphics[width=\textwidth]{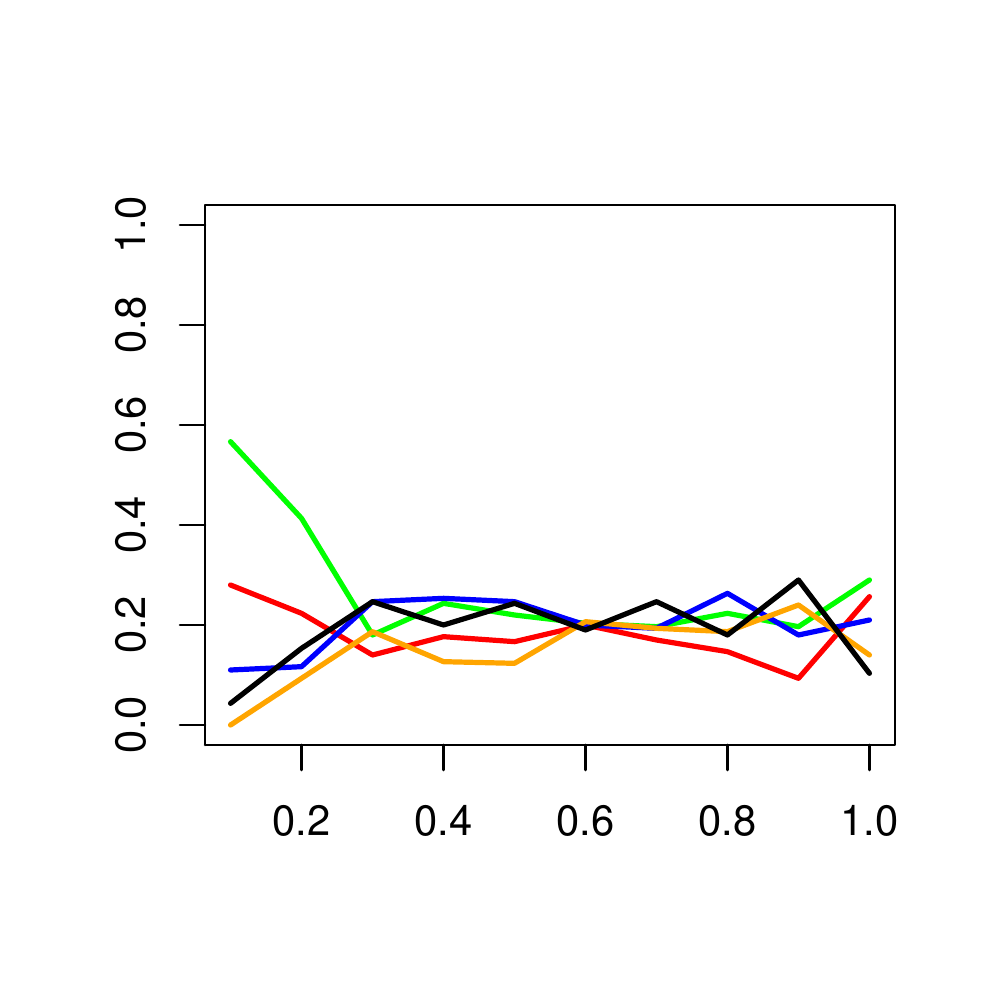}
		\vspace{-10mm} 
  	\caption{\scriptsize{Ackley}} \label{fig:gr100Ac}
	\end{subfigure}%
	\hspace{-6mm} 
	\begin{subfigure}[t]{.2\textwidth}
		\includegraphics[width=\textwidth]{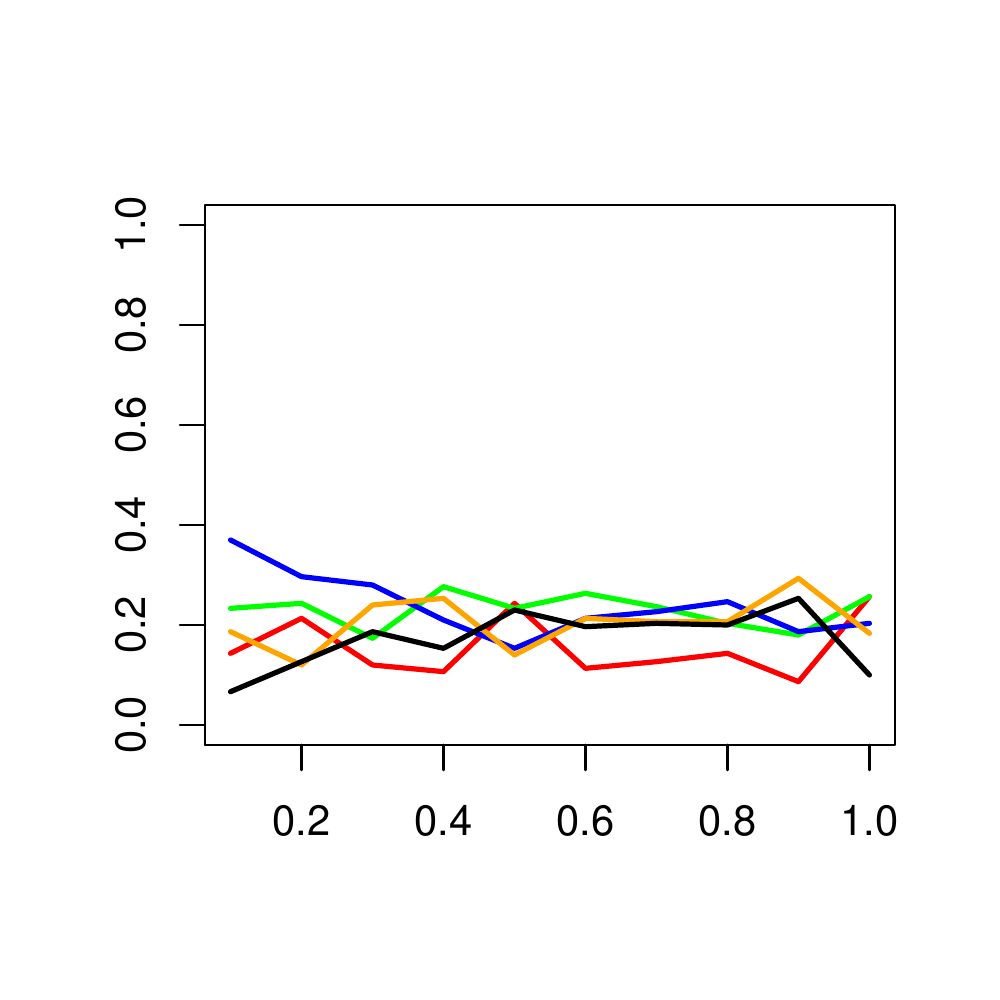}
		\vspace{-10mm} 
  	\caption{\scriptsize{Sphere}} \label{fig:gr100Sp}
	\end{subfigure}%
	\hspace{-6mm} 
	\begin{subfigure}[t]{.2\textwidth}
		\includegraphics[width=\textwidth]{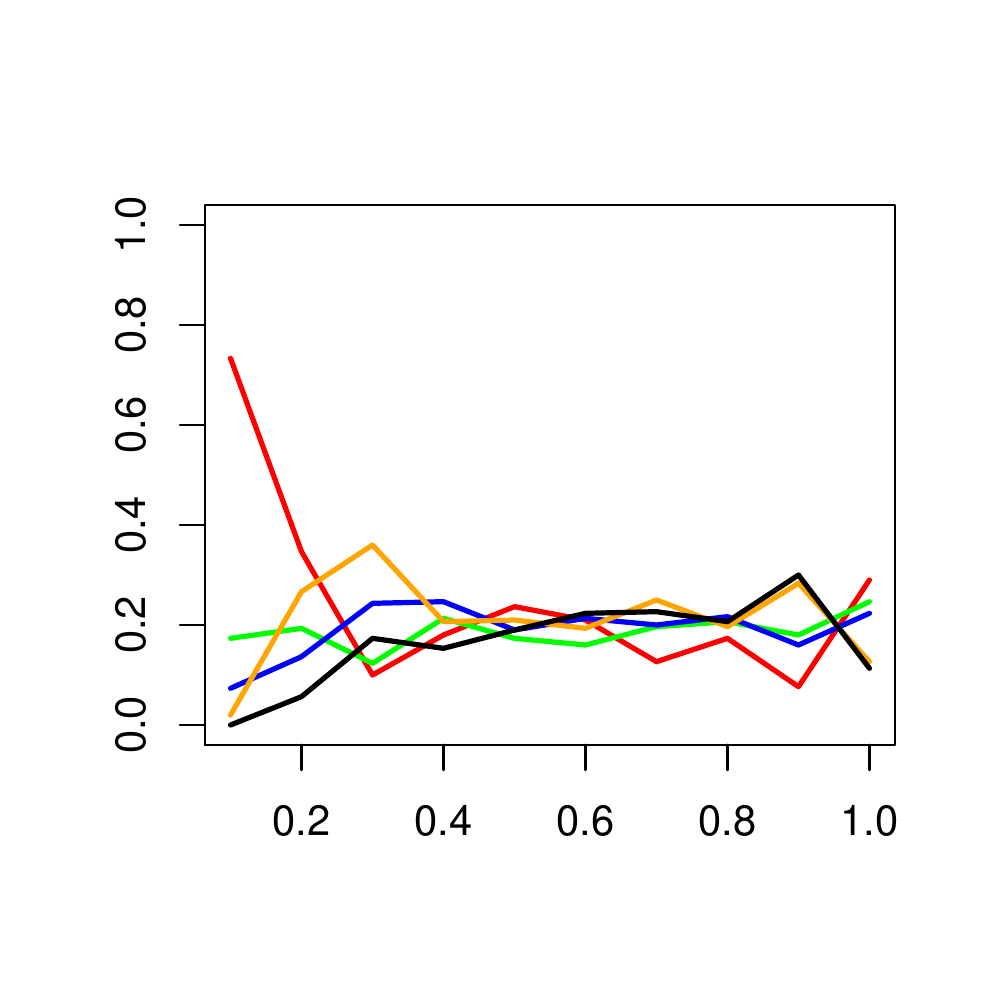}
		\vspace{-10mm} 
  	\caption{\scriptsize{Rosenbrock}} \label{fig:gr100Ro}
	\end{subfigure}
		
	\caption{Component -- decile breakdowns for the group (swarm) size, across all GGGP heuristics at 100D. Components: [2 -- 10] (red), [11 -- 20] (green), [21 -- 30] (blue), [31 -- 40] (orange), \textgreater \ 40 (black).}
	\label{fig:group100Comp}
	
\end{figure}

\begin{figure}[H]
	\centering
	
	\vspace{-5mm} 

	\begin{subfigure}[t]{.24\textwidth}
		\includegraphics[width=\textwidth]{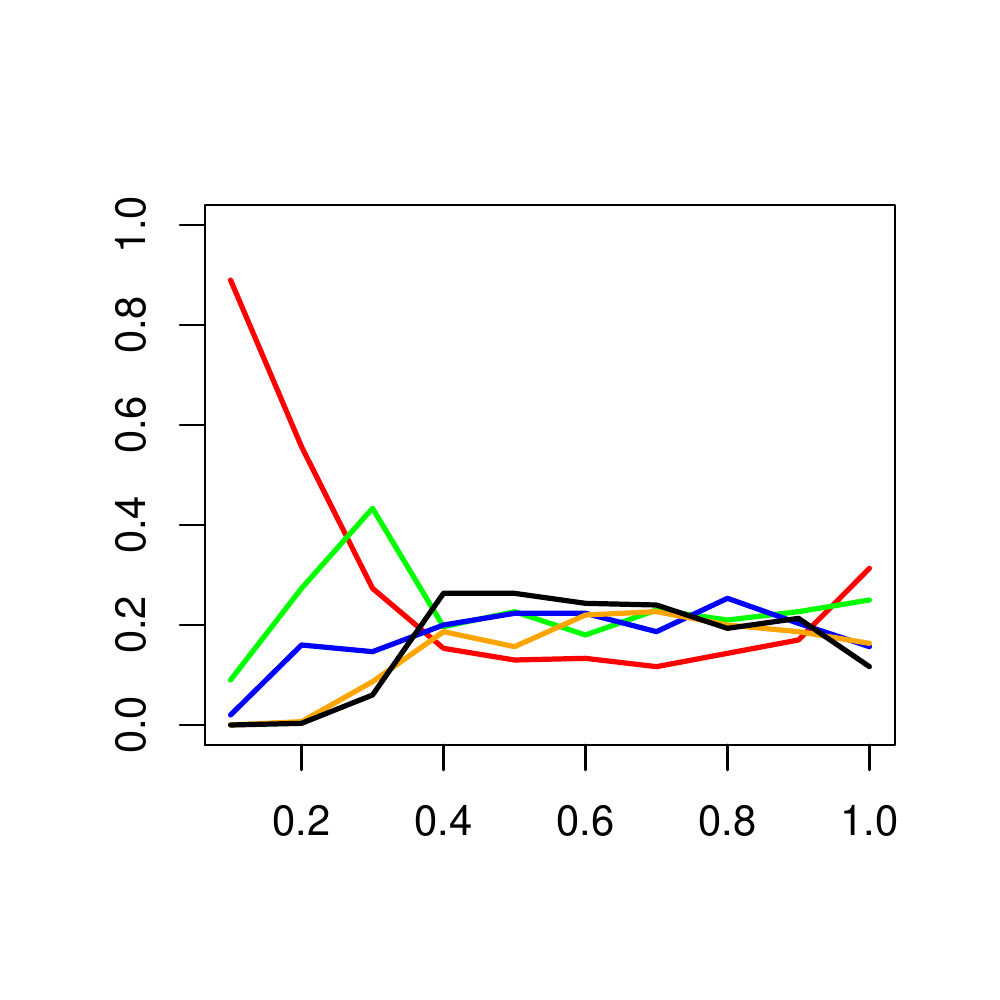}
		\vspace{-10mm} 
  	\caption{\scriptsize{30D}} \label{fig:gr30}
	\end{subfigure}%
	\hspace{-6mm} 
	\begin{subfigure}[t]{.24\textwidth}
		\includegraphics[width=\textwidth]{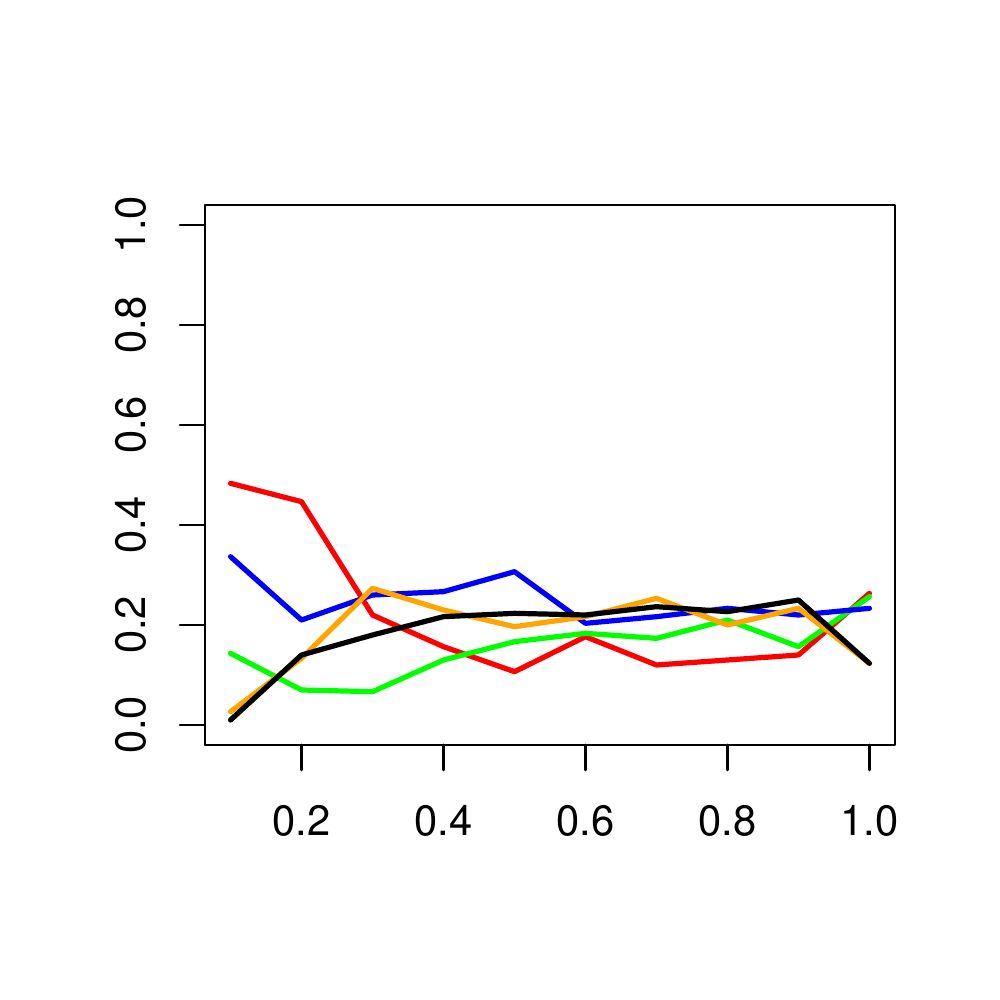}
		\vspace{-10mm} 
  	\caption{\scriptsize{100D}} \label{fig:gr100}
	\end{subfigure}

	\caption{Component -- decile breakdowns for the group (swarm) size, across all GGGP heuristics at 30D and 100D for the general heuristics. Components: [2 -- 10] (red), [11 -- 20] (green), [21 -- 30] (blue), [31 -- 40] (orange), \textgreater \ 40 (black).}
	\label{fig:groupComp}
	
\end{figure}

For the group (outer PSO swarm) size, in Table~\ref{fig:proportions} all of the different categories are quite well represented considering all heuristics generated in the GP runs. Most common is the lowest range, 2-10 (red), and in the best performing third of heuristics this range stands out somewhat, increasing to 43\% and 34\% for 30D and 100D respectively. A similar pattern is observed for the individual case heuristics decile analysis in Figures~\ref{fig:group30Comp} and~\ref{fig:group100Comp}, although the next group size range, 11-20 (green), is also favoured in the lowest deciles for a few problems and in particular at 100D. For the general case  at 100D, Figure~\ref{fig:groupComp}, the lower decile results are well distributed across the group size ranges.

\vspace{5mm} 

\subsubsection{Use of a stopping condition}
\label{sec:compStop}

Considering the inclusion of a stopping condition in the inner maximisation, from Table~\ref{fig:proportions} it can be seen that this is preferred in 63\% and 59\% of all heuristics from the GP runs, for 30D and 100D respectively. In the top third of results this increases to 84\% and 73\%. A decile level analysis confirms this dominance in the best performing heuristics across all GP runs.

\subsubsection{Use of neighbourhood information for dormancy}
\label{sec:compDormancy}

When a heuristic uses the dormancy-relocation capability, (+LEH or +DD+LEH), the assessment of dormancy may (Yes) or may not (No) make use of existing uncertainty neighbourhood information from the history set. Table~\ref{fig:proportions} indicates that the use or non-use of this information is quite evenly apportioned, particularly for 30D. These high level results are reflected at the individual case and general case decile level analysis. In a few cases the use of information performs better in the lowest deciles.

\subsubsection{Use of neighbourhood information for $\pb^j$}
\label{sec:compPbest}

For the use of uncertainty neighbourhood information in the history set to update a particle's robust value on completion of the inner maximisation search, the use of such information (Yes) versus non-use (No) is somewhat evenly apportioned -- although there is some limited preference for using the information. Again these results are reflected at the more detailed decile level analysis.

\subsubsection{Particle level mutation}
\label{sec:compMute}

\begin{figure}[H]
	\centering
	

	\begin{subfigure}[t]{.2\textwidth}
		\includegraphics[width=\textwidth]{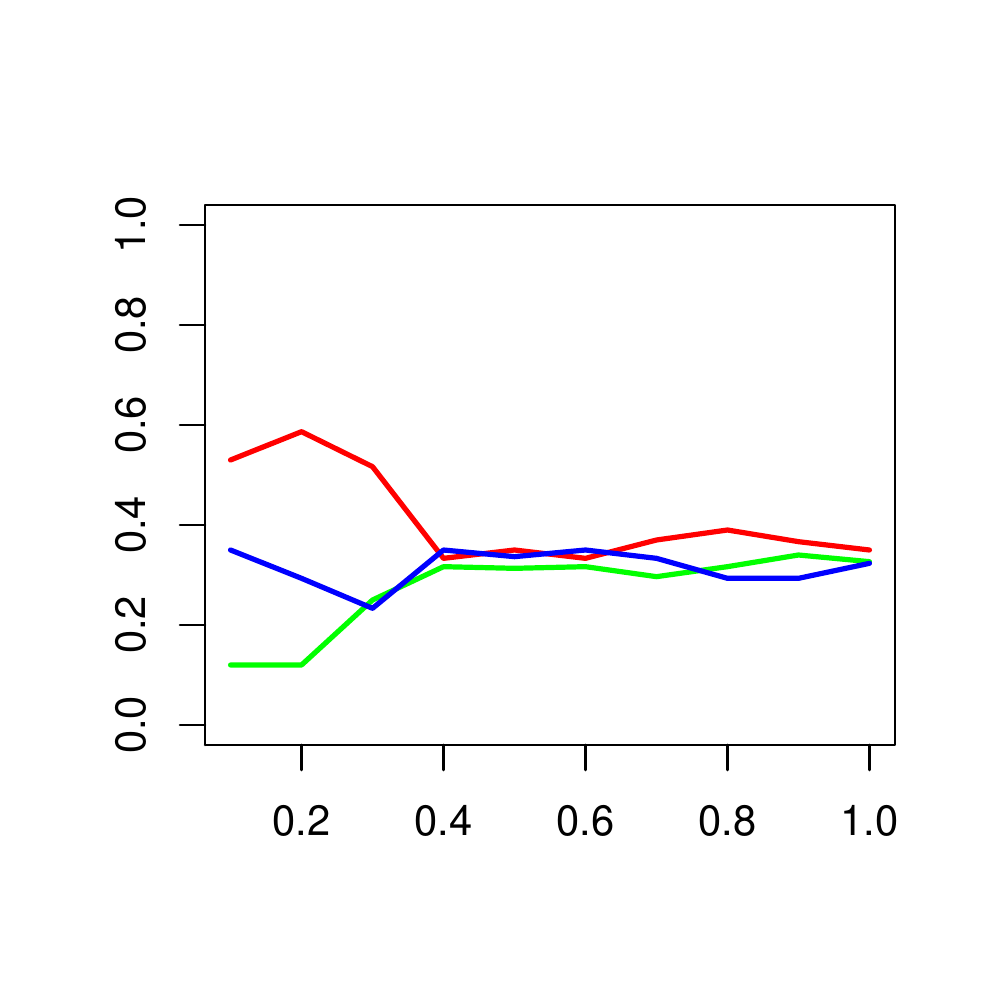}
		\vspace{-10mm} 
  	\caption{\scriptsize{Rastrigin}} \label{fig:mu30Ra}
	\end{subfigure}%
	\hspace{-6mm} 
	\begin{subfigure}[t]{.2\textwidth}
		\includegraphics[width=\textwidth]{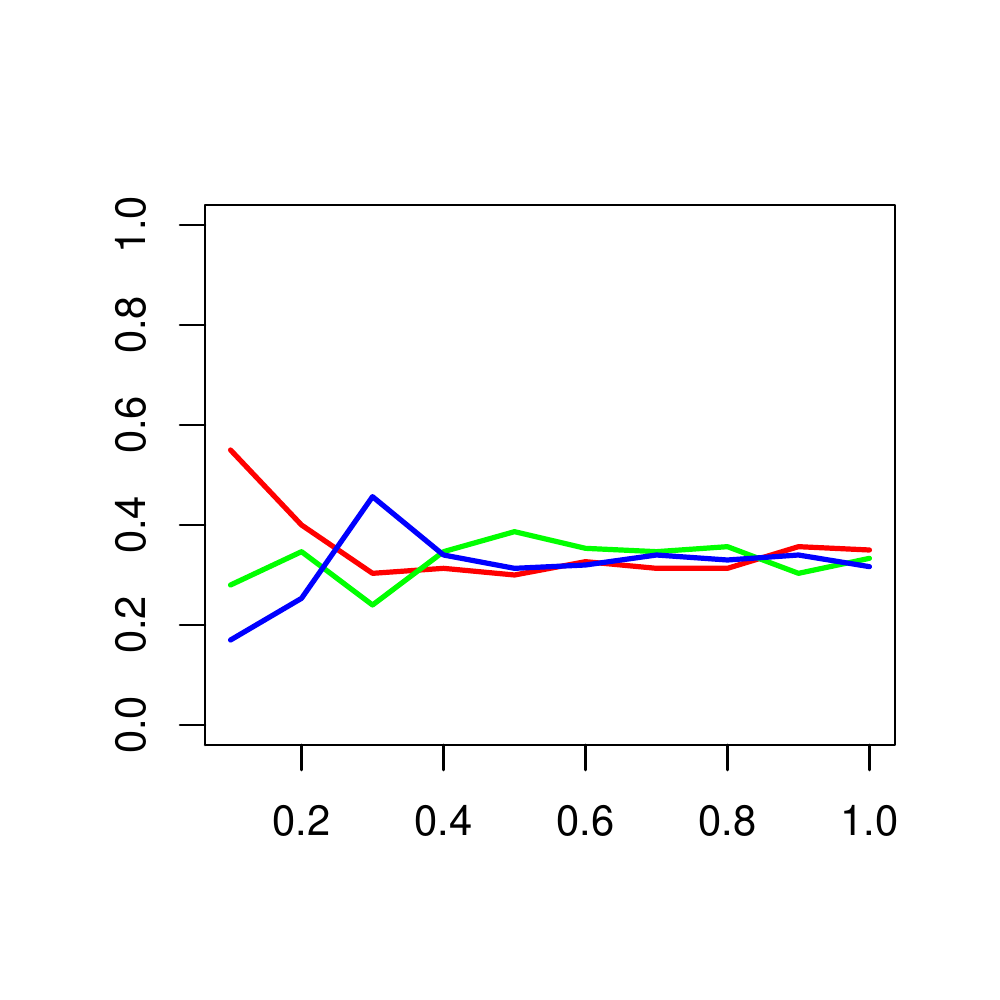}
		\vspace{-10mm} 
  	\caption{\scriptsize{Multipeak F1}} \label{fig:mu30M1}
	\end{subfigure}%
	\hspace{-6mm} 
	\begin{subfigure}[t]{.2\textwidth}
		\includegraphics[width=\textwidth]{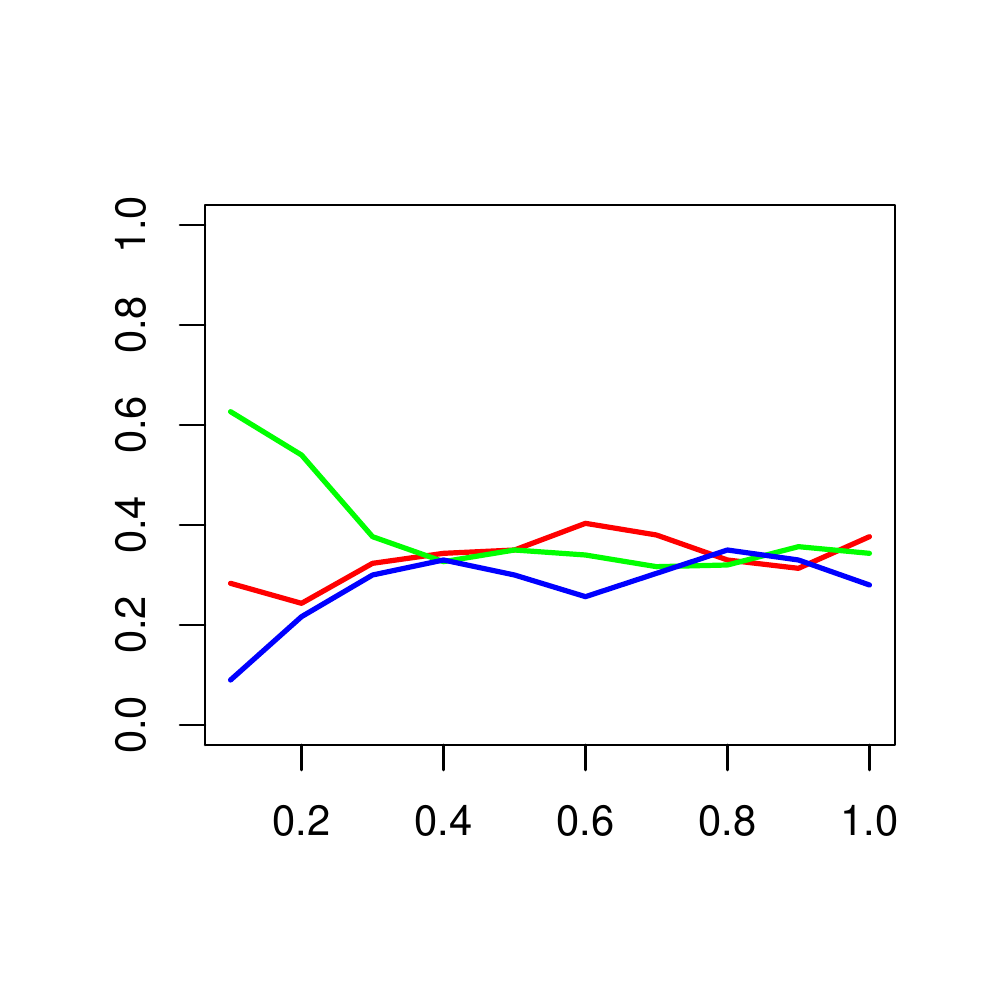}
		\vspace{-10mm} 
  	\caption{\scriptsize{Multipeak F2}} \label{fig:mu30M2}
	\end{subfigure}%
	\hspace{-6mm} 
	\begin{subfigure}[t]{.2\textwidth}
		\includegraphics[width=\textwidth]{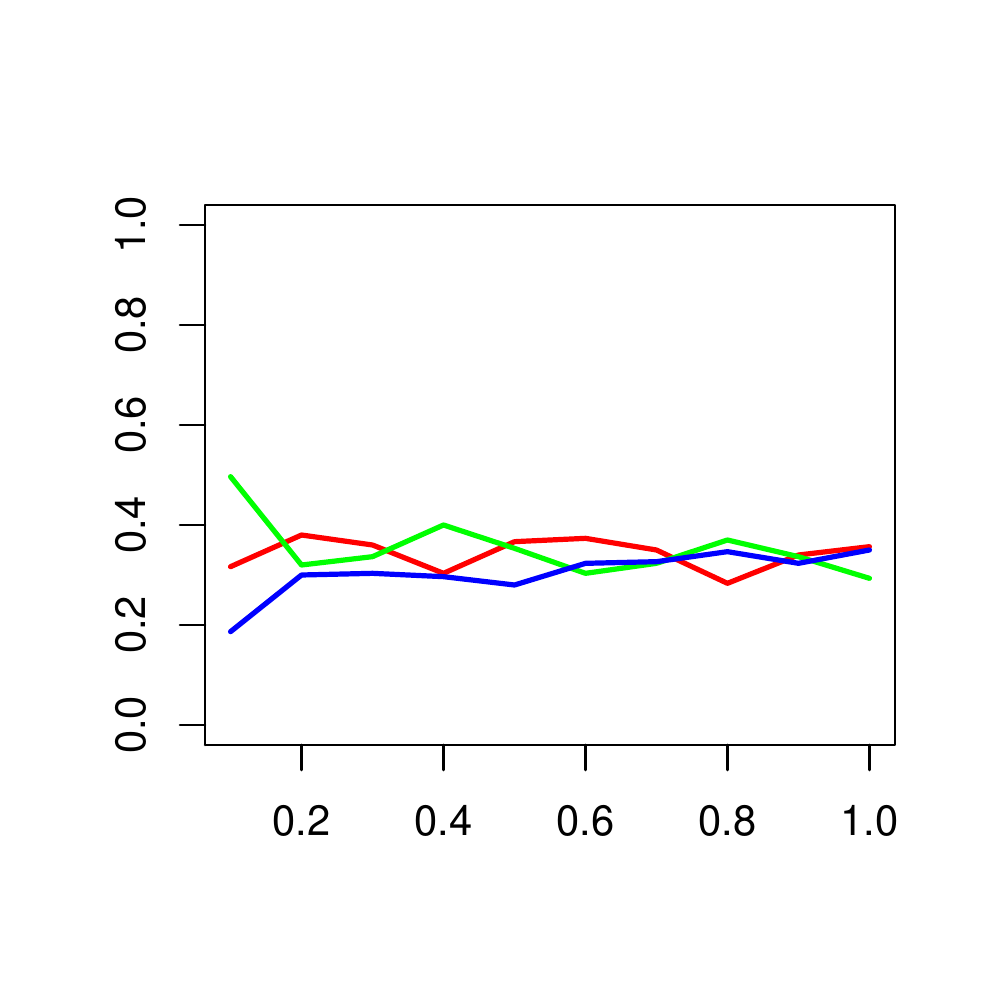}
		\vspace{-10mm} 
  	\caption{\scriptsize{Brankes}} \label{fig:mu30Br}
	\end{subfigure}%
	\hspace{-6mm} 
	\begin{subfigure}[t]{.2\textwidth}
		\includegraphics[width=\textwidth]{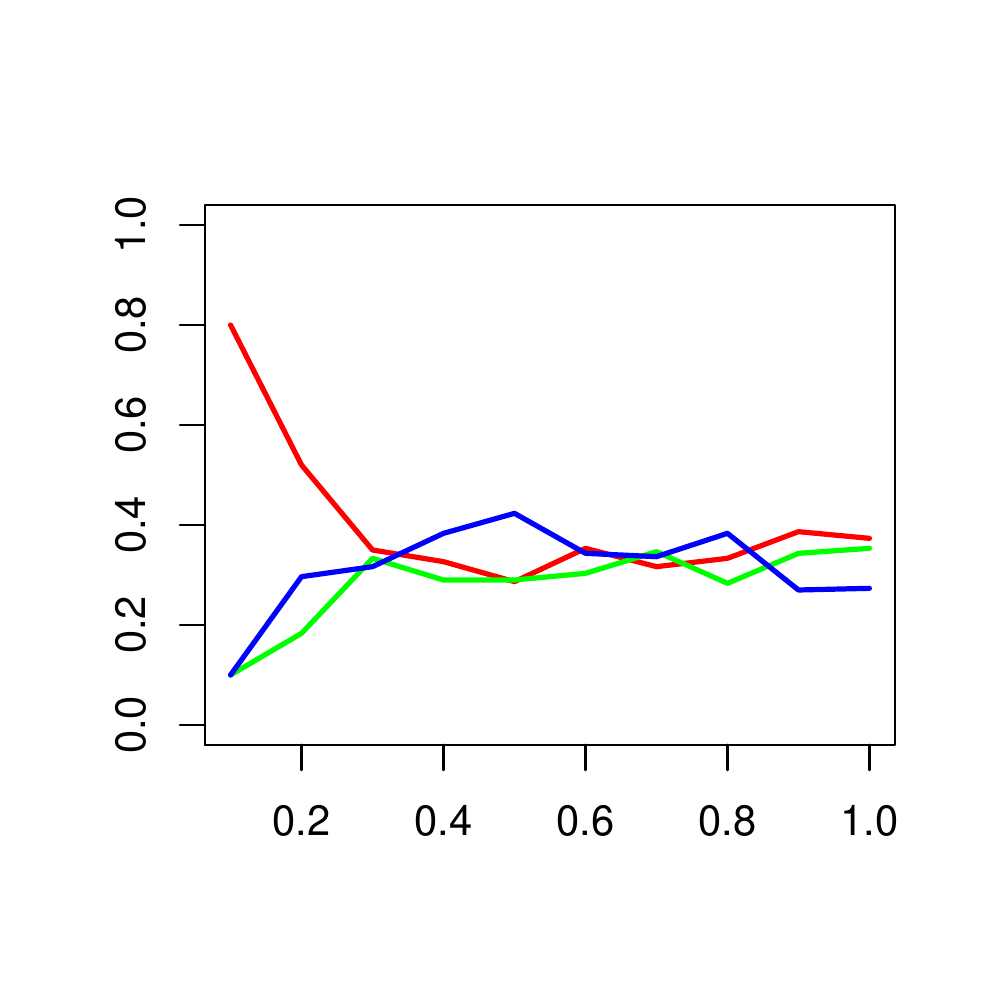}
		\vspace{-10mm} 
  	\caption{\scriptsize{Pickelhaube}} \label{fig:mu30Pi}
	\end{subfigure}
	
	\vspace{-2mm} 
		
	\begin{subfigure}[t]{.2\textwidth}
		\includegraphics[width=\textwidth]{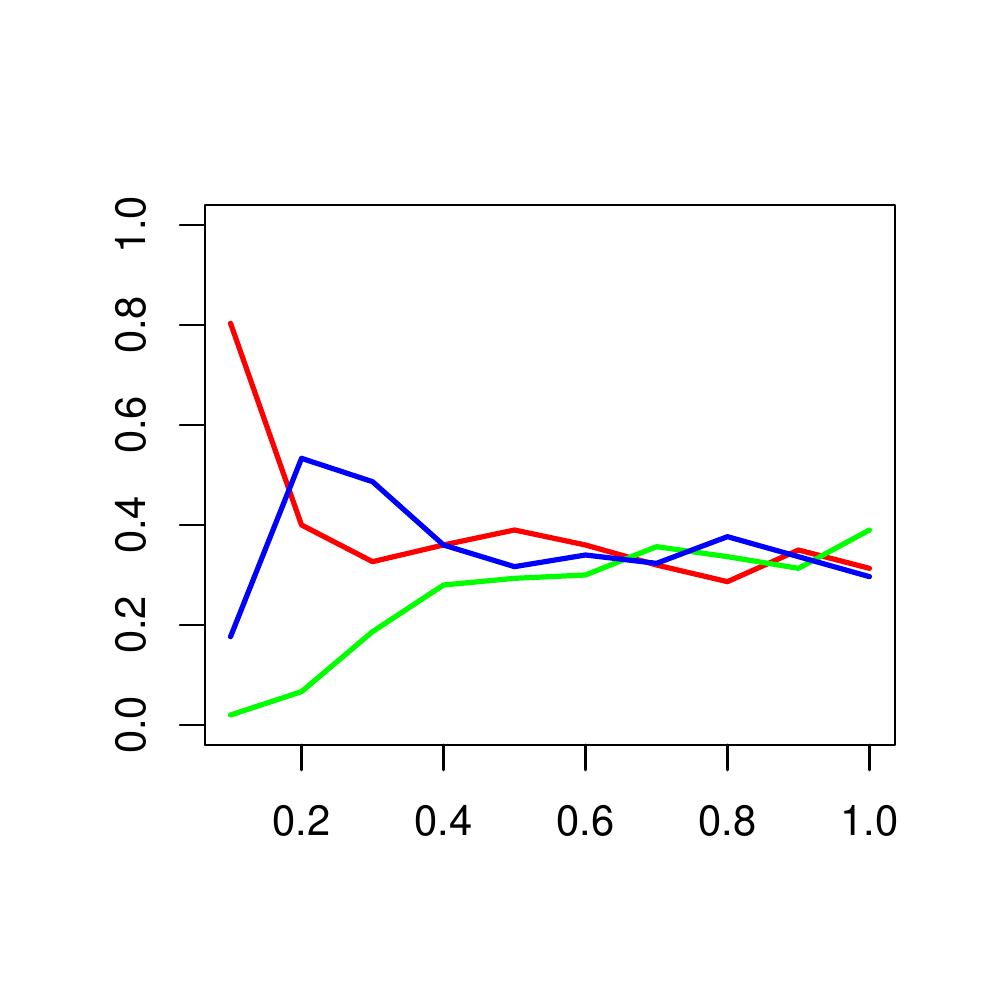}
		\vspace{-10mm} 
  	\caption{\scriptsize{Heaviside}} \label{fig:mu30Hv}
	\end{subfigure}%
	\hspace{-6mm} 
	\begin{subfigure}[t]{.2\textwidth}
		\includegraphics[width=\textwidth]{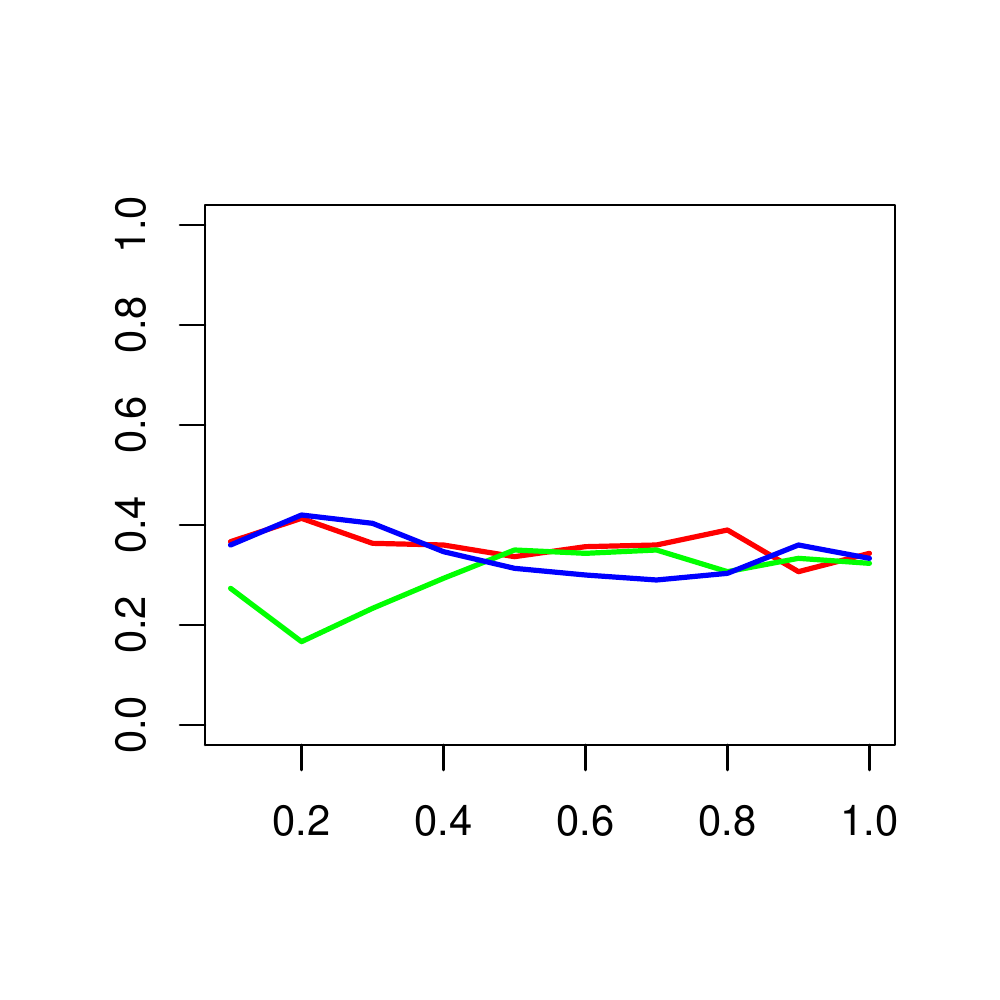}
		\vspace{-10mm} 
  	\caption{\scriptsize{Sawtooth}} \label{fig:mu30Sa}
	\end{subfigure}%
	\hspace{-6mm} 
	\begin{subfigure}[t]{.2\textwidth}
		\includegraphics[width=\textwidth]{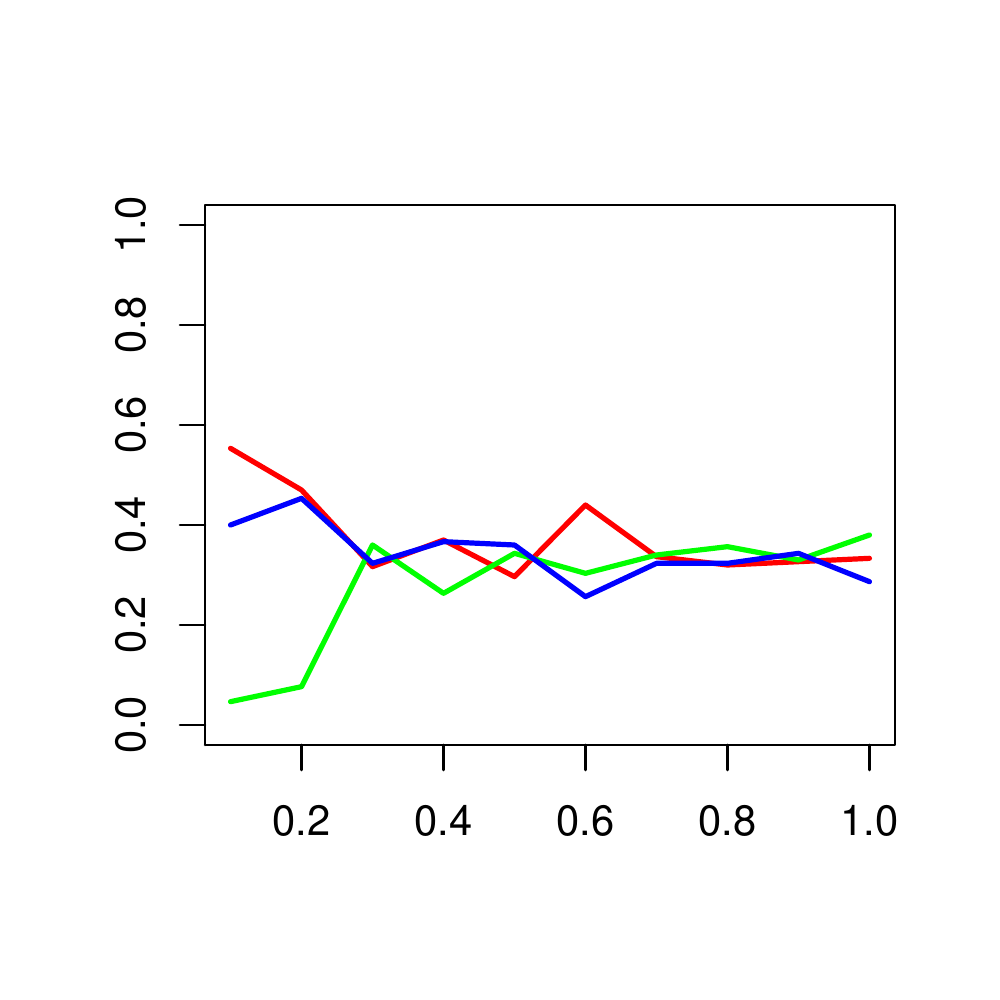}
		\vspace{-10mm} 
  	\caption{\scriptsize{Ackley}} \label{fig:mu30Ac}
	\end{subfigure}%
	\hspace{-6mm} 
	\begin{subfigure}[t]{.2\textwidth}
		\includegraphics[width=\textwidth]{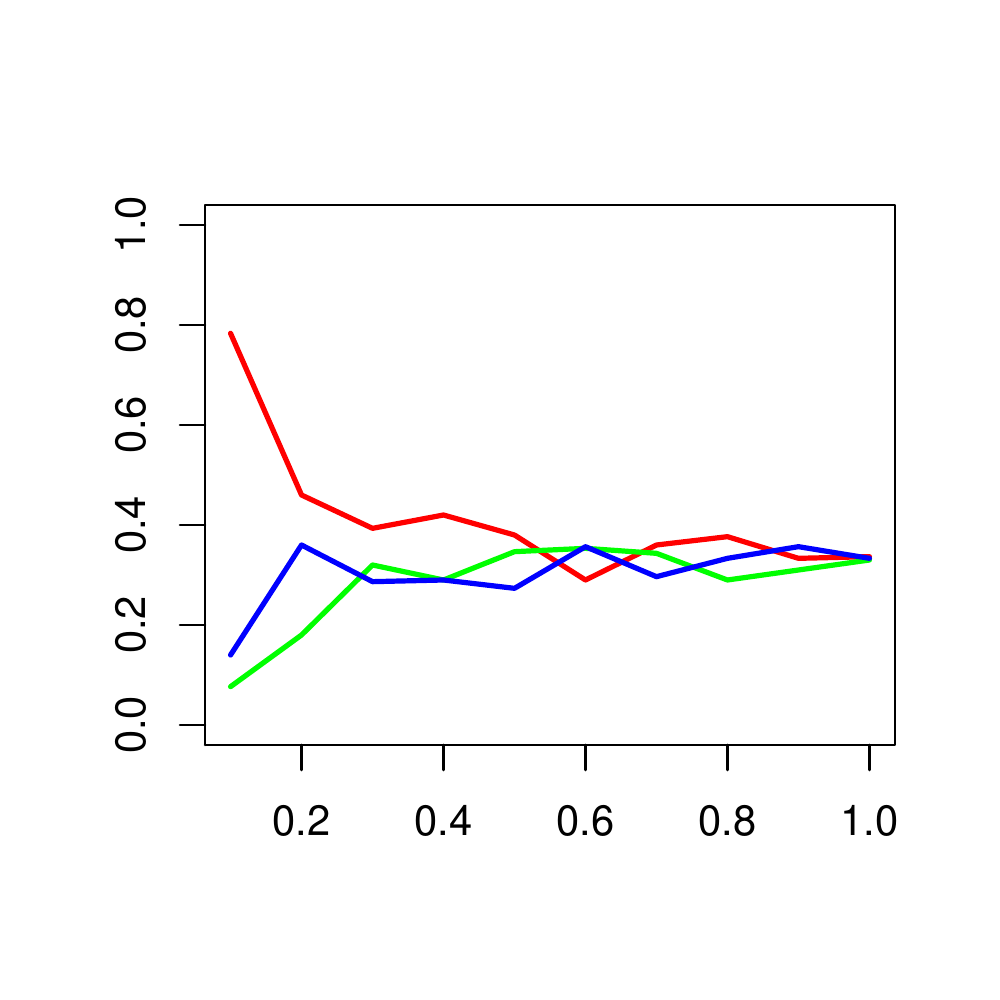}
		\vspace{-10mm} 
  	\caption{\scriptsize{Sphere}} \label{fig:mu30Sp}
	\end{subfigure}%
	\hspace{-6mm} 
	\begin{subfigure}[t]{.2\textwidth}
		\includegraphics[width=\textwidth]{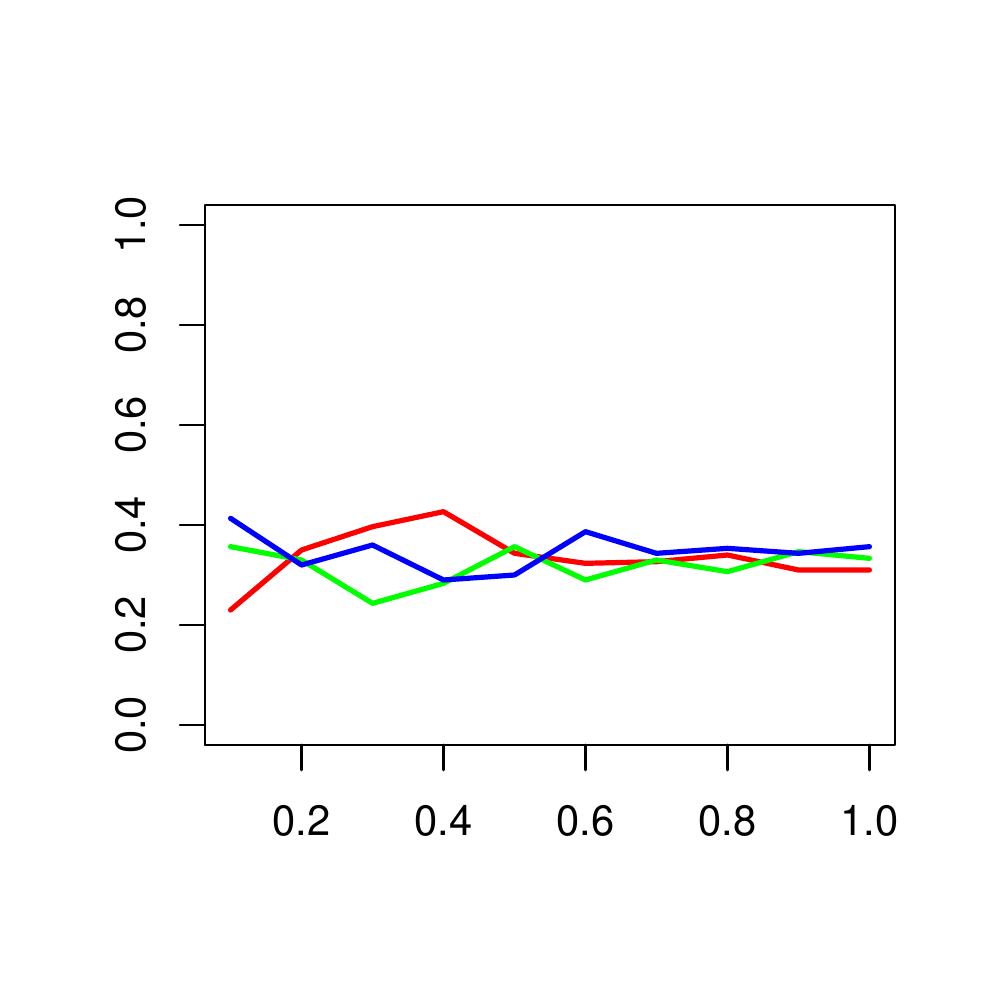}
		\vspace{-10mm} 
  	\caption{\scriptsize{Rosenbrock}} \label{fig:mu30Ro}
	\end{subfigure}
		
	\caption{Component -- decile breakdowns for the form of PSO mutation, across all GGGP heuristics at 30D. Components: None (red), Uniform (green), Gaussian (blue).}
	\label{fig:mutation30Comp}

\end{figure}

\begin{figure}[H]
	\centering
	
	\vspace{-5mm} 

	\begin{subfigure}[t]{.2\textwidth}
		\includegraphics[width=\textwidth]{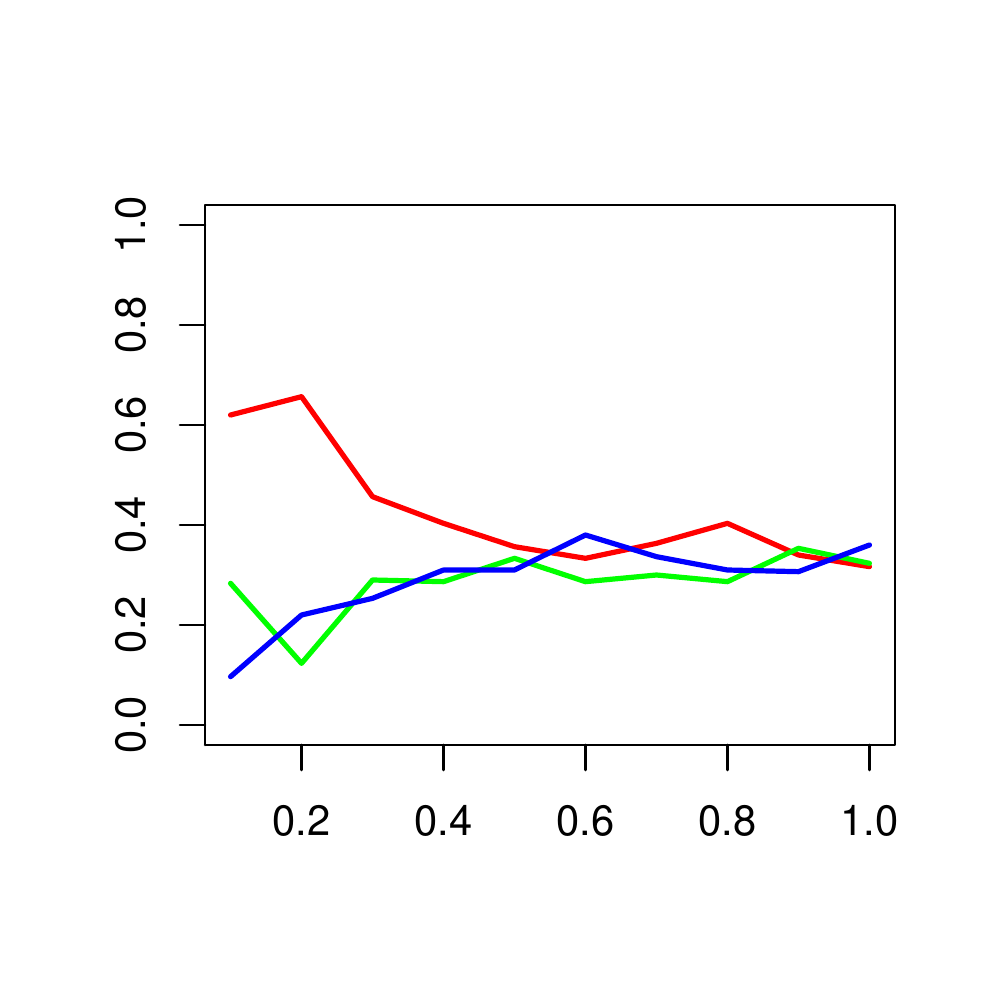}
		\vspace{-10mm} 
  	\caption{\scriptsize{Rastrigin}} \label{fig:mu100Ra}
	\end{subfigure}%
	\hspace{-6mm} 
	\begin{subfigure}[t]{.2\textwidth}
		\includegraphics[width=\textwidth]{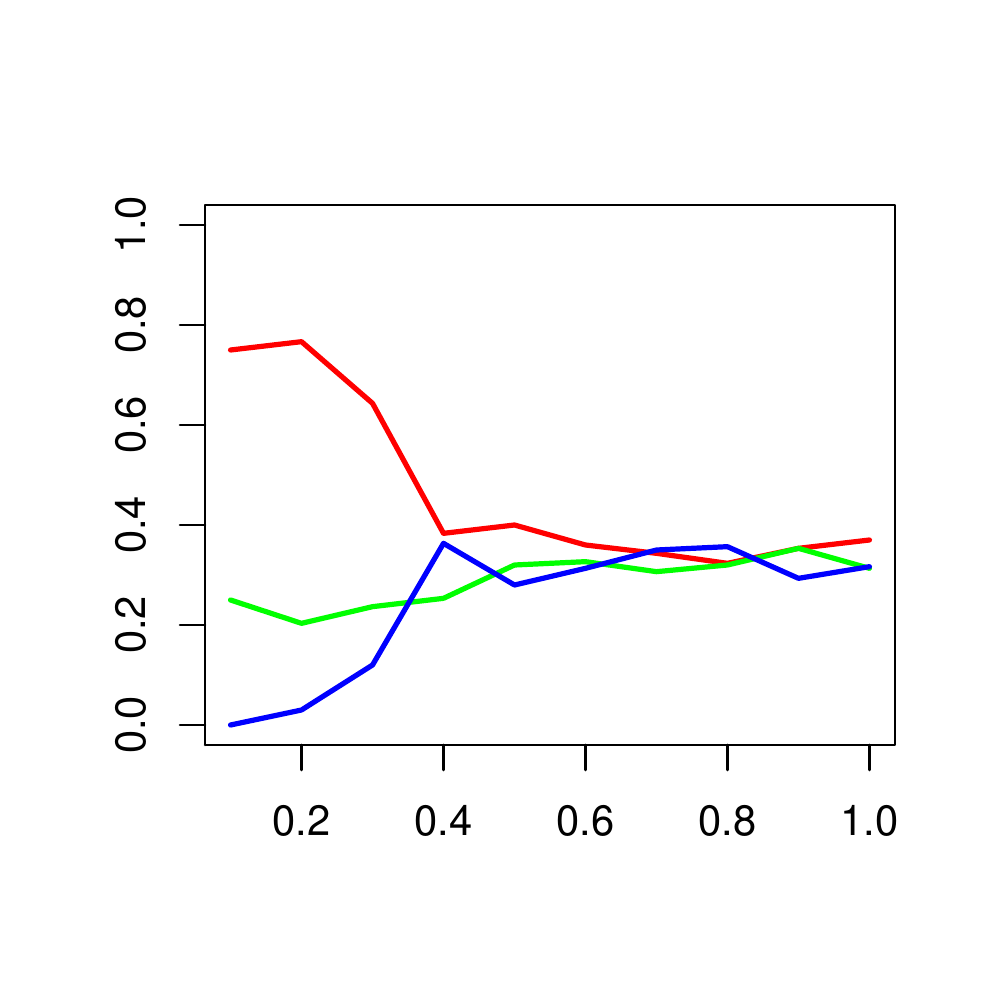}
		\vspace{-10mm} 
  	\caption{\scriptsize{Multipeak F1}} \label{fig:mu100M1}
	\end{subfigure}%
	\hspace{-6mm} 
	\begin{subfigure}[t]{.2\textwidth}
		\includegraphics[width=\textwidth]{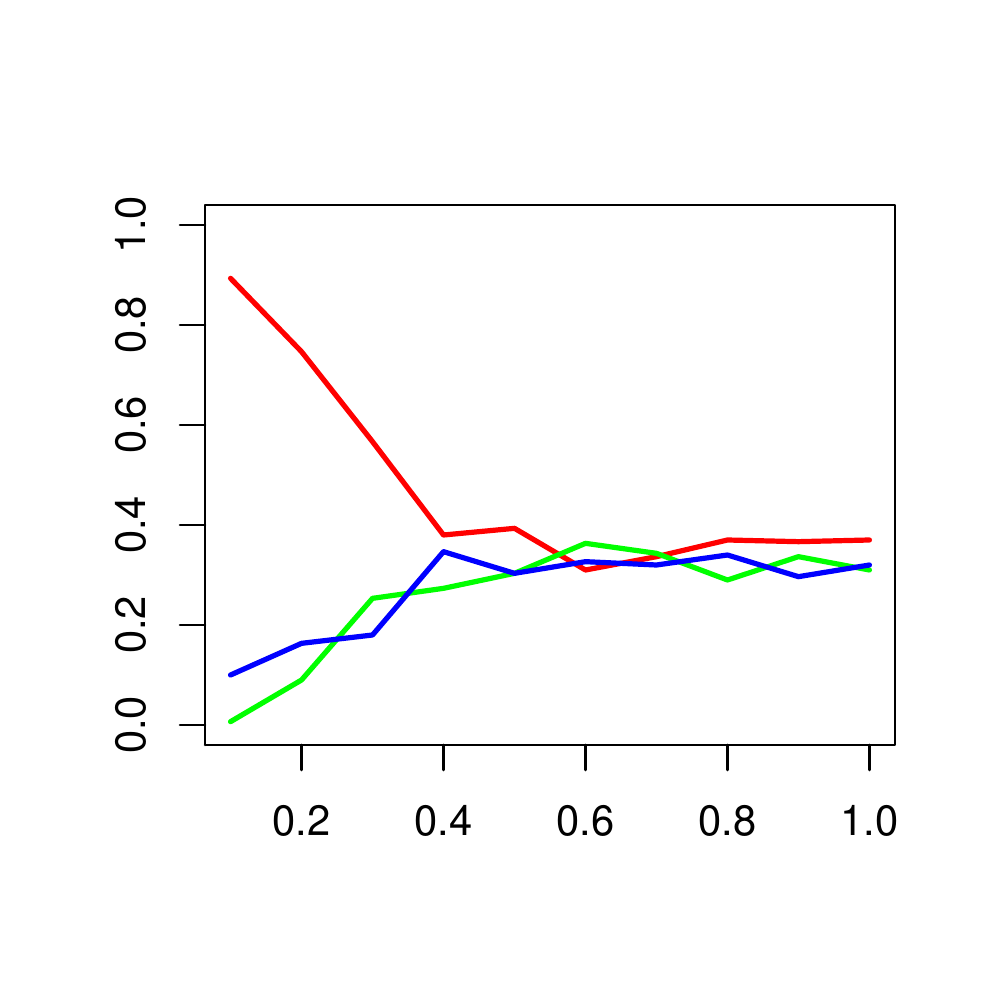}
		\vspace{-10mm} 
  	\caption{\scriptsize{Multipeak F2}} \label{fig:mu100M2}
	\end{subfigure}%
	\hspace{-6mm} 
	\begin{subfigure}[t]{.2\textwidth}
		\includegraphics[width=\textwidth]{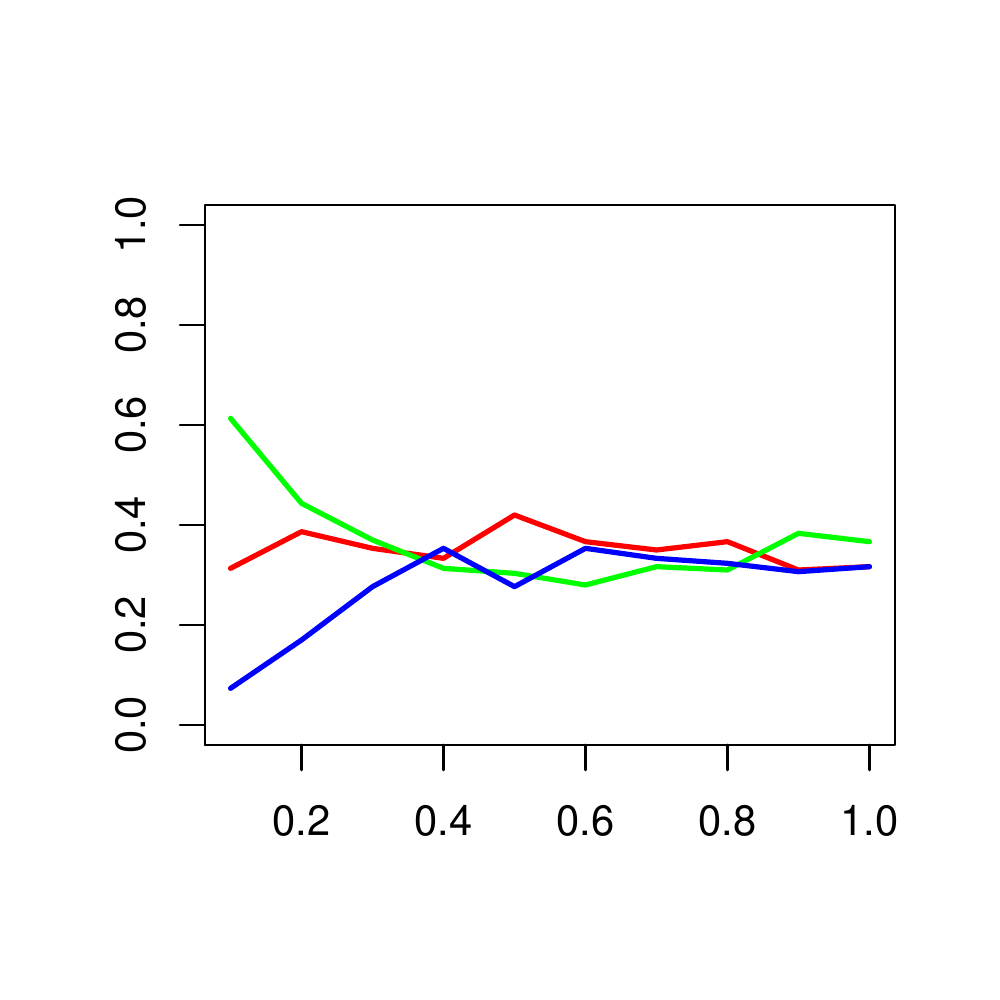}
		\vspace{-10mm} 
  	\caption{\scriptsize{Brankes}} \label{fig:mu100Br}
	\end{subfigure}%
	\hspace{-6mm} 
	\begin{subfigure}[t]{.2\textwidth}
		\includegraphics[width=\textwidth]{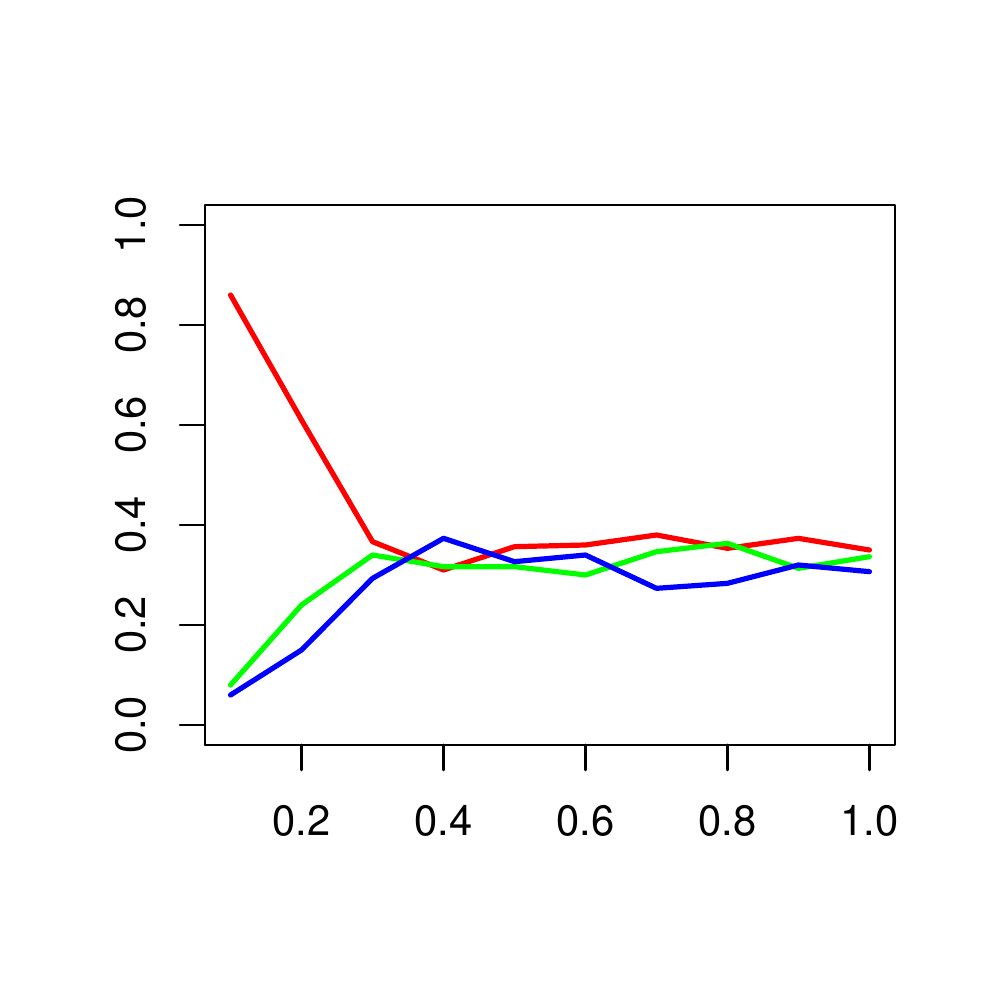}
		\vspace{-10mm} 
  	\caption{\scriptsize{Pickelhaube}} \label{fig:mu100Pi}
	\end{subfigure}
	
	\vspace{-2mm} 
		
	\begin{subfigure}[t]{.2\textwidth}
		\includegraphics[width=\textwidth]{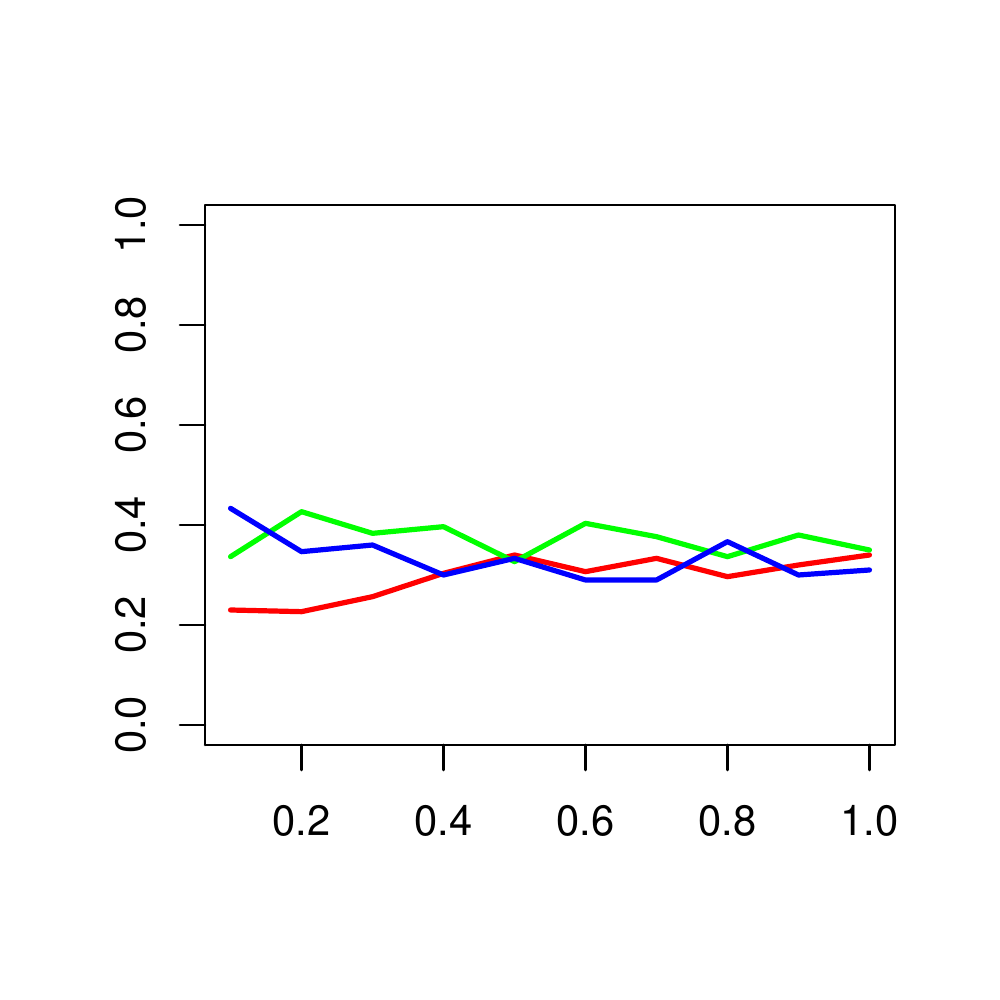}
		\vspace{-10mm} 
  	\caption{\scriptsize{Heaviside}} \label{fig:mu100Hv}
	\end{subfigure}%
	\hspace{-6mm} 
	\begin{subfigure}[t]{.2\textwidth}
		\includegraphics[width=\textwidth]{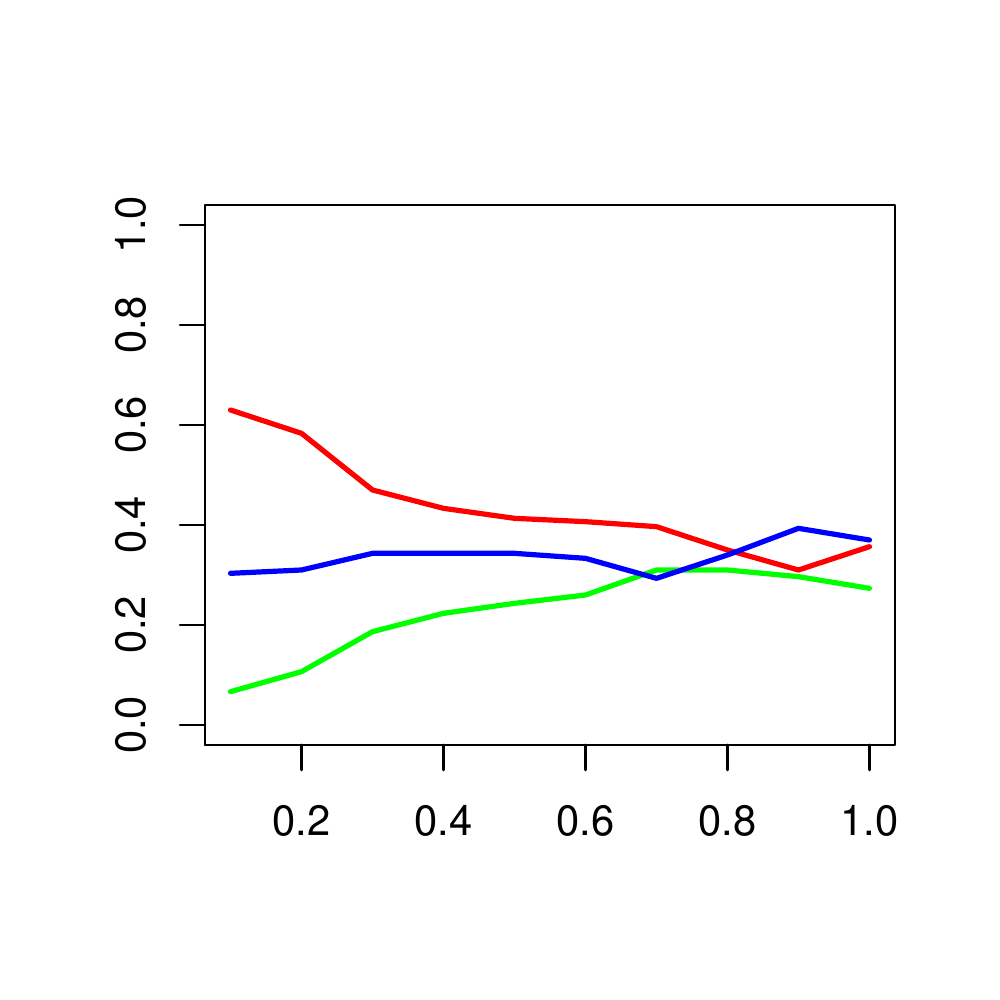}
		\vspace{-10mm} 
  	\caption{\scriptsize{Sawtooth}} \label{fig:mu100Sa}
	\end{subfigure}%
	\hspace{-6mm} 
	\begin{subfigure}[t]{.2\textwidth}
		\includegraphics[width=\textwidth]{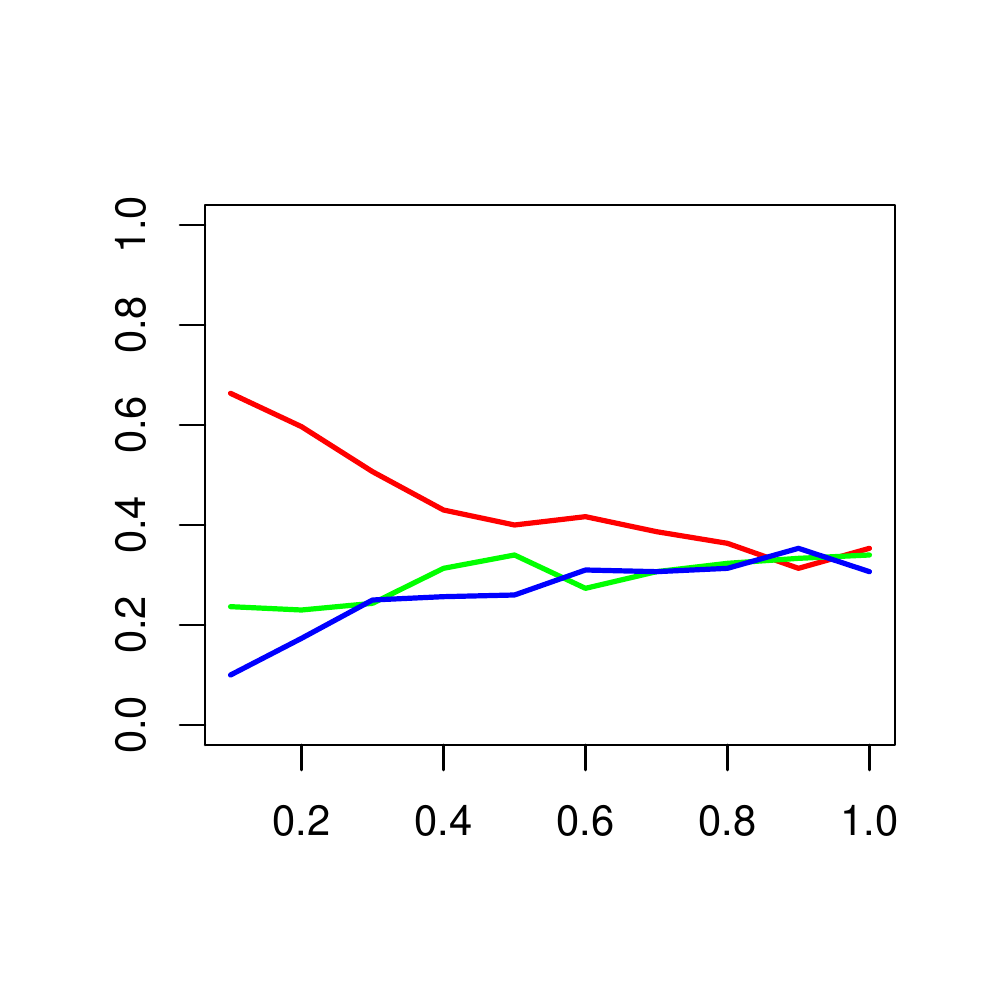}
		\vspace{-10mm} 
  	\caption{\scriptsize{Ackley}} \label{fig:mu100Ac}
	\end{subfigure}%
	\hspace{-6mm} 
	\begin{subfigure}[t]{.2\textwidth}
		\includegraphics[width=\textwidth]{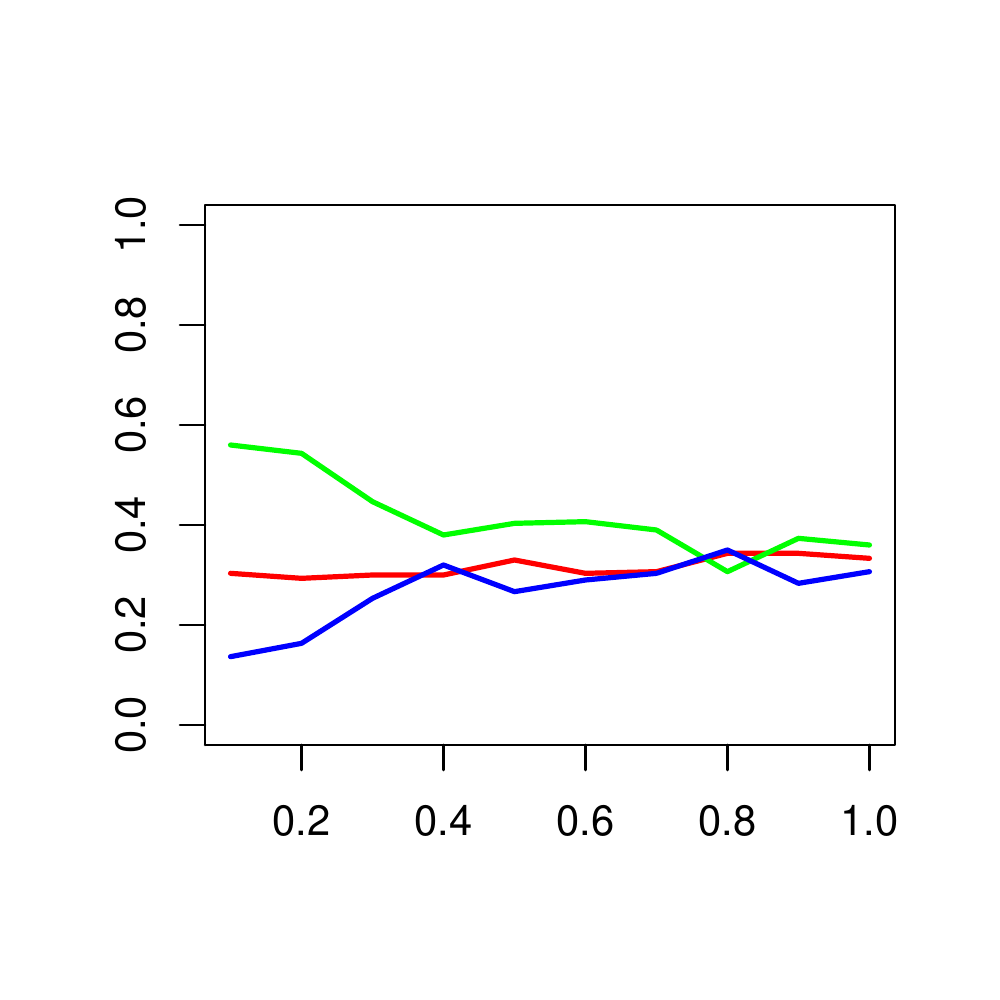}
		\vspace{-10mm} 
  	\caption{\scriptsize{Sphere}} \label{fig:mu100Sp}
	\end{subfigure}%
	\hspace{-6mm} 
	\begin{subfigure}[t]{.2\textwidth}
		\includegraphics[width=\textwidth]{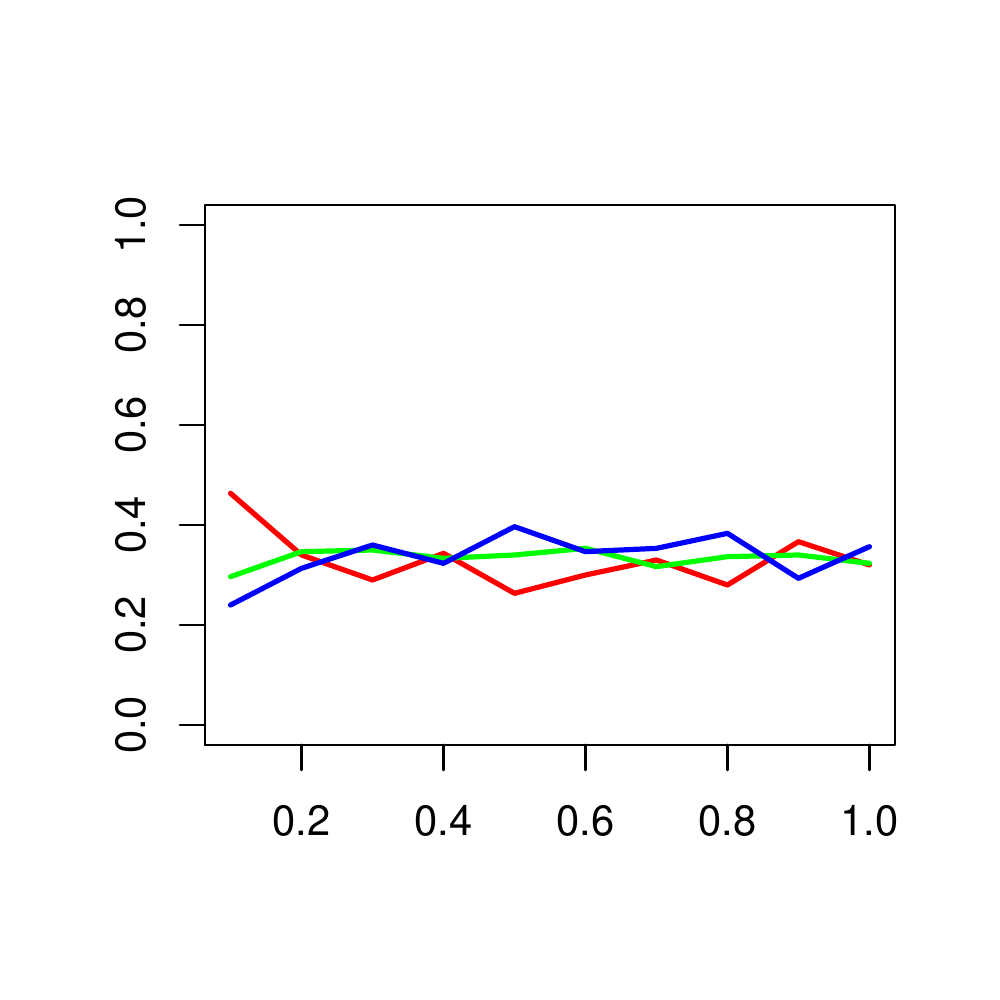}
		\vspace{-10mm} 
  	\caption{\scriptsize{Rosenbrock}} \label{fig:mu100Ro}
	\end{subfigure}
		
	\caption{Component -- decile breakdowns for the form of PSO mutation, across all GGGP heuristics at 100D. Components: None (red), Uniform (green), Gaussian (blue).}
	\label{fig:mutation100Comp}
	
\end{figure}

\begin{figure}[H]
	\centering
	
	\vspace{-5mm} 

	\begin{subfigure}[t]{.24\textwidth}
		\includegraphics[width=\textwidth]{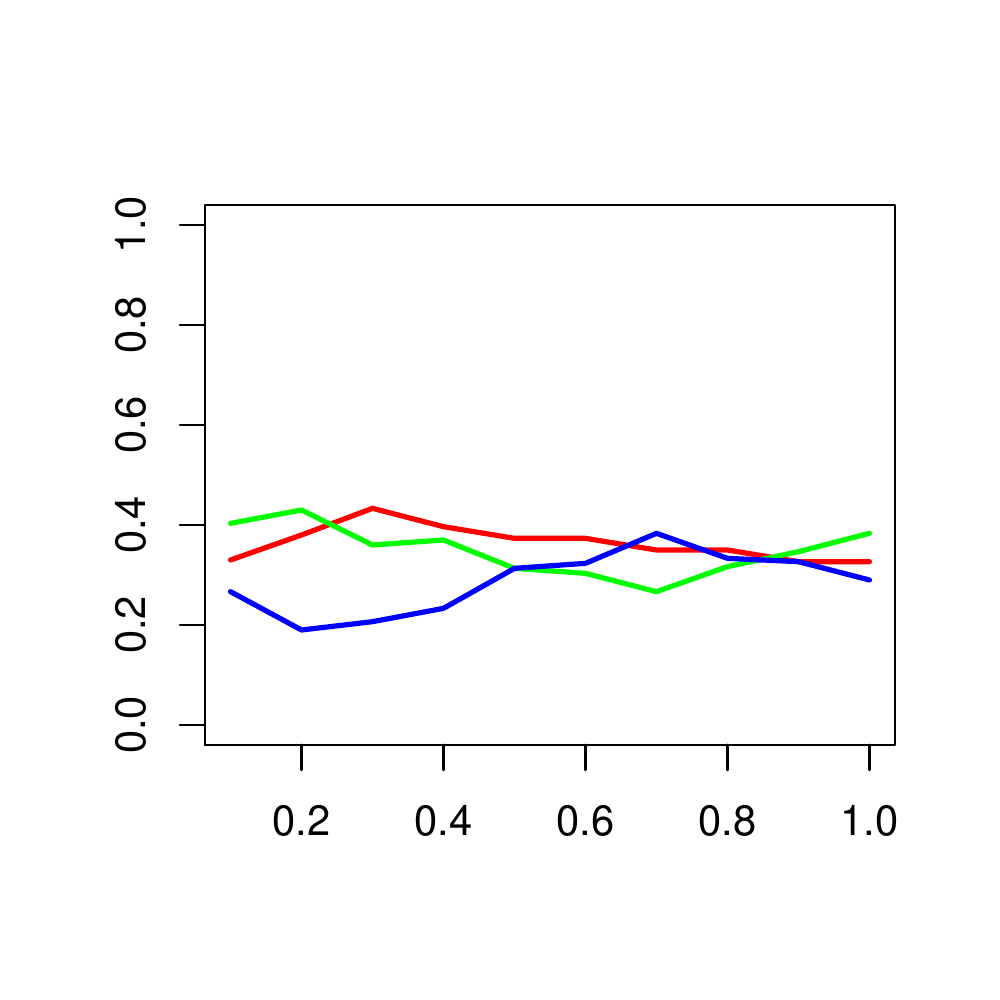}
		\vspace{-10mm} 
  	\caption{\scriptsize{30D}} \label{fig:mu30}
	\end{subfigure}%
	\hspace{-6mm} 
	\begin{subfigure}[t]{.24\textwidth}
		\includegraphics[width=\textwidth]{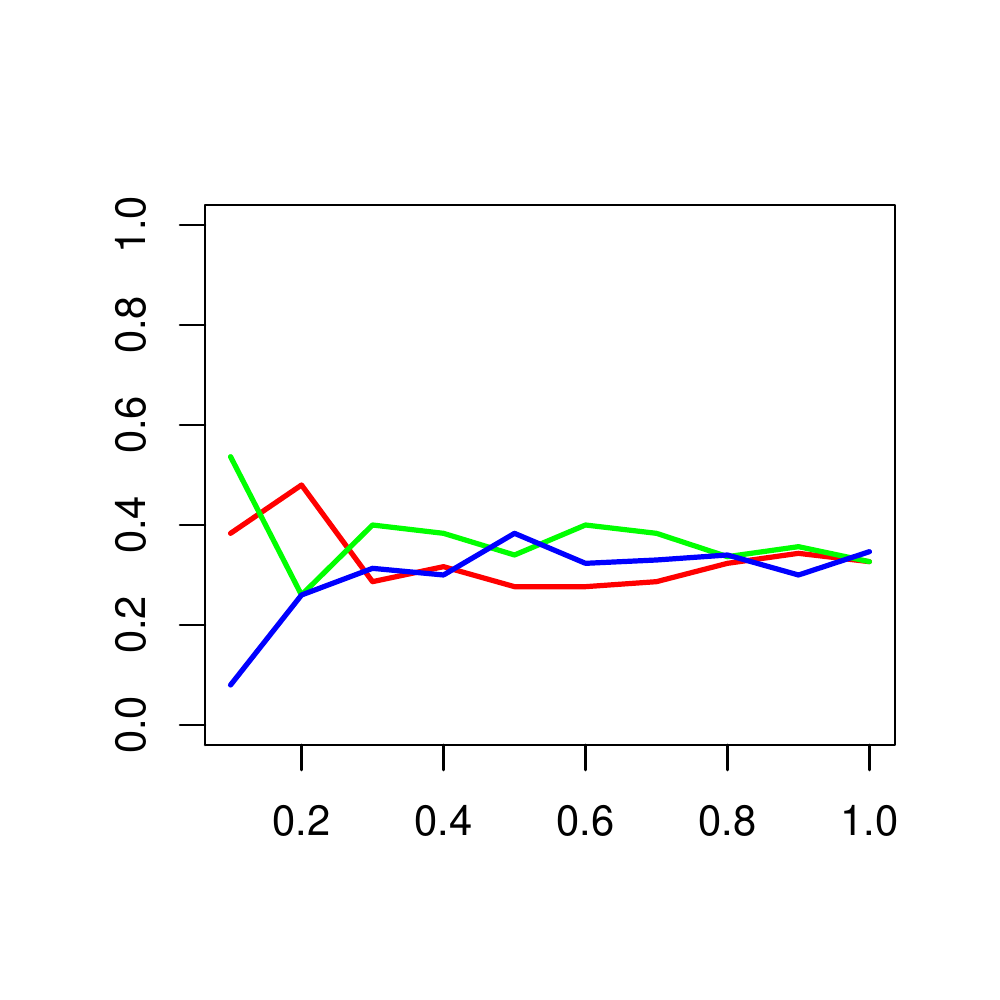}
		\vspace{-10mm} 
  	\caption{\scriptsize{100D}} \label{fig:mu100}
	\end{subfigure}

	\caption{Component -- decile breakdowns for the form of PSO mutation, across all GGGP heuristics at 30D and 100D for the general heuristics. Components: None (red), Uniform (green), Gaussian (blue).}
	\label{fig:mutationComp}
	
\end{figure}

Our grammar includes the ability to mutate a particle's intended position. Analysis of this component is shown in Table~\ref{fig:proportions} and Figures~\ref{fig:mutation30Comp} to~\ref{fig:mutationComp}. A choice of no mutation (red), or mutation due to either Uniform (green) or Gaussian (blue) random sampling, is available. In Table~\ref{fig:proportions} all three alternatives are well represented, with no mutation performing best, appearing in over 40\% of the top third performing heuristics. In the decile analysis, for 30D individual cases preference is quite even, whilst at 100D the non-use of mutation is more preferred at the lowest deciles. In the general cases apportionment is evenly distributed.

\subsubsection{Extended movement capability features: relocation by LEH}
\label{sec:compLEH}

If a heuristic employs the dormancy-relocation capability (+LEH or +DD+LEH), there is a choice of how dormant particles are relocated. This is by the calculation of the largest empty hypersphere devoid of poor points, or completely randomly. From Table~\ref{fig:proportions} the use of relocation using LEH is seen to dominate when dormancy is employed. At 30D 88\% of all heuristics use LEH-relocation, rising to 100\% of the top performers. At 100D 89\% rises to 100\%. This complete dominance is also observed in the decile level analysis.

\subsubsection{Extended movement capability features: descent directions $\pmb{r}_3$ vector}
\label{sec:compR3}

When a heuristic use the descent directions capability (+DD or +DD+LEH), the $\pmb{r}_3$ vector may be generated randomly or set to the unit vector. Table~\ref{fig:proportions} indicates a reasonably even use of the randomised or unit $\pmb{r}_3$ vectors. Where d.d.\ is employed, the use of a unit vector rises from 56\% across all heuristics to 61\% in the top performers, for 30D. Whereas at 100D this figure remains static at 48\%. The decile level analysis reflects these high level patterns. At 30D for several cases the use of a unit vector performs better in the lowest deciles, whilst at 100D the preference is fairly evenly distributed across cases.

\subsection{Summary of experimental analysis}
\label{sec:ExperimentSummary}

The analysis of the results due to the best heuristics generated in the GP runs shows a strong performance. For individual case performance the indications against comparator results is encouraging. For the general cases the newly developed heuristics show an improvement over the best comparators, in some cases significant, despite using a budget 60\% lower.

For the component level analysis, inner maximisation using random sampling on a small number of points performed best, with a particle level stopping condition strongly preferred. For the outer minimisation the best heuristic performance is (separately) related to a relatively small swarm size, communication using a Global typology, and a particle movement formulation consisting of an inertia based velocity equation plus d.d.\ and LEH extended capabilities.

In addition to the decile level component analysis reported here, consideration was given to potential correlations between alternatives across different components. No such correlation was observed.

\section{Conclusions and further work}
\label{sec:concusionsFurtherWork}

We have used grammar-guided genetic programming to automatically generate particle swarm based metaheuristics for robust problems, in order to determine improved search algorithms and assess the effectiveness of various algorithmic sub-components. This has involved the generation of a grammar consisting of a number of heuristic building blocks, the design rules for constructing heuristics from these components, and an evolutionary GP process. We have searched a heuristic sub-algorithm space not previously investigated, encompassing specialised robust-focussed capabilities alongside more standard elements such as network topologies and alternatives for the inner maximisation. 

Using a suite of 10 test problems at 30D and 100D, the best evolved heuristics were identified at individual and general (all problems simultaneously) test case levels. Using comparators, significant improvements are observed for the best new general heuristics, whilst indicative individual case results are highly promising.

The GP process generates substantial numbers of heuristics, enabling an assessment of algorithmic sub-components against heuristic performance. In the context of a budget of 2,000 function evaluations, this identifies an inner maximisation by random sampling on a small number of points as most effective, including the use of a particle level stopping condition. For the outer minimisation small numbers of particles are preferred, sharing information through a Global topology. Other topologies exhibit some good performance. The preferred particle movement uses a baseline inertia velocity equation plus some largest empty hypersphere \cite{HughesGoerigkWright2019, HughesGoerigkDokka2020a} and descent directions \cite{BertsimasNohadaniTeo2010, HughesGoerigkDokka2020a} heuristic capabilities. This includes the assessment of particle dormancy and relocation to the centre of the LEH, or the use of a d.d.\ vector component in the velocity formulation, or both.

There are a number of ways in which this work can be built upon, most obviously in terms of extending the sub-algorithmic space over which the GP operates. Moving away from a PSO structure for all of the heuristics to a more general agent based setting, using other population-based metaheuristics, would introduce alternative movement and information sharing capabilities into our grammar for the outer minimisation layer.

As the use of random sampling for the inner maximisation layer has proven effective here, the inclusion in our grammar of some efficient sampling techniques such as the specialised Latin hypercube approach described in \cite{FeiBrankeGulpinar2019}, would seem appropriate.

The potential efficiencies offered by emulation in either the outer minimisation or inner maximisation layers, warrants investigation. The introduction of emulation based components into the grammar, including sub-elements of specific emulation approaches, could significantly extend the heuristic solution space.

A final consideration might be the use of alternatives to the GGGP approach, to automatically generate heuristics for robust problems.


\newpage

\appendix

\section*{Appendix}

\section{Test functions}
\label{sec:testFunctionFormulae}

\begin{figure}[htbp]
	\centering
	
	\vspace{-5mm} 
		
	\begin{subfigure}[t]{.24\textwidth}
		\includegraphics[width=\textwidth]{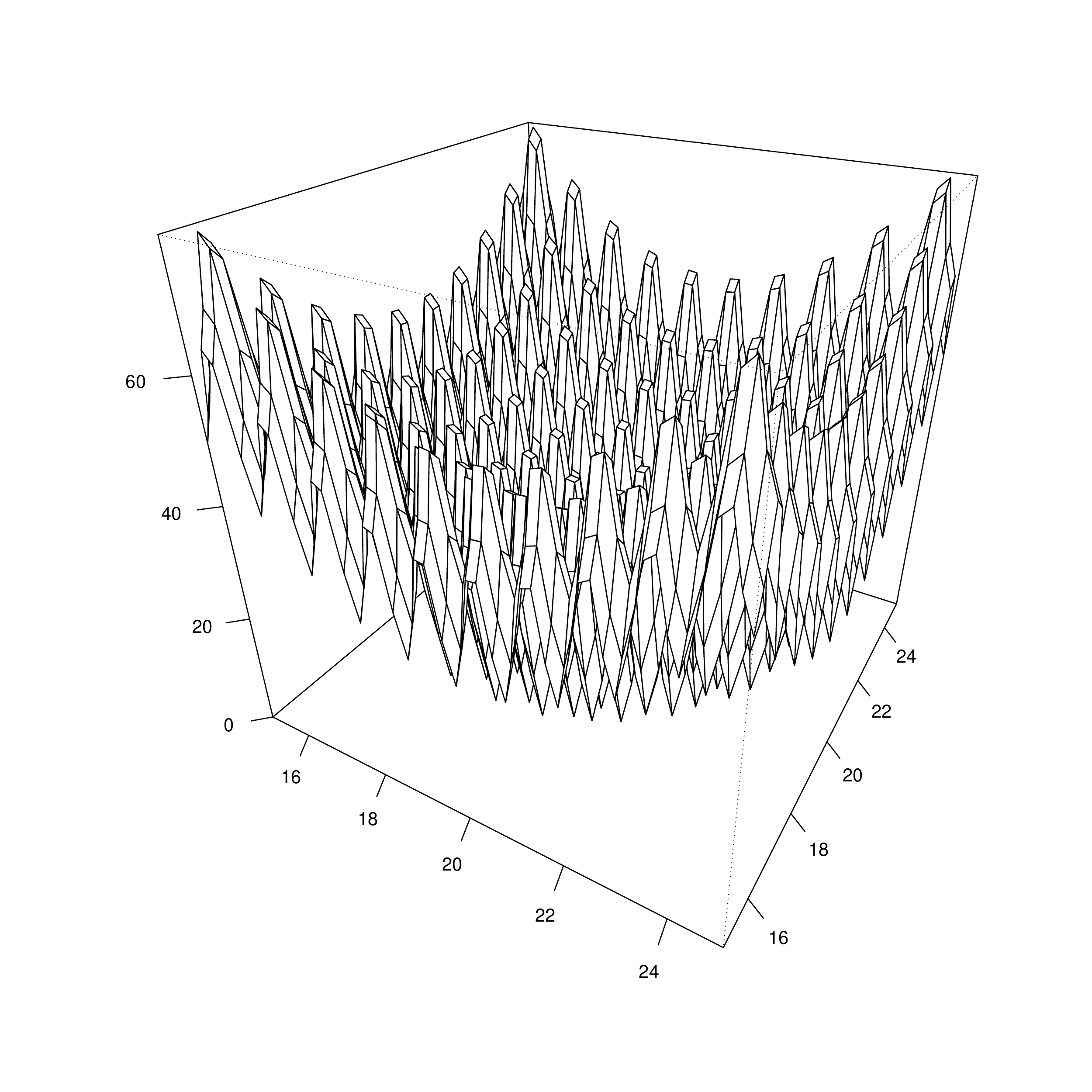}		
  	\caption{\scriptsize{Rastrigin Nom}} \label{fig:RastriginNom}
	\end{subfigure}%
	\begin{subfigure}[t]{.24\textwidth}
		\includegraphics[width=\textwidth]{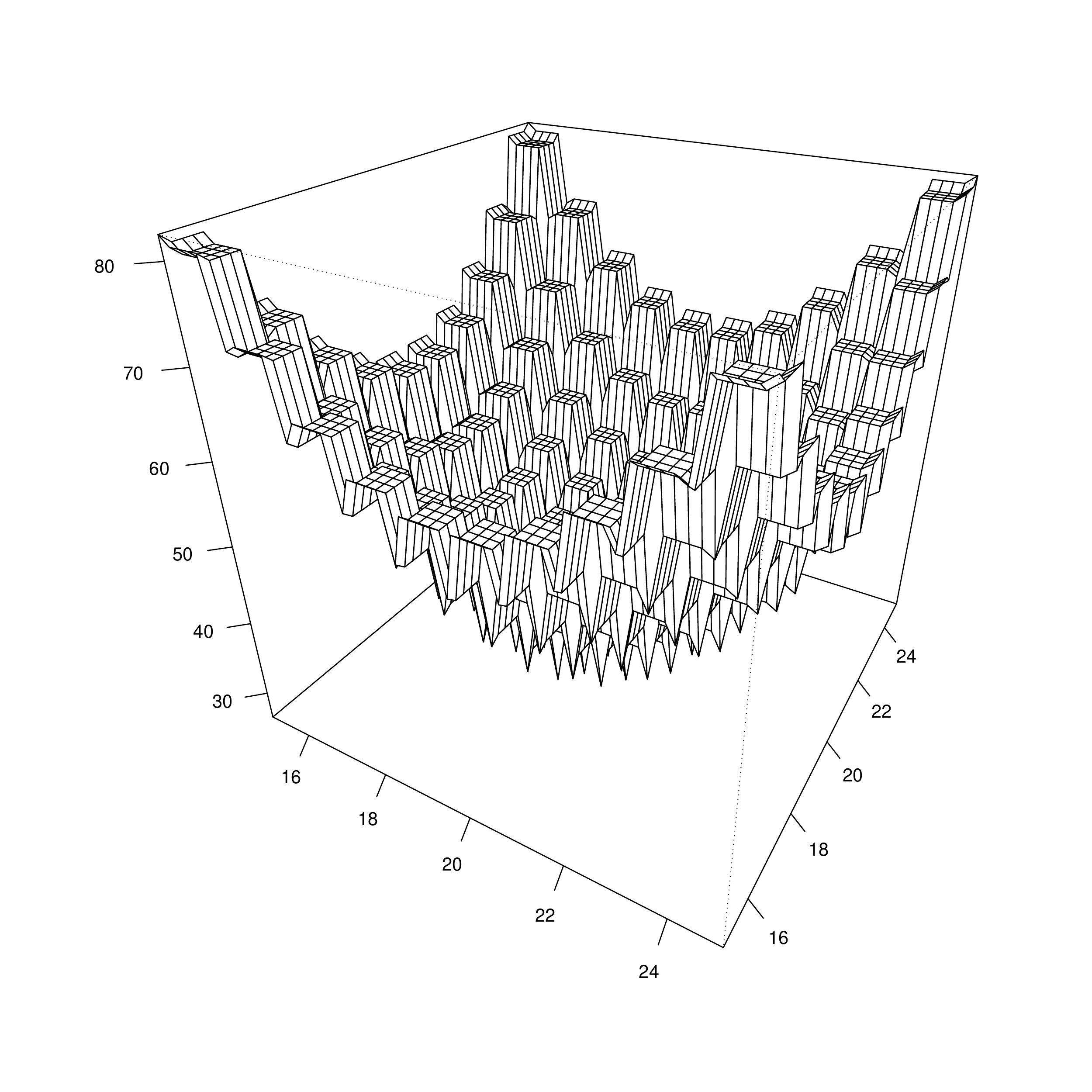}	
  	\caption{\scriptsize{Rastrigin Worst}} \label{fig:RastriginWorst}
	\end{subfigure}		
	\begin{subfigure}[t]{.24\textwidth}
		\includegraphics[width=\textwidth]{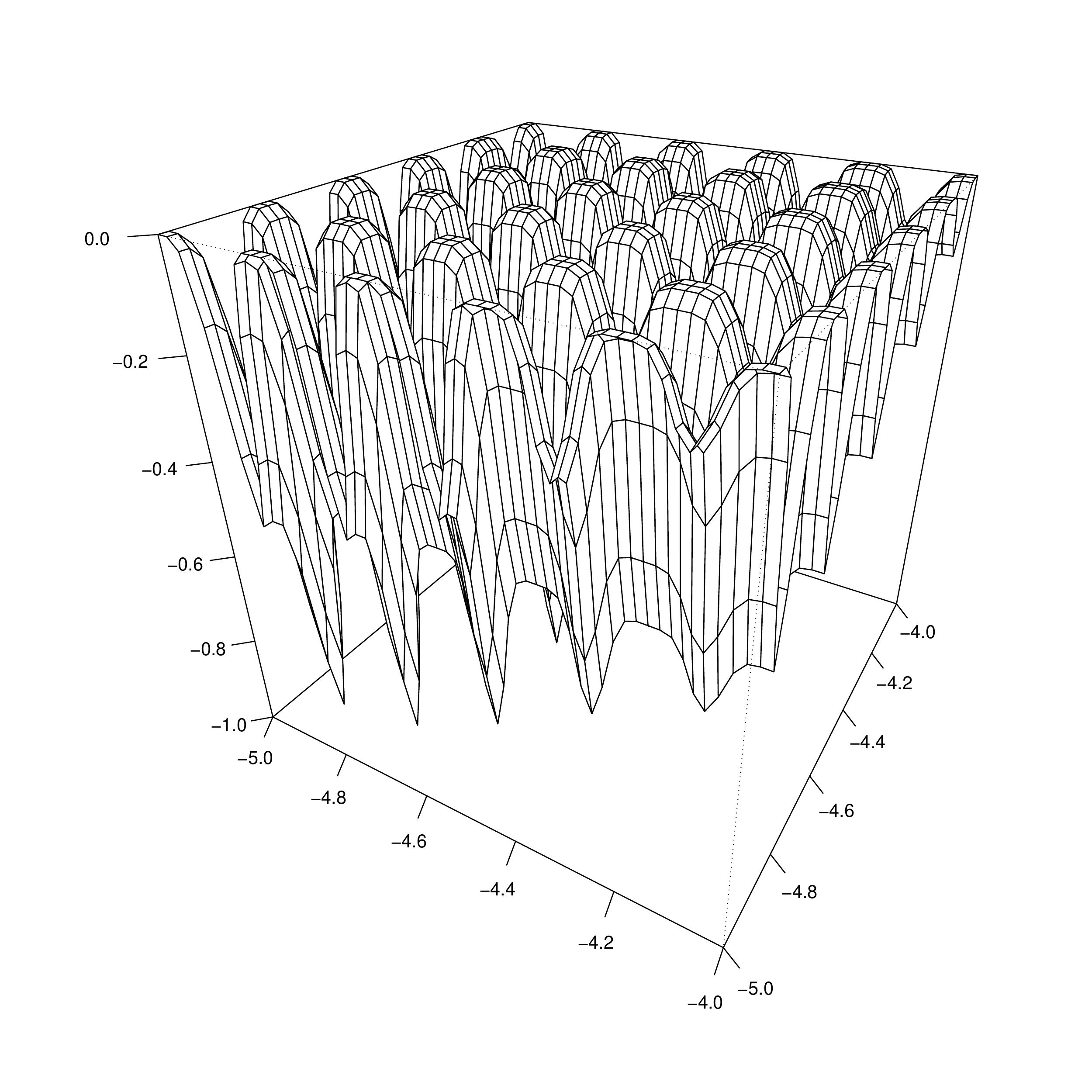}		
  	\caption{\scriptsize{Multipeak F1 Nom}} \label{fig:MultipeakF1Nom}
	\end{subfigure}%
	\begin{subfigure}[t]{.24\textwidth}
		\includegraphics[width=\textwidth]{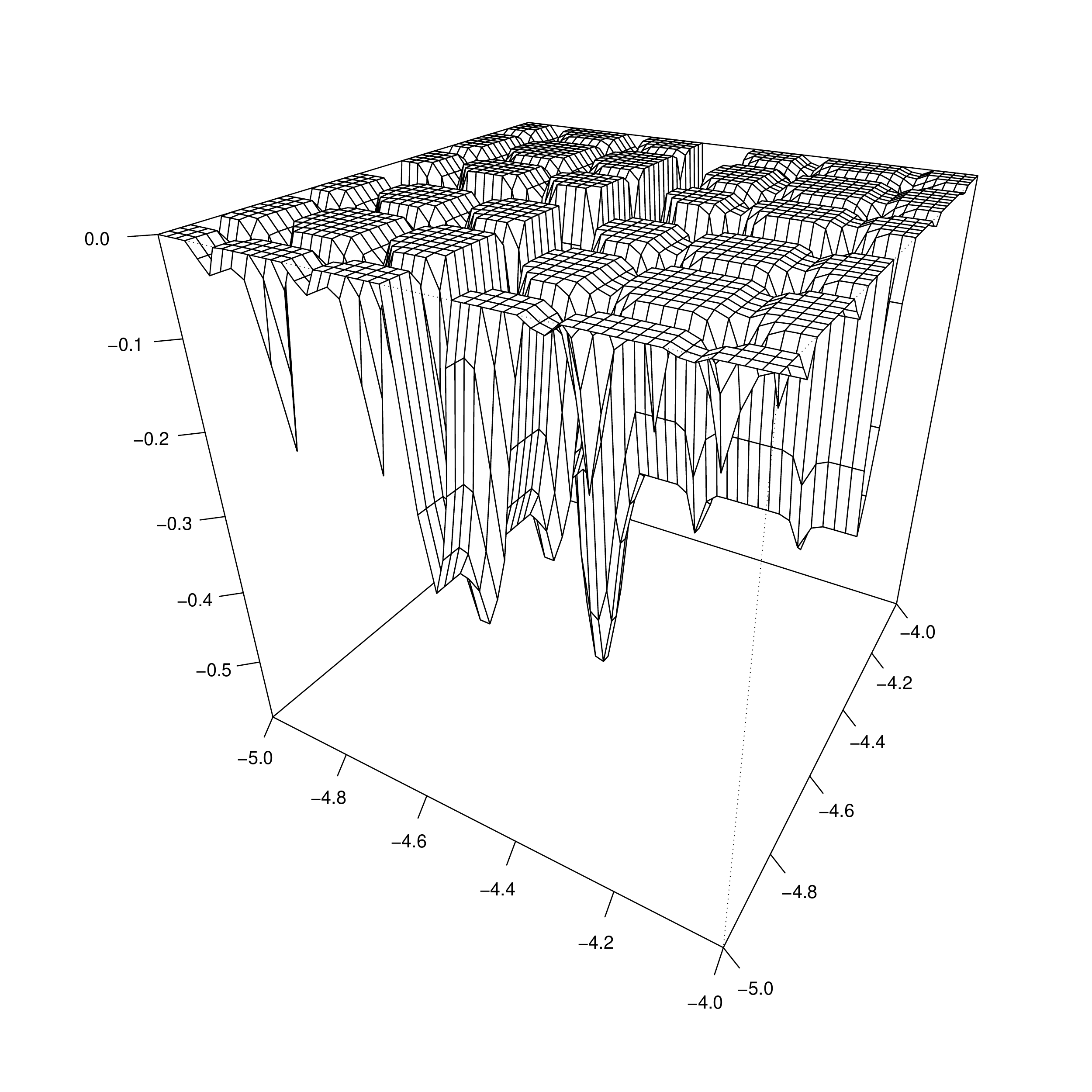}	
  	\caption{\scriptsize{Multipeak F1 Worst}} \label{fig:MultipeakF1Worst}
	\end{subfigure}
	\begin{subfigure}[t]{.24\textwidth}
		\includegraphics[width=\textwidth]{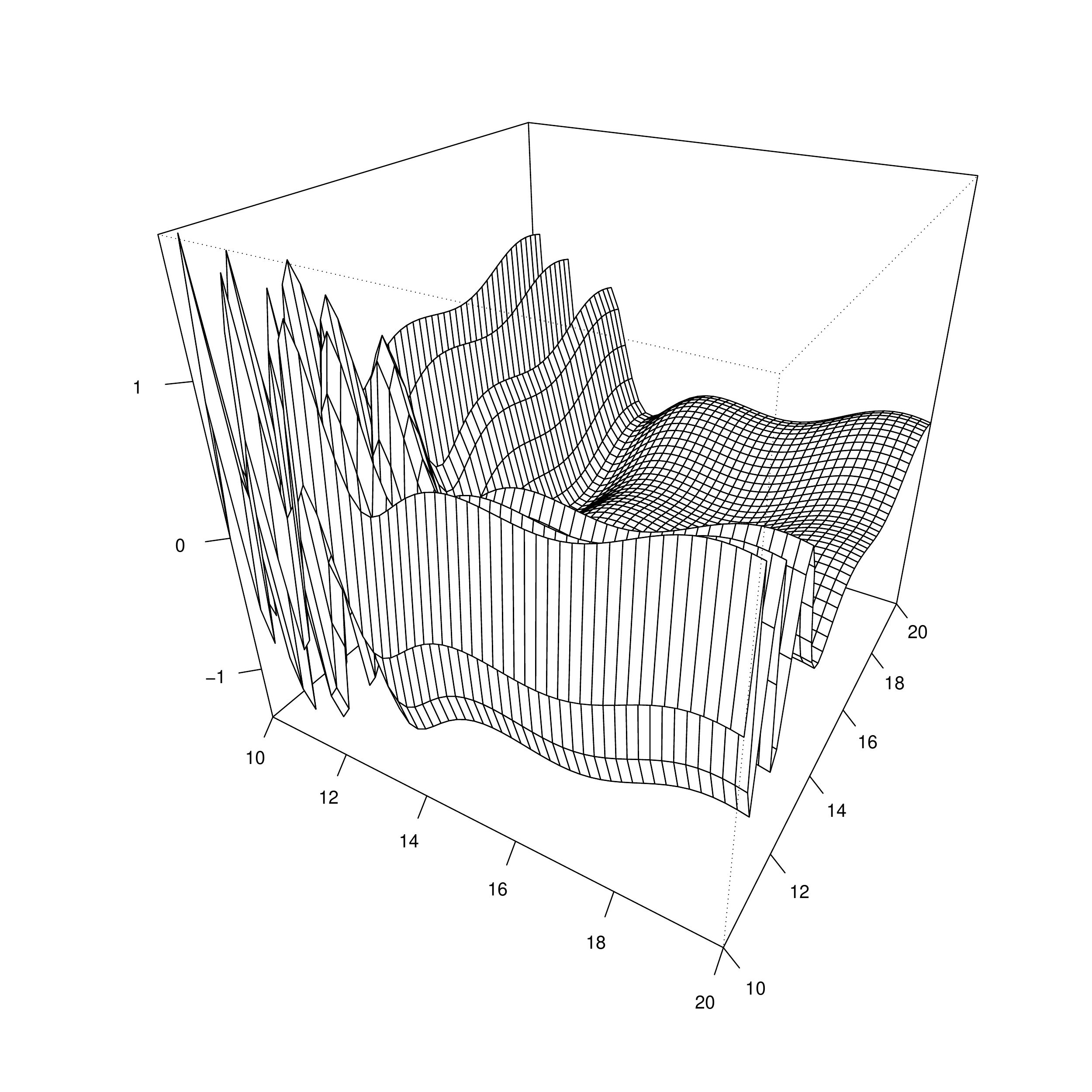}	
  	\caption{\scriptsize{Multipeak F2 Nom}} \label{fig:MultipeakF2Nom}
	\end{subfigure}%
	\begin{subfigure}[t]{.24\textwidth}
		\includegraphics[width=\textwidth]{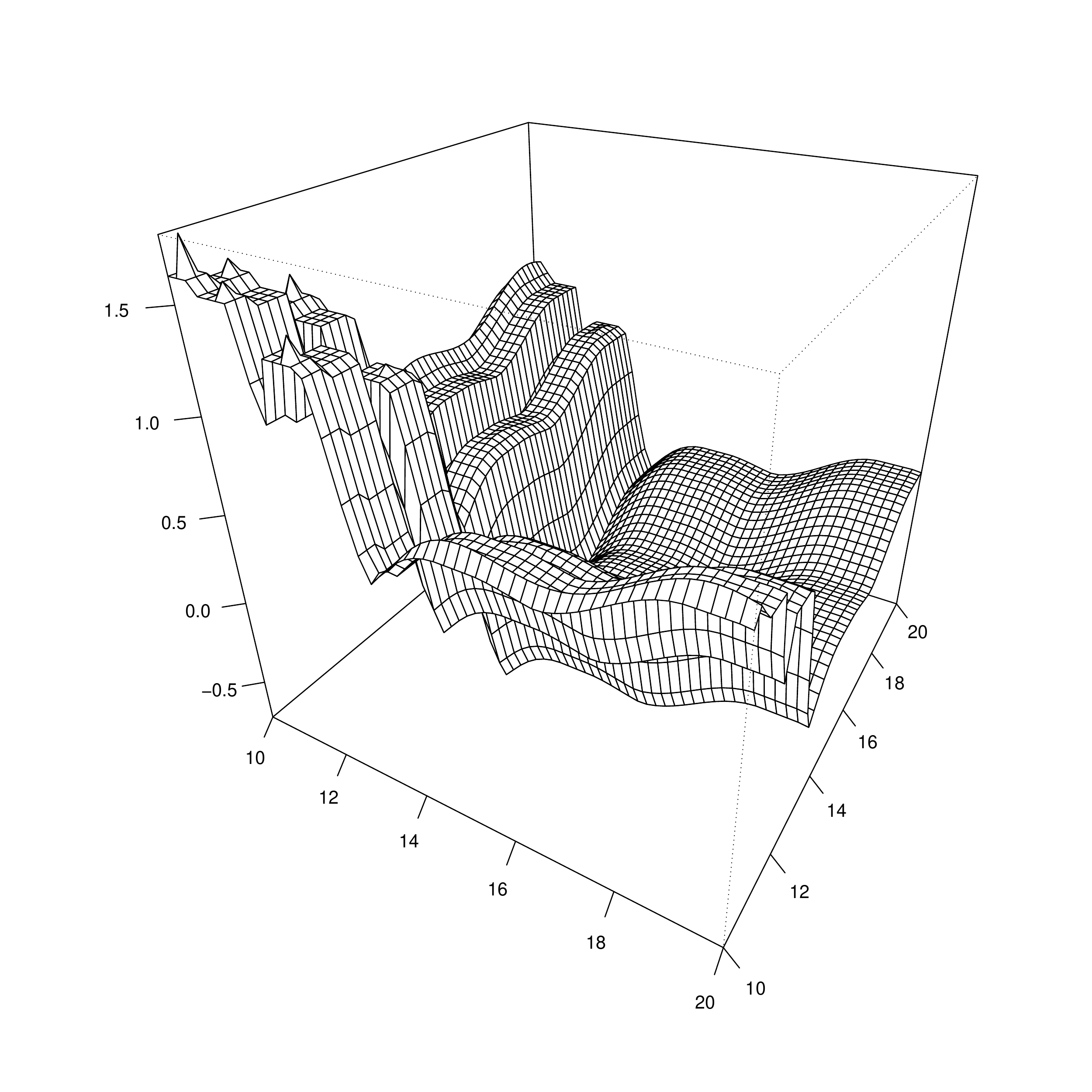}	
  	\caption{\scriptsize{Multipeak F2 Worst}} \label{fig:MultipeakF2Worst}
	\end{subfigure}
  \begin{subfigure}[t]{.24\textwidth}
		\includegraphics[width=\textwidth]{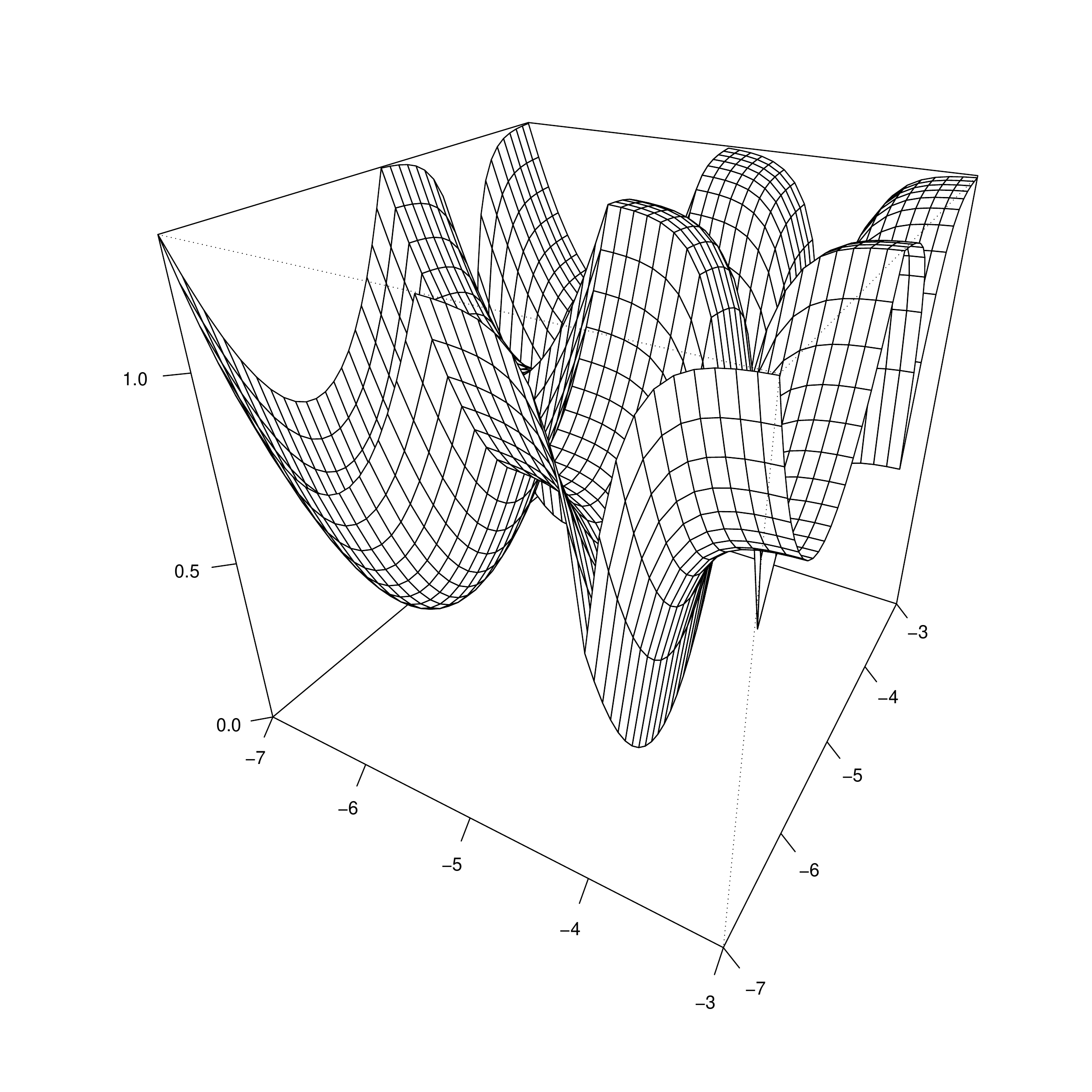}
  	\caption{\scriptsize{Brankes Multi Nom}} \label{fig:BrankesMultipeakNom}
	\end{subfigure}%
	\begin{subfigure}[t]{.24\textwidth}
		\includegraphics[width=\textwidth]{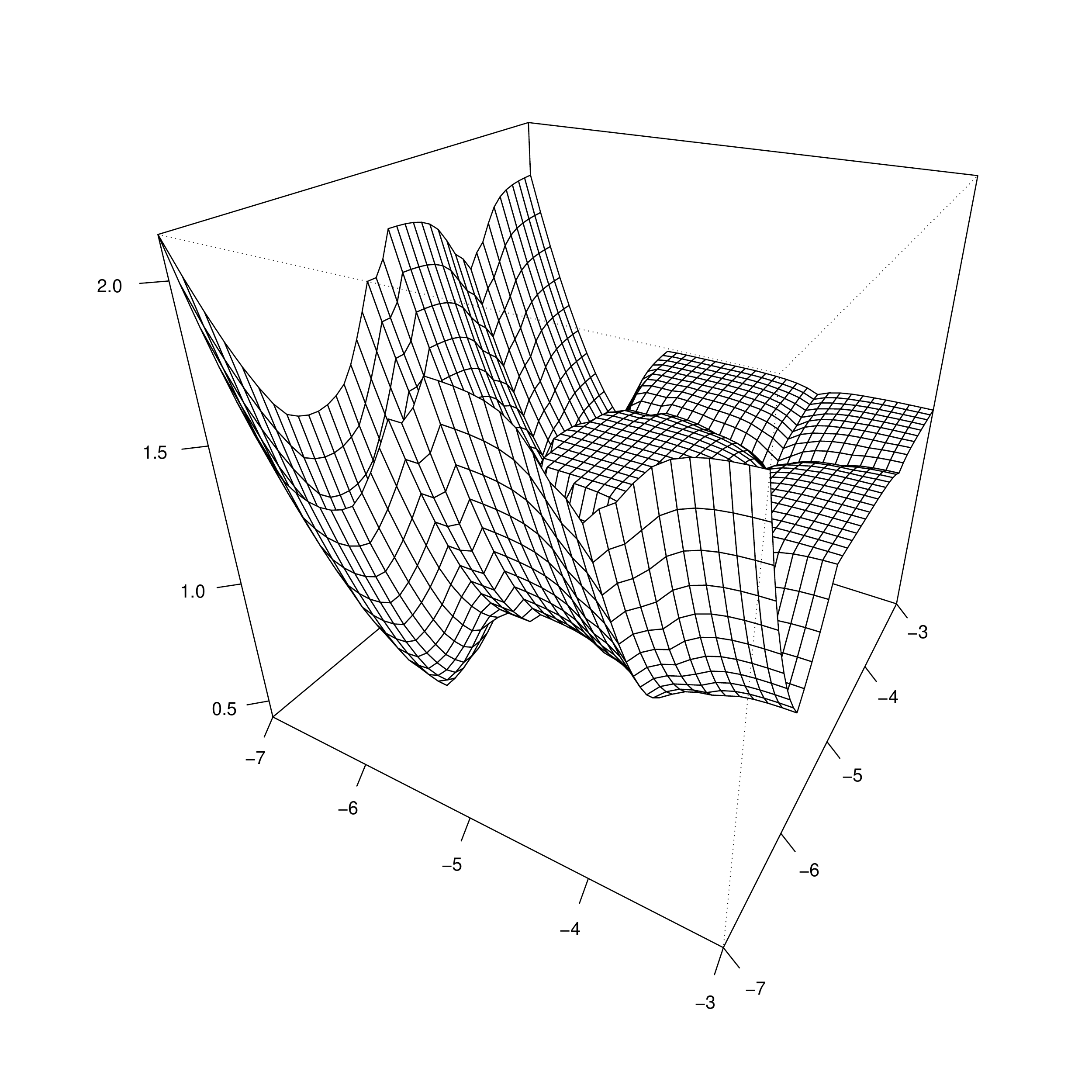}	
  	\caption{\scriptsize{Brankes Multi Worst}} \label{fig:BrankesMultipeakWorst}
	\end{subfigure}
	\begin{subfigure}[t]{.24\textwidth}
		\includegraphics[width=\textwidth]{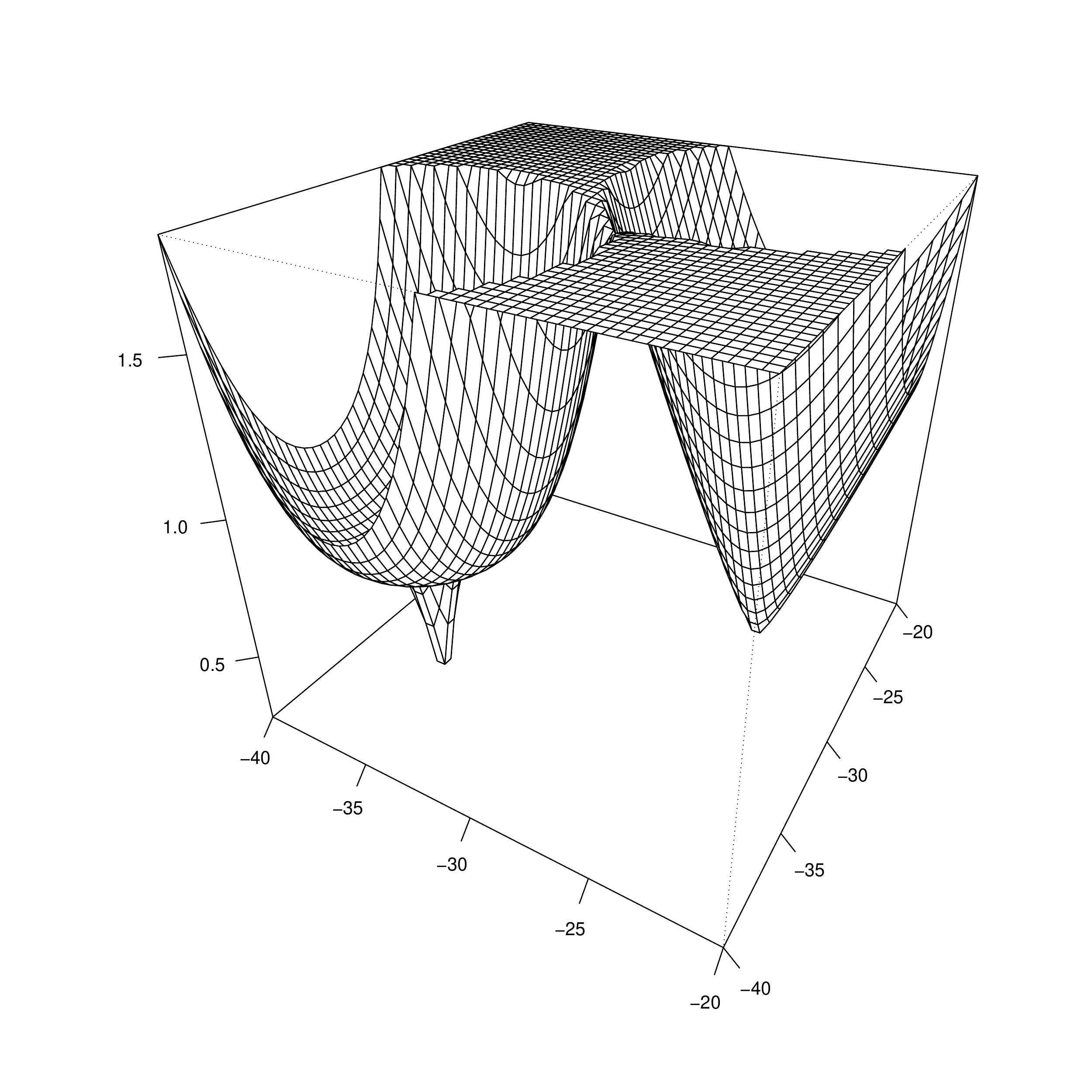}		
	  \caption{\scriptsize{Pickelhaube Nom}}	\label{fig:PickelhaubeNom}
	\end{subfigure}%
	\begin{subfigure}[t]{.24\textwidth}
		\includegraphics[width=\textwidth]{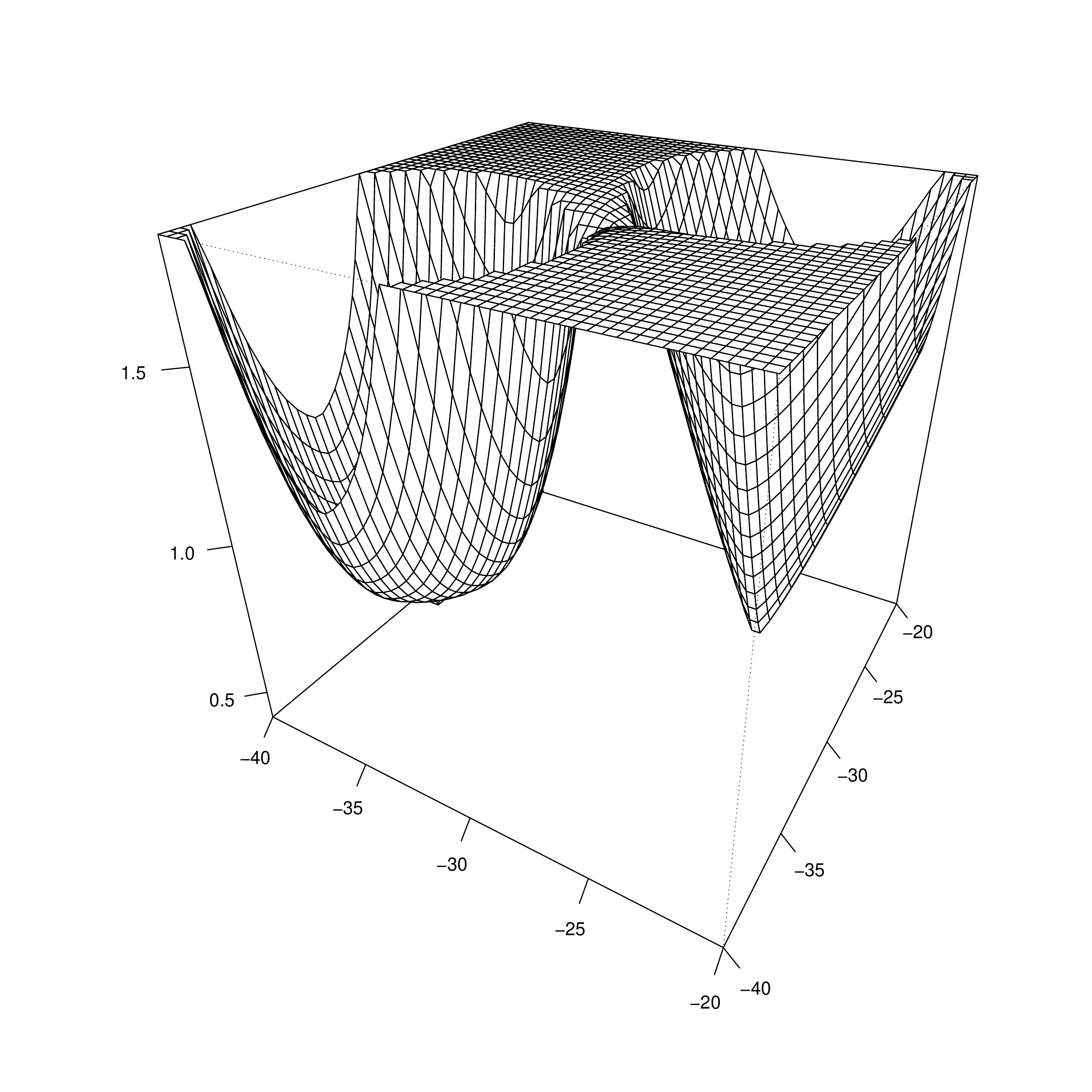}	
  	\caption{\scriptsize{Pickelhaube Worst}} \label{fig:PickelhaubeWorst}
	\end{subfigure}	
	\begin{subfigure}[t]{.24\textwidth}
		\includegraphics[width=\textwidth]{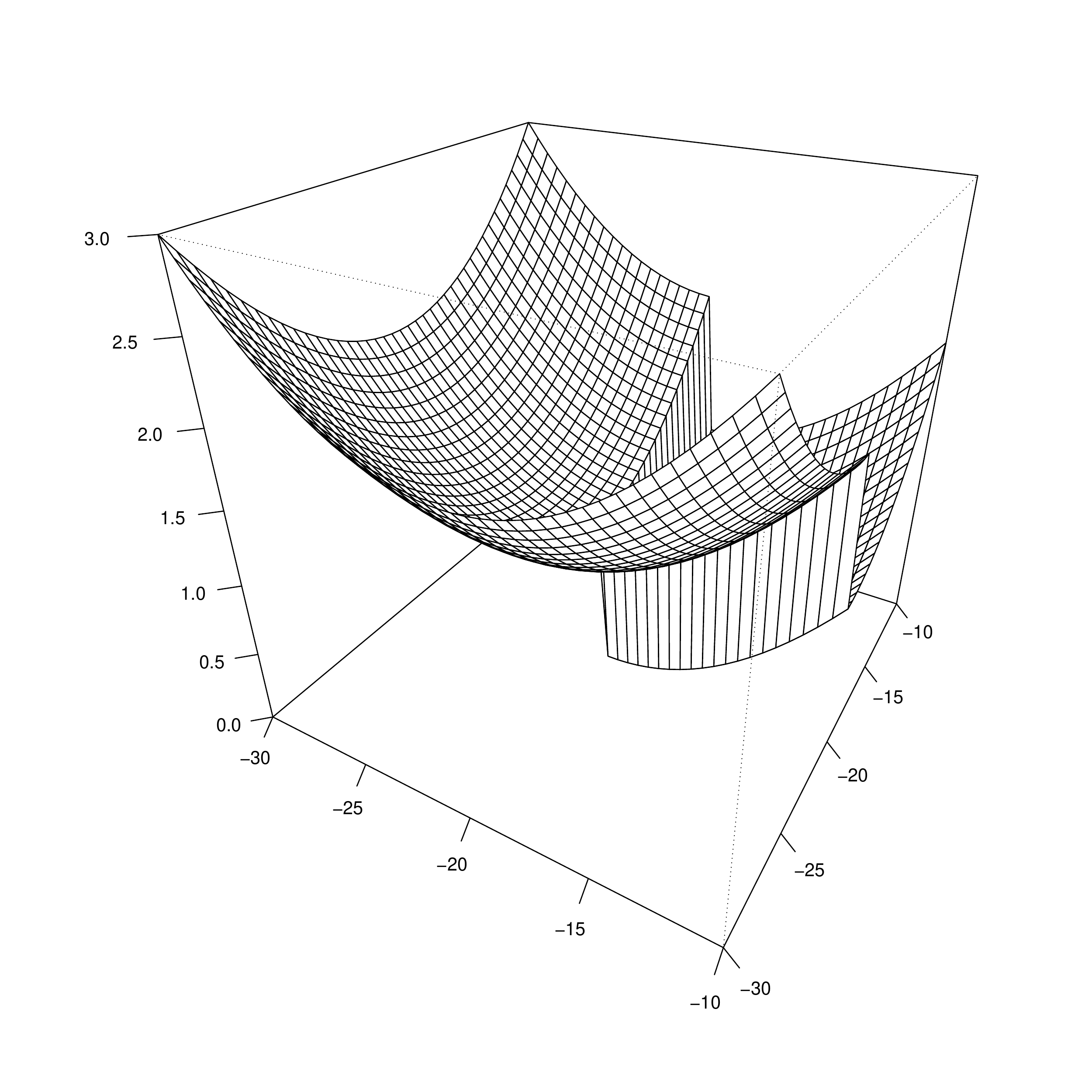}	
		\caption{\scriptsize{Heaviside S Nom}} \label{fig:HeavisideSphereNom}
	\end{subfigure}%
	\begin{subfigure}[t]{.24\textwidth}	
		\includegraphics[width=\textwidth]{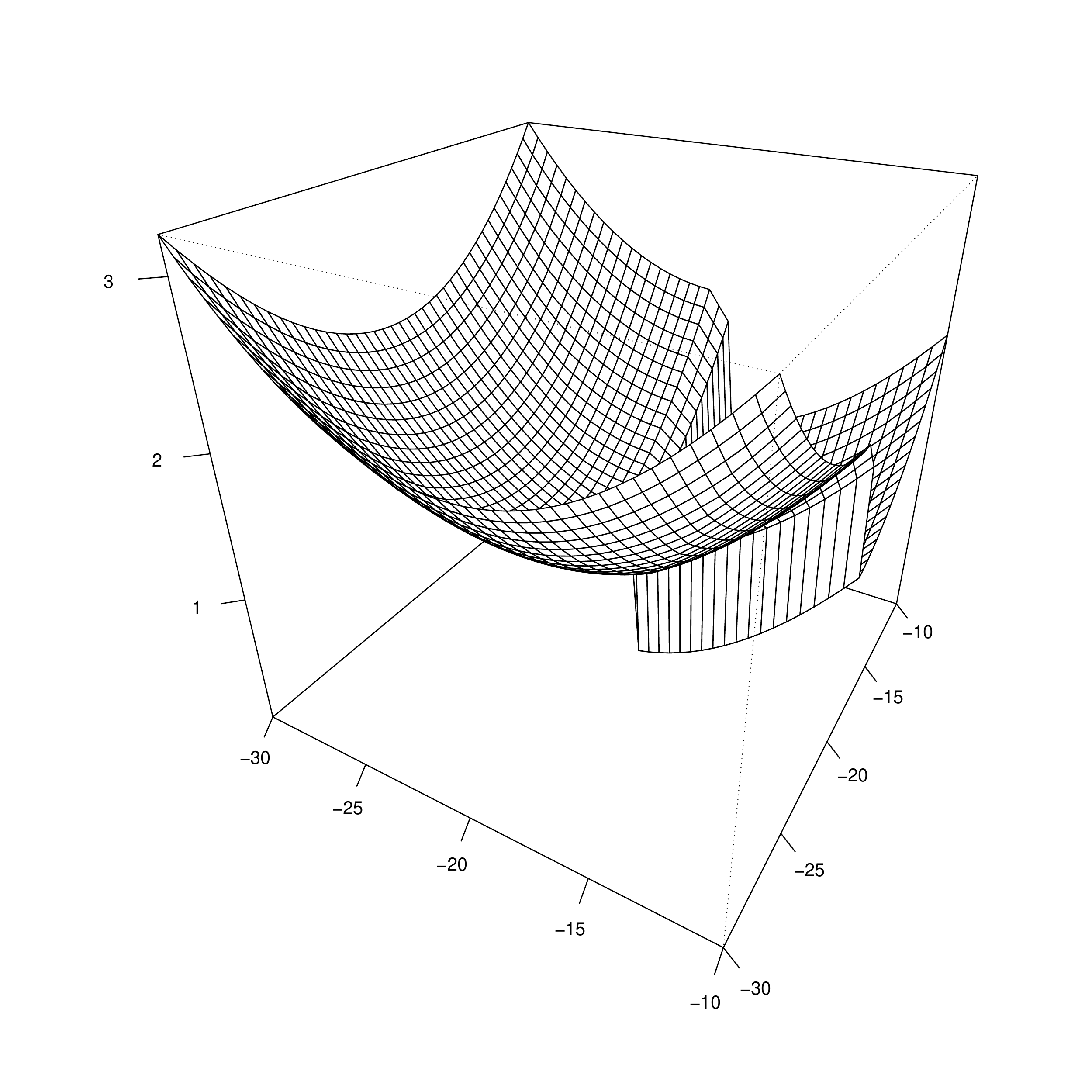}	
  	\caption{\scriptsize{Heaviside S Worst}} \label{fig:HeavisideSphereWorst}
	\end{subfigure}	
  \begin{subfigure}[t]{.24\textwidth}
		\includegraphics[width=\textwidth]{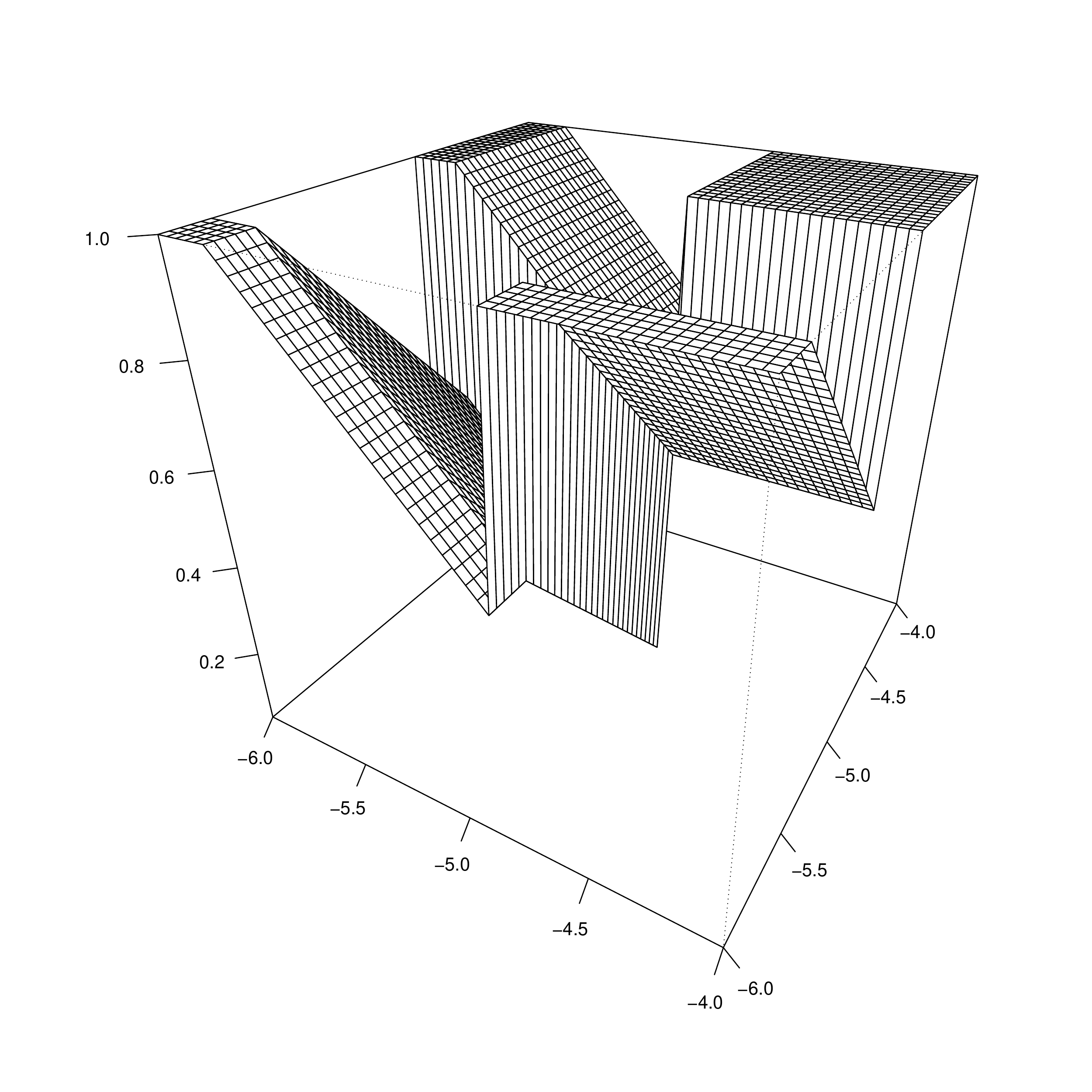}	
	  \caption{\scriptsize{Sawtooth Nom}} \label{fig:SawtoothNom}
	\end{subfigure}%
	\begin{subfigure}[t]{.24\textwidth}
		\includegraphics[width=\textwidth]{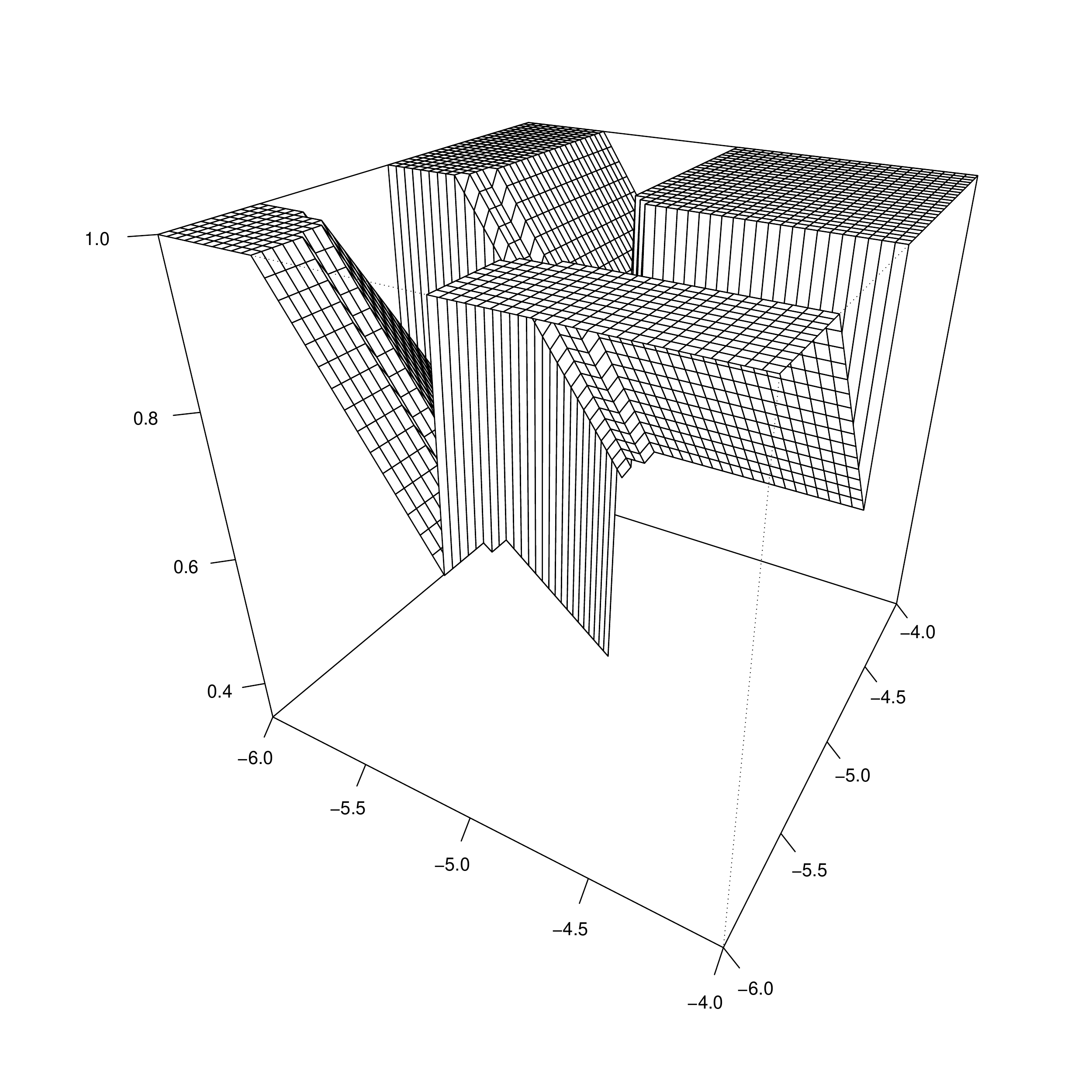}	
  	\caption{\scriptsize{Sawtooth Worst}} \label{fig:SawtoothWorst}
	\end{subfigure}	
	\begin{subfigure}[t]{.24\textwidth}
		\includegraphics[width=\textwidth]{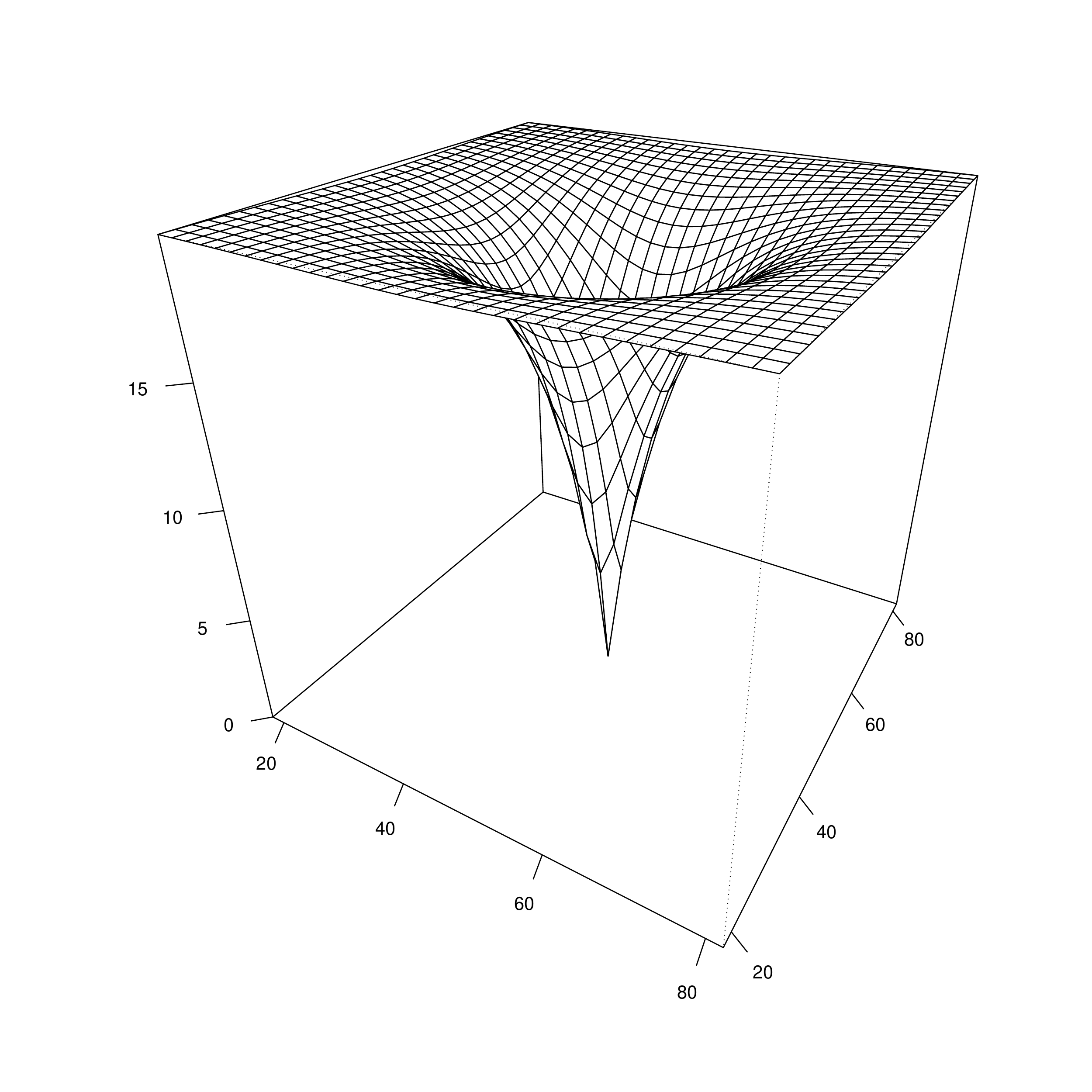}	
  	\caption{\scriptsize{Ackley Nom}} \label{fig:AckleyNom}
	\end{subfigure}%
	\begin{subfigure}[t]{.24\textwidth}
		\includegraphics[width=\textwidth]{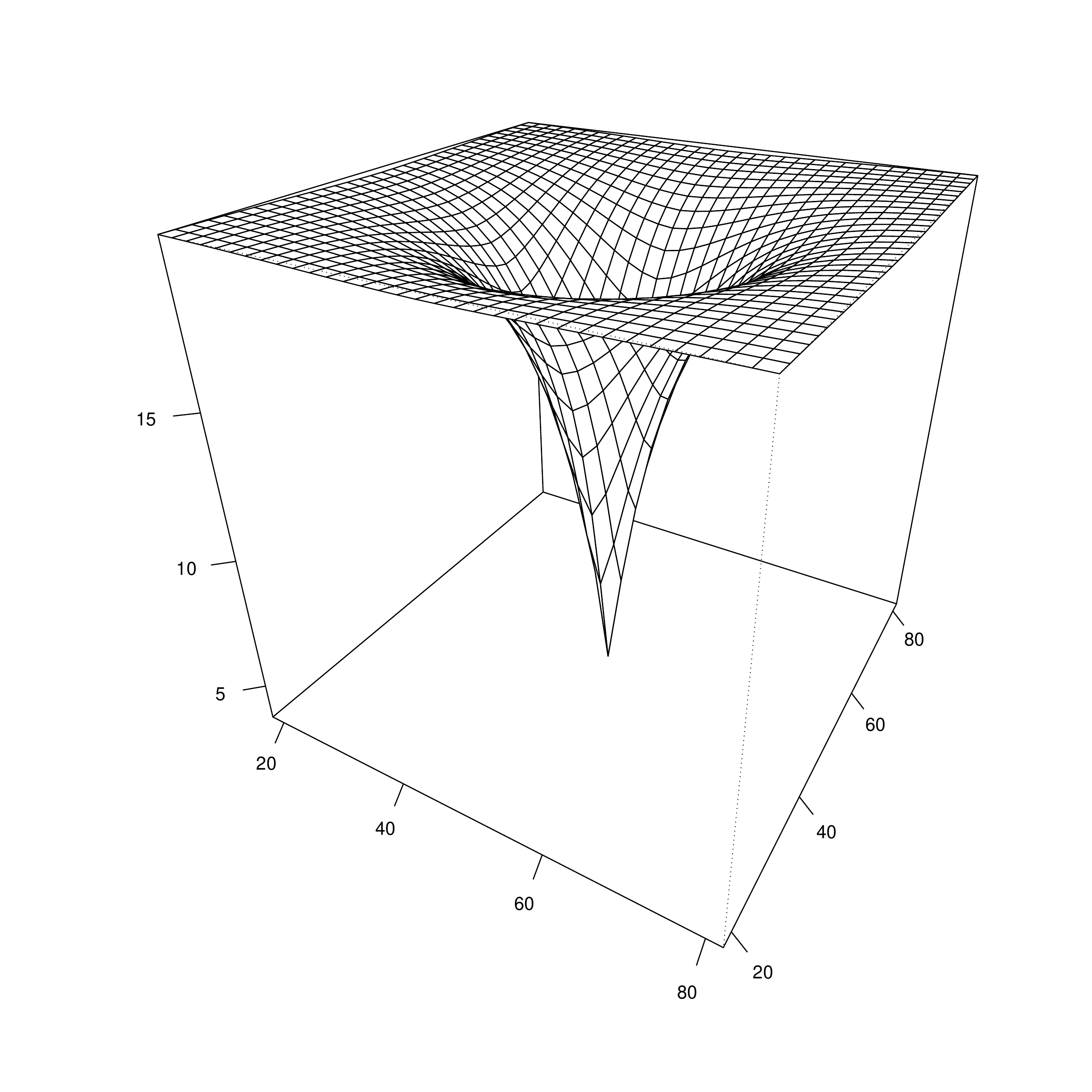}	
  	\caption{\scriptsize{Ackley Worst}} \label{fig:AckleyWorst}
	\end{subfigure}
	\begin{subfigure}[t]{.24\textwidth}
		\includegraphics[width=\textwidth]{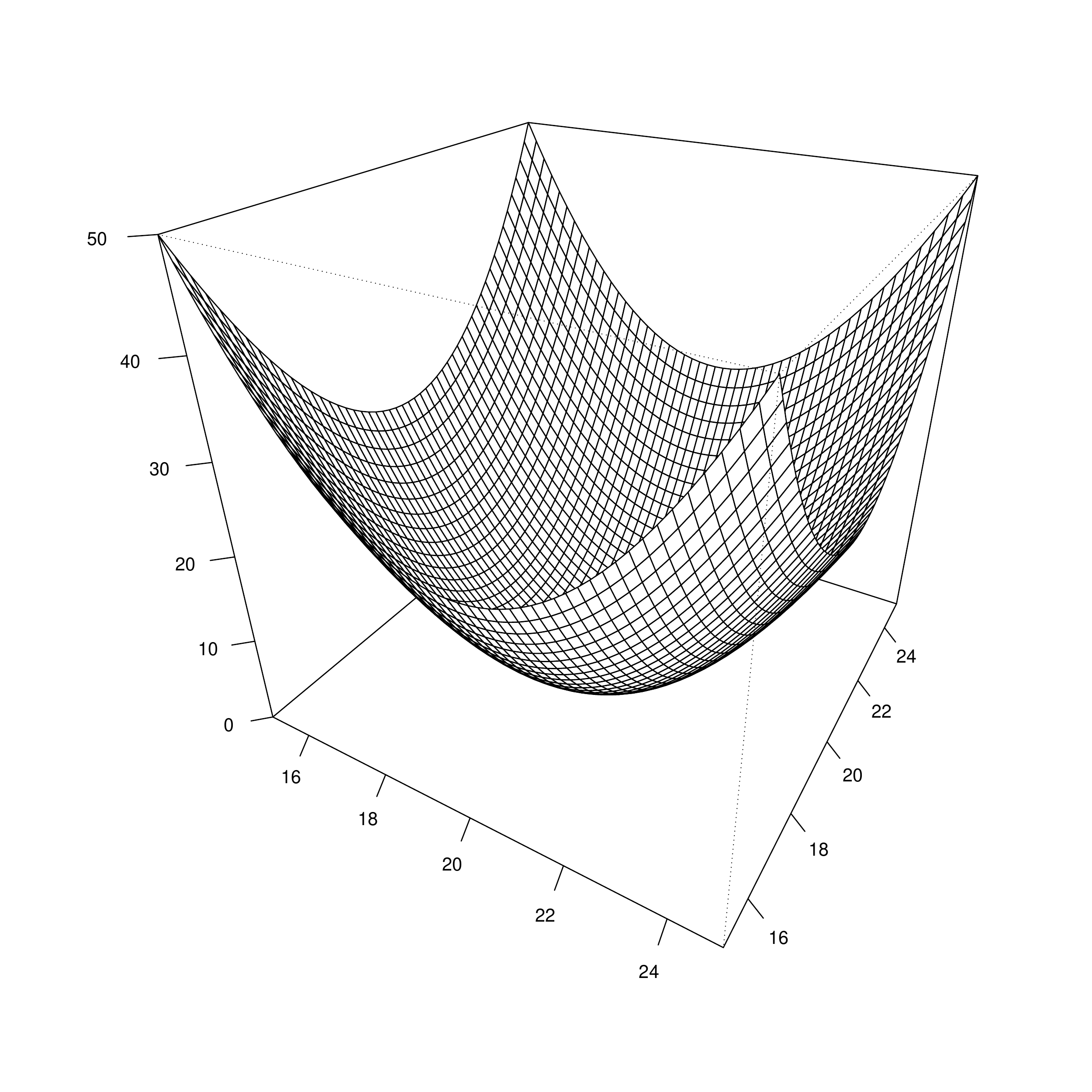}
  	\caption{\scriptsize{Sphere Nom}} \label{fig:SphereNom}
	\end{subfigure}%
	\begin{subfigure}[t]{.24\textwidth} 
		\includegraphics[width=\textwidth]{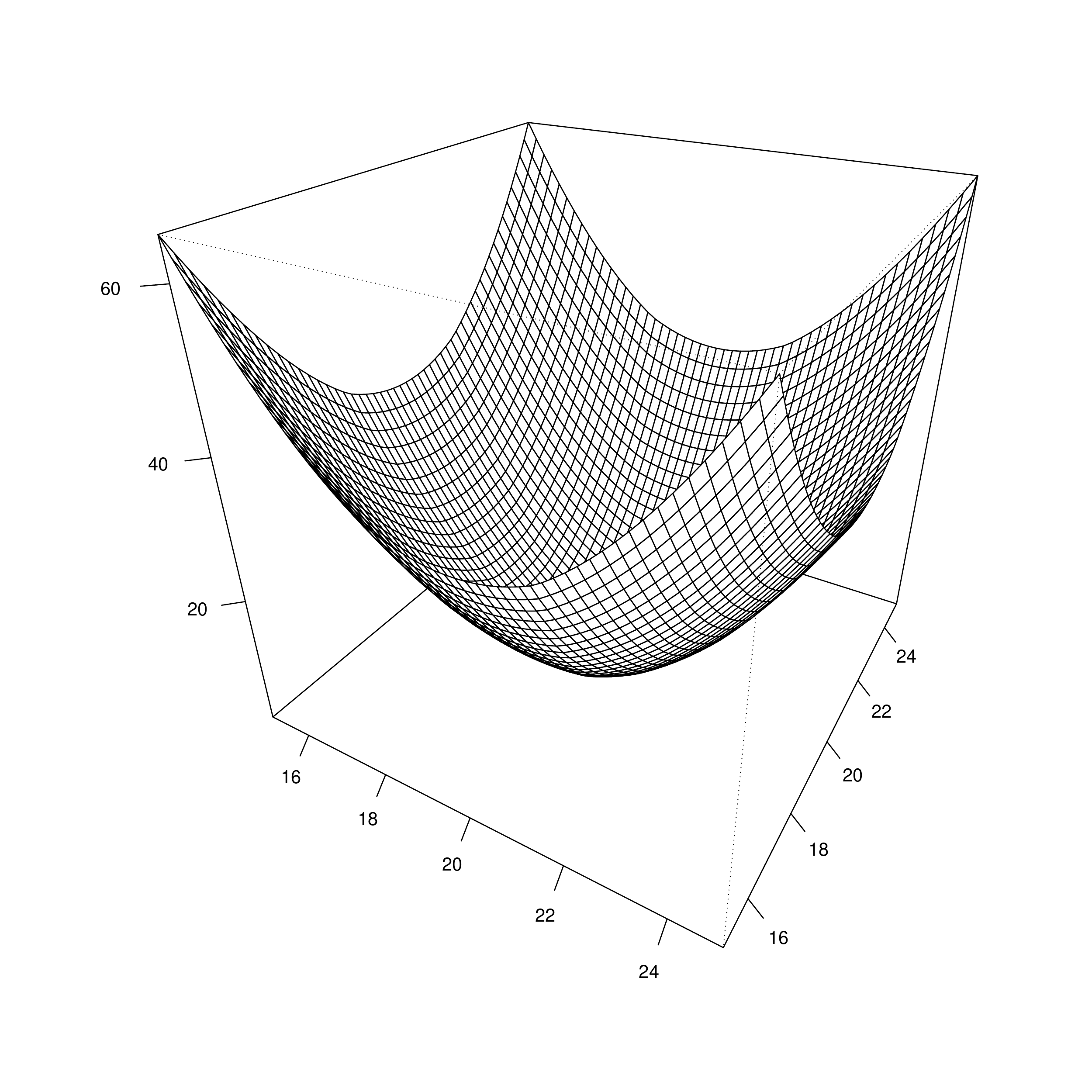}	
  	\caption{\scriptsize{Sphere Worst}} \label{fig:SphereWorst}
	\end{subfigure}		
	\begin{subfigure}[t]{.24\textwidth}
		\includegraphics[width=\textwidth]{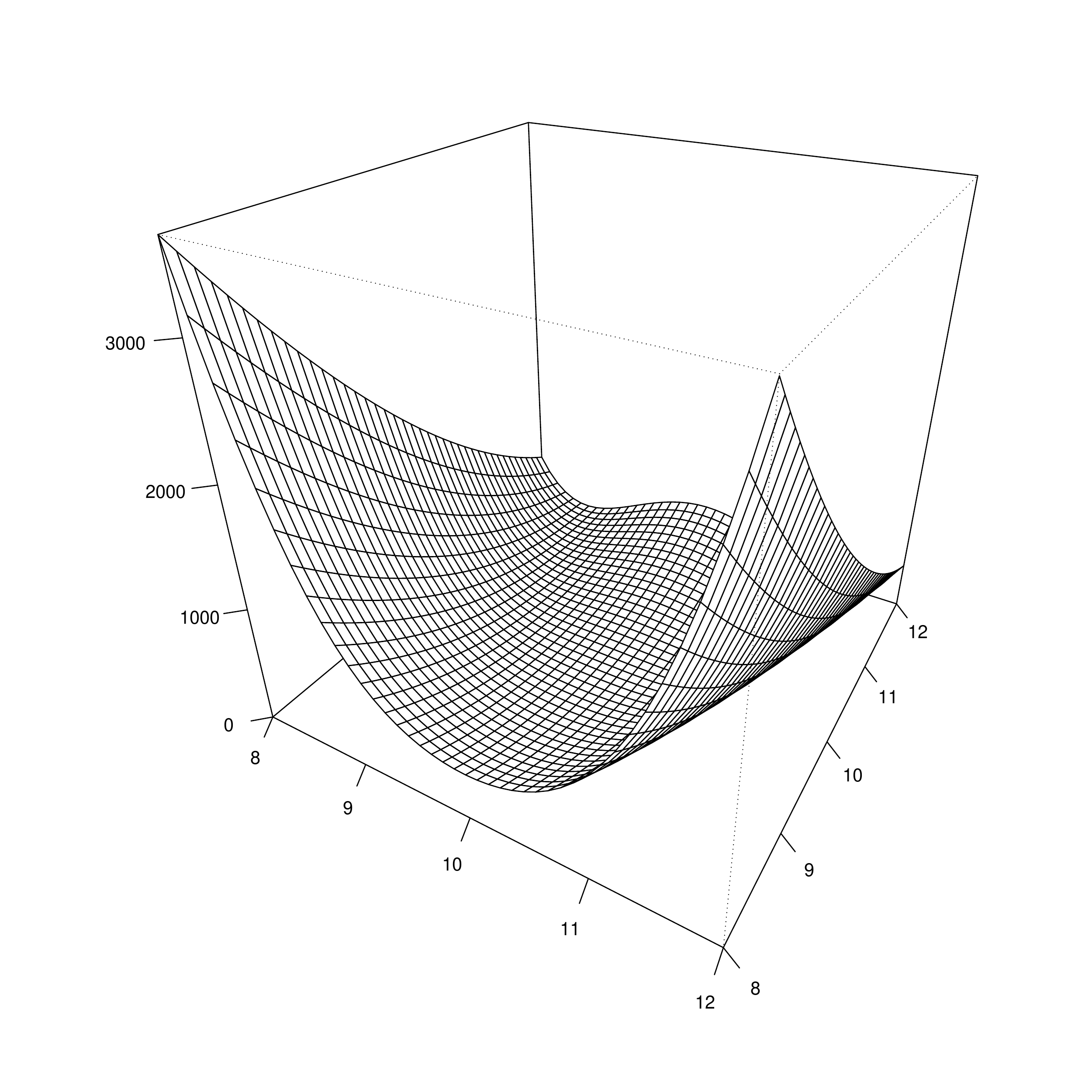}
	  \caption{\scriptsize{Rosenbrock Nom}} \label{fig:RosenbrockNom}
	\end{subfigure}%
	\begin{subfigure}[t]{.24\textwidth}
		\includegraphics[width=\textwidth]{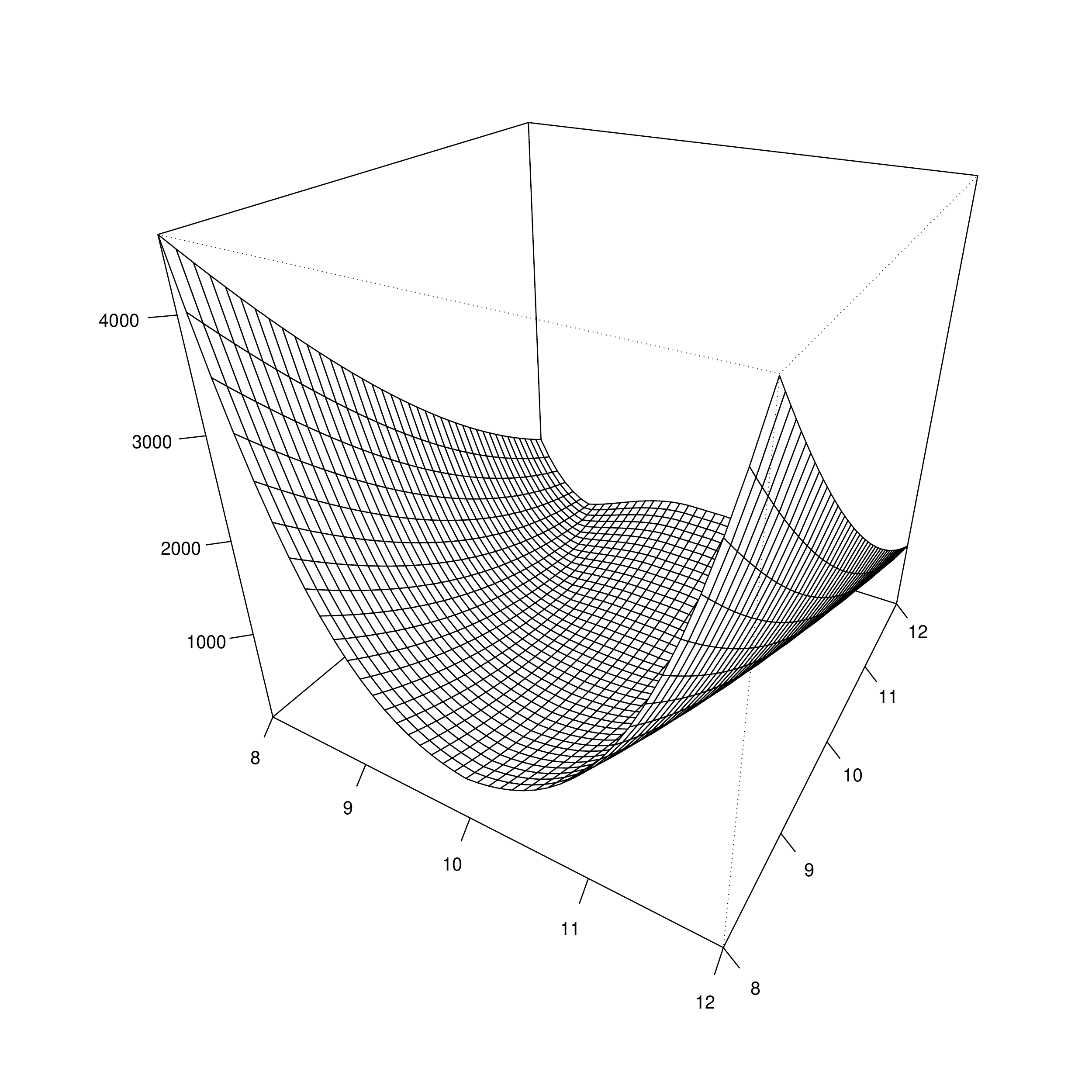}		
  	\caption{\scriptsize{Rosenbrock Worst}} \label{fig:RosenbrockWorst}
	\end{subfigure}	
	
	\vspace{2mm} 
	
	\caption{Plots of 2D versions of the functions used in our experimental testing.}
	\label{fig:rPSOtestSuite}
\end{figure}

The mathematical descriptions for the 10 test functions used in the experimental testing are given below, with 3D plots of their 2D versions shown in Figure~\ref{fig:rPSOtestSuite}. All functions are taken from the literature: \cite{Branke1998, KruisselbrinkEmmerichBack2010, KruisselbrinkReehuisDeutzBackEmmerich2011, Kruisselbrink2012, JamilYang2013}.

\begingroup
\allowdisplaybreaks
\begin{align*}
\text{Rastrigin:}\quad &  f(\pmb{x}) = 10 n + \sum_{i=1}^{n}{[{(x_i-20)}^{2} - 10 \cos(2\pi (x_i-20))]} \\
& \X = [14.88, 25.12]^n  \\
\text{Multipeak F1:}\quad & f(\pmb{x}) = - \frac{1}{n} \sum_{i=1}^n g(x_i) \\
& g(x_i) =
\begin{cases}
	e^{-2 \ln 2 ( \frac{(x_i+5)-0.1}{0.8} )^2} \sqrt{ \left| \sin(5\pi (x_i+5)) \right| }  &\quad\text{if } 0.4 < x_i+5 \le 0.6 \text{ ,} \\
  e^{-2 \ln 2 ( \frac{(x_i+5)-0.1}{0.8} )^2} \sin^{6}(5\pi (x_i+5))  &\quad\text{otherwise} \\
\end{cases}\\
& \X = [-5, -4]^n \\
\text{Multipeak F2:}\quad & f(\pmb{x}) = \frac{1}{n} \sum_{i=1}^n g(x_i) \\
& g(x_i) = 2\sin(10 \exp (-0.2(x_i-10)) (x_i-10)) \exp (-0.25(x_i-10))\\
& \X = [10, 20]^n \\
\text{Branke's Multipeak:}\quad &f(\pmb{x}) = \max \lbrace c_1, c_2 \rbrace - \frac{1}{n} \sum_{i-1}^n g(x_i) \\
& g(x_i) =
\begin{cases}
	c_1 \Bigg(1-{\frac{4((x_i+5)+\frac{b_1}{2})^2}{b_1^2}} \Bigg)  &\quad\text{if } -b_1 \le (x_i+5) < 0 \text{ ,}
	\\
	c_2 \cdot 16^{\frac{-2 \left| b_2 - 2(x_i+5) \right|}{b_2}} &\quad\text{if } 0 \le (x_i+5) \le b_2 \text{ ,}
	\\
  0  &\quad\text{otherwise} \\
\end{cases} \\
& b_1=2, b_2=2, c_1=1, c_2=1.3 \\ 
& \X = [-7, -3]^n\\
\text{Pickelhaube}\quad & f(\pmb{x}) = \frac{5}{5-\sqrt{5}} - \max \lbrace g_0(\pmb{x}), g_{1a}(\pmb{x}), g_{1b}(\pmb{x}), g_2(\pmb{x}) \rbrace \\
& g_0(\pmb{x}) = \frac{1}{10} e^{-\frac{1}{2} \| \pmb{x}+30 \|} \\
&g_{1a}(\pmb{x}) = \frac{5}{5-\sqrt{5}} \Bigg( 1-\sqrt{\frac{\| \pmb{x}+30+5 \|}{5\sqrt{n}}} \Bigg) \\
&g_{1b}(\pmb{x}) = c_1 \Bigg( 1-\Bigg(\frac{\| \pmb{x}+30+5 \|}{5\sqrt{n}} \Bigg)^4 \Bigg) \\
&g_2(\pmb{x}) =  c_2 \Bigg( 1-\Bigg(\frac{\| \pmb{x}+30-5 \|}{5\sqrt{n}} \Bigg)^{d_2} \Bigg) \\
& c_1=625/624, c_2=1.5975, d_2=1.1513 \\
& \X = [-40, -20]^n \\
\text{Heaviside Sphere}\quad & f(\pmb{x}) = \Bigg( 1 - \prod_{i=1}^n g(x_i) \Bigg) +  \sum_{i=1}^n \Bigg(\frac{(x_i+20)}{10} \Bigg)^2 \\
& g(x_i) =
\begin{cases}
	0  &\quad\text{if } 0 < (x_i+20) \text{ ,} \\
  1  &\quad\text{otherwise} \\
\end{cases} \\
& \X =  [-30, -10]^n\\
\text{Sawtooth}\quad & f(\pmb{x}) = 1 - \frac{1}{n} \sum_{i=1}^n g(x_i) \\
& g(x_i) =
\begin{cases}
	(x_i+5) + 0.8  &\quad\text{if } -0.8 \le (x_i+5) < 0.2 \text{ ,} \\
  0  &\quad\text{otherwise} \\
\end{cases} \\
& \X =  [-6, -4]^n \\
\text{Ackleys}\quad & f(\pmb{x}) = -20 \exp \Bigg( -0.2 \sqrt{\frac{1}{n} \sum_{i=1}^n (x_i-50)^2} \Bigg) \\ 
& \phantom{f(\pmb{x}) = }- \exp \Bigg( \frac{1}{n} \sum_{i=1}^n \cos(2\pi (x_i-50)) \Bigg) + 20 + \exp(1) \\
& \X = [17.232, 82.768]^n \\
\text{Sphere}\quad & f(\pmb{x}) = \sum_{i=1}^n (x_i-20)^2 \\
& \X = [15, 25]^n \\
\text{Rosenbrock}\quad & f(\pmb{x}) = \sum_{i=1}^{n-1} [100((x_{i+1}-10) - (x_i-10)^2)^2 + ((x_i-10)-1)^2]\\
& \X = [7.952, 12.048]
\end{align*}

\endgroup

\clearpage 
\newpage


\bibliographystyle{alpha}

\newcommand{\etalchar}[1]{$^{#1}$}

\end{document}